\documentclass[11pt,twoside,leqno]{aomamlt2e}
  \pageno{1}
\received{April 14, 2006; Accepted on October 3, 2006}
\newtheorem{theorem}{Theorem}[section]
\newtheorem{lemma}{Lemma}[section]
\newtheorem{proposition}{Proposition}[section]

\newtheorem{definition}{Definition}[section]
\newtheorem{remark}{Remark}[section]

\newcommand{\bR}{ {\bf R}} %%%% real
\newcommand{\R}{ {\bf R}} %%%% real n-dim
\newcommand{\dist}{ \mbox{dist}}
\newcommand{\divg}{ \mbox{div}}

\newcommand{\PtUpL}{{P_1}}
\newcommand{\PtLwL}{{P_2}}
\newcommand{\PtLwR}{{P_3}}
\newcommand{\PtUpR}{{P_4}}

\newcommand{\setK}{{\mathcal K}}
\newcommand{\setS}{{\mathcal S}(N)}
\newcommand{\Kphi}{\phi}

\newcommand{\defd}{:=}
\newcommand{\Nl}{{\mathcal N}}
\newcommand{\Ml}{{\mathcal M}}
\newcommand{\epsP}{\sigma}
\newcommand{\vr}{\varrho}
\newcommand{\Dom}{{\mathcal D}}
\newcommand{\grad}{{D}}
\newcommand{\mxx}{{\xi}}
\newcommand{\mxy}{{\eta}}

\newcommand{\ElDom}{\Omega}
\newcommand{\ElDomS}{\Omega'}
\newcommand{\hElDomS}{\hat\Omega'}
\newcommand{\ElDomU}{\Omega''}
\newcommand{\DomS}{{\mathcal D}'}
\newcommand{\DomU}{{\mathcal D}''}
\newcommand{\sonic}{\Gamma_{sonic}}
\newcommand{\shock}{\Gamma_{shock}}
\newcommand{\wedgeB}{\Gamma_{wedge}}
\newcommand{\corners}{{\mathcal P}}

\newcommand{\UreflAbs}{q_2}

\newcommand{\mX}{\mathcal X}

\newcommand{\Ex}{{{\mathcal E}_2}}

\newcommand{\ExCFn}{{\mathcal {P}}}

\newcommand{\Dnu}{\mu}

\newcommand{\bea}{{\begin{eqnarray}}}
\newcommand{\eea}{{\end{eqnarray}}}

\begin{document}
\currannalsline{0}{2006}
\title{Global Solutions of Shock Reflection by Large-Angle Wedges for
Potential Flow} \shorttitle{Global Solutions of Shock Reflection by
Large-Angle Wedges} \twoauthors{Gui-Qiang Chen}{Mikhail Feldman}
\institution{Gui-Qiang Chen\\
 Department of Mathematics, Northwestern University,
 Evanston, IL 60208-2730, USA;\\
School of Mathematical Sciences, Fudan University,
 Shanghai 200433, PRC;
\email{gqchen@math.northwestern.edu}\\
Mikhail Feldman\\
Department of Mathematics, University of Wisconsin, Madison, WI
53706-1388, USA;
         \email{feldman@math.wisc.edu}}
%\author{Gui-Qiang Chen\\
%  Department of Mathematics\\
%         Northwestern University \\
%         Evanston, IL 60208-2730, USA; School of Mathematical Sciences,
% Fudan University, Shanghai 200433,
%PRC \\
%gqchen@math.northwestern.edu}
%\author{Mikhail Feldman\\
% Department of Mathematics\\
%         University of Wisconsin\\
%         Madison, WI 53706-1388, USA\\
%         feldman@math.wisc.edu}
%\keywords{Regular shock reflection,
%existence, stability, global solutions, self-similar,
%elliptic-hyperbolic, nonlinear equations, second-order, mixed type,
%transonic shocks, free boundary problems, degenerate elliptic,
%corner singularity, Euler equations, compressible flow, analytical
%approach, iteration methods}
%\subjclass{35M10,35J65,35R35,35J70,76H05,76L05,35B60,35B35,35B65}
\date{March 28, 2006}

\begin{abstract}
When a plane shock hits a wedge head on, it experiences a
reflection-diffraction process and then a self-similar reflected
shock moves outward as the original shock moves forward in time.
Experimental, computational, and asymptotic analysis has shown that
various patterns of shock reflection may occur, including regular
and Mach reflection. However, most of the fundamental issues for
shock reflection have not been understood, including the global
structure, stability, and transition of the different patterns of
shock reflection. Therefore, it is essential to establish the global
existence and structural stability of solutions of shock reflection
in order to understand fully the phenomena of shock reflection. On
the other hand, there has been no rigorous mathematical result on
the global existence and structural stability of shock reflection,
including the case of potential flow which is widely used in
aerodynamics. Such problems involve several challenging difficulties
in the analysis of nonlinear partial differential equations such as
mixed equations of elliptic-hyperbolic type, free boundary problems,
and corner singularity where an elliptic degenerate curve meets a
free boundary. In this paper we develop a rigorous mathematical
approach to overcome these difficulties involved and establish a
global theory of existence and stability for shock reflection by
large-angle wedges for potential flow. The techniques and ideas
developed here will be useful for other nonlinear problems involving
similar difficulties.
\end{abstract}
%\maketitle

\section{Introduction}
\numberwithin{equation}{section}
%\numberwithin{equation}{section}

We are concerned with the problems of shock reflection by wedges.
These problems arise not only in many important physical situations
but also are fundamental in the mathematical theory of
multidimensional conservation laws since their solutions are
building blocks and asymptotic attractors of general solutions to
the multidimensional Euler equations for compressible fluids (cf.
Courant-Friedrichs \cite{CF}, von Neumann \cite{Neumann}, and
Glimm-Majda \cite{GlimmMajda}; also see
\cite{BD,ChangChen,GlimmK,LaxLiu,Morawetz2,Serre,VD}). When a plane
shock hits a wedge head on, it experiences a reflection-diffraction
process and then a self-similar reflected shock moves outward as the
original shock moves forward in time. The complexity of reflection
picture was first reported by Ernst Mach \cite{Mach} in 1878, and
experimental, computational, and asymptotic analysis has shown that
various patterns of shock reflection may occur, including regular
and Mach reflection (cf.
\cite{BD,GRT,GlimmMajda,hunter1,HK,KB,Morawetz2,VD,Neumann}).
However, most of the fundamental issues for shock reflection have
not been understood, including the global structure, stability, and
transition of the different patterns of shock reflection. Therefore,
it is essential to establish the global existence and structural
stability of solutions of shock reflection in order to understand
fully the phenomena of shock reflection. On the other hand, there
has been no rigorous mathematical result on the global existence and
structural stability of shock reflection, including the case of
potential flow which is widely used in aerodynamics (cf.
\cite{Bers1,CC,GlimmMajda,MajdaTh,Morawetz2}). One of the main
reasons is that the problems involve several challenging
difficulties in the analysis of nonlinear partial differential
equations such as mixed equations of elliptic-hyperbolic type, free
boundary problems, and corner singularity where an elliptic
degenerate curve meets a free boundary. In this paper we develop a
rigorous mathematical approach to overcome these difficulties
involved and establish a global theory of existence and stability
for shock reflection by large-angle wedges for potential flow. The
techniques and ideas developed here will be useful for other
nonlinear problems involving similar difficulties.

\medskip
The Euler equations for potential flow consist of the conservation
law of mass and the Bernoulli law for the density $\rho$ and
velocity potential $\Phi$:
\begin{eqnarray}
&&\partial_t\rho + \divg_{\bf x}(\rho\nabla_{\bf x}\Phi)=0, \label{1.1.1} \\
&&\partial_t\Phi +\frac{1}{2}|\nabla_{\bf x}\Phi|^2+i(\rho)=K,
\label{1.1.2}
\end{eqnarray}
where $K$ is the Bernoulli constant determined by the incoming flow
and/or boundary conditions, and
$$
i'(\rho)=p'(\rho)/\rho=c^2(\rho)/\rho
$$
with $c(\rho)$ being the sound speed. For polytropic gas,
$$
p(\rho)=\kappa\rho^\gamma,\qquad
c^2(\rho)=\kappa\gamma\rho^{\gamma-1}, \qquad \gamma>1, \,\,
\kappa>0.
$$
Without loss of generality, we choose $\kappa=(\gamma-1)/\gamma$ so
that
$$
i(\rho)=\rho^{\gamma-1}, \qquad c(\rho)^2=(\gamma-1)\rho^{\gamma-1},
$$
which can be achieved by the  following scaling:
$$
({\bf x},t,K)\to  (\alpha {\bf x}, \alpha^2 t, \alpha^{-2} K), \quad
\alpha^2=\kappa\gamma/(\gamma-1).
$$
Equations \eqref{1.1.1}--\eqref{1.1.2} can written as the following
nonlinear equation of second order:
\begin{equation}\label{nonlinear-second-order:1}
\partial_t\hat{\rho}\big(K-\partial_t\Phi-\frac{1}{2}|\nabla_{\bf x}\Phi|^2\big)
+\divg_{\bf
x}\big(\hat{\rho}(K-\partial_t\Phi-\frac{1}{2}|\nabla_{\bf
x}\Phi|^2)
        \nabla_{\bf x}\Phi\big)=0,
\end{equation}
where $\hat{\rho}(s)=s^{1/(\gamma-1)}=i^{-1}(s)$ for $s\ge 0$.

\medskip
When a plane shock in the $({\bf x},t)$--coordinates, ${\bf
x}=(x_1,x_2)\in\R^2$, with left state $ (\rho,\nabla_{\bf
x}\Psi)=(\rho_1, u_1,0) $ and right state $(\rho_0, 0,0), u_1>0,
\rho_0<\rho_1$,
hits a symmetric wedge
$$
W:=\{|x_2|< x_1 \tan\theta_w, x_1>0\}
$$
head on,
it experiences a reflection-diffraction process, and the reflection
problem can be formulated as the following mathematical problem.

\medskip
{\bf Problem 1 (Initial-Boundary Value Problem)}. {\it Seek a
solution of system \eqref{1.1.1}--\eqref{1.1.2} with
$K=\rho_0^{\gamma-1}$, the initial condition at $t=0$:
\begin{equation}\label{initial-condition}
(\rho,\Phi)|_{t=0} =\begin{cases}
(\rho_0, 0) \qquad&  \mbox{for}\,\, |x_2|>x_1\tan\theta_w, x_1>0,\\
(\rho_1, u_1 x_1) \qquad &\mbox{for}\,\, x_1<0,
\end{cases}
\end{equation}
and the slip boundary condition along the wedge boundary $\partial
W$:
\begin{equation}\label{boundary-condition}
\nabla\Phi\cdot \nu|_{\partial W}=0,
\end{equation}
where $\nu$ is the exterior unit normal to $\partial W$ (see Fig.
{\rm 1}).}

\begin{figure}[h]
\centering
\includegraphics[height=2.5in,width=2.7in]{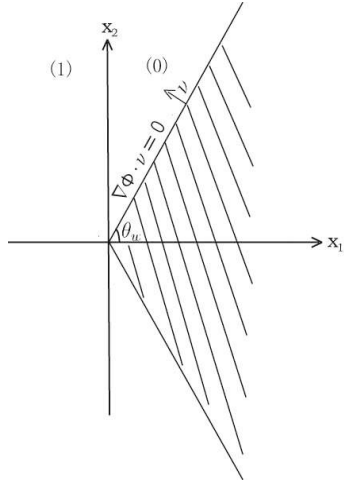}  %%% for latex (non-pdf)
\caption[]{Initial-boundary value problem} \label{fig:IBVP}
\end{figure}

\medskip
Notice that the initial-boundary value problem
\eqref{1.1.1}--\eqref{boundary-condition} is invariant under the
self-similar scaling:
$$
({\bf x}, t)\to (\alpha {\bf x}, \alpha t), \quad (\rho, \Phi)\to
(\rho, \Phi/\alpha) \qquad \quad\mbox{for}\quad \alpha\ne 0.
$$
Thus, we seek self-similar solutions with the form
$$
\rho({\bf x},t)=\rho(\xi,\eta), \quad \Phi({\bf
x},t)=t\,\psi(\xi,\eta) \qquad\quad \mbox{for}\quad (\xi,\eta)={\bf
x}/t.
$$
Then the pseudo-potential function
$\varphi=\psi-\frac{1}{2}(\xi^2+\eta^2)$ satisfies the following
Euler equations for self-similar solutions:
\begin{eqnarray}
&&\divg\,(\rho\, D\varphi)+2\rho=0, \label{1.1.3}\\
&&\frac{1}{2}|D\varphi|^2+\varphi+ \rho^{\gamma-1}
  =\rho_0^{\gamma-1},
   \label{1.1.4}
\end{eqnarray}
where the divergence $\divg$ and gradient $D$ are with respect to
the self-similar variables $(\xi,\eta)$. This implies that the
pseudo-potential function $\varphi(\xi,\eta)$ is governed by the
following potential flow equation of second order:
\begin{equation}
\divg\, \big(\rho(|D\varphi|^2, \varphi)D\varphi\big)
+2\rho(|D\varphi|^2, \varphi)=0 \label{1.1.5}
\end{equation}
with
\begin{equation}
\rho(|D\varphi|^2, \varphi)
=\hat{\rho}(\rho_0^{\gamma-1}-\varphi-\frac{1}{2}|D\varphi|^2).
\label{1.1.6}
\end{equation}
Then we have
\begin{equation}\label{c-through-density-function}
c^2=c^2(|D\varphi|^2,\varphi,\rho_0^{\gamma-1})
=(\gamma-1)(\rho_0^{\gamma-1}-\frac{1}{2}|D\varphi|^2-\varphi).
\end{equation}

Equation \eqref{1.1.5} is a mixed equation of elliptic-hyperbolic
type. It is elliptic if and only if
\begin{equation}
|D\varphi| < c(|D\varphi|^2,\varphi,\rho_0^{\gamma-1}),
\label{1.1.8}
\end{equation}
which is equivalent to
\begin{equation}
|D \varphi| <c_*(\varphi, \rho_0, \gamma)
:=\sqrt{\frac{2(\gamma-1)}{\gamma+1}(\rho_0^{\gamma-1}-\varphi)}.
\label{1.1.8a}
\end{equation}
Shocks are discontinuities in the pseudo-velocity $D\varphi$. That
is, if $\Omega^+$ and $\Omega^-:=\Omega\setminus\overline{\Omega^+}$
are two nonempty open subsets of $\Omega\subset\bR^2$ and
$S:=\partial\Omega^+\cap\Omega$ is a $C^1$ curve where $D\varphi$
has a jump,  then $\varphi\in W^{1,1}_{loc}(\Omega)\cap
C^1(\Omega^\pm\cup S)\cap C^2(\Omega^\pm)$ is a global weak solution
of (\ref{1.1.5}) in $\Omega$ if and only if $\varphi$ is in
$W^{1,\infty}_{loc}(\Omega)$ and satisfies equation \eqref{1.1.5} in
$\Omega^\pm$ and the Rankine-Hugoniot condition on $S$:
\begin{equation}\label{FBConditionSelfSim-0}
\left[\rho(|D\varphi|^2,\varphi)D\varphi\cdot\nu\right]_S=0.
\end{equation}
The continuity of $\varphi$ is followed by the continuity of the
tangential derivative of $\varphi$ across $S$, which is a direct
corollary of irrotationality of the pseudo-velocity. The
discontinuity $S$ of $D\varphi$ is called a shock if $\varphi$
further satisfies the physical entropy condition that the
corresponding density function $\rho(|D\varphi|^2,\varphi)$
increases across $S$ in the pseudo-flow direction. We remark that
the Rankine-Hugoniot condition \eqref{FBConditionSelfSim-0} with the
continuity of $\varphi$ across a shock for \eqref{1.1.5} is also
fairly good approximation to the corresponding Rankine-Hugoniot
conditions for the full Euler equations for shocks of small
strength, since the errors are third-order in strength of the shock.

\medskip
The plane incident shock solution in the $({\bf x},t)$--coordinates
with states $(\rho, \nabla_{\bf x}\Psi)=(\rho_0, 0,0)$ and $(\rho_1,
u_1,0)$
corresponds to a continuous weak solution $\varphi$ of (\ref{1.1.5})
in the self-similar coordinates $(\xi,\eta)$ with the following
form:
\begin{eqnarray}
&&\varphi_0(\xi,\eta)=-\frac{1}{2}(\xi^2+\eta^2) \qquad
 \hbox{for } \,\, \xi>\xi_0,
 \label{flatOrthSelfSimShock1} \\
&&\varphi_1(\xi,\eta)=-\frac{1}{2}(\xi^2+\eta^2)+ u_1(\xi-\xi_0)
\qquad
 \hbox{for } \,\, \xi<\xi_0,
 \label{flatOrthSelfSimShock2}
\end{eqnarray}
respectively, where
\begin{equation}\label{shocklocation}
\xi_0=\rho_1
\sqrt{\frac{2(\rho_1^{\gamma-1}-\rho_0^{\gamma-1})}{\rho_1^2-\rho_0^2}}
=\frac{\rho_1u_1}{\rho_1-\rho_0}>0
\end{equation}
is the location of the incident shock, uniquely determined by
$(\rho_0,\rho_1,\gamma)$ through (\ref{FBConditionSelfSim-0}).
Since the problem is symmetric with respect to the axis $\eta=0$, it
suffices to consider the problem in the half-plane $\eta>0$ outside
the half-wedge
$$
\Lambda:=\{\xi\le 0,\eta>0\}\cup\{\eta>\xi \tan\theta_w,\, \xi>0\}.
$$
%
%Old version 2-21-07
%$$
%\Lambda:=\{\xi<0,\eta>0\}\cup\{\eta>\xi \tan\theta_w,\, \xi>0\}.
%$$
Then the initial-boundary value problem
\eqref{1.1.1}--\eqref{boundary-condition} in the $({\bf x},
t)$--coordinates can be formulated as the following
boundary value problem in the self-similar coordinates $(\xi,\eta)$.

\medskip
{\bf Problem 2 (Boundary Value Problem)} (see Fig. 2). {\it Seek a
solution $\varphi$ of equation \eqref{1.1.5} in the self-similar
domain $\Lambda$ with the slip boundary condition on the wedge
boundary $\partial\Lambda$:
\begin{equation}\label{boundary-condition-3}
D\varphi\cdot\nu|_{\partial\Lambda}=0
\end{equation}
and the asymptotic boundary condition at infinity:
\begin{equation}\label{boundary-condition-2}
\varphi\to\bar{\varphi}=
\begin{cases} \varphi_0 \qquad\mbox{for}\,\,\,
                         \xi>\xi_0, \eta>\xi \tan\theta_w,\\
              \varphi_1 \qquad \mbox{for}\,\,\,
                          \xi<\xi_0, \eta>0,
\end{cases}
\qquad \mbox{when $\xi^2+\eta^2\to \infty$},
\end{equation}}
where (\ref{boundary-condition-2}) holds in the sense that $
\displaystyle
\lim_{R\to\infty}\|\varphi-\overline{\varphi}\|_{C(\Lambda\setminus
B_R(0))}=0. $

\begin{figure}[h]
\centering
\includegraphics[height=2.6in,width=2.8in]{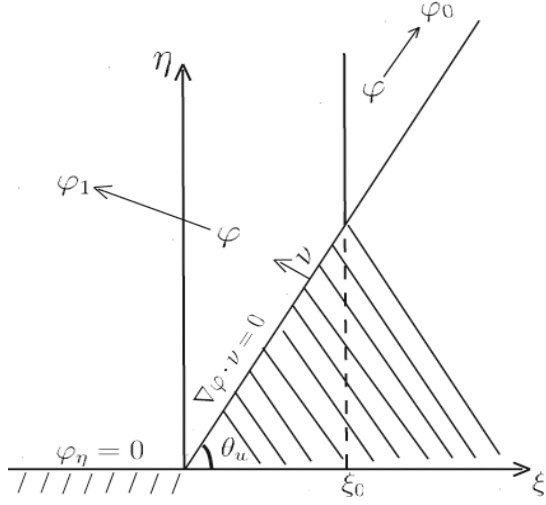}  %%% for latex (non-pdf)
\caption[]{Boundary value problem in the unbounded domain}
\label{fig:BVP}
\end{figure}

\begin{figure}[h]
\centering
\includegraphics[height=2.6in,width=3.4in]{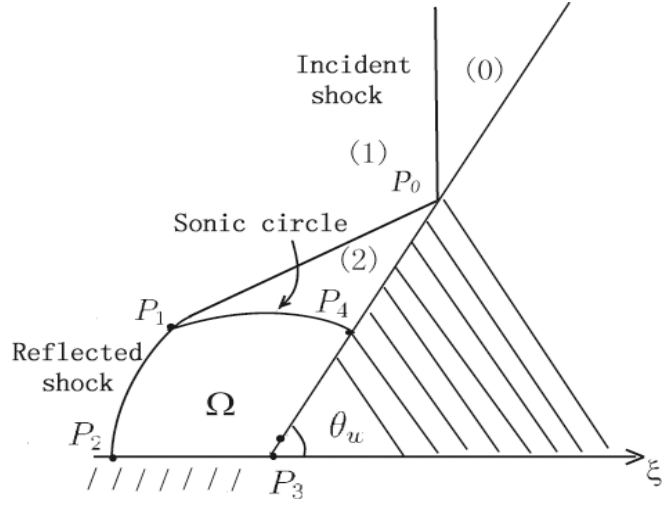}  %%% for latex (non-pdf)
\caption[]{Regular reflection} \label{fig:RegularReflection}
\end{figure}

\medskip
Since $\varphi_1$ does not satisfy the slip boundary condition
\eqref{boundary-condition-3}, the solution must differ from
$\varphi_1$ in $\{\xi<\xi_0\}\cap\Lambda$, thus a shock diffraction
by the wedge occurs. In this paper, we first follow the von Neumann
criterion to establish a local existence theory of regular shock
reflection near the reflection point and show that the structure of
solution is as in Fig. \ref{fig:RegularReflection}, when the wedge
angle is large and close to $\pi/2$, in which the vertical line is
the incident shock $S=\{\xi=\xi_0\}$ that hits the wedge at the
point $P_0=(\xi_0, \xi_0 \tan\theta_w)$,
and state (0) and state (1) ahead of and behind $S$ are given by
$\varphi_0$ and $\varphi_1$ defined in \eqref{flatOrthSelfSimShock1}
and \eqref{flatOrthSelfSimShock2}, respectively. The solutions
$\varphi$ and $\varphi_1$ differ only in the domain $P_0\PtUpL
\PtLwL \PtLwR$ because of shock diffraction by the wedge vertex,
where the curve $P_0\PtUpL \PtLwL$ is the reflected shock with the
straight segment $P_0\PtUpL$. State (2) behind $P_0\PtUpL$ can be
computed explicitly with the form:
\begin{equation}\label{state2a}
\varphi_2(\xi,\eta)=-\frac{1}{2}(\xi^2+\eta^2)+u_2(\xi-\xi_0)+
(\eta-\xi_0\tan\theta_w)u_2\tan\theta_w,
\end{equation}
which satisfies
$$
D\varphi\cdot\nu=0 \qquad \hbox{on }\, \partial\Lambda\cap
\{\xi>0\};
$$
the constant velocity $u_2$ and the angle $\theta_s$ between
$P_0\PtUpL$ and the $\xi$--axis are determined by
$(\theta_w,\rho_0,\rho_1,\gamma)$ from the two algebraic equations
expressing (\ref{FBConditionSelfSim-0}) and continuous  matching of
state (1) and state (2) across $P_0\PtUpL$, whose existence is
exactly guaranteed by the condition on
$(\theta_w,\rho_0,\rho_1,\gamma)$ under which regular shock
reflection is expected to occur.

\medskip
We develop a rigorous mathematical approach to extend the local
theory to a global theory for solutions of regular shock reflection,
which converge to the unique solution of the normal shock reflection
when $\theta_w$ tends to $\pi/2$. The solution $\varphi$ is
pseudo-subsonic within the sonic circle for state (2) with center
$(u_2, u_2\tan\theta_w)$ and radius $c_2>0$ (the sonic speed) and is
pseudo-supersonic outside this circle containing the arc
$\PtUpL\PtUpR$ in Fig. \ref{fig:RegularReflection}, so that
$\varphi_2$ is the unique solution in the domain $P_0\PtUpL \PtUpR$,
as argued in \cite{ChangChen,Serre}. In the domain $\Omega$, the
solution is expected to be pseudo-subsonic, smooth, and
$C^1$-smoothly matching with state (2) across $\PtUpL\PtUpR$ and to
satisfy $\varphi_\eta=0$ on $\PtLwL\PtLwR$; the transonic shock
curve $\PtUpL\PtLwL$ matches up to second-order with $P_0\PtUpL$ and
is orthogonal to the $\xi$-axis at the point $\PtLwL$ so that the
standard reflection about the $\xi$--axis yields a global solution
in the whole plane. Then the solution of Problem 2 can be shown to
be the solution of Problem 1.

\medskip
{\bf Main Theorem.} There exist
$\theta_c=\theta_c(\rho_0,\rho_1,\gamma) \in (0,\pi/2)$ and
$\alpha=\alpha(\rho_0,\rho_1,\gamma)\in (0, 1/2)$ such that, when
$\theta_w\in [\theta_c,\pi/2)$, there exists a global self-similar
solution
$$
\Phi({\bf x}, t) =t\,\varphi(\frac{\bf x}{t}) +\frac{|\bf x|^2}{2t}
\qquad\mbox{for}\,\, \frac{\bf x}{t}\in \Lambda,\, t>0
$$
with
$$
\rho({\bf x}, t)=(\rho_0^{\gamma-1}-\Phi_t
      -\frac{1}{2}|\nabla_{\bf x}\Phi|^2)^{\frac{1}{\gamma-1}}
$$
of Problem 1 (equivalently, Problem 2) for shock reflection by the
wedge, which satisfies that, for $(\xi,\eta)={\bf x}/t$,
$$
\varphi\in C^{\infty}(\Omega)\cap C^{1,\alpha}(\bar{\Omega}),
$$
\begin{equation}\label{phi-states-0-1-2}
\varphi=\left\{\begin{array}{ll}
\varphi_0 \qquad\mbox{for}\,\, \xi>\xi_0 \mbox{ and } \eta>\xi\tan\theta_w,\\
\varphi_1 \qquad\mbox{for}\,\, \xi<\xi_0
  \mbox{ and above the reflection shock} \,\, P_0\PtUpL\PtLwL,\\
\varphi_2 \qquad \mbox{in}\,\, P_0\PtUpL\PtUpR,
\end{array}\right.
\end{equation}
$\varphi$ is $C^{1,1}$ across the part $\PtUpL\PtUpR$ of the sonic
circle including the endpoints $\PtUpL$ and $\PtUpR$, and the
reflected shock $P_0\PtUpL\PtLwL$ is $C^2$ at $\PtUpL$ and
$C^\infty$ except $\PtUpL$. Moreover, the solution $\varphi$ is
stable with respect to the wedge angle in $W^{1,1}_{loc}(\overline
\Lambda)$ and converges in $W^{1,1}_{loc}(\overline\Lambda)$ to the
solution of the normal reflection described in \S\ref{section:4} as
$\theta_w\to\pi/2$.

\medskip
One of the main difficulties for the global existence is that the
ellipticity condition \eqref{1.1.8a} for (\ref{1.1.5}) is hard to
control, in comparison to our earlier work on steady flow
\cite{ChenFeldman1,ChenFeldman2}. The second difficulty is that the
ellipticity degenerates at the sonic circle $\PtUpL\PtUpR$ (the
boundary of the pseudo-subsonic flow). The third difficulty is that,
on $\PtUpL\PtUpR$, we need to match the solution in $\Omega$ with
$\varphi_2$ at least in $C^1$, that is, the two conditions on the
fixed boundary $\PtUpL\PtUpR$: the Dirichlet and conormal
conditions, which are generically overdetermined for an elliptic
equation since the conditions on the other parts of boundary have
been prescribed. Thus we have to prove that, if $\varphi$ satisfies
(\ref{1.1.5}) in  $\Omega$, the Dirichlet continuity condition on
the sonic circle, and the appropriate conditions on the other parts
of $\partial\Omega$ derived from Problem 2, then the normal
derivative $D\varphi\cdot\nu$ automatically matches with
$D\varphi_2\cdot\nu$ along $\PtUpL\PtUpR$. We show that, in fact,
this follows from the structure of elliptic degeneracy of
(\ref{1.1.5}) on $\PtUpL\PtUpR$ for the solution $\varphi$. Indeed,
equation (\ref{1.1.5}), written in terms of the function
$u=\varphi-\varphi_2$ in the $(x,y)$--coordinates defined near
$\PtUpL\PtUpR$ such that $\PtUpL\PtUpR$ becomes a segment on
$\{x=0\}$, has the form:
\begin{equation}\label{degenerate-equation}
\big(2x-(\gamma+1)u_x\big)u_{xx}+\frac{1}{c_2^2}u_{yy}-u_x=0
 \qquad\,\, \mbox{in } x>0 \mbox{  and near } x=0,
\end{equation}
plus the ``small'' terms that are controlled by $\pi/2-\theta_w$ in
appropriate norms. Equation \eqref{degenerate-equation} is elliptic
if $u_x<2x/(\gamma+1)$. Thus, we need to obtain the $C^{1,1}$
estimates near $\PtUpL\PtUpR$ to ensure $|u_x|<2x/(\gamma+1)$ which
in turn implies both the ellipticity of the equation in $\Omega$ and
the match of normal derivatives
$D\varphi\cdot\nu=D\varphi_2\cdot\nu$ along $\PtUpL\PtUpR$. Taking
into account the ``small'' terms to be added to equation
\eqref{degenerate-equation}, we need to make the stronger estimate
$|u_x|\le 4x/\big(3(\gamma+1)\big)$ and assume that $\pi/2-\theta_w$
is appropriately small to control these additional terms. Another
issue is the non-variational structure and nonlinearity of this
problem which makes it hard to apply directly the approaches of
Caffarelli \cite{Ca} and Alt-Caffarelli-Friedman \cite{AC,ACF}.
Moreover, the elliptic degeneracy and geometry of the problem makes
it difficult to apply the hodograph transform approach in
Kinderlehrer-Nirenberg \cite{KinderlehrerNirenberg} and Chen-Feldman
\cite{ChenFeldman4} to fix the free boundary.

\medskip
For these reasons, one of the new ingredients in
our approach is to further develop the iteration scheme in
\cite{ChenFeldman1,ChenFeldman2} to a partially modified equation.
We modify equation  (\ref{1.1.5}) in $\Omega$ by a proper cutoff
that depends on the distance to the sonic circle, so that the
original and modified equations coincide for $\varphi$ satisfying
$|u_x| \le 4x/\big(3(\gamma+1)\big)$, and the modified equation
${\mathcal N}\varphi=0$ is elliptic
%so that the expected global solution $\varphi$ of Problem 2 is also
%a solution of the modified equation ${\mathcal N}\varphi=0$, which
%is elliptic
in $\Omega$ with elliptic degeneracy on $\PtUpL\PtUpR$. Then we
solve a free boundary problem for this modified equation: The free
boundary is the curve $\PtUpL\PtLwL$, and the free boundary
conditions on $\PtUpL\PtLwL$ are $\varphi=\varphi_1$ and the
Rankine-Hugoniot condition (\ref{FBConditionSelfSim-0}).

On each step, an ``iteration free boundary'' curve $\PtUpL\PtLwL$ is
given,
and  a solution of the modified equation ${\mathcal N}\varphi=0$ is
constructed in $\Omega$ with the boundary condition
(\ref{FBConditionSelfSim-0}) on $\PtUpL\PtLwL$, the Dirichlet
condition $\varphi=\varphi_2$ on the degenerate circle
$\PtUpL\PtUpR$, and $D\varphi\cdot\nu=0$ on $\PtLwL\PtLwR$ and
$\PtLwR\PtUpR$.
Then we prove that $\varphi$ is in fact $C^{1,1}$ up to the boundary
$\PtUpL\PtUpR$, especially $|\grad(\varphi-\varphi_2)|\le Cx$, by
using the nonlinear structure of elliptic degeneracy near
$\PtUpL\PtUpR$ which is modeled by equation
(\ref{degenerate-equation}) and a scaling technique similar to
Daskalopoulos-Hamilton \cite{DG} and Lin-Wang \cite{LW}.
Furthermore, we modify the ``iteration free boundary'' curve
$\PtUpL\PtLwL$ by using the Dirichlet condition $\varphi=\varphi_1$
on $\PtUpL\PtLwL$.
 A fixed point $\varphi$ of this iteration procedure is a solution
of the free boundary problem for the modified equation. Moreover, we
prove the precise gradient estimate:
$|u_x|<4x/\big(3(\gamma+1)\big)$, which implies that  $\varphi$
satisfies the original equation (\ref{1.1.5}).

\smallskip Some efforts have been made mathematically for the
reflection problem via simplified models. One of these models, the
unsteady transonic small-disturbance (UTSD) equation, was derived
and used in Keller-Blank \cite{KB}, Hunter-Keller \cite{HK}, Hunter
\cite{hunter1}, Morawetz \cite{Morawetz2}, and the references cited
therein for asymptotic analysis of shock reflection. Also see Zheng
\cite{Zheng1} for the pressure gradient equation and
Canic-Keyfitz-Kim \cite{CKK1} for the UTSD equation and the
nonlinear wave system.
%when the wedge angle is close to $\pi/2$.
On the other hand, in order to deal with the reflection problem,
some asymptotic methods have been also developed. Lighthill
\cite{Lighthill1,Lighthill2} studied shock reflection  under the
assumption that the wedge angle is either very small or close to
$\pi/2$. Keller-Blank \cite{KB}, Hunter-Keller \cite{HK}, and
Harabetian \cite{Harabetian} considered the problem under the
assumption that the shock is so weak that its motion can be
approximated by an acoustic wave. For a weak incident shock and a
wedge with small angle in the context of potential flow, by taking
the jump of the incident shock as a small parameter, the nature of
the shock reflection pattern was explored in Morawetz
\cite{Morawetz2} by a number of different scalings, a study of mixed
equations, and matching the asymptotics for the different scalings.
Also see Chen \cite{Sxchen} for a linear approximation of shock
reflection when the wedge angle is close to $\pi/2$
 and Serre
\cite{Serre} for an apriori analysis of solutions of shock
reflection and related discussions in the context of the Euler
equations for isentropic and adiabatic fluids.

\smallskip
The organization of this paper is the following. In \S 2, we present
the potential flow equation in self-similar coordinates and exhibit
some basic properties of solutions to the potential flow equation.
In \S 3, we discuss the normal reflection solution and then follow
the von Neumann criterion to derive the necessary condition for the
existence of regular reflection and show that the shock reflection
can be regular locally when the wedge angle is large. In \S 4, the
shock reflection problem is reformulated and reduced to a free
boundary problem for a second-order nonlinear equation of mixed type
in a convenient form. In \S 5, we develop an iteration scheme, along
with an elliptic cutoff technique, to solve the free boundary
problem and set up the ten detailed steps of the iteration
procedure.

Finally, we complete the remaining steps in our iteration procedure
in \S6--\S9: Step 2 for the existence of solutions of the boundary
value problem to the degenerate elliptic equation via the vanishing
viscosity approximation in \S6; Steps 3--8 for the existence of the
iteration map and its fixed point in \S7; and Step 9 for the removal
of the ellipticity cutoff in the iteration scheme by using
appropriate comparison functions and deriving careful global
estimates for some directional derivatives of the solution in \S8.
We complete the proof of Main Theorem in \S9. Careful estimates of
the solutions to both the ``almost tangential derivative" and
oblique derivative boundary value problems for elliptic equations
are made in Appendix, which are applied in \S6--\S7.

\section{Self-Similar Solutions of the Potential Flow Equation}
\renewcommand{\theequation}{\thesection.\arabic{equation}}

In this section we present the potential flow equation in
self-similar coordinates and exhibit some basic properties of
solutions of the potential flow equation (also see Morawetz
\cite{Morawetz2}).

\Subsec{\bf\large The potential flow equation for self-similar
solutions}
%\renewcommand{\theequation}{\thesubsection.\arabic{equation}}
%\medskip
Equation \eqref{1.1.5}
%or \eqref{potent-flow-nondiv-phi}
is a mixed equation of elliptic-hyperbolic type. It is elliptic if
and only if
\eqref{1.1.8a} holds. The hyperbolic-elliptic boundary is the
pseudo-sonic curve: $|D\varphi|=c_*(\varphi,\rho_0,\gamma)$.

\medskip
We first define the notion of weak solutions of
\eqref{1.1.5}--\eqref{1.1.6}. Essentially, we require the equation
to be satisfied in the distributional sense.

\begin{definition}[Weak Solutions] \label{def2.1}
A function $\varphi\in W^{1,1}_{loc}(\Lambda)$ is called a weak
solution of \eqref{1.1.5}--\eqref{1.1.6} in a self-similar domain
$\Lambda$ if
\begin{enumerate}\renewcommand{\theenumi}{\roman{enumi}}
\item
$\rho_0^{\gamma-1}-\varphi-\frac{1}{2}|D\varphi|^2\ge 0$ a.e. in
$\Lambda$;
\item
$(\rho(|D\varphi|^2, \varphi),
\rho(|D\varphi|^2, \varphi)|\grad\varphi|)\in
(L^1_{loc}(\Lambda))^2$;
\item
For every $\zeta\in C^\infty_c(\Lambda)$,
$$
\int_\Lambda\Big(\rho(|D\varphi|^2,
\varphi)\grad\varphi\cdot\grad\zeta -2\rho(|D\varphi|^2,
\varphi)\zeta\Big) \,d\mxx d\mxy=0.
$$
\end{enumerate}
\end{definition}

\smallskip
It is straightforward to verify the equivalence between
time-dependent self-similar solutions and weak solutions of
\eqref{1.1.5} defined in Definition \ref{def2.1} in the weak sense.
It can also be verified that, if $\varphi\in C^{1,1}(\Lambda)$ (and
thus $\varphi$ is twice differentiable a.e. in $\Lambda$), then
$\varphi$ is a weak solution of (\ref{1.1.5}) in $\Lambda$ if and
only if $\varphi$ satisfies  equation \eqref{1.1.5} a.e. in
$\Lambda$. Finally, it is easy to see that, if $\Lambda^+$ and
$\Lambda^-=\Lambda\setminus\overline{\Lambda^+}$ are two nonempty
open subsets of $\Lambda\subset\bR^2$ and $S=\partial \Lambda^+\cap
\Lambda$ is a $C^1$ curve where $D\varphi$ has a jump, then
$\varphi\in W^{1,1}_{loc}(D)\cap C^1(\Lambda^\pm\cup S)\cap
C^{1,1}(\Lambda^\pm)$ is a weak solution of (\ref{1.1.5}) in
$\Lambda$ if and only if $\varphi$ is in
$W^{1,\infty}_{loc}(\Lambda)$ and satisfies equation \eqref{1.1.5}
a.e. in $\Lambda^\pm$ and the Rankine-Hugoniot condition
(\ref{FBConditionSelfSim-0}) on $S$.

\smallskip
Note that, for  $\varphi\in C^1(\Lambda^\pm\cup S)$, the condition
$\varphi\in W^{1,\infty}_{loc}(\Lambda)$ implies
\begin{equation}
[\varphi]_{S}=0. \label{1.1.14}
\end{equation}
Furthermore, the Rankine-Hugoniot conditions imply
\begin{equation}
[\varphi_\xi][\rho\varphi_\xi]-[\varphi_\eta][\rho\varphi_\eta]=0
\qquad \text{on } S \label{1.1.16}
\end{equation}
which is a useful identity.

\medskip
A discontinuity of $\grad \varphi$ satisfying the Rankine-Hugoniot
conditions \eqref{1.1.14} and \eqref{FBConditionSelfSim-0} is called
a shock if it satisfies the physical entropy condition: {\it The
density function $\rho$ increases across a shock in the pseudo-flow
direction}. The entropy condition indicates that {\it the normal
derivative function $\varphi_\nu$ on a shock always decreases across
the shock in the pseudo-flow direction}.

\medskip\noindent
\Subsec{\bf \large The states with constant density}
 When the density $\rho$ is
constant, \eqref{1.1.5}--\eqref{1.1.6} imply that $\varphi$
satisfies
\begin{equation*}
\Delta \varphi +2=0,\qquad \frac{1}{2}|D\varphi|^2 +\varphi=const.
\end{equation*}
This implies $(\Delta\varphi)_\xi=0, (\Delta\varphi)_\eta=0,$ and
$(\varphi_{\xi\xi}+1)^2+\varphi_{\xi\eta}^2=0$.
%\end{eqnarray*}
Thus, we have
$$
\varphi_{\xi\xi}=-1, \quad \varphi_{\xi\eta}=0, \quad
\varphi_{\eta\eta}=-1,
$$
which yields
\begin{equation}
\varphi(\xi,\eta) =-\frac{1}{2}(\xi^2+\eta^2)
  +a\xi +b\eta +c,
\label{1.1.11}
\end{equation}
where $a, b$, and $c$ are constants.

\smallskip
\Subsec{\bf \large Location of the incident shock}
%
%\medskip
Consider state $(0)$: $(\rho_0,u_0,v_0)=(\rho_0,0,0)$ with
$\rho_0>0$ and state $(1)$: $(\rho_1,u_1,v_1)$ $=(\rho_1,u_1,0)$
with $\rho_1>\rho_0>0$ and $u_1>0$. The plane incident shock
solution with state (0) and state (1) corresponds to a continuous
weak solution $\varphi$ of (\ref{1.1.5}) in the self-similar
coordinates $(\xi,\eta)$ with form \eqref{flatOrthSelfSimShock1} and
\eqref{flatOrthSelfSimShock2} for state (0) and state (1)
respectively, where $\xi=\xi_0>0$ is the location of the incident
shock.

The unit normal to the shock line is $\nu=(1,0)$. Using
\eqref{1.1.16}, we have
$$
u_1=\frac{\rho_1-\rho_0}{\rho_1}\xi_0>0.
$$
Then \eqref{1.1.6} implies
$$
\rho_1^{\gamma-1}-\rho_0^{\gamma-1}
=-\frac{1}{2}|D\varphi_1|^2-\varphi_1
=\frac{1}{2}\frac{\rho_1^2-\rho_0^2}{\rho_1^2}\xi_0^2.
$$
Therefore, we have
\begin{equation}
u_1=(\rho_1-\rho_0)
\sqrt{\frac{2(\rho_1^{\gamma-1}-\rho_0^{\gamma-1})}{\rho_1^2-\rho_0^2}},
          \label{1.2.4}
\end{equation}
and the location of the incident shock in the self-similar
coordinates is $\xi=\xi_0>u_1$ determined by \eqref{shocklocation}.

\section{The von Neumann Criterion and Local Theory for Shock Reflection}
\label{section:3}

In this section, we first discuss the normal reflection solution.
Then we follow the von Neumann criterion to derive the necessary
condition for the existence of regular reflection and show that the
shock reflection can be regular locally when the wedge angle is
large, that is, when $\theta_w$ is close to $\pi/2$ and,
equivalently, the angle between the incident shock and the wedge
\begin{equation}\label{angleCloseToPiOver2}
\sigma:=\pi/2-\theta_w
\end{equation}
tends to zero.

\medskip
\Subsec{\bf\large Normal shock reflection}\label{section:4} In this
case, the wedge angle is $\pi/2$, i.e., $\sigma=0$, and the incident
shock normally reflects (see Fig. 4). The reflected shock is also a
plane at $\xi=\bar{\xi}<0$, which will be defined below. Then
$\bar{u}_2=\bar{v}_2=0$, state (1) has form
\eqref{flatOrthSelfSimShock2}, and state (2) has the form:
\begin{eqnarray}
\varphi_2(\xi,\eta)=-\frac{1}{2}(\xi^2+\eta^2)+u_1(\bar\xi-\xi_0)\qquad\mbox{for}\,\,
\xi\in (\bar{\xi},0), \label{phi-2-a}
\end{eqnarray}
where
$\xi_0=\rho_1 u_1/(\rho_1-\rho_0)>0$ may be regarded to be the
position of the incident shock.

\begin{figure}[h]
\centering
\includegraphics[height=1.8in,width=3.0in]{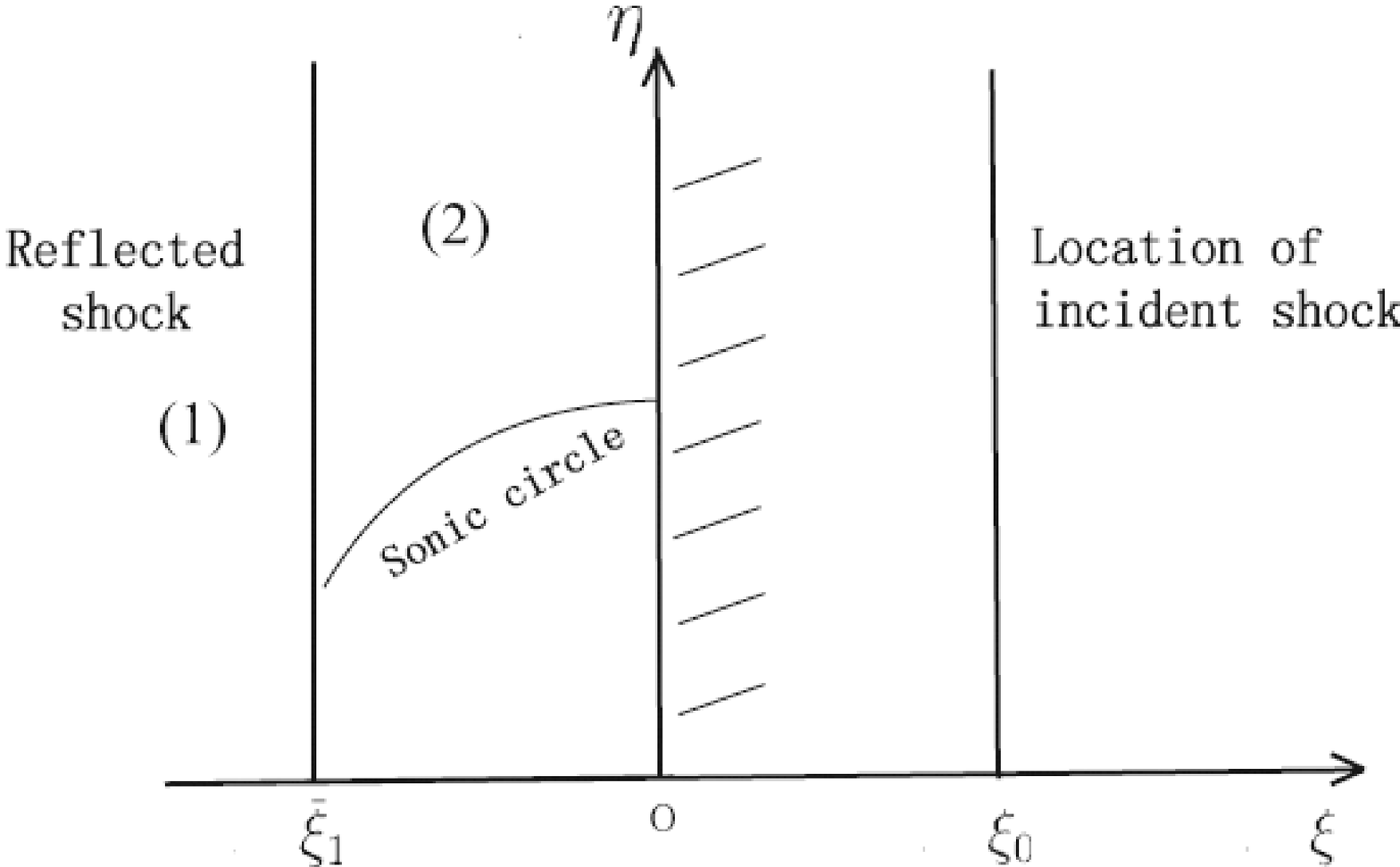}  %%% for latex (non-pdf)
\caption[]{Normal reflection} \label{fig:NF}
\end{figure}

At the reflected shock $\xi=\bar{\xi}<0$, the Rankine-Hugoniot
condition \eqref{1.1.16} implies
\begin{equation} \label{3.14}
\bar{\xi}=-\frac{\rho_1 u_1}{\bar{\rho}_2-\rho_1}<0.
\end{equation}

We use the Bernoulli law (\ref{1.1.4}):
$$
\rho_0^{\gamma-1}=\rho_1^{\gamma-1}+\frac{1}{2}u_1^2-u_1 \xi_0
=\bar{\rho}_2^{\gamma-1}+u_1(\bar{\xi}-\xi_0)
$$
to obtain
\begin{equation} \label{3.15}
\bar{\rho}_2^{\gamma-1}=\rho_1^{\gamma-1}+\frac{1}{2}u_1^2
 +\frac{\rho_1u_1^2}{\bar{\rho}_2-\rho_1}.
\end{equation}

It can be shown that there is a unique solution $\bar{\rho}_2$ of
\eqref{3.15} such that
$$
\bar{\rho}_2>\rho_1.
$$
Indeed, for fixed $\gamma>1$ and $\rho_1, u_1>0$ and for
$F(\bar{\rho}_2)$ that is the right-hand side of (\ref{3.15}), we
have
\begin{eqnarray*}
&&\lim_{s\to\infty}F(s)=\rho_1^{\gamma-1}+\frac{1}{2}u_1^2>\rho_1^{\gamma-1},
\quad \lim_{s\to\rho_1+}F(s)=\infty,\\
&&F'(s)=-\frac{\rho_1u_1^2}{(s-\rho_1)^2}<0 \qquad\, \mbox{ for}
\,\, s>\rho_1.
\end{eqnarray*}
Thus there exists a unique
$\bar{\rho}_2\in(\rho_1, \infty)$ satisfying
$\bar{\rho}_2^{\gamma-1}=F(\bar{\rho}_2)$, i.e., \eqref{3.15}. Then
the position of the reflected shock $\xi=\bar{\xi}<0$ is uniquely
determined by \eqref{3.14}.

Moreover, for the sonic speed
$\bar{c}_2=\sqrt{(\gamma-1)\bar{\rho}_2^{\gamma-1}}$ of state (2),
we have
\begin{equation}\label{sonic-intersect-shock-normal}
|\bar{\xi}|< \bar{c}_2.
\end{equation}
This can be seen as follows. First note that
\begin{equation}\label{meanValTh}
\bar{\rho}_2^{\gamma-1}-\rho_1^{\gamma-1}=\beta(\bar{\rho}_2-\rho_1),
\end{equation}
where $\beta=(\gamma-1)\rho_*^{\gamma-2}>0$ for some $\rho_*\in
(\rho_1,\bar{\rho}_2)$. We consider two cases, respectively.

\medskip
\noindent {\bf Case 1.} $\gamma\ge 2$. Then
\begin{equation}\label{meanValTh-A-1}
0<(\gamma-1)\rho_1^{\gamma-2}\le \beta\le
(\gamma-1)\bar{\rho}_2^{\gamma-2}.
\end{equation}
Since $\beta>0$ and $\bar{\rho}_2>\rho_1$, we use \eqref{3.15} and
\eqref{meanValTh} to find
$$
\bar{\rho}_2=\rho_1+\frac{u_1}{4\beta}\big(u_1+\sqrt{u_1^2+ 16\beta
\rho_1}\big),
$$
and hence
\begin{equation}\label{meanValTh-A-xi}
\bar{\xi}=-\frac{4\beta \rho_1}{u_1+\sqrt{u_1^2+16\beta \rho_1}}.
\end{equation}
Then using \eqref{meanValTh-A-1}--(\ref{meanValTh-A-xi}),
$\bar{\rho}_2>\rho_1>0$, and $u_1>0$ yields
$$
|\bar{\xi}|=\frac{4\beta \rho_1}{u_1+\sqrt{u_1^2+16\beta \rho_1}}
<\sqrt{\beta \rho_1}\le
\sqrt{(\gamma-1)\bar{\rho}_2^{\gamma-2}\bar{\rho}_2}=\bar{c}_2.
$$

\medskip
\noindent {\bf Case 2.} $1<\gamma<2$. Then, since
$\bar{\rho}_2>\rho_1>0$,
\begin{equation}\label{meanValTh-A-2}
0<(\gamma-1)\bar{\rho}_2^{\gamma-2}\le \beta\le
(\gamma-1)\rho_1^{\gamma-2}.
\end{equation}
Since $\beta>0$, (\ref{meanValTh-A-xi}) holds by the calculation as
in Case 1. Now we use (\ref{meanValTh-A-xi})--\eqref{meanValTh-A-2},
$\bar{\rho}_2>\rho_1>0$, $u_1>0$, and $1<\gamma<2$ to find again
$$
|\bar{\xi}|<\sqrt{\beta \rho_1}
%\le\sqrt{(\gamma-1)\rho_1^{\gamma-2}\rho_1}
\le \sqrt{(\gamma-1)\rho_1^{\gamma-1}}\le
\sqrt{(\gamma-1)\bar{\rho}_2^{\gamma-1}} =\bar{c}_2.
$$
This shows that (\ref{sonic-intersect-shock-normal}) holds in
general.

\Subsec{\bf\large The von Neumann criterion and local theory for
regular reflection}\label{section:3.3} In this subsection, we first
follow the von Neumann criterion to derive the necessary condition
for the existence of regular reflection and show that, when the
wedge angle is large, there exists a unique state (2) with two-shock
structure at the reflected point, which is close to the solution
$(\bar{\rho}_2,\bar{u}_2,\bar{v}_2)=(\bar{\rho}_2,0,0)$ of normal
reflection for which $\theta_w=\pi/2$ in \S 3.1.

For a possible two-shock configuration satisfying the corresponding
boundary condition on the wedge $\eta=\xi\tan\theta_w$, the three
state functions $\varphi_j, j=0,1,2$, must be of form
\eqref{flatOrthSelfSimShock1}, \eqref{flatOrthSelfSimShock2}, and
\eqref{state2a} (cf. \eqref{1.1.11}).

Set the reflected point $P_0=(\xi_0, \xi_0 \tan\theta_w)$ and assume
that the line that coincides with the reflected shock in state (2)
will intersect with the axis $\eta=0$ at the point $(\tilde{\xi},
0)$ with the angle $\theta_s$ between the line and $\eta=0$.

Note that $\varphi_1(\xi,\eta)$ is defined by
(\ref{flatOrthSelfSimShock2}). The continuity of $\varphi$ at
$(\tilde{\xi}, 0)$ yields
\begin{equation}\label{state2}
\varphi_2(\xi,\eta) =-\frac{1}{2}(\xi^2+\eta^2)+u_2\xi+v_2\eta
  +\big(u_1(\tilde\xi-\xi_0)-u_2\tilde{\xi}\big).
\end{equation}
Furthermore, $\varphi_2$ must satisfy the slip boundary condition at
$P_0$:
\begin{equation}\label{3.11a}
v_2=u_2 \tan\theta_w.
\end{equation}

Also we have
\begin{equation}\label{state2.5}
\tilde{\xi}=\xi_0-\xi_0\frac{\tan \theta_w}{\tan \theta_s}.
\end{equation}
The Bernoulli law (\ref{1.1.4}) becomes
\begin{equation}\label{state2.2}
\rho_0^{\gamma-1}=\rho_2^{\gamma-1}+\frac{1}{2}(u_2^2+v_2^2)
+(u_1-u_2)\tilde{\xi}-u_1\xi_0.
\end{equation}
Moreover, the continuity of $\varphi$ on the shock implies that
$D(\varphi_2-\varphi_1)$ is orthogonal to the tangent direction of
the reflected shock:
\begin{equation}\label{thetaS-uv}
(u_2-u_1, v_2)\cdot (\cos\theta_s, \sin\theta_s)=0,
\end{equation}
that is,
\begin{equation}\label{state2.3}
u_2=u_1\frac{\cos\theta_w\cos\theta_s}{\cos(\theta_w-\theta_s)}.
\end{equation}
The Rankine-Hugoniot condition (\ref{FBConditionSelfSim-0}) along
the reflected shock is
$$
[\rho\,D\varphi]\cdot (\sin\theta_s, -\cos\theta_s)=0,
$$
that is,
\begin{equation}\label{state2.4}
\rho_1(u_1-\tilde\xi)\sin\theta_s
=\rho_2\big(u_2\frac{\sin(\theta_s-\theta_w)}{\cos\theta_w}-\tilde\xi\sin\theta_s\big).
\end{equation}

Combining \eqref{state2.5}--\eqref{state2.4}, we obtain the
following system for $(\rho_2, \theta_s, \tilde{\xi})$:
\begin{eqnarray}
&&(\tilde{\xi}-\xi_0)\cos\theta_w
 +\xi_0\sin\theta_w\cot\theta_s=0, \label{suf:1}\\
&&\rho_2^{\gamma-1}+\frac{u_1^2\cos^2\theta_s}{2\cos^2(\theta_w-\theta_s)}
 +\frac{u_1\sin\theta_w\sin\theta_s}{\cos(\theta_w-\theta_s)}\tilde{\xi}
 -u_1\xi_0-\rho_0^{\gamma-1}=0, \label{suf:2}\\
&&\big(u_1\cos\theta_s\tan (\theta_s-\theta_w)
  -\tilde{\xi}\sin\theta_s\big)\rho_2
  -\rho_1(u_1-\tilde{\xi})\sin\theta_s =0. \label{suf:3}
\end{eqnarray}
The condition for solvability of this system is the necessary
condition for the existence of regular shock reflection.

Now we compute the Jacobian $J$ in terms of $(\rho_2,\theta_s,
\tilde{\xi})$ at the normal reflection solution state
$(\bar{\rho}_2,\frac{\pi}{2}, \bar{\xi})$ in \S 3.1 for state $(2)$
when $\theta_w=\pi/2$ to obtain
$$
J=-\xi_0\big((\gamma-1)\bar{\rho}_2^{\gamma-2}(\bar{\rho}_2-\rho_1)
  -u_1 \bar{\xi}\big)< 0,
$$
since $\bar{\rho}_2>\rho_1$ and $\bar{\xi}<0$.
Then, by the Implicit Function Theorem, when $\theta_w$ is near
$\pi/2$, there exists a unique solution $(\rho_2,\theta_s,
\tilde{\xi})$ close to $(\bar{\rho}_2,\frac{\pi}{2}, \bar{\xi})$ of
system \eqref{suf:1}--\eqref{suf:3}. Moreover, $(\rho_2,\theta_s,
\tilde{\xi})$ are smooth functions of $\sigma=\pi/2-\theta_w\in
(0,\sigma_1)$ for $\epsP_1>0$ depending only on $\rho_0, \rho_1$,
and $\gamma$. In particular,
\begin{eqnarray}
\label{theta_s-close-w} |\rho_2-\bar{\rho}_2|+|\pi/2-\theta_s|+
|\tilde{\xi}-\bar{\xi}|+|c_2-\bar{c}_2| \le C\sigma,
\end{eqnarray}
where $\displaystyle c_2=\sqrt{(\gamma-1)\rho_2^{\gamma-1}}$ is the
sonic speed of state (2).

Reducing $\epsP_1>0$ if necessary, we find that, for any $\sigma\in
(0,\sigma_1)$,
\begin{equation} \label{xi-1-negative}
\tilde{\xi}<0
\end{equation}
from (\ref{3.14}) and (\ref{theta_s-close-w}). Since
$\theta_w\in(\pi/2-\epsP_1, \pi/2)$, then $\theta_s\in(\pi/4,
3\pi/4)$ if $\epsP_1$ is small, which implies $\sin\theta_s>0$. We
conclude from (\ref{suf:1}), (\ref{xi-1-negative}), and $\xi_0>0$
that $ \tan\theta_w>\tan\theta_s>0$. Thus,
\begin{equation} \label{theta_s<w}
\pi/4<\theta_s<\theta_w<\pi/2.
\end{equation}

Now, given $\theta_w$, we define  $\varphi_2$ as follows: We have
shown that there exists a unique solution $(\rho_2,\theta_s,
\tilde{\xi})$ close to $(\bar{\rho}_2,\frac{\pi}{2}, \bar{\xi})$ of
system \eqref{suf:1}--\eqref{suf:3}. Define $u_2$ by
(\ref{state2.3}), $v_2$ by (\ref{3.11a}), and $\varphi_2$ by
(\ref{state2}). Then  the shock connecting state (1) with state (2)
is the straight line $S_{12}=\{(\mxx,\mxy)\;:\,
\varphi_1(\mxx,\mxy)=\varphi_2(\mxx,\mxy)\}$, which is
$\xi=\mxy\cot\theta_s +\tilde{\xi}$ by (\ref{flatOrthSelfSimShock2}),
(\ref{state2}), and (\ref{state2.3}). Now (\ref{suf:3}) implies that
the Rankine-Hugoniot condition (\ref{FBConditionSelfSim-0}) holds on
$S_{12}$. Moreover, (\ref{3.11a}) and (\ref{state2.3}) imply
(\ref{thetaS-uv}). Thus the solution $(\theta_s, \rho_2, u_2, v_2)$
satisfies (\ref{3.11a})--(\ref{suf:3}). Furthermore, (\ref{suf:1})
implies that the point $P_0$ lies on $S_{12}$, and (\ref{suf:2})
implies (\ref{state2.2}) that is the Bernoulli law:
\begin{equation}\label{bernLawState2}
\rho_2^{\gamma-1}+{1\over 2}|\grad\varphi_2|^2+\varphi_2=
\rho_0^{\gamma-1}.
\end{equation}
Thus we have established the local existence of the two-shock
configuration near the reflected point so that, behind the straight
reflected shock emanating from the reflection point, state (2) is
pseudo-supersonic up to the sonic circle of state (2). Furthermore,
this local structure is stable in the limit $\theta_w\to\pi/2$,
i.e., $\sigma\to 0$.

We also notice from (\ref{3.11a}) and (\ref{state2.3}) with the use
of (\ref{theta_s-close-w}) and (\ref{theta_s<w}) that
\begin{equation}\label{u2-v2-bound}
|u_2|+ |v_2|\le C \sigma.
\end{equation}
Furthermore, from (\ref{sonic-intersect-shock-normal}) and the
continuity of $\rho_2$ and $\tilde{\xi}$ with respect to $\theta_w$
on $(\pi/2-\epsP_1, \pi/2]$, it follows that, if $\sigma>0$ is
small,
\begin{equation}\label{sonic-intersect-shock}
|\tilde{\xi}|< c_2.
\end{equation}

\medskip
In \S\ref{reformulProbSection}--\S\ref{proofSection}, we prove that
this local theory for the existence of two shock configuration can
be extended to a global theory for regular shock reflection.

\section{Reformulation of the Shock Reflection Problem}
\label{reformulProbSection}

We first assume that $\varphi$ is a solution of the shock reflection
problem in the elliptic domain $\ElDom$ in Fig.
\ref{fig:RegularReflection} and that $\varphi-\varphi_2$ is small in
$C^1(\overline\ElDom)$. Under such assumptions, we rewrite the
equation and boundary conditions for solutions of the shock
reflection problem in the elliptic region.

\Subsec{\bf\large Shifting coordinates} \label{shiftCoordSection} It
is more convenient to change the coordinates in the self-similar
plane by shifting the origin to the center of sonic circle of state
(2). Thus we define
$$
(\mxx,\mxy)_{new}:=(\mxx, \mxy)-(u_2,v_2).
$$
For simplicity of notations, throughout this paper below, we will
always work in the new coordinates without changing the notation
$(\mxx, \mxy)$, and we will not emphasize this again later.

In the new shifted coordinates, the domain $\ElDom$ is expressed as
\begin{equation}\label{ellipticDomain}
\ElDom=B_{c_2}(0)\cap\{\mxy>-v_2\}\cap \{f(\mxy)<\mxx<
\mxy\cot\theta_w\},
\end{equation}
where $f$ is the position function of the free boundary, i.e., the
curved part of the reflected shock $\shock :=\{\mxx=f(\mxy)\}$. The
function $f$ in \eqref{ellipticDomain} will be determined below so
that
\begin{equation}\label{FBfunct-estimate}
\|f-l\|\le C \sigma
\end{equation}
in an appropriate norm, specified later. Here $\xi=l(\mxy)$ is the
location of the reflected shock of state (2) which is a straight
line, that is,
\begin{equation}\label{reflected-shock-s2}
l(\mxy)=\mxy\cot\theta_s +\hat\mxx
\end{equation}
and
\begin{equation}\label{x1-in-shifed}
\hat\mxx=\tilde\mxx- u_2 +v_2\cot\theta_s<0,
\end{equation}
if $\sigma=\pi/2-\theta_w>0$ is sufficiently small, since $u_2$ and
$v_2$ are small and $\tilde\mxx<0$ by (\ref{3.14}) in this case. Also
note that, since $u_2 =v_2\cot\theta_w>0$, it follows from
(\ref{theta_s<w}) that
\begin{equation}\label{x1-in-shifed-2}
\hat\mxx>\tilde\mxx.
\end{equation}

Another condition on $f$ comes from the fact that the curved part
and straight part of the reflected shock should match at least up to
first-order.
Denote by $P_1=(\mxx_1,\mxy_1)$ with $\mxy_1>0$ the intersection point
of the line $\mxx=l(\mxy)$ and the sonic circle $\mxx^2+\mxy^2=c_2^2$,
i.e., $(\mxx_1,\mxy_1)$ is the unique point for small $\sigma>0$
satisfying
\begin{equation}\label{coord-P4}
l(\mxy_1)^2+\mxy_1^2=c_2^2,\qquad \mxx_1=l(\mxy_1), \qquad \mxy_1>0.
\end{equation}
The existence and uniqueness of such point $(\mxx_1,\mxy_1)$ follows
from $-c_2<\tilde{\xi}<0$, which holds from (\ref{theta_s<w}),
(\ref{sonic-intersect-shock}), (\ref{x1-in-shifed}), and the
smallness of $u_2$ and $v_2$. Then $f$ satisfies
\begin{equation}\label{curved-straight-shock-match}
f(\mxy_1)=l(\mxy_1), \qquad f'(\mxy_1)=l'(\mxy_1)=\cot\theta_s.
\end{equation}
Note also that, for small $\sigma>0$, we obtain from
(\ref{sonic-intersect-shock}),
(\ref{x1-in-shifed})--(\ref{x1-in-shifed-2}), and
$l'(\eta)=\cot\theta_s>0$ that
\begin{equation}\label{inSonicRegion-in-shifed}
-c_2<\tilde\mxx<\hat\mxx< \mxx_1<0, \qquad c_2-|\tilde{\xi}|\ge
\frac{\bar{c}_2-|\bar{\xi}|}{2}>0.
\end{equation}

Furthermore, equations (\ref{1.1.5})--(\ref{1.1.6}) and the
Rankine-Hugoniot conditions (\ref{FBConditionSelfSim-0}) and
(\ref{1.1.14}) on $\shock$ do not change under the shift of
coordinates. That is, we seek $\varphi$ satisfying
(\ref{1.1.5})--(\ref{1.1.6}) in $\ElDom$ so that the equation is
elliptic on $\varphi$ and satisfying the following boundary
conditions on $\shock$: The continuity of the pseudo-potential
function across the shock:
\begin{equation}
\varphi=\varphi_1\qquad\mbox{on }\;\shock
\label{cont-accross-shock-mod-phi}
\end{equation}
and the gradient jump condition:
\begin{equation}
\rho(|\grad\varphi|^2,\varphi)\grad\varphi\cdot\nu_s= \rho_1
\grad\varphi_1\cdot\nu_s \qquad\mbox{on }\;\shock,
\label{RH-mod-phi}
\end{equation}
where $\nu_s$ is the interior unit normal to $\Omega$ on $\shock$.

The boundary conditions on the other parts of $\partial\ElDom$ are
\begin{eqnarray}
&&\varphi=\varphi_2\qquad\mbox{on }\;\sonic=\partial \ElDom\cap
\partial B_{c_2}(0),
\label{condOnSonicLinePhi}
\\
&& \varphi_\nu=0\qquad\mbox{on }\;\wedgeB=\partial \ElDom\cap
\{\mxy=\mxx\tan\theta_w\}, \label{condOnWedgePhi}
\\
&&\varphi_\nu=0\qquad\mbox{on }\;\partial \ElDom\cap \{\mxy=-v_2\}.
\label{condOnSymmtryLinePhi}
\end{eqnarray}

\medskip
Rewriting the background solutions in the shifted coordinates, we
find
\begin{eqnarray}
&&\qquad
\varphi_0(\mxx,\mxy)=-\frac{1}{2}(\xi^2+\eta^2)-(u_2\mxx+v_2\mxy)
-{1\over 2}\UreflAbs^2,
\label{phi-0-shifted} \\
&&\qquad \varphi_1(\mxx,\mxy)=-\frac{1}{2}(\xi^2+\eta^2)
+(u_1-u_2)\mxx-v_2\mxy -{1\over 2}\UreflAbs^2+u_1(u_2-\mxx_0),
\label{phi-1-shifted} \\
&&\qquad \varphi_2(\mxx,\mxy)=-\frac{1}{2}(\xi^2+\eta^2) -{1\over
2}\UreflAbs^2 +(u_1-u_2)\hat{\xi}+u_1(u_2-\xi_0),
\label{phi-2-shifted}
\end{eqnarray}
where $q_2^2=u_2^2+v_2^2$.

Furthermore, substituting  $\tilde\mxx$ in (\ref{x1-in-shifed}) into
equation (\ref{suf:1}) and using (\ref{3.11a}) and
(\ref{thetaS-uv}), we find
\begin{equation}
\rho_2\hat\mxx=\rho_1\big(\hat\mxx-{(u_1-u_2)^2+v_2^2\over
u_1-u_2}\big), \label{RH-states-1-2'}
\end{equation}
which expresses the  Rankine-Hugoniot conditions on the reflected
shock of state (2) in terms of $\hat\mxx$. We use this equality
below.

\begin{figure}[h]
\centering
\includegraphics[height=2.9in,width=2.8in]{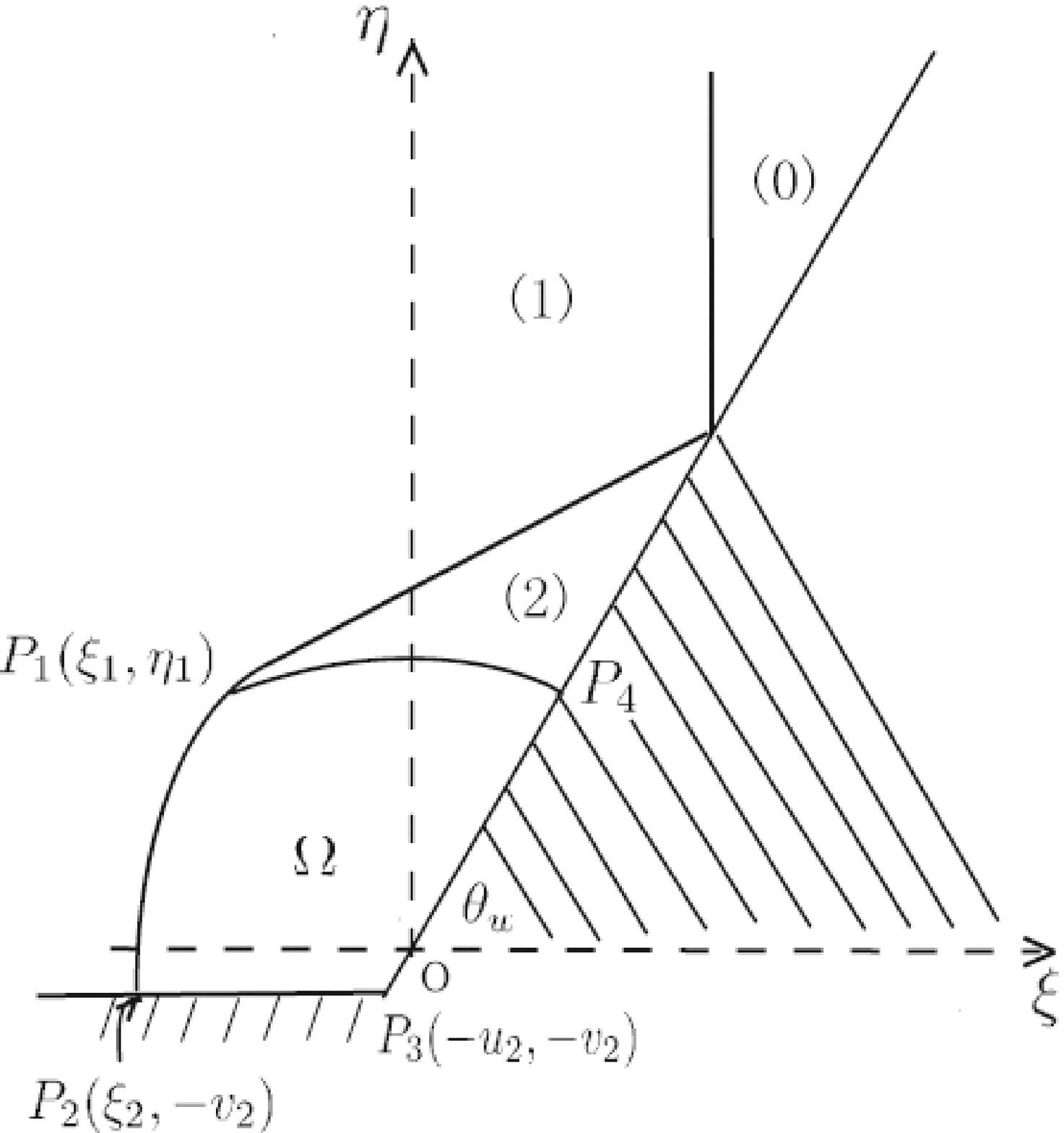}  %%% for latex (non-pdf)
\caption[]{Regular reflection in the new coordinates} \label{fig:RR}
\end{figure}

\Subsec{\bf\large   The equations and boundary conditions in terms
of $\psi=\varphi-\varphi_2$} \label{equationForPsiSection} It is
convenient to study the problem in terms of the difference between
our solution $\varphi$ and  the function $\varphi_2$ that is a
solution for state (2) given by \eqref{phi-2-shifted}. Thus we
introduce a function
\begin{equation}\label{psi-definition}
\psi=\varphi-\varphi_2\qquad\mbox{in }\;\ElDom.
\end{equation}
Then it follows from
(\ref{1.1.5})--(\ref{c-through-density-function}),
(\ref{bernLawState2}), and \eqref{phi-2-shifted} by explicit
calculation that $\psi$ satisfies the following equation in
$\ElDom$:
\begin{eqnarray}\label{potent-flow-nondiv-psi-1}
&&\quad
\big(c^2(D\psi,\psi,\xi,\eta)-(\psi_\mxx-\mxx)^2\big)\psi_{\mxx\mxx}
+\big(c^2(D\psi,\psi,\xi,\eta)-(\psi_\mxy-\mxy)^2\big)\psi_{\mxy\mxy}
 \\
&&\qquad\quad -2(\psi_\mxx-\mxx)(\psi_\mxy-\mxy)\psi_{\mxx\mxy}
=0,\nonumber
\end{eqnarray}
and the expressions of the density and sound speed in $\ElDom$ in
terms of $\psi$ are
\begin{eqnarray}
&&\rho(\grad\psi,\psi,\xi,\eta) =\Big(\rho_2^{\gamma-1}
+\mxx\psi_\mxx+\mxy\psi_\mxy-\frac{1}{2}|\grad\psi|^2-\psi
\Big)^\frac{1}{\gamma-1},
\label{density-psi} \\
&&c^2(\grad\psi, \psi,\xi,\eta) =c_2^2
+({\gamma-1})\Big(\mxx\psi_\mxx+\mxy\psi_\mxy-\frac{1}{2}|\grad\psi|^2
  -\psi \Big).
\label{speedOfsound-psi}
\end{eqnarray}
where $\rho_2$ is the density of state (2). In the polar coordinates
$(r, \theta)$ with $r=\sqrt{\xi^2+\eta^2}$,
 $\psi$ satisfies
\begin{equation}\label{potent-flow-nondiv-psi-Polar}
\big(c^2-(\psi_r-r)^2\big)\psi_{rr}
-\frac{2}{r^2}(\psi_r-r)\psi_\theta\psi_{r\theta}
+\frac{1}{r^2}(c^2-\frac{1}{r^2}\psi_\theta^2)\psi_{\theta\theta}
+\frac{c^2}{r}\psi_r +\frac{1}{r^3}(\psi_r-2r)\psi_\theta^2 =0
\end{equation}
with
\begin{equation}
\label{speedSound-expression-psi-Polar-1}
c^2=(\gamma-1)\Big(\rho_2^{\gamma-1}-\psi
+r\psi_r-\frac{1}{2}\big(\psi_r^2+\frac{1}{r^2}\psi_\theta^2\big)\Big).
\end{equation}
Also, from (\ref{condOnSonicLinePhi})--(\ref{condOnWedgePhi}) and
(\ref{phi-2-shifted})--(\ref{psi-definition}), we obtain
\begin{eqnarray}
&&\psi=0\qquad\mbox{on }\;\sonic=\partial \ElDom\cap
\partial B_{c_2}(0),
\label{condOnSonicLine-Psi}
\\
&&\psi_\nu=0\qquad\mbox{on }\;\wedgeB=\partial \ElDom\cap
\{\mxy=\mxx\tan\theta_w\}, \label{condOnWedge-Psi}
\\
&&\psi_\mxy=-v_2\qquad\mbox{on }\;\partial\ElDom\cap\{\mxy=-v_2\}.
\label{condOnSymmtryLine-Psi}
\end{eqnarray}

Using (\ref{phi-1-shifted})--(\ref{phi-2-shifted}), the
Rankine-Hugoniot conditions in terms of $\psi$ take the following
form: The continuity of the pseudo-potential function across
(\ref{cont-accross-shock-mod-phi}) is written as
\begin{equation}
\psi-{1\over 2}\UreflAbs^2+\hat{\mxx}(u_1-u_2)+u_1(u_2-\xi_0)
=\mxx(u_1-u_2)-\mxy v_2 -{1\over 2}\UreflAbs^2+u_1(u_2-\xi_0)
\quad\mbox{on }\;\shock, \label{cont-accross-shock-psi}
\end{equation}
that is,
\begin{equation}
\mxx=\frac{\psi(\mxx,\mxy)+v_2\mxy}{u_1-u_2}+\hat\mxx,
\label{cont-accross-shock-psi-resolved}
\end{equation}
where $\hat\mxx$ is defined by (\ref{x1-in-shifed}); and the gradient
jump condition (\ref{RH-mod-phi}) is
\begin{equation}
\rho(\grad\psi,\psi)\left(\grad\psi-
\left(\mxx,\mxy\right)\right)\cdot\nu_s =\rho_1 \left(u_1-u_2-\mxx,
-v_2-\mxy\right)\cdot\nu_s \qquad\mbox{on }\;\shock, \label{RH-psi}
\end{equation}
where $\rho(\grad\psi,\psi)$ is defined by (\ref{density-psi}) and
$\nu_s$ is the interior unit normal to  $\ElDom$ on $\shock$. If
$|(u_2, v_2, \grad\psi)|<u_1/50$, the unit normal $\nu_s$ can be
expressed as
\begin{equation}
\nu_s=\frac{\grad(\varphi_1-\varphi)}{|\grad(\varphi_1-\varphi)|} =
\frac{(u_1-u_2-\psi_\mxx,-v_2-\psi_\mxy)}
{\sqrt{(u_1-u_2-\psi_\mxx)^2+(v_2+\psi_\mxy)^2}},
\label{norm-to-Shock}
\end{equation}
where we have used (\ref{phi-1-shifted})--(\ref{phi-2-shifted}) and
(\ref{psi-definition}) to obtain the last expression.

Now we rewrite the jump condition (\ref{RH-psi}) in a more
convenient form for $\psi$ satisfying
(\ref{cont-accross-shock-mod-phi}) when $\sigma>0$ and
$\|\psi\|_{C^1(\bar{\ElDom})}$ are sufficiently small.

We first discuss the smallness assumptions for $\sigma>0$ and
$\|\psi\|_{C^1(\bar{\ElDom})}$. By (\ref{1.2.4}),
(\ref{theta_s-close-w}), and (\ref{u2-v2-bound}), it follows that,
if $\sigma$ is small depending only on the data, then
\begin{equation}\label{condRewritingRH-0}
\frac{5 \bar c_2}{6}\le c_2\le \frac{6\bar c_2}{5}, \quad
\frac{5\bar \rho_2}{6}\le \rho_2\le \frac{6\bar\rho_2}{5}, \quad
\sqrt{u_2^2+v_2^2}\le {u_1\over 50}.
\end{equation}
We also require that $\|\psi\|_{C^1(\bar{\ElDom})}$ is sufficiently
small so that, if (\ref{condRewritingRH-0}) holds, the expressions
(\ref{density-psi}) and (\ref{norm-to-Shock}) are well-defined in
$\ElDom$, and $\mxx$ defined by the right-hand side of
(\ref{cont-accross-shock-psi-resolved}) satisfies $|\mxx|\le 7\bar
c_2/5$ for $\mxy\in (-v_2, c_2)$, which is the range of $\mxy$ on
$\shock$. Since (\ref{condRewritingRH-0}) holds and $\ElDom\subset
B_{c_2}(0)$ by (\ref{ellipticDomain}), it suffices to assume
\begin{equation}\label{condRewritingRH}
\|\psi\|_{C^1(\bar{\ElDom})}\le \min\big(
{\bar\rho_2^{\gamma-1}\over 50(1+4\bar c_2)}, \min(1,\bar
c_2){u_1\over 50} \big)=:\delta^*.
\end{equation}

For the rest of this section, we assume that
(\ref{condRewritingRH-0}) and (\ref{condRewritingRH}) hold.

Under these conditions, we can substitute  the right-hand side of
 (\ref{norm-to-Shock}) for  $\nu_s$ into (\ref{RH-psi}).
Thus,  we rewrite (\ref{RH-psi}) as
\begin{equation}\label{RH-psi-int1}
F(\grad\psi, \psi, u_2, v_2, \mxx, \mxy)=0 \qquad\mbox{on}\;\shock,
\end{equation}
where, denoting $p=(p_1,p_2)\in\bR^2$ and $z\in\bR$, %$U_2=(u_2, v_2)$,
% and $\mX=(z, p_1,p_2, u_2, v_2, \mxx, \mxy)$,
\begin{equation}\label{RH-psi-func1}
F(p, z, u_2, v_2, \mxx, \mxy)=\big(\tilde\rho\left(p-
\left(\mxx,\mxy\right)\right)-\rho_1 \left(u_1-u_2-\mxx,
-v_2-\mxy\right)\big)\cdot\hat\nu
\end{equation}
with $\tilde\rho:=\tilde\rho(p, z, \mxx, \mxy)$ and
$\hat\nu:=\hat\nu(p, u_2, v_2)$ defined by
\begin{eqnarray}
&&\tilde\rho(p, z, \mxx, \mxy) =\left(\rho_2^{\gamma-1} +\mxx
p_1+\mxy p_2-\frac{|p|^2}{2}-z \right)^{\frac{1}{\gamma-1}},
\label{RH-psi-func2}
\\
&&\hat\nu(p,  u_2, v_2)=\frac{(u_1-u_2-p_1,-v_2-p_2)}
{\sqrt{(u_1-u_2-p_1)^2+(v_2+p_2)^2}}. \label{RH-psi-func3}
\end{eqnarray}
{}From the explicit definitions of $\tilde\rho$ and $\hat\nu$, it
follows from (\ref{condRewritingRH-0}) that
$$
\tilde\rho\in C^\infty(\overline{B_{\delta^*}(0)\times (-\delta^*,
\delta^*)\times B_{2\bar c_2}(0)}), \quad \hat\nu\in
C^\infty(\overline{B_{\delta^*}(0)\times B_{u_1/50}(0)}),
$$
where $B_{R}(0)$ denotes the ball in $\bR^2$ with center $0$ and
radius $R$ and, for $k\in {\bf N}$ (the set of nonnegative
integers), the $C^k$--norms of $\tilde\rho$ and $\hat\nu$ over the
regions specified above are bounded by the constants depending only
on $\gamma, u_1, \bar\rho_2, \bar c_2$, and $k$, that is, by
\S\ref{section:3}, the $C^k$--norms depend only on the data and $k$.
Thus,
\begin{equation}\label{rewriteRH-reg-1}
F\in C^\infty(\overline{ B_{\delta^*}(0)\times (-\delta^*, \delta^*)
\times B_{u_1/50}(0)\times B_{2\bar c_2}(0)}),
\end{equation}
with its $C^k$--norm depending only on the data and $k$.

Furthermore, since $\psi$ satisfies
(\ref{cont-accross-shock-mod-phi}) and hence
(\ref{cont-accross-shock-psi-resolved}), we can substitute  the
right-hand side
 of (\ref{cont-accross-shock-psi-resolved}) for $\mxx$ into
(\ref{RH-psi-int1}). Thus we rewrite
 (\ref{RH-psi}) as
\begin{equation}\label{RH-psi-int2}
\Psi(\grad\psi, \psi, u_2, v_2, \mxy)=0 \qquad\mbox{on}\;\shock,
\end{equation}
where
\begin{equation}\label{RH-psi-func4}
\Psi(p, z, u_2, v_2, \mxy)=F(p, z, u_2, v_2,
(z+v_2\mxy)/(u_1-u_2)+\hat\mxx, \mxy).
\end{equation}
If $\eta\in (-6\bar c_2/5, 6\bar c_2/5)$ and $|z|\le \delta^*$,
then, from (\ref{inSonicRegion-in-shifed}) and
(\ref{condRewritingRH-0})--(\ref{condRewritingRH}), it follows that
$\big|(z+v_2\mxy)/(u_1-u_2)+\hat\mxx\big|\le 7\bar c_2/5$. That is,
$((z+v_2\mxy)/(u_1-u_2)+\hat\mxx,\,\eta)\in B_{2\bar c_2}(0)$ if
$\eta\in (-6\bar c_2/5, 6\bar c_2/5)$ and $|z|\le \delta^*$. Thus,
from (\ref{rewriteRH-reg-1}) and (\ref{RH-psi-func4}), $\Psi\in
C^\infty(\overline{\mathcal A})$ with
$\|\Psi\|_{C^k(\overline{\mathcal A})}$ depending only on the data
and $k\in{\bf N}$, where ${\mathcal A}=B_{\delta^*}(0)\times
(-\delta^*, \delta^*)\times B_{u_1/50}(0)\times (-6\bar c_2/5, 6\bar
c_2/5)$.

Using the explicit expression of $\Psi$ given by
(\ref{RH-psi-func1})--(\ref{RH-psi-func3}) and (\ref{RH-psi-func4}),
we calculate
\begin{eqnarray*}
&&\Psi((0,0),0, u_2,v_2,\mxy)\\
&& =- {(u_1-u_2)\rho_2\hat\mxx\over\sqrt{(u_1-u_2)^2+v_2^2}}
-\rho_1\big(
\sqrt{(u_1-u_2)^2+v_2^2}-{(u_1-u_2)\hat\mxx\over\sqrt{(u_1-u_2)^2+v_2^2}}
\big).
\end{eqnarray*}
Now, using (\ref{RH-states-1-2'}), we have
$$
\Psi((0,0), 0,u_2,v_2,\mxy) =0 \qquad\text{for any } (u_2,v_2,\mxy)\in
B_{u_1/50}(0)\times (-6\bar c_2/5, 6\bar c_2/5).
$$
Then, denoting $p_0=z$ and $\mX=((p_1,p_2), p_0, u_2,v_2,\mxy)\in
{\mathcal A}$, we have
\begin{equation}
\label{Psi-Function-1} \Psi(\mX)=\sum_{i=0}^2p_i D_{p_i}\Psi((0,0),
0, u_2,v_2,\mxy) +\sum_{i,j=0}^2 p_ip_j g_{ij}(\mX),
\end{equation}
where $ g_{ij}(\mX)=\int_0^1(1-t)D^2_{p_i p_j}\Psi((tp_1,tp_2),
tp_0, u_2,v_2,\mxy)dt$ for $i,j=0,1,2$. Thus, $g_{ij}\in
C^\infty(\overline{\mathcal A})$ and
$\|g_{ij}\|_{C^k(\overline{\mathcal
A})}\le\|\Psi\|_{C^{k+2}(\overline{\mathcal A})}$ depending only on
the data and $k\in {\bf N}$.

Next, denoting
$\rho_2'\defd\hat\rho'(\rho_2^{\gamma-1})=\rho_2/c_2^2>0,$ we
compute from the explicit expression of $\Psi$ given by
(\ref{RH-psi-func1})--(\ref{RH-psi-func3}) and (\ref{RH-psi-func4}):
\begin{eqnarray*}
D_{(p,z)}\Psi((0,0), 0, 0,0,\mxy)= \big(\rho_2'(c_2^2-{\hat\mxx}^2),\;
\big(\frac{\rho_2-\rho_1}{u_1}-\rho_2'\hat\mxx \big)\mxy,\,\;
\rho_2'\hat\mxx-\frac{\rho_2-\rho_1}{u_1}\big).
\end{eqnarray*}
Note that, for $i=0,1,2$,
$$
\partial_{p_i} \Psi((0,0), 0, u_2, v_2,\mxy)=\partial_{p_i} \Psi((0,0),0,0,0,\mxy)
+h_i(u_2, v_2,\mxy)
$$
with $\|h_i\|_{C^k(\overline{B_{u_1/50}(0)\times (-6\bar c_2/5,
6\bar c_2/5)})} \le\|\Psi\|_{C^{k+2}(\overline{\mathcal A})}$ for
$k\in {\bf N}$, and $|h_i(u_2, v_2,\mxy)|\le C(|u_2|+|v_2|)$ with
$C=\|D^2\Psi\|_{C(\overline{\mathcal A})}$. Then we obtain from
(\ref{Psi-Function-1}) that, for all $\mX=(p, z, u_2,
v_2,\mxy)\in{\mathcal A}$,
\begin{equation}
\label{Psi-Function-2} \Psi(\mX)=\rho_2'(c_2^2-\hat{\mxx}^2)p_1
+\big(\frac{\rho_2-\rho_1}{u_1}-\rho_2'\hat{\mxx} \big)(\mxy p_2-z)
+\hat E_1(\mX)\cdot p +\hat E_2(\mX)z,
\end{equation}
where $\hat E_1\in C^\infty(\overline{\mathcal A}; \bR^2)$ and $\hat
E_2\in C^\infty(\overline{\mathcal A})$ with
\begin{eqnarray*}
&&\|\hat E_i\|_{C^k(\overline{\mathcal A})}\le
\|\Psi\|_{C^{k+2}(\overline{\mathcal A})},
\qquad i=1,2, \quad k\in {\bf N}, \\
&&|\hat E_i(p,z, u_2, v_2, \mxy)| \le C(|p|+|z|+|u_2|+|v_2|)
\qquad\mbox{for all }(p, z, u_2, v_2, \mxy)\in {\mathcal A},
\end{eqnarray*}
for $C$ depending only on $\|D^2\Psi\|_{C(\overline{\mathcal A})}$.

{}From now on, we fix $(u_2, v_2)$ to be equal to the velocity of
state (2) obtained in \S \ref{section:3.3} and write $E_i(p, z,
\mxy)$ for $\hat E_i(p,z,u_2, v_2, \mxy)$. We conclude that, if
(\ref{condRewritingRH-0}) holds and  $\psi\in C^1(\ElDom)$ satisfies
(\ref{condRewritingRH}), then $\psi=\varphi-\varphi_2$ satisfies
(\ref{cont-accross-shock-mod-phi})--(\ref{RH-mod-phi}) on $\shock$
if and only if $\psi$ satisfies  conditions
(\ref{cont-accross-shock-psi-resolved}) on $\shock$,
\begin{equation}
\rho_2'(c_2^2-\hat{\mxx}^2)\psi_\mxx
+\big(\frac{\rho_2-\rho_1}{u_1}-\rho_2'\hat{\mxx}
\big)(\mxy\psi_\mxy-\psi) +E_1(\grad\psi,\psi, \mxy)\cdot \grad\psi
+E_2(\grad\psi,\psi,\mxy)\psi =0, \label{RH-psi-2}
\end{equation}
and the functions $E_i(p, z, \mxy), i=1,2,$ are smooth on
$$
\overline{B_{\delta^*}(0)\times(-\delta^*, \delta^*) \times (-6\bar
c_2/5, 6\bar c_2/5)}
$$
and satisfy that, for all $(p,z,\mxy)\in
B_{\delta^*}(0)\times(-\delta^*, \delta^*)\times (-6\bar c_2/5,
6\bar c_2/5)$,
\begin{equation}
|E_i(p,z,\mxy)|\le C\left(|p|+|z|+ \sigma\right)
\label{RH-psi-2-error-term1}
\end{equation}
and, for all $(p,z,\mxy)\in B_{\delta^*}(0)\times(-\delta^*,
\delta^*)\times (-6\bar c_2/5, 6\bar c_2/5)$,
\begin{equation}
|(\grad_{(p,z,\mxy)}E_i, \; D^2_{(p,z,\mxy)}E_i)|\le C,
\label{RH-psi-2-error-term2}
\end{equation}
where we have used (\ref{u2-v2-bound}) in the derivation of
(\ref{RH-psi-2-error-term1}) and $C$ depends only on the data.

\medskip
Denote by $\nu_0$ the unit normal on the reflected shock to the
region of state (2). Then $\nu_0=(\sin\theta_s, -\cos\theta_s)$ from
the definition of $\theta_s$. We compute
\begin{eqnarray}\label{obliquenessRH}
&&
\big(\rho_2'(c_2^2-\hat{\mxx}^2),(\frac{\rho_2-\rho_1}{u_1}-\rho_2'\hat{\mxx}
)\mxy \big) \cdot\nu_0\\
 &&=\rho_2'(c_2^2-\hat{\mxx}^2)\sin\theta_s-
\big(\frac{\rho_2-\rho_1}{u_1}-\rho_2'\hat\mxx
\big)\mxy\cos\theta_s\nonumber
\\
&&\ge \frac{1}{2}\rho_2'(c_2^2-\hat{\mxx}^2)>0,\nonumber
\end{eqnarray}
if $\pi/2-\theta_s$ is small and $\mxy\in Proj_\mxy(\shock)$. {}From
(\ref{thetaS-uv}) and (\ref{norm-to-Shock}), we obtain
$\|\nu_s-\nu_0\|_{L^\infty(\shock)}\le
C\|\grad\psi\|_{C(\overline\ElDom)}$. Thus, if $\sigma>0$ and
$\|\grad\psi\|_{C(\overline\ElDom)}$ are small depending only on the
data, then (\ref{RH-psi-2}) is an oblique derivative condition on
$\shock$.

\Subsec{\bf\large  The equation and boundary conditions near the
sonic circle} \label{sectionEqNearSonicLine} For the shock
reflection solution, equation (\ref{1.1.5}) is expected to be
elliptic in the domain $\ElDom$ and degenerate on the sonic circle
of state (2) which is the curve $\sonic=\partial \ElDom\cap \partial
B_{c_2}(0)$. Thus we consider the subdomains:
\begin{equation}\label{defOfSubdomains}
\begin{array}{l}
\ElDomS:=\ElDom\cap\{(\mxx,\mxy)\; : \;\dist((\mxx,\mxy),\sonic)<2\varepsilon\},\\
\ElDomU:=\ElDom\cap \{(\mxx,\mxy)\; :
\;\dist((\mxx,\mxy),\sonic)>\varepsilon\},
\end{array}
\end{equation}
where the small constant $\varepsilon>0$ will be chosen later.
Obviously, $\ElDomS$ and $\ElDomU$ are open subsets of $\ElDom$, and
$\ElDom=\ElDomS\cup\ElDomU$. Equation  (\ref{1.1.5}) is expected to
be degenerate elliptic in $\ElDomS$ and uniformly elliptic in
$\ElDomU$ on the solution of the shock reflection problem.

In order to display the structure of the equation near the sonic
circle where the ellipticity degenerates, we introduce the new
coordinates in $\ElDomS$ which flatten $\sonic$ and rewrite equation
(\ref{1.1.5}) in these new coordinates. Specifically, denoting $(r,
\theta)$ the polar coordinates in the $(\mxx, \mxy)$--plane, i.e.,
$(\mxx,\mxy)=(r\cos\theta, r\sin\theta)$, we consider the coordinates:
\begin{equation}\label{coordNearSonic}
x=c_2-r, \quad y=\theta-\theta_w \qquad \,\,
 \hbox{on } \ElDomS.
\end{equation}
By \S \ref{section:3.3}, the domain  $\DomS$ does not contain the
point $(\mxx,\mxy)=(0,0)$ if $\varepsilon$ is small. Thus, the
change of coordinates $(\mxx, \mxy)\to (x,y)$ is smooth and smoothly
invertible on $\ElDomS$. Moreover, it follows from the geometry of
domain $\ElDom$ especially from
(\ref{FBfunct-estimate})--(\ref{curved-straight-shock-match}) that,
if $\sigma>0$ is small, then, in the $(x,y)$--coordinates,
$$
\ElDomS=\{(x,y)\: : \; 0<x<2\varepsilon,\,\,
0<y<\pi+\mbox{arctan}\left({\mxy(x)/f(\mxy(x))}\right)-\theta_w \},
$$
where $\mxy(x)$ is the unique solution, close to $\mxy_1$, of the
equation $\mxy^2+f(\mxy)^2=(c_2-x)^2$.

We write the equation for $\psi$ in the $(x,y)$--coordinates. As
discussed in \S \ref{equationForPsiSection}, $\psi$ satisfies
equation
(\ref{potent-flow-nondiv-psi-Polar})--(\ref{speedSound-expression-psi-Polar-1})
in the polar coordinates. Thus, in the $(x,y)$--coordinates in
$\ElDomS$, the equation for $\psi$ is
\begin{equation}
\big(2x-(\gamma+1)\psi_x+O_1 \big)\psi_{xx} +O_2\psi_{xy} + ({1\over
c_2}+O_3 )\psi_{yy} -(1+O_4)\psi_{x} +O_5\psi_{y}=0,
\label{equationForPsi-sonicStruct}
\end{equation}
where
\begin{eqnarray}
\\
O_1(\grad\psi,\psi,x) &=& -\frac{x^2}{c_2}+{\gamma+1\over
2c_2}(2x-\psi_x)\psi_x -{\gamma-1\over c_2}\big(\psi+{1\over
2(c_2-x)^2}\psi_y^2\big),\nonumber
\\
O_2(\grad\psi,\psi, x)&=&-{2\over c_2(c_2-x)^2}(\psi_x+c_2-x)\psi_y,
\nonumber
\\
O_3(\grad\psi,\psi, x) &=&{1\over c_2(c_2-x)^2}\Big(x(2c_2-x)-
(\gamma-1)(\psi+(c_2-x)\psi_x+{1\over 2}\psi_x^2) \nonumber
\\
\label{erTerms-xy-nontrunc} && \qquad\qquad\qquad
-\frac{\gamma+1}{2(c_2-x)^2}\psi_y^2\Big), \nonumber
\\
O_4(\grad\psi,\psi, x) &=&\frac{1}{c_2-x}\Big(x- {\gamma-1\over
c_2}\big(\psi+(c_2-x)\psi_x+{1\over 2}\psi_x^2 +\frac{\psi_y^2}{2
(c_2-x)^2}\big)\Big), \nonumber
\\
O_5(\grad\psi,\psi, x)&=&
-\frac{1}{c_2(c_2-x)^3}\big(\psi_x+2c_2-2x\big)\psi_y. \nonumber
\end{eqnarray}
The terms $O_k(\grad\psi,\psi, x)$ are small perturbations of the
leading terms of equation (\ref{equationForPsi-sonicStruct}) if the
function $\psi$ is small in an appropriate norm considered below. In
order to see this, we note the following properties: For any
$(p,z,x)\in \bR^2\times\bR\times (0, c_2/2)$ with $|p|<1$,
\begin{eqnarray}\label{estSmallterms}
&&|O_1(p,z,x)|\le C(|p|^2+|z|+|x|^2), \nonumber\\
&&|O_3(p,z,x)|+|O_4(p,z,x)|\le C(|p|+|z|+|x|),\\
&&|O_2(p,z,x)|+|O_5(p,z,x)|\le C(|p|+|x|+1)|p|.\nonumber
\end{eqnarray}

In particular, dropping the terms $O_k$, $k=1,\dots, 5$, from
equation (\ref{equationForPsi-sonicStruct}), we obtain the {\bf
transonic small disturbance equation} (cf. \cite{Morawetz2}):
\begin{equation}
\big(2x-(\gamma+1)\psi_x \big)\psi_{xx} + \frac{1}{c_2}\psi_{yy}
-\psi_{x} =0. \label{equation-TSD}
\end{equation}

Now we write the boundary conditions on $\sonic$, $\shock$, and
$\wedgeB$ in the $(x,y)$--coordinates. Conditions
(\ref{condOnSonicLine-Psi}) and (\ref{condOnWedge-Psi}) become
\begin{eqnarray}
&&\psi=0\qquad\mbox{on }\;\sonic=\partial \ElDom\cap \{x=0\},
\label{condOnSonicLine-Psi-xy}
\\
&&\psi_\nu\equiv\psi_y=0\qquad\mbox{on }\;\wedgeB=\partial
\ElDom\cap \{y=0\}. \label{condOnWedge-Psi-xy}
\end{eqnarray}

It remains to write condition (\ref{RH-psi-2}) on $\shock$ in the
$(x,y)$--coordinates. Expressing $\psi_\mxx$ and $\psi_\mxy$ in the
polar coordinates $(r,\theta)$ and using (\ref{coordNearSonic}), we
write (\ref{RH-psi-2}) on $\shock\cap\{x<2\varepsilon\}$ in the
form:
\begin{equation}
\begin{array}{l}
\left(-\rho_2'(c_2^2-\hat{\mxx}^2) \cos(y+\theta_w) -
(\frac{\rho_2-\rho_1}{u_1}
 -\rho_2'\hat{\mxx})(c_2-x)\sin^2(y+\theta_w)\right)\psi_x
\displaystyle\\
\,\,+ \sin(y+\theta_w)
\left(-\frac{\rho_2'}{c_2-x}(c_2^2-\hat{\mxx}^2)
+(\frac{\rho_2-\rho_1}{u_1}
 -\rho_2'\hat{\mxx})\cos(y+\theta_w)\right)\psi_y
\displaystyle\\
\,\, -\left(\frac{\rho_2-\rho_1}{u_1}-\rho_2'\hat{\mxx} \right)\psi
+\tilde E_1(\grad_{(x,y)}\psi, \psi, x,y)\cdot \grad_{(x,y)}\psi+
\tilde E_2(\grad_{(x,y)}\psi, \psi, x,y)\psi=0, \displaystyle
\end{array}
\label{RH-psi-2-xy}
\end{equation}
where $\tilde E_i(p,z,x,y), i=1,2,$ are smooth functions of
$(p,z,x,y)\in\bR^2\times\bR\times\bR^2$ satisfying
$$
|\tilde E_i(p,z,x,y)|\le C\left(|p|+|z|+ \sigma\right)
\qquad\mbox{for }\, |p|+|z|+x\le\varepsilon_0(u_1,\bar{\rho}_2).
$$

We now rewrite (\ref{RH-psi-2-xy}). We note first that, in the
$(\mxx, \mxy)$--coordinates, the point $\PtUpL=\sonic\cap\shock$ has
the coordinates $(\mxx_1,\mxy_1)$ defined by (\ref{coord-P4}). Using
(\ref{theta_s-close-w}), (\ref{theta_s<w}),
(\ref{reflected-shock-s2}), and (\ref{coord-P4}), we find
$$
0\le |\hat\mxx|- |\mxx_1|\le C \sigma.
$$
In the $(x,y)$--coordinates, the point $\PtUpL$ is $(0, y_1)$, where
$y_1$ satisfies
\begin{equation}\label{xy-xieta-at-P1}
c_2\cos(y_1+\theta_w)=\mxx_1,\qquad c_2\sin(y_1+\theta_w)=\mxy_1,
\end{equation}
from (\ref{coord-P4}) and  (\ref{coordNearSonic}). Using this and
noting that the leading terms of the coefficients of
(\ref{RH-psi-2-xy}) near $\PtUpL=(0, y_1)$ are the coefficients at
$(x,y)=(0, y_1)$, we rewrite (\ref{RH-psi-2-xy}) as follows:
 \begin{equation}
\begin{array}{l}
-\frac{\rho_2-\rho_1}{u_1c_2} \mxy^2_1\psi_x
-\left(\rho_2'-\frac{\rho_2-\rho_1}{u_1
c_2^2}\mxx_1\right)\mxy_1\psi_y -\left(\frac{\rho_2-\rho_1}{u_1}-
\rho_2'\mxx_1\right)\psi
\displaystyle\\
+\hat E_1(\grad_{(x,y)}\psi, \psi, x,y)\cdot \grad_{(x,y)}\psi+ \hat
E_2(\grad_{(x,y)}\psi,\psi, x,y)\psi=0\,\,\mbox{on
}\shock\cap\{x<2\varepsilon\}, \displaystyle
\end{array}
\label{RH-psi-3-xy}
\end{equation}
where the terms $\hat E_i(p,z,x,y), i=1,2,$  satisfy
\begin{equation}
|\hat E_i(p,z,x,y)|\le C\left(|p|+|z|+x+|y-y_1|+\sigma\right) \qquad
\label{RH-psi-2-error-term-xy-1}
\end{equation}
for  $(p,z,x,y)\in {\mathcal
T}:=\{(p,z,x,y)\in\bR^2\times\bR\times\bR^2 :\;
|p|+|z|\le\varepsilon_0(u_1,\bar{\rho}_2)\}$ and
\begin{equation} \label{RH-psi-2-error-term-xy-2}
\|(\grad_{(p,z,x,y)}\hat E_i,\; D^2_{(p,z,x,y)}\hat
E_i)\|_{L^\infty( {\mathcal T})} \le C.
\end{equation}

We note that the left-hand side of (\ref{RH-psi-3-xy}) is obtained
by expressing the left-hand side of
 (\ref{RH-psi-2}) on $\shock\cap\{c_2-r<2\varepsilon\}$ in the $(x,y)$--coordinates.
Assume $\varepsilon<\bar c_2/4$. In this case, transformation
(\ref{coordNearSonic})
 is smooth on $\{0<c_2-r<2\varepsilon\}$ and has nonzero Jacobian.
Thus, condition (\ref{RH-psi-3-xy}) is equivalent to
(\ref{RH-psi-2}) and hence to
 (\ref{RH-psi}) on $\shock\cap\{x<2\varepsilon\}$
if $\sigma>0$ is small so that (\ref{condRewritingRH-0}) holds and
if $\|\psi\|_{C^1(\overline\ElDom)}$ is small depending only on the
data such that (\ref{condRewritingRH}) is satisfied.

\section{Iteration Scheme}
\label{iterSchemeSection}

In this section, we develop an iteration scheme to solve the free
boundary problem and set up the detailed steps of the iteration
procedure in the shifted coordinates.

\Subsec{\bf\large  Iteration domains} \label{iteration-dom-subsect}
Fix $\theta_w<\pi/2$ close to $\pi/2$. Since our problem is a free
boundary problem, the elliptic domain $\ElDom$ of the solution is
apriori unknown and thus we perform the iteration in a larger domain
\begin{equation}\label{ellipticDomainFull}
\Dom\equiv\Dom_{\theta_w}\defd B_{c_2}(0)\cap\{\mxy>-v_2\}\cap \{l(\mxy)<\mxx<
\mxy\cos\theta_w\},
\end{equation}
where $l(\mxy)$ is defined by (\ref{reflected-shock-s2}). We will
construct a solution with $\Omega\subset\Dom$. Moreover, the
reflected shock for this solution coincides  with $\{\mxx=l(\mxy)\}$
outside the sonic circle, which implies $\partial \Dom\cap\partial
B_{c_2}(0) =\partial \ElDom\cap\partial B_{c_2}(0)=:\sonic$.
Then we
decompose $\Dom$ similar to (\ref{defOfSubdomains}):
\begin{equation}\label{defOfSubdomains-iteration}
\begin{array}{l}
\displaystyle
\DomS:=\Dom\cap \{(\mxx,\mxy)\; : \;\dist((\mxx,\mxy),\sonic)<2\varepsilon\},\\
\displaystyle \DomU:=\Dom\cap \{(\mxx,\mxy)\; :
\;\dist((\mxx,\mxy),\sonic)>{\varepsilon/2}\}.
\end{array}
\end{equation}
The universal constant $C>0$ in the estimates of this section
depends only on the data and is independent on $\theta_w$.

\medskip
We will work in the $(x,y)$--coordinates (\ref{coordNearSonic}) in
the domain $\Dom\cap\{c_2-r<\kappa_0\}$, where $\kappa_0\in (0, \bar
c_2)$ will be determined depending only on the data for the sonic
speed $\bar c_2$  of state (2) for normal reflection (see
\S\ref{section:4}). Now we determine $\kappa_0$ so that
$\varphi_1-\varphi_2$  in the $(x,y)$--coordinates satisfies certain
bounds independent of $\theta_w$ in $\Dom\cap\{c_2-r<\kappa_0\}$ if
$\sigma=\pi/2-\theta_w$ is small.

We first consider the case of normal reflection $\theta_w=\pi/2$.
Then, from (\ref{flatOrthSelfSimShock2}) and (\ref{phi-2-a}) in the
$(x,y)$--coordinates (\ref{coordNearSonic}) with $c_2=\bar c_2$ and
$\theta_w=\pi/2$, we obtain
$$
\varphi_1-\varphi_2=-u_1 (\bar c_2-x)\sin y -u_1 \bar\xi \qquad
\mbox{for }\;0<x<\bar c_2, \;0<y<\pi/2.
$$
Recall $\bar\xi<0$ and $|\bar\xi|<\bar c_2$ by
(\ref{sonic-intersect-shock}). Then, in the region
$\Dom_0:=\{0<x<\bar c_2, \;0<y<\pi/2\}$, we have
$\varphi_1-\varphi_2=0$ only on the line
$$
y=\hat f_{0,0}(x)\defd\arcsin\big(\frac{|\bar\xi|}{\bar c_2-x}\big)
\qquad \mbox{for }x\in (0, \bar c_2-|\bar\xi|).
$$

Denote $\kappa_0\defd {(\bar c_2-|\bar\xi|)/2}$. Then $\kappa_0\in
(0, \bar c_2)$ by (\ref{sonic-intersect-shock-normal}) and  depends
only on the data. Now we show that there exists $\sigma_0>0$ small,
depending only on the data, such that, if
$\theta_w\in(\pi/2-\epsP_0, \pi/2)$, then
\begin{eqnarray}
&&C^{-1}\le\partial_x(\varphi_1-\varphi_2),
-\partial_y(\varphi_1-\varphi_2)\le C\;\label{nondegenPolar-1}\\
&& \qquad\qquad\qquad\qquad\qquad\mbox{ on }[0,\kappa_0]\times[{\hat
f_{0,0}(0)\over 2}, {\hat f_{0,0}(\kappa_0)+{\pi/2}\over 2}],
\nonumber
\\
&&\varphi_1-\varphi_2\ge C^{-1}>0 \qquad\mbox{on }
[0,\kappa_0]\times[0, \frac{\hat f_{0,0}(0)}{2}],
\label{nondegenPolar-2}
\\
&&\varphi_1-\varphi_2\le -C^{-1}<0 \qquad\mbox{on }
[0,\kappa_0]\times\{\frac{\hat f_{0,0}(\kappa_0)+\pi/2}{2}\},
\label{nondegenPolar-3}
 \end{eqnarray}
where ${\hat f_{0,0}(\kappa_0)+{\pi/2}\over 2}<{\pi/2}$.

\medskip
We first prove (\ref{nondegenPolar-1})--(\ref{nondegenPolar-3}) in
the case of normal reflection $\theta_w=\pi/2$. We compute from the
explicit expressions of $\varphi_1-\varphi_2$ and $\hat f_{0,0}$
given above to obtain
\begin{eqnarray*}
0<\arcsin\big({2|\bar\xi|\over \bar c_2+|\bar\xi|}\big)<\hat
f_{0,0}(x)< \arcsin\big({|\bar\xi|\over \bar c_2}\big)<{\pi\over
2},\quad && C^{-1}\le \hat f_{0,0}'(x)\le C\\
&&\qquad\quad\mbox{for }x\in [0,\kappa_0],
\end{eqnarray*}
$
\partial_x(\varphi_1-\varphi_2)=u_1\sin y$, and
$\partial_y(\varphi_1-\varphi_2)=-u_1(\bar c_2-x)\cos y$, which
imply (\ref{nondegenPolar-1}). Now, (\ref{nondegenPolar-2}) is true
since $\bar\xi=-\bar c_2\sin(\hat f_{0,0}(0))$ and thus
$\varphi_1-\varphi_2=u_1\big(\bar c_2\sin(\hat f_{0,0}(0))-(\bar
c_2-x)\sin y\big)$, and (\ref{nondegenPolar-3}) follows from
(\ref{nondegenPolar-1}) since $(\varphi_1-\varphi_2)(\kappa_0, \hat
f_{0,0}(\kappa_0))=0$ and ${(\hat f_{0,0}(\kappa_0)+{\pi/ 2})/
2}-\hat f_{0,0}(\kappa_0)\ge C^{-1}$.

\medskip
Now let $\theta_w<\pi/2$. Then, from
(\ref{thetaS-uv})--(\ref{phi-2-shifted})
and (\ref{coordNearSonic}), we have
$$
\varphi_1-\varphi_2=-(c_2-x)\sin(y+\theta_w-\theta_s)\sqrt{(u_1-u_2)^2+v_2^2}
-(u_1-u_2)\hat\xi.
$$
By \S\ref{section:3.3}, when $\theta_w\to \pi/2$, we know that
$(u_2, v_2)\to (0,0)$, $\theta_s\to\pi/2$, $\tilde\xi\to\bar\xi$,
and thus, by (\ref{x1-in-shifed}), we also have $\hat\xi\to\bar\xi$.
This shows that, if $\sigma_0>0$ is small depending only on the
data, then, for all $\theta_w\in(\pi/2-\epsP_0, \pi/2)$, estimates
(\ref{nondegenPolar-1})--(\ref{nondegenPolar-3}) hold
 with $C$ that is equal to twice the constant $C$ from the respective estimates
(\ref{nondegenPolar-1})--(\ref{nondegenPolar-3})
for $\theta_w=\pi/2$.

{}From (\ref{nondegenPolar-1})--(\ref{nondegenPolar-3}) for
$\theta_w\in(\pi/2-\epsP_0, \pi/2)$ and since
$$
\Dom\cap\{c_2-r<\kappa_0\} =\{\varphi_1>\varphi_2\}\cap\{0\le
x\le\kappa_0, 0\le y\le {\hat f_{0,0}(\kappa_0)+{\pi/2}\over 2} \},
$$
there exists $\hat f_0:=\hat f_{0,\pi/2-\theta_w}\in
C^\infty(\overline{\bR_+})$ such that
\begin{eqnarray}
\label{domain-in-xy-0}
&&\Dom\cap\{c_2-r<\kappa_0\}
=\{0< x< \kappa_0,\quad 0<y<\hat f_0(x)\},
\\
\label{domain-in-xy-funct-0} &&\hat f_0(0)=y_{\PtUpL},\qquad
C^{-1}\le \hat f_0'(x)\le C \,\,\,\, \mbox{ on } [0,\kappa_0],
\\
&&{\hat f_{0,0}(0)/2}\le \hat f_0(0)<\hat f_0(\kappa_0)\le
({\hat f_{0,0}(\kappa_0)+{\pi/2})/2}.
\label{fbFUnctionCloseToNormalPolar}
\end{eqnarray}
In fact, the line $y=\hat f_0(x)$ is the line $\mxx=l(\mxy)$ expressed
in the $(x,y)$--coordinates, and thus we obtain explicitly with the
use of (\ref{thetaS-uv}) that
\begin{equation}\label{referenceFB-polar} \hat
f_0(x)=\arcsin\big(\frac{|\hat\xi|\sin\theta_s}{(c_2-x)}\big)-\theta_w+\theta_s
\qquad\mbox{on }[0,\kappa_0].
\end{equation}

\Subsec{\bf\large  H\"{o}lder norms in $\ElDom$} For the elliptic
estimates, we need the H\"{o}lder norms in $\ElDom$ weighted by the
distance to the corners $\PtLwL=\shock\cap\{\mxy=-v_2\}$ and
$\PtLwR=(-u_2, -v_2)$, and with a ``parabolic'' scaling near the
sonic circle.

More generally, we consider a subdomain $\ElDom\subset \Dom$ of the
form $\ElDom:= \Dom\cap \{\mxx\ge f(\mxy)\}$ with $f\in C^1(\bR)$ and
set the subdomains
$\ElDomS:=\ElDom\cap\DomS$ and $\ElDomU:=\ElDom\cap\DomU$ defined by
(\ref{defOfSubdomains}). Let $\Sigma\subset\partial\ElDomU$ be
closed. We now introduce the H\"{o}lder norms in $\ElDomU$ weighted
by the distance to $\Sigma$. Denote by $X=(\mxx,\mxy)$ the points of
$\ElDomU$ and set
$$
\delta_X:=\dist(X,\Sigma),\quad \delta_{X,Y}:=\min(\delta_X,
\delta_Y) \qquad\mbox{for }\, X, Y\in \ElDomU.
$$
Then, for $k\in \bR$, $\alpha\in (0, 1)$, and $m\in {\bf N}$,
define
\begin{eqnarray}
&&\qquad\,\,|u\|^{(k,\Sigma)}_{m,0,\ElDomU} :=\sum_{0\le |\beta|\le
m} \sup_{X\in\ElDomU}
   \left(\delta_X^{\max(|\beta|+k,0)}|D^\beta u(X)|\right),
\nonumber \\
&&\qquad\,\, [u]^{(k,\Sigma)}_{m,\alpha,\ElDomU}
:=\sum_{|\beta|=m}\sup_{X,Y\in\ElDomU, X\ne Y}
 \left(\delta_{X,Y}^{\max(m+\alpha+k,0)}
  {\frac{|D^\beta u(X)-D^\beta u(Y)|}{|X-Y|^\alpha}}\right),
\label{weightNormsApp}
\\
&&\qquad\,\, \|u\|^{(k,\Sigma)}_{m,\alpha,\ElDomU}
:=\|u\|^{(k,\Sigma)}_{m,0,\ElDomU}
+[u]^{(k,\Sigma)}_{m,\alpha,\ElDomU}, \nonumber
\end{eqnarray}
where $D^\beta=\partial_{\mxx}^{\beta_1}\partial_{\mxy}^{\beta_2}$,
and $\beta=(\beta_1,\beta_2)$ is a multi-index with $\beta_j\in {\bf
N}$ and $|\beta|=\beta_1+\beta_2$.
We denote by $C^{(k,\Sigma)}_{m,\alpha,\ElDomU}$ the space of
functions with finite norm
$\|\cdot\|^{(k,\Sigma)}_{m,\alpha,\ElDomU}$.

\begin{remark}
If $m\ge -k\ge 1$ and $k$ is an integer, then any function $u\in
C^{(k,\Sigma)}_{m,\alpha,\ElDomU}$ is $C^{|k|-1,1}$  up to $\Sigma$,
but not necessarily $C^{|k|}$ up to $\Sigma$.
\end{remark}

\medskip
In $\ElDomS$, the equation is degenerate elliptic, for which the
H\"{o}lder norms with parabolic scaling are natural. We define the
norm $\|\psi\|_{2,\alpha,\ElDomS}^{(par)}$ as follows: Denoting
$z=(x,y)$ and $\tilde z=(\tilde x,\tilde y)$ with $x, \tilde x\in(0,
2\varepsilon)$ and
$$
\delta^{(par)}_\alpha(z, \tilde z)\defd \left(|x-\tilde x|^2+
\min(x,\tilde x)|y-\tilde y|^2\right)^{\alpha/2},
$$
then, for $u\in C^2(\ElDomS)\cap C^{1,1}(\overline\ElDomS)$ written
in the $(x,y)$--coordinates (\ref{coordNearSonic}), we define
\begin{eqnarray}
&&\qquad\,\, \|u\|^{(par)}_{2,0,\ElDomS} :=\sum_{0\le k+l\le 2}
\sup_{z\in\ElDomS}\left(x^{k+l/2-2}|\partial_x^k\partial_y^lu(z)|\right),
\nonumber\\
&&\qquad\,\, [u]^{(par)}_{2,\alpha,\ElDomS} :=\sum_{k+l=2}\sup_{z,
\tilde z\in\ElDomS, z\ne \tilde z}
 \bigg(\min(x,\tilde x)^{\alpha-l/2}
 \frac{|\partial_x^k\partial_y^lu(z)-\partial_x^k\partial_y^lu(\tilde z)|}
 {\delta^{(par)}_\alpha(z,\tilde z)}\bigg),
\label{parabNormsApp}  \\
&&\qquad\,\,\|u\|^{(par)}_{2,\alpha,\ElDomS}
:=\|u\|^{(par)}_{2,0,\ElDomS} +[u]^{(par)}_{2,\alpha,\ElDomS}.
\nonumber
\end{eqnarray}
To motivate this definition, especially
the parabolic scaling, we consider a scaled version of the function
$u(x,y)$ in the parabolic rectangles:
\begin{equation}\label{parabRectangles}
R_{(x,y)}=\Big\{(s,t)\;\;:\;\; |s-x|<\frac{x}{4},
|t-y|<\frac{\sqrt{x}}{4}\Big\}\cap \ElDom \qquad\mbox{for }\,
z=(x,y)\in \ElDomS.
\end{equation}
Denote  $Q_1\defd (-1, 1)^2$. Then the rescaled rectangle
(\ref{parabRectangles}) is
\begin{equation}\label{rescaled-parabRectangles}
Q_1^{(z)}:=\Big\{(S,T)\in Q_1\; : \; (x+\frac{x}{4}S,
y+\frac{\sqrt{x}}{4}T)\in \ElDom\Big\}.
\end{equation}
Denote by $u^{(z)}(S, T)$ the following function in $Q_1^{(z)}$:
\begin{equation}\label{parabRescaling}
u^{(z)}(S, T):=\frac{1}{x^2}u(x+\frac{x}{4}S, y+\frac{\sqrt{x}}{4}T)
\qquad\mbox{for } (S, T)\in Q_1^{(z)}.
\end{equation}
Then we have
$$
C^{-1}\sup_{z\in\ElDomS\cap\{x<3\varepsilon/2\}}
\|u^{(z)}\|_{C^{2,\alpha}\big(\overline{Q_1^{(z)}}\big)}\leq
\|u\|^{(par)}_{2,\alpha,\ElDomS}\leq
C\sup_{z\in\ElDomS}\|u^{(z)}\|_{C^{2,\alpha}\big(\overline
{Q_1^{(z)}}\big)},
$$
where $C$ depends only on the domain $\ElDom$ and is independent of
$\varepsilon\in (0, \kappa_0/2)$.

\Subsec{\bf\large  Iteration set} \label{iterSet-Section} We
consider the wedge angle close to $\pi/2$, that is,
$\epsP=\frac{\pi}{2}-\theta_w> 0$ is small which will be chosen
below. Set $\Sigma_0:=\partial\Dom\cap\{\mxy=-v_2\}$. Let
$\varepsilon, \epsP> 0$ be the constants from
(\ref{defOfSubdomains-iteration}) and (\ref{angleCloseToPiOver2}).
Let $M_1, M_2 \geq 1$. We define $\setK\equiv \setK(\epsP,
\varepsilon, M_1, M_2)$ by
\begin{equation}\label{defSetK_R}
\setK \defd
\bigg\{
 \Kphi \in C^{1,\alpha}(\overline\Dom)\cap C^{2}(\Dom)
\, : \, \|\Kphi\|_{2,\alpha,\DomS}^{(par)} \leq M_1,
    \|\Kphi\|_{2,\alpha,\DomU}^{(-1-\alpha, \Sigma_0)} \leq M_2\epsP,
    \Kphi\ge 0 \mbox{ in }\Dom
\bigg\} \quad
\end{equation}
for $\alpha\in (0, 1/2)$. Then $\setK$
is convex. Also,  $\Kphi\in\setK$ implies that
$$
\|\Kphi\|_{C^{1,1}(\overline\DomS)}\le M_1,\qquad
\|\Kphi\|_{C^{1,\alpha}(\overline\DomU)}\le M_2\epsP,
$$
so that $\setK$ is a bounded subset in
$C^{1,\alpha}(\overline\Dom)$. Thus, $\setK$ is a compact and convex
subset of $C^{1,\alpha/2}(\overline\Dom)$.

We note that the choice of  constants $M_1, M_2\ge 1$ and
$\varepsilon,\epsP>0$ below will guarantee the following property:
\begin{equation}
\epsP\max(M_1,M_2)+\varepsilon^{1/4} M_1+ \epsP
M_2/\varepsilon^2 \le {\hat C}^{-1} \label{condConst-00}
\end{equation}
for some sufficiently large $\hat C>1$ depending only on the data.
In particular,  (\ref{condConst-00}) implies that $\epsP\le {\hat
C}^{-1}$ since $\max(M_1,M_2)\ge 1$, which implies
$\pi/2-\theta_w\le{\hat C}^{-1}$ from (\ref{angleCloseToPiOver2}).
Thus, if we choose $\hat C$ large depending only on the data, then
(\ref{condRewritingRH-0}) holds. Also, for $\psi\in \setK$, we have
$$
|(\grad\psi, \psi)(x,y)|\le M_1x^2+M_1x \,\,\mbox{ in }\, \DomS,
\qquad \|\psi\|_{C^1(\bar{\DomU})}\le M_2\epsP.
$$
Furthermore, $0<x<2\varepsilon$ in $\DomS$ by (\ref{coordNearSonic})
and (\ref{defOfSubdomains-iteration}). Now it follows from
(\ref{condConst-00}) that $\|\psi\|_{C^1}\le 2/{\hat C}$. Then
(\ref{condRewritingRH}) holds if $\hat C$ is large depending only on
the data. Thus, in the rest of this paper, we always assume that
(\ref{condRewritingRH-0}) holds and that $\psi\in \setK$ implies
(\ref{condRewritingRH}). Therefore, (\ref{RH-psi}) is equivalent to
(\ref{RH-psi-2-error-term1})--(\ref{RH-psi-2-error-term2}) for
$\psi\in \setK$.

We also note the following fact.
%*****************BEGIN LEMMA*************
\begin{lemma}\label{relatingNorms}
%Let $\alpha\in(0,1/2)$.
There exist $\hat C$ and $C$ depending only on the data such that,
if $\epsP, \varepsilon>0$ and $M_1, M_2\ge 1$ in {\rm
(\ref{defSetK_R})} satisfy {\rm (\ref{condConst-00})}, then, for
every $\Kphi\in\setK$,
\begin{equation}\label{relateHolderNorms}
\|\Kphi\|_{2,\alpha,\Dom}^{(-1-\alpha, \Sigma_0\cup\sonic)} \le
C(M_1\varepsilon^{1-\alpha}+M_2\epsP).
\end{equation}
\end{lemma}
%*****************END LEMMA*************
\Proof
%*****************BEGIN PROOF LEMMA*************
In this proof, $C$ denotes a universal constant depending only on
the data. We use definitions
(\ref{weightNormsApp})--(\ref{parabNormsApp}) for the norms. We
first show that
\begin{equation}\label{relateHolderNorms-Sing}
\|\Kphi\|_{2,\alpha,\DomS}^{(-1-\alpha, \sonic)} \le
CM_1\varepsilon^{1-\alpha},
\end{equation}
where $\delta_{(x,y)}:=\dist((x,y), \sonic)$ in
(\ref{weightNormsApp}). First we show (\ref{relateHolderNorms-Sing})
in the $(x,y)$--coordinates. Using (\ref{domain-in-xy-0}), we have
$\DomS=\{0< x< 2\varepsilon,\, 0<y<\hat f_0(x)\}$ with
$\sonic=\{x=0,\; 0<y<\hat f_0(x)\}$, where $\|f_0'\|_{L^\infty((0,
2\varepsilon))}$ depends only the data, and thus $\dist((x,y),
\sonic)\le Cx$ in $\DomS$. Then, since
$\|\Kphi\|_{2,\alpha,\DomS}^{(par)} \leq M_1$, we obtain that, for
$(x,y)\in\DomS$,
\begin{eqnarray*}
&&|\Kphi(x,y)|\le M_1x^2\le  M_1\varepsilon^2, \qquad
|D\Kphi(x,y)|\le M_1x\le  M_1\varepsilon, \\
&&\delta_{(x,y)}^{1-\alpha}|D^2\Kphi(x,y)|
=x^{1-\alpha}|D^2\Kphi(x,y)|\le \varepsilon^{1-\alpha}M_1.
\end{eqnarray*}
Furthermore, from (\ref{condConst-00}) with $\hat C\ge 16$, we
obtain $\varepsilon\le 1/2$. Thus, denoting $z=(x,y)$ and $\tilde
z=(\tilde x,\tilde y)$ with $x, \tilde x\in(0, 2\varepsilon)$, we
have
\begin{eqnarray*}
\delta^{(par)}_\alpha(z, \tilde z)&\defd& \left(|x-\tilde x|^2+
\min(x,\tilde x)|y-\tilde y|^2\right)^{\alpha/2}\\
&\le& \left(|x-\tilde x|^2+ 2\varepsilon|y-\tilde
y|^2\right)^{\alpha/2}\le |z- \tilde z|^\alpha,
\end{eqnarray*}
and $\min(\delta_{z},\delta_{\tilde z})=\min(x, \tilde x)$, which
implies
\begin{eqnarray*}
\min(\delta_{z},\delta_{\tilde z}) {|D^2\Kphi(z)-D^2\Kphi(\tilde
z)|\over |z- \tilde z|^\alpha} &\le& C\varepsilon^{1-\alpha}\min(x,
\tilde x)^\alpha {|D^2\Kphi(z)-D^2\Kphi(\tilde z)|\over
\delta^{(par)}_\alpha(z, \tilde z)}\\
&\le& C\varepsilon^{1-\alpha}M_1.
\end{eqnarray*}
Thus we have
proved (\ref{relateHolderNorms-Sing}) in the $(x,y)$--coordinates.
By (\ref{condRewritingRH-0}) and (\ref{condConst-00}), we have
$\varepsilon\le c_2/50$ if  $\hat C$ is large depending only on the
data. Then the change $(\mxx, \mxy)\to (x,y)$ in $\DomS$ and its
inverse have bounded $C^3$--norms in terms of the data. Thus,
(\ref{relateHolderNorms-Sing}) holds in the $(\mxx,
\mxy)$--coordinates.

Since $\Kphi\in\setK$, then $\|\Kphi\|_{2,\alpha,\DomU}^{(-1-\alpha,
\Sigma_0)} \leq M_2\epsP$. Thus, in order to complete the proof of
(\ref{relatingNorms}), it suffices to estimate
$\{\min(\delta_{z},\delta_{\tilde z}) {|D^2\Kphi(z)-D^2 \Kphi(\tilde
z)|\over |z- \tilde z|^\alpha}\}$ in the case
$z\in\DomS\setminus\DomU$ and $\tilde z\in\DomU\setminus\DomS$ for
$\delta_z=\dist(z, \sonic\cup\Sigma_0)$. {}From
$z\in\DomS\setminus\DomU$ and $\tilde z\in\DomU\setminus\DomS$, we
obtain $0<c_2-|z|<\varepsilon/2$ and $c_2-|\tilde z|\ge
2\varepsilon$, which implies that $|z-\tilde z|\ge 3\varepsilon/2$.
We have $c_2-|z|\le\dist(z,\sonic)\le C(c_2-|z|)$, where we have
used (\ref{condRewritingRH-0}) and (\ref{ellipticDomainFull}). Thus,
$\min(\delta_z, \delta_{\tilde z})\le C(c_2-|z|)\le C\varepsilon$.
Also we have $|D^2\Kphi(z)|\le M_1$ by (\ref{parabNormsApp}). If
$\delta_{\tilde z}\ge \delta_z$, then $\delta_{\tilde z}\ge
\varepsilon/2$ and thus $|D^2\Kphi(\tilde z)|\le
(\varepsilon/2)^{-1+\alpha}M_2\epsP$ by (\ref{weightNormsApp}). Then
we have
$$
\min(\delta_{z},\delta_{\tilde z}) {|D^2\Kphi(z)-D^2\Kphi(\tilde
z)|\over |z- \tilde z|^\alpha} \le
C\varepsilon{M_1+(2\varepsilon)^{-1+\alpha}M_2\epsP\over
(3\varepsilon/2)^\alpha} \le
C\big(\varepsilon^{1-\alpha}M_1+M_2\epsP\big).
$$
If $\delta_{\tilde z}\le \delta_z$, then $\dist(\tilde z,
\Sigma_0)\le \dist(\tilde z, \sonic)$, which implies by
(\ref{inSonicRegion-in-shifed}) that $|z-\tilde z|\ge 1/C$ if
$\varepsilon$ is sufficiently small, depending only on the data.
Then $|D^2\Kphi(\tilde z)|\le \delta_{\tilde z}^{-1+\alpha}M_2\epsP$
and
$$
\min(\delta_{z},\delta_{\tilde z}) {|D^2 \Kphi(z)-D^2 \Kphi(\tilde
z)|\over |z- \tilde z|^\alpha} \le
C\big(\delta_{z}M_1+\delta_{\tilde z}\delta_{\tilde
z}^{-1+\alpha}M_2\epsP\big) \le C\big(\varepsilon M_1+M_2\epsP\big).
$$
\Endproof

\Subsec{\bf\large  Construction of the iteration scheme and choice
of $\alpha$} \label{Constr-iter-section} In this section, for
simplicity of notations, the universal constant $C$ depends only on
the data and may be different at each occurrence.

By (\ref{u2-v2-bound}), it follows that, if $\epsP$ is sufficiently
small depending on the data, then
\begin{equation}\label{q2-u1}
q_2\le u_1/10,
\end{equation}
where $q_2=\sqrt{u_2^2+v_2^2}$.
Let $\Kphi\in\setK$. From
(\ref{phi-1-shifted})--(\ref{phi-2-shifted}) and (\ref{q2-u1}), it follows that
\begin{equation}\label{nondegeneracy}
(\varphi_1-\varphi_2-\Kphi)_{\mxx}(\mxx,\mxy)\geq u_1/2>0\quad
\mbox{in}\;\;\Dom.
\end{equation}
Since $\varphi_1-\varphi_2=0$ on $\{\mxx=l(\mxy)\}$ and $\Kphi\ge 0$
in $\Dom$, we have $\Kphi\ge \varphi_1-\varphi_2$ on
$\{\mxx=l(\mxy)\}\cap\partial\Dom$, where $l(\mxy)$ is defined by
(\ref{reflected-shock-s2}).
Then there exists $f_\Kphi\in
C^{1,\alpha}(\bR)$ such that
\begin{equation}\label{shockPL}
\{\Kphi=\varphi_1-\varphi_2\}\cap\Dom =\{(f_\Kphi(\mxy),\mxy)\; :
\;\mxy\in(-v_2,\eta_2)\}.
\end{equation}
It follows that $f_\Kphi(\mxy)\ge l(\mxy)$ for all $\mxy\in[-
v_2,\eta_2)$ and
\begin{equation}\label{OmegaPL}
\Omega^+(\Kphi):=\{\mxx>f_\Kphi(\mxy)\}\cap \Dom =
\{\Kphi<\varphi_1-\varphi_2\}\cap \Dom.
\end{equation}
Moreover,
$\partial\Omega^+(\Kphi)=\shock\cup\sonic\cup\wedgeB\cup\Sigma_0$,
where
\begin{equation}\label{shockIterDef}
\begin{array}{l}
\displaystyle
\shock(\Kphi):=\{\mxx=f_\Kphi(\mxy)\}\cap\partial\Omega^+(\Kphi),\qquad
\displaystyle \sonic:=\partial \Dom\cap
\partial B_{c_2}(0),\\
\displaystyle \wedgeB:=\partial \Dom\cap
\{\mxy=\mxx\tan\theta_w\},\qquad \displaystyle
\Sigma_0(\Kphi):=\partial \Omega^+(\Kphi)\cap \{\mxy=-v_2\}.
\end{array}
\end{equation}
We denote by $P_j, 1\le j\le 4$, the corner points of
$\Omega^+(\Kphi)$. Specifically, $\PtLwL=\shock(\Kphi)\cap \Sigma_0(\Kphi)$
and $\PtLwR=(-u_2, -v_2)$ are the corners on the symmetry line $\{\mxy=-v_2\}$,
and $\PtUpL=\sonic\cap\shock(\Kphi)$ and
$\PtUpR=\sonic\cap\wedgeB$ are the corners on the sonic circle. Note
that, since $\Kphi\in\setK$ implies $\Kphi=0$ on $\sonic$, %%
it follows that $\PtUpL$ is the intersection point $(\mxx_1,\mxy_1)$
of the line $\mxx=l(\mxy)$ and the sonic circle $\mxx^2+\mxy^2=c_2^2$,
where $(\mxx_1,\mxy_1)$ is determined by (\ref{coord-P4}).

We also note that $f_0=l$ for $0\in\setK$. From $\Kphi\in\setK$ and
Lemma \ref{relatingNorms} with $\alpha\in (0, 1/2)$, we obtain the
following estimate of $f_\Kphi$ on the interval $(-v_2, \mxy_1)$:
\begin{eqnarray}
&&\|f_\Kphi-l\|_{2,\alpha,(-v_2, \mxy_1)}^{(-1-\alpha, \{-v_2,
\mxy_1\})} \leq C\big(M_1\varepsilon^{1/2}+M_2\epsP\big)\le
\varepsilon^{1/4}, \label{OmegaPL-f-higher}
\end{eqnarray}
where the second inequality in (\ref{OmegaPL-f-higher}) follows from
(\ref{condConst-00}) with sufficiently large $\hat C$.

\medskip
We also work in the $(x,y)$--coordinates. Denote $\kappa\defd\kappa_0/2$.
Choosing $\hat C$ in (\ref{condConst-00}) large
depending only on the data, we conclude from
(\ref{nondegenPolar-1})--(\ref{nondegenPolar-3}) that, for every
$\Kphi\in\setK$, there exists a function $\hat f\equiv \hat
f_\Kphi\in C^{(-2, \{0\})}_{2,\alpha,(0,\kappa)}$ such that
\begin{equation}\label{domain-in-rescaled-lemma}
\Omega^+(\Kphi)\cap\{c_2-r<\kappa\} =\{0< x< \kappa,\quad 0<y<\hat
f_\Kphi(x)\},
\end{equation}
with
\begin{equation}\label{holder-hat-f}
\hat f_\Kphi(0)=\hat f_0(0)>0, \quad \hat f_\Kphi'>0 \mbox{ on }
(0,\kappa), \quad \|\hat f_\Kphi-\hat f_0\|^{(-1-\alpha,
\{0\})}_{2,\alpha,(0,\kappa)}\le
C\big(M_1\varepsilon^{1-\alpha}+M_2\epsP\big),
\end{equation}
where we have used Lemma \ref{relatingNorms}. More precisely,
\begin{equation}
\begin{array}{l}
\displaystyle \sum_{k=0}^2
\sup_{x\in(0,2\varepsilon)}\big(x^{k-2}|D^k(\hat f_\Kphi-\hat
f_0)(x)|\big)
\\
\displaystyle \qquad + \sup_{x_1\ne x_2\in(0,2\varepsilon)}
\Big((\min(x_1,x_2))^\alpha\, {|(\hat f_\Kphi''-\hat f_0'')(x_1)-
(\hat f_\Kphi''-\hat f_0'')(x_2)| \over |x_1-x_2|^\alpha}\Big) \le
CM_1,
\\
\displaystyle
\|\hat f_\Kphi-\hat f_0\|_{2,\alpha,({\varepsilon/2},
\kappa)}\le CM_2\epsP.
\end{array}
\label{holder-hat-f-S}
\end{equation}

\smallskip
Note that, in the $(\xi,\eta)$--coordinates, the angles
$\theta_{\PtLwL}$ and $\theta_{\PtLwR}$ at the corners $\PtLwL$ and
$\PtLwR$ of $\Omega^+(\Kphi)$ respectively satisfy
\begin{equation}\label{anglesCloseToPi2}
|\theta_{P_i}-\frac{\pi}{2}|\le \frac{\pi}{16}\qquad\mbox{for }i=2,3.
\end{equation}
Indeed,
$\theta_{\PtLwR}=\pi/2-\theta_w$.
The estimate for $\theta_{\PtLwL}$ follows
from (\ref{OmegaPL-f-higher}) with (\ref{condConst-00}) for large $\hat C$.

We now consider the following problem in the domain
$\Omega^+(\Kphi)$:
\begin{eqnarray}
&&\Nl(\psi)\defd A_{11}\psi_{\mxx\mxx}+
2A_{12}\psi_{\mxx\mxy}
+A_{22}\psi_{\mxy\mxy}=0 \qquad
\mbox{ in }\;\;\Omega^+(\Kphi),
\label{iterationEquation} \\
&&\Ml(\psi)\defd \rho_2'(c_2^2-\hat{\mxx}^2)\psi_\mxx +
\big(\frac{\rho_2-\rho_1}{u_1}-\rho_2'\hat\mxx
\big)(\mxy\psi_\mxy-\psi)
\label{iterationRH}\\
&&\qquad\qquad\,\, +E_1^\Kphi(\mxx,\mxy)\cdot
D\psi+E_2^\Kphi(\mxx,\mxy)\psi=0 \qquad\mbox{on }\;\shock(\Kphi),
\nonumber\\
&&\psi=0\qquad\,\,\mbox{on }\;\sonic,
\label{iterationCondOnSonicLine}
\\
&&\psi_\nu=0\qquad\,\, \mbox{on }\;\wedgeB,
\label{iterationCondOnWedge}
\\
&&\psi_\mxy=-v_2\qquad\mbox{on }\;\partial \Omega^+(\Kphi)\cap
\{\mxy=-v_2\}, \label{iterationCondOnSymmtryLine}
\end{eqnarray}
where $A_{ij}=A_{ij}(\grad\psi,\mxx,\mxy)$ will be defined below, and
equation (\ref{iterationRH}) is obtained from (\ref{RH-psi-2}) by
substituting $\Kphi$ into $E_i, i=1,2,$ i.e.,
\begin{equation}
E_i^\Kphi(\mxx, \mxy)=E_i(\grad\Kphi(\mxx, \mxy),\Kphi(\mxx, \mxy), \mxy).
\label{iteration-RH-error-term}
\end{equation}
Note that, for $\Kphi\in\setK$ and $(\mxx, \mxy)\in \Dom$, we have
$(\grad\Kphi(\mxx, \mxy),\Kphi(\mxx, \mxy), \mxy)\in
B_{\delta^*}(0)\times(-\delta^*, \delta^*)\times (-6\bar c_2/5,
6\bar c_2/5)$ by (\ref{condRewritingRH-0})--(\ref{condRewritingRH}).
Thus, the right-hand side of (\ref{iteration-RH-error-term}) is
well-defined.

Also, we now fix $\alpha$ in the definition of $\setK$. Note that
the angles $\theta_{\PtLwL}$ and $\theta_{\PtLwR}$ at the corners
$\PtLwL$ and $\PtLwR$ of $\Omega^+(\Kphi)$ satisfy
(\ref{anglesCloseToPi2}). Near these corners, equation
(\ref{iterationEquation}) is linear and its ellipticity constants
near the corners are uniformly bounded in terms of the data.
Moreover, the directions in the oblique derivative conditions on the
arcs meeting at the corner  $\PtLwR$ (resp. $\PtLwL$) are at the
angles within the range $(7\pi/16, 9\pi/16)$, since
(\ref{iterationRH}) can be written in the form
$\psi_\mxx+e\psi_\mxy-d\psi=0$, where $|e|\le C\epsP$ near $\PtLwL$
from $\mxy(\PtLwL)=-v_2$, (\ref{u2-v2-bound}),
(\ref{RH-psi-2-error-term1})--(\ref{RH-psi-2-error-term2}), and
(\ref{condConst-00}). Then, by \cite{Lieberman88}, there exists
$\alpha_0\in(0, 1)$ such that, for any  $\alpha\in (0,\alpha_0)$,
the solution of
(\ref{iterationEquation})--(\ref{iterationCondOnSymmtryLine}) is in
$C^{1,\alpha}$ near and up to $\PtLwL$ and $\PtLwR$ if  the arcs are
in $C^{1,\alpha}$ and the coefficients of the equation and the
boundary conditions are in the appropriate H\"{o}lder spaces with
exponent $\alpha$. We use $\alpha=\alpha_0/2$ in the definition of
$\setK$ for $\alpha_0=\alpha_0(9\pi/16, 1/2)$, where
$\alpha_0(\theta_0, \varepsilon)$ is defined in \cite[Lemma
1.3]{Lieberman88}. Note that $\alpha\in(0, 1/2)$ since
$\alpha_0\in(0, 1)$.

\Subsec{\bf\large  An elliptic cutoff and the equation for the
iteration} \label{eqForIterationSection} In this subsection, we fix
$\Kphi\in\setK$ and define equation (\ref{iterationEquation}) such
that

(i) It is strictly elliptic inside the domain $\Omega^+(\Kphi)$ with
elliptic degeneracy at the sonic circle $\sonic=\partial
\Omega^+(\Kphi)\cap\partial B_{c_2}(0)$;

(ii) For a fixed point $\psi=\Kphi$ satisfying an appropriate
smallness condition of $|\grad\psi|$, equation
(\ref{iterationEquation}) coincides with the original equation
\eqref{potent-flow-nondiv-psi-1}.

\smallskip
We define the coefficients $A_{ij}$ of equation
(\ref{iterationEquation}) in the larger domain $\Dom$. More
precisely, we define the coefficients separately in the domains
$\DomS$ and $\DomU$ and then combine them.
%by using a partition of unity in $\Dom$.

\medskip
In $\DomU$, we define  the coefficients of (\ref{iterationEquation}) by
substituting $\Kphi$ into the coefficients
of (\ref{potent-flow-nondiv-psi-1}), i.e.,
\begin{equation}\label{iterationUniforDomEquation}
\begin{array}{ll}
A^1_{11}(\mxx,\mxy)= c^2(\grad\Kphi, \Kphi, \xi, \eta)
-(\Kphi_\mxx-\mxx)^2, \,\,
 A^1_{22}(\mxx,\mxy)=c^2(\grad\Kphi, \Kphi,
\xi, \eta) -(\Kphi_\mxy-\mxy)^2, \\
A^1_{12}(\mxx,\mxy)=A^1_{21}(\mxx,\mxy)=
-(\Kphi_\mxx-\mxx)(\Kphi_\mxy-\mxy),
\end{array}
\end{equation}
where $\Kphi, \Kphi_\mxx$, and $\Kphi_\mxy$ are
evaluated at $(\mxx,\mxy)$. Thus, (\ref{iterationEquation}) in
$\Omega^+(\Kphi)\cap\DomU$ is a linear equation
$$
A^1_{11}\psi_{\mxx\mxx}+
2A^1_{12}\psi_{\mxx\mxy}
+A^1_{22}\psi_{\mxy\mxy}=0 \qquad
\mbox{ in }\;\;\Omega^+(\Kphi)\cap\DomU.
$$
{}From the definition of $\DomU$, it follows that
$\sqrt{\mxx^2+\mxy^2}\le c_2-\varepsilon$ in $\DomU$. Then
calculating explicitly the eigenvalues of matrix $(A^1_{ij})_{1\le
i,j\le 2}$ defined by (\ref{iterationUniforDomEquation}) and using
(\ref{condRewritingRH-0}) yield that there exists $C=C(\gamma, \bar
c_2)$ such that, if $\varepsilon<\min(1, \bar c_2)/10$ and
$\|\Kphi\|_{C^1}\le \varepsilon/C$, then
\begin{equation}
\label{ellipticityInUniformDomain} {\varepsilon\bar c_2\over
8}|\mu|^2\le \sum_{i,j=1}^2A^1_{ij}(\mxx,\mxy)\mu_i\mu_j\le 4\bar
c_2^2|\mu|^2\qquad \mbox{for any $(\mxx, \mxy)\in\DomU$ and
$\mu\in\bR^2$.}
\end{equation}
The required smallness of  $\varepsilon$ and $\|\Kphi\|_{C^1}$ is
achieved by choosing sufficiently large $\hat C$ in
(\ref{condConst-00}), since $\Kphi\in\setK$.

\medskip
In $\DomS$, we use (\ref{equationForPsi-sonicStruct}) and substitute
$\Kphi$ into the terms $O_1,\dots,O_5$. However, it is essential
that we do not substitute $\Kphi$ into the term $(\gamma+1)\psi_x$
of the coefficient of $\psi_{xx}$ in
(\ref{equationForPsi-sonicStruct}), since this nonlinearity allows
us to obtain some crucial estimates (see Lemma
\ref{quadraticGrowthPsi-Lemma} and Proposition
\ref{boundPsiXfromAbove-Prop}). Thus, we make an elliptic cutoff of
this term. In order to motivate our construction, we note that, if
$$
|O_k|\le \frac{x}{10\max(c_2,1)(\gamma+1)}, \qquad
\psi_x<\frac{4x}{3(\gamma+1)}\qquad\quad \mbox{ in }\, \DomS,
$$
then equation (\ref{equationForPsi-sonicStruct}) is strictly
elliptic in $\DomS$. Thus we want to replace the term
$(\gamma+1)\psi_x$ in the coefficient of $\psi_{xx}$ in
(\ref{equationForPsi-sonicStruct}) by $\displaystyle
(\gamma+1)x\zeta_1\big(\frac{\psi_x}{x} \big)$, where
$\zeta_1(\cdot)$ is a cutoff function. On the other hand, we also
need to keep form (\ref{iterationEquation}) for the modified
equation in the $(\xi,\eta)$--coordinates, i.e., the form without
lower-order terms. This form is used in Lemma
\ref{negativeDerivPsiLemma}. Thus we perform a cutoff in equation
(\ref{potent-flow-nondiv-psi-1}) in the $(\xi,\eta)$--coordinates
such that the modified equation satisfies the following two
properties:

(i) Form (\ref{iterationEquation}) is preserved;

(ii) When written in the $(x,y)$--coordinates, the modified equation
has the main terms as in (\ref{equationForPsi-sonicStruct}) with the
cutoff described above and corresponding modifications in the terms
$O_1,\dots, O_5$ of (\ref{equationForPsi-sonicStruct}).

Also, since the equations in $\DomS$ and $\DomU$ will be combined
and the specific form of the equation is more important in $\DomS$,
we define our equation in a larger domain
$\DomS_{4\varepsilon}\defd\Dom\cap \{c_2-r<4\varepsilon\}$.

We first rewrite equation (\ref{potent-flow-nondiv-psi-1})
in the form  %%, related to polar coordinates:
$$
I_1+I_2+I_3+I_4=0,
$$
where
\begin{eqnarray*}
&&I_1\defd \big(c^2(D\psi,\psi,\xi,\eta)-(\xi^2+\eta^2)\big)
\Delta\psi,\\
&&I_2\defd
\mxy^2\psi_{\mxx\mxx}+\mxx^2\psi_{\mxy\mxy}-2\mxx\mxy\psi_{\mxx\mxy},
\\
&& I_3\defd 2\big(\mxx\psi_\mxx\psi_{\mxx\mxx}
+(\mxx\psi_\mxy+\mxy\psi_\mxx)\psi_{\mxx\mxy}+\mxy\psi_\mxy\psi_{\mxy\mxy}\big),\\
&&I_4\defd -\frac{1}{2}\left(\psi_\mxx(|\grad\psi|^2)_\mxx
+\psi_\mxy(|\grad\psi|^2)_\mxy\right).
\end{eqnarray*}
Note that, in
the polar coordinates, $I_1,\dots,I_4$ have the following
expressions:
\begin{eqnarray*}
&&I_1=
\big(c_2^2-r^2+(\gamma-1)(r\psi_r-\frac{1}{2}|\grad\psi|^2-\psi)\big)
\Delta\psi,\\
&&I_2= \psi_{\theta\theta}+r\psi_r,\\
&&I_3=r(|\grad\psi|^2)_r=2r\psi_r\psi_{rr}
   +\frac{2}{r^2}\psi_\theta\psi_{r\theta}
-\frac{2}{r^2}\psi_{\theta}^2,\\
&& I_4=-\frac{1}{2}\big(\psi_r(|\grad\psi|^2)_r
+\frac{1}{r^2}\psi_\theta(|\grad\psi|^2)_\theta\big)
\end{eqnarray*}
with
$|\grad\psi|^2=\psi_r^2+\frac{1}{r^2}\psi_\theta^2$ and
$\Delta\psi=\psi_{rr}+\frac{1}{r^2}\psi_{\theta\theta}+\frac{1}{r}\psi_r$.

{}From this, by (\ref{coordNearSonic}), we see that the dominating
terms of (\ref{equationForPsi-sonicStruct}) come only from $I_1,
I_2$, and the term $2r\psi_r\psi_{rr}$ of $I_3$, i.e., the remaining
terms of $I_3$ and $I_4$ affect only the terms $O_1,\dots, O_5$ in
(\ref{equationForPsi-sonicStruct}). Moreover, the term
$(\gamma+1)\psi_x$ in the coefficient of $\psi_{xx}$ in
(\ref{equationForPsi-sonicStruct}) is obtained as the leading term
in the sum of the coefficient $(\gamma-1)r\psi_r$ of $\psi_{rr}$ in
$I_1$ and the coefficient $2r\psi_r$ of $\psi_{rr}$ in $I_3$. Thus
we modify the terms $I_1$ and $I_3$ by cutting off the
$\psi_r$-component of first derivatives in the coefficients of
second-order terms as follows.  Let $\zeta_1\in C^\infty(\bR)$
satisfy
\begin{equation}\label{defZeta-1}
\zeta_1(s)=\left\{
\begin{array}{ll}
s,\quad&\displaystyle \mbox{if }\;|s|<4/\big(3(\gamma+1)\big),\\
\displaystyle 5\,\mbox{sign}(s)/[3(\gamma+1)], \quad&\displaystyle
\mbox{if }\; |s|>2/(\gamma+1),
\end{array}
\right.
\end{equation}
so that
\begin{eqnarray}
&&\zeta_1'(s)\ge 0,\quad
\zeta_1(-s)=-\zeta_1(s)\qquad \mbox{on }\;\bR;\label{defZeta-1a} \\
&& \zeta_1''(s)\le 0\qquad \mbox{on }\;\{s\ge
0\}.\label{zeta-1-concave}
\end{eqnarray}
Obviously, such a smooth function $\zeta_1\in C^\infty(\bR)$ exits.
Property (\ref{zeta-1-concave}) will be used only in Proposition
\ref{boundPsiXfromAbove-Prop}. Now we note that
$\psi_\mxx=\frac{\mxx}{r}\psi_r-\frac{\mxy}{r^2}\psi_\theta$ and
$\psi_\mxy=\frac{\mxy}{r}\psi_r+\frac{\mxx}{r^2}\psi_\theta$,
and define
\begin{eqnarray*}
\hat I_1&\defd & \Big(c_2^2-r^2+ (\gamma-1)r(c_2-r)
\zeta_1(\frac{\mxx\psi_\mxx+\mxy\psi_\mxy}{r(c_2-r)})
-(\gamma-1)(\frac{1}{2}|\grad\psi|^2+\psi)\Big)
\Delta\psi,\\
\hat I_3&\defd & \;\; 2\Big( \frac{\mxx}{r}
(c_2-r)\zeta_1(\frac{\mxx\psi_\mxx+\mxy\psi_\mxy}{r(c_2-r)})
-\frac{\mxy}{r^2}(\mxx\psi_\mxy-\mxy\psi_\mxx)\Big)
 (\mxx\psi_{\mxx\mxx}+\mxy\psi_{\mxx\mxy})\\
&& + 2\Big( \frac{\mxy}{r}
(c_2-r)\zeta_1(\frac{\mxx\psi_\mxx+\mxy\psi_\mxy}{r(c_2-r)})
+\frac{\mxx}{r^2}(\mxx\psi_\mxy-\mxy\psi_\mxx)\Big)(\mxx\psi_{\mxx\mxy}
+\mxy\psi_{\mxy\mxy}).
\end{eqnarray*}
The modified equation in the domain $\DomS_{4\varepsilon}$ is
\begin{equation}
\hat I_1+I_2+\hat I_3+I_4=0.
\label{iteration-equation-sonic-cartesian}
\end{equation}
By (\ref{defZeta-1}), the modified equation
(\ref{iteration-equation-sonic-cartesian}) coincides with the
original equation (\ref{potent-flow-nondiv-psi-1}) if
$$
\left|\frac{\mxx}{r}\psi_\mxx+\frac{\mxy}{r}\psi_\mxy
\right|<\frac{4(c_2-r)}{3(\gamma+1)},
$$
i.e., if $\displaystyle
\left|\psi_x\right|<4x/\big(3(\gamma+1)\big)$ in the
$(x,y)$--coordinates. Also, equation
(\ref{iteration-equation-sonic-cartesian}) is of form
(\ref{iterationEquation}) in the $(\xi,\eta)$--coordinates.

Now we define (\ref{iterationEquation}) in $\DomS_{4\varepsilon}$ by
substituting $\Kphi$ into the coefficients of
(\ref{iteration-equation-sonic-cartesian}) except for the terms
involving $\displaystyle
\zeta_1(\frac{\mxx\psi_\mxx+\mxy\psi_\mxy}{r(c_2-r)})$. Thus, we obtain
an equation of form (\ref{iterationEquation}) with the coefficients:
\begin{equation}\label{iterationSonicDomEquation}
\begin{array}{ll}
A^2_{11}(\grad\psi,\mxx,\mxy)=& c_2^2- (\gamma-1)\Big(r(c_2-r)
\zeta_1(\frac{\mxx\psi_\mxx+\mxy\psi_\mxy}{r(c_2-r)})
 +\frac{1}{2}|\grad\Kphi|^2+\Kphi\Big)
%\nonumber
\\
&-(\Kphi_\mxx^2+\mxx^2)+2\mxx \Big(\frac{\mxx}{r}
(c_2-r)\zeta_1(\frac{\mxx\psi_\mxx+\mxy\psi_\mxy}{r(c_2-r)})-
\frac{\mxy}{r^2}(\mxx\Kphi_\mxy-\mxy\Kphi_\mxx)\Big),
%\nonumber
\\
A^2_{22}(\grad\psi,\mxx,\mxy)=& c_2^2- (\gamma-1)\Big(r(c_2-r)
\zeta_1(\frac{\mxx\psi_\mxx+\mxy\psi_\mxy}{r(c_2-r)})
+\frac{1}{2}|\grad\Kphi|^2+\Kphi\Big)
%\nonumber
\\
&-(\Kphi_\mxy^2+\mxy^2)+2\mxy \Big( \frac{\mxy}{r}
(c_2-r)\zeta_1(\frac{\mxx\psi_\mxx+\mxy\psi_\mxy}{r(c_2-r)})+
\frac{\mxx}{r^2}(\mxx\Kphi_\mxy-\mxy\Kphi_\mxx)\Big),
\\
A^2_{12}(\grad\psi,\mxx,\mxy)=&
-(\Kphi_\mxx\Kphi_\mxy+\mxx\mxy)+2\Big( \frac{\mxx\mxy}{r}
(c_2-r)\zeta_1(\frac{\mxx\psi_\mxx+\mxy\psi_\mxy}{r(c_2-r)})+
\frac{\mxx^2-\mxy^2}{r^2}(\mxx\Kphi_\mxy-\mxy\Kphi_\mxx)\Big),
%\nonumber
\\
A^2_{21}(\grad\psi,\mxx,\mxy)=&A^2_{12}(\grad\psi, \mxx,\mxy),
%\nonumber
\end{array}
\end{equation}
where $\Kphi, \Kphi_\mxx$, and $\Kphi_\mxy$ are
evaluated at $(\mxx,\mxy)$.

Now we write  (\ref{iteration-equation-sonic-cartesian}) in the
$(x,y)$--coordinates. By calculation, the terms $\hat I_1$ and $\hat
I_3$ in the polar coordinates are
\begin{eqnarray*}
&&\hat I_1=
\Big(c_2-r^2+(\gamma-1)\big(r(c_2-r)\zeta_1(\frac{\psi_r}{c_2-r})
-\frac{1}{2}|\grad\psi|^2-\psi\big)\Big)\Delta\psi,\\
&&\hat I_3=2r(c_2-r)\zeta_1(\frac{\psi_r}{c_2-r})\psi_{rr}
  +\frac{2}{r^2}\psi_\theta\psi_{r\theta}
-\frac{2}{r^2}\psi_{\theta}^2.
\end{eqnarray*}
Thus, equation (\ref{iteration-equation-sonic-cartesian}) in the
$(x,y)$--coordinates in $\DomS_{4\varepsilon}$  has the form
\iffalse
\begin{equation}
\left(2x-(\gamma+1)x\zeta_1(\frac{\psi_x}{x})
  +\tilde O_1 \right)\psi_{xx}
%
+\tilde O_2\psi_{xy}
%
+ \left({1\over c_2}+\tilde O_3\right)\psi_{yy}
%
-(1+\tilde O_4)\psi_{x} +\tilde O_5\psi_{y}=0
\label{cutOff-equation-sonicStruct}
\end{equation}
with
\begin{eqnarray}
&&\tilde O_i(p,z,x)= O_i(p,z,x)\qquad\mbox{ for }\; i=2, 5,
\nonumber
\\
&&\tilde O_1(\grad\psi,\psi,x)= -\frac{x^2}{2c_2}+{\gamma+1\over
2c_2} \left(2x^2\zeta_1(\frac{\psi_x}{x})-\psi_x^2\right)
-(\gamma-1)\left(\psi+{1\over 2c_2(c_2-x)^2}\psi_y^2\right),
\nonumber
\\
&&\tilde O_3(\grad\psi,\psi,x)={1\over c_2(c_2-x)^2}\bigg(
x(2c_2-x)-\frac{\gamma+1}{2(c_2-x)^2}\psi_y^2
\label{erTerms-xy-trunc} \\
&& \qquad\qquad\qquad\qquad\qquad\qquad\quad +(\gamma-1)\big(\psi+
(c_2-x)x\zeta_1(\frac{\psi_x}{x})+ \frac{\psi_x^2}{2} \big) \bigg),
\nonumber\\
&&\tilde O_4(\grad\psi,\psi, x)=\frac{1}{c_2-x}\left(x  -
{\gamma-1\over c_2}\big(\psi+(c_2-x)x\zeta_1(\frac{\psi_x}{x})
+\frac{\psi_x^2}{2} +\frac{\psi_y^2}{2 (c_2-x)^2}\big) \right),
\nonumber
\end{eqnarray}
where $O_i(p,z,x), i=2,5,$ are given by (\ref{erTerms-xy-nontrunc}).
It follows that $\tilde O_1(p,z,x),\dots, \tilde O_5(p,z,x)$ satisfy
estimates (\ref{estSmallterms}). In the $(x,y)$--coordinates, this
equation has the form
\fi
%
\begin{equation}
\big( 2x-(\gamma+1)x\zeta_1(\frac{\psi_x}{x})
+O_1^\Kphi\big)\psi_{xx}
+O_2^\Kphi\psi_{xy}
+
\left({1\over c_2}+O_3^\Kphi\right)\psi_{yy}
-(1+O_4^\Kphi)\psi_{x}
+O_5^\Kphi\psi_{y}=0,
\label{iteration-equation-sonicStruct}
\end{equation}
with $\tilde O_k^\Kphi(p,x,y)$ defined by
\begin{equation}
\begin{array}{ll}
&\tilde O_1^\Kphi(p,x, y)= -\frac{x^2}{c_2}+{\gamma+1\over 2c_2}
\big(2x^2\zeta_1(\frac{p_1}{x})-\Kphi_x^2\big) -{\gamma-1\over
c_2}\left(\Kphi+{1\over 2(c_2-x)^2}\Kphi_y^2\right),
\\
&\tilde O_k^\Kphi(x,y)=\tilde O_k(\grad\Kphi(x,y),\Kphi(x,y), x)
\qquad\mbox{ for }\; i= 2, 5,
%\label{5.31a}
\\
&\tilde O_3^\Kphi(p,x,y)={1\over
c_2(c_2-x)^2}\Big(x(2c_2-x)-\frac{\gamma+1}{2(c_2-x)^2}\Kphi_y^2\\
&\qquad\qquad\qquad\qquad\qquad \,\,-(\gamma-1)\big(\Kphi+
(c_2-x)x\zeta_1\left(\frac{p_1}{x}\right)+\frac{1}{2}\Kphi_x^2\big)
\Big),
\\
&\tilde O_4^\Kphi(p,x, y)=\frac{1}{c_2-x}\left(x - {\gamma-1\over
c_2}\big(\Kphi+(c_2-x)x\zeta_1\left(\frac{p_1}{x}\right)
+\frac{\Kphi_x^2}{2} +\frac{\Kphi_y^2}{2 (c_2-x)^2}\big) \right),
\qquad
\end{array}
\label{iteration-equation-sonicStruct-coef}
\end{equation}
where $p=(p_1, p_2)$, and $(\grad\Kphi, \Kphi)$ are evaluated at
$(x,y)$. The estimates in (\ref{estSmallterms}), the definition of
the cutoff function $\zeta_1$, and $\Kphi\in\setK$ with
(\ref{condConst-00}) imply
\begin{equation}
|\tilde O_1^\Kphi(p,x, y)|\le  C|x|^{3/2}, \qquad |\tilde
O_k^\Kphi(x, y)|\le C|x|\qquad \mbox{for }\; k=2,\dots,5,
\label{estSmallterms-iter}
\end{equation}
for all $p\in \bR^2$ and $(x,y)\in \DomS_{4\varepsilon}$. Indeed,
using that $\Kphi\in\setK$ implies
$\|\Kphi\|_{2,\alpha,\DomS}^{(par)} \leq M_1$, we find that, for all
$p\in \bR^2$ and $(x,y)\in \DomS \equiv \DomS_{2\varepsilon}$,
\begin{eqnarray}
&&|\tilde O_1^\Kphi(p,x, y)|\le C (M_1^2+1)|x|^2\le C|x|^{3/2},
\nonumber
\\
&&|\tilde O_k^\Kphi(x, y)|\le C(1+M_1|x|)M_1|x|^{3/2}\le C|x|\qquad
\mbox{for }\; k=2,5, \label{estSmallterms-iter-S}
\\
&&|\tilde O_k^\Kphi(p,x,y)|\le C(|x|+M_1^2|x|^2)\le C|x|
\qquad \mbox{for }\; k=3,4.
\nonumber
\end{eqnarray}
In order to obtain the corresponding estimates in the domain
$\DomS_{4\varepsilon}\setminus \DomS_{2\varepsilon}$, we note that
$\DomS_{4\varepsilon}\setminus \DomS_{2\varepsilon}\subset\DomU$.
Since $2\varepsilon\le x\le 4\varepsilon$ in
$\DomS_{4\varepsilon}\setminus \DomS_{2\varepsilon}$ and
$\Kphi\in\setK$ implies $\|\Kphi\|_{2,\alpha,\DomU}^{(-1-\alpha,
\Sigma_0)} \leq M_2\epsP$, we find that, for any $p\in \bR^2$ and
$(x,y)\in \DomS_{4\varepsilon}\setminus \DomS_{2\varepsilon}$,
\begin{eqnarray}
&&\qquad\quad |\tilde O_1^\Kphi(p,x, y)|\le C
(1+M_2^2\epsP^2+M_2\epsP)\varepsilon^2 \le C\varepsilon^2\le C|x|^2,
\nonumber
\\
&&\qquad\quad |\tilde O_k^\Kphi(x, y)|\le C(1+M_2\epsP)M_2\epsP\le
C\varepsilon^2\le C|x|^2\qquad \mbox{for }\; k=2,5,
\label{estSmallterms-iter-U}
\\
&&\qquad\quad |\tilde O_k^\Kphi(p,x,y)|\le
C(\varepsilon+M_2^2\epsP^2+M_2\epsP)\le C\varepsilon\le C|x| \qquad
\mbox{for }\; k=3,4. \nonumber
\end{eqnarray}
Estimates (\ref{estSmallterms-iter-S})--(\ref{estSmallterms-iter-U})
imply (\ref{estSmallterms-iter}).

The estimates in (\ref{estSmallterms-iter}) imply that,  if
$\Kphi\in\setK$ and $\varepsilon$ is sufficiently small depending
only on the data (which is guaranteed by (\ref{condConst-00}) with
sufficiently large $\hat C$), equation
(\ref{iteration-equation-sonicStruct}) is nonuniformly elliptic in
$\DomS$. First, in the $(x,y)$--coordinates, writing
(\ref{iteration-equation-sonicStruct}) as
$$
a_{11}\psi_{xx}
+2a_{12}\psi_{xy}
+
a_{22}\psi_{yy}
+ a_1\psi_{x} +a_2\psi_{y}=0,
$$
with $a_{ij}=a_{ij}(\grad\psi, x,y)=a_{ji}$ and
$a_{i}=a_{i}(\grad\psi, x,y)$, and using (\ref{condRewritingRH-0}),
we have
$$
\frac{x}{6}|\mu|^2\le \sum_{i,j=1}^2a_{ij}(p,x,y)\mu_i\mu_j \le
\frac{2}{\bar c_2}|\mu|^2 \quad\mbox{for any }
(p,x,y)\in\bR^2\times\DomS_{4\varepsilon} \,\mbox{and}
\,\mu\in\bR^2.
$$
In order to show similar ellipticity in the
$(\mxx,\mxy)$--coordinates, we note that,
 by (\ref{condRewritingRH-0}), the change of coordinates $(\mxx,
\mxy)$ to $(x,y)$ in $\DomS_{4\varepsilon}$ and its inverse have
$C^1$ norms bounded by a constant depending only on the data if
$\varepsilon<\bar c_2/10$. Then there exists $\tilde\lambda>0$
depending only on the data such that, for any
$(p,\mxx,\mxy)\in\bR^2\times\DomS_{4\varepsilon}$ and $\mu\in\bR^2$,
\begin{equation}\label{ellipticityInDegenerateDomain}
\tilde\lambda(c_2-r)|\mu|^2\le \sum_{i,j=1}^2A^2_{ij}(p,\mxx,\mxy)\mu_i\mu_j \le
\tilde{\lambda}^{-1}|\mu|^2,
\end{equation}
where $A^2_{ij}(p,\mxx,\mxy), i,j=1,2,$ are defined by
(\ref{iterationSonicDomEquation}), and  $r=\sqrt{\mxx^2+\mxy^2}$.

Next, we combine the equations introduced above by defining the
coefficients of (\ref{iterationEquation}) in $\Dom$ as follows. Let
$\zeta_2\in C^\infty(\bR)$ satisfy
$$
\zeta_2(s)=\left\{
\begin{array}{ll}
0,\quad&\displaystyle \mbox{if }\;s\le2\varepsilon,\\ \displaystyle
1,\quad&\displaystyle \mbox{if }\;s\ge 4\varepsilon,
\end{array}
\right.
\qquad
\mbox{and }\;\,\,
0\le\zeta_2'(s) \le 10/\varepsilon \quad \mbox{on }\;\bR.
$$
Then we define that, for $p\in\bR^2$ and $(\mxx,\mxy)\in\Dom$,
\begin{equation}\label{combineCoeffs}
A_{ij}(p,\mxx,\mxy)=\zeta_2(c_2-r)A_{ij}^1(\mxx,\mxy)
  +\big(1-\zeta_2(c_2-r)\big)A_{ij}^2(p,\mxx,\mxy).
\end{equation}
Then (\ref{iterationEquation}) is strictly elliptic in $\Dom$ and
uniformly elliptic in  $\DomU$ with ellipticity constant $\lambda>0$
depending only on the data and $\varepsilon$. We state this and
other properties of $A_{ij}$ in the following lemma.

\begin{lemma}\label{propertiesNonlinCoeffs}
There exist  constants $\lambda>0$, $C$, and $\hat C$ depending only
on the data such that, if $M_1,M_2, \varepsilon$, and $\epsP$
satisfy {\rm (\ref{condConst-00})}, then, for any $\Kphi \in\setK$,
the coefficients $A_{ij}(p,\mxx, \mxy)$ defined by {\rm
(\ref{combineCoeffs})}, $i,j=1,2$, satisfy
\begin{enumerate}\renewcommand{\theenumi}{\roman{enumi}}
\item \label{propertiesNonlinCoeffs-i1}
For any $(\mxx, \mxy)\in\Dom$ and $p,\mu\in\bR^2$,
\begin{equation*}
\displaystyle \lambda(c_2-r)|\mu|^2 \le
\sum_{i,j=1}^2A_{ij}(p,\mxx,\mxy)\mu_i\mu_j
\le \lambda^{-1}|\mu|^2 \qquad \mbox{with}\,\, r=\sqrt{\mxx^2+\mxy^2};
\end{equation*}

\item\label{propertiesNonlinCoeffs-i2}
$A_{ij}(p,\mxx,\mxy)=A^1_{ij}(\mxx, \mxy)$ for any $(\mxx, \mxy)\in
\Dom\cap\{c_2-r>4\varepsilon\}$ and $p\in\bR^2$, where
$A^1_{ij}(\mxx, \mxy)$ are defined by {\rm
(\ref{iterationUniforDomEquation})}. Moreover,
$$
A^1_{ij}\in
C^{1,\alpha}(\overline{\Dom\cap\{c_2-r>4\varepsilon\}})
$$ with
$
\|A^1_{ij}\|_{{1,\alpha}
(\overline{\Dom\cap\{c_2-r>4\varepsilon\}})}\le C;
$

\item\label{propertiesNonlinCoeffs-i3}
%The functions $A_{ij}(p,\mxx,\mxy)$ satisfy
$|A_{ij}|+|D_{(p,\mxx,\mxy)}A_{ij}|
%%% , |p\cdot D_p A_{ij}|
\le C$  for any $(\mxx,\mxy)\in \Dom\cap \{0<c_2-r<12\varepsilon\}$
and $p\in\bR^2$.
\end{enumerate}
\end{lemma}

\Proof
Property (\ref{propertiesNonlinCoeffs-i1}) follows from
(\ref{ellipticityInUniformDomain}) and
(\ref{ellipticityInDegenerateDomain})--(\ref{combineCoeffs}).
Properties
(\ref{propertiesNonlinCoeffs-i2})--(\ref{propertiesNonlinCoeffs-i3})
follow from the explicit expressions
(\ref{iterationUniforDomEquation}) and
(\ref{iterationSonicDomEquation}) with $\Kphi\in\setK$. In
estimating these expressions in property
(\ref{propertiesNonlinCoeffs-i3}), we use that $|s\zeta_1'(s)|\le C$
%on $\bR$,
which follows from the smoothness of $\zeta_1$ and
(\ref{defZeta-1}).
\Endproof

Also, equation (\ref{iterationEquation}) coincides with equation
(\ref{iteration-equation-sonicStruct}) in the domain $\DomS$. Assume
that $\varepsilon<\kappa_0/24$, which can be achieved by choosing
$\hat C$ large in (\ref{condConst-00}). Then, in the larger domain
$\Dom\cap\{c_2-r<12\varepsilon\}$, equation
(\ref{iterationEquation}) written in the $(x,y)$--coordinates has
form (\ref{iteration-equation-sonicStruct}) with the only difference
that the term $x\zeta_1(\frac{\psi_x}{x})$ in the coefficient of
$\psi_{xx}$ of (\ref{iteration-equation-sonicStruct}) and in the
terms $\tilde O_1^\Kphi$,  $\tilde O_3^\Kphi$, and $\tilde
O_4^\Kphi$ given by (\ref{iteration-equation-sonicStruct-coef}) is
replaced by
$$
x\Big(\zeta_2(x)\zeta_1(\frac{\Kphi_x}{x}) +
(1-\zeta_2(x))\zeta_1(\frac{\psi_x}{x}) \Big).
$$
{}From this, we have

\begin{lemma}\label{propertiesNonlinCoeffs-xy}
There exist $C$ and $\hat C$ depending only on the data such that
the following holds. Assume that $M_1,M_2, \varepsilon$, and $\epsP$
satisfy {\rm (\ref{condConst-00})}. Let $\Kphi\in\setK$. Then
equation {\rm (\ref{iterationEquation})} written in the
$(x,y)$--coordinates in $\Dom\cap\{c_2-r<12\varepsilon\}$ has the
form
\begin{equation}\label{nonlinIterEq-xy-lg}
\hat A_{11}\psi_{xx}+2\hat A_{12}\psi_{xy}+\hat A_{22}\psi_{yy}
+\hat A_{1}\psi_{x}+\hat A_{2}\psi_{y}=0,
\end{equation}
where $\hat A_{ij}=\hat A_{ij}(\psi_x,x,y)$, $\hat A_{i}=\hat
A_{i}(\psi_x,x,y)$, and  $\hat A_{21}=\hat A_{12}$. Moreover, the
coefficients $\hat A_{ij}(p,x,y)$ and $\hat A_{i}(p,x,y)$ with
$p=(p_1,p_2)\in \bR^2$ satisfy

\begin{enumerate}\renewcommand{\theenumi}{\roman{enumi}}

\item \label{propertiesNonlinCoeffs-xy-i1} For any
$(x,y)\in\Dom\cap\{x<12\varepsilon\}$ and $p,\mu\in\bR^2$,
\begin{equation}
\displaystyle {x\over 6}|\mu|^2\le \sum_{i,j=1}^2\hat
A_{ij}(p,x,y)\mu_i\mu_j
\le \frac{2}{\bar c_2}|\mu|^2;
\end{equation}

\item\label{propertiesNonlinCoeffs-xy-i3}
For any $(x,y)\in \Dom\cap \{x<12\varepsilon\}$ and $p\in\bR^2$,
$$
|(\hat A_{ij}, D_{(p,x,y)}\hat A_{ij})|+|(\hat A_{i},
D_{(p,x,y)}\hat A_{i})| \le C;
$$

\item \label{propertiesNonlinCoeffs-xy-i4-0}
$\hat A_{11}$, $\hat A_{22}$, and $\hat A_1$ are independent of
$p_2$;

\item \label{propertiesNonlinCoeffs-xy-i4}
$\hat A_{12}$, $\hat A_{21}$, and $\hat A_2$ are independent of $p$,
and
$$|(\hat A_{12},\hat A_{21},\hat A_{2})(x,y)|\le
C|x|, \quad |D(\hat A_{12},\hat A_{21},\hat A_{2})(x,y)| \le
C|x|^{1/2}.
$$
\end{enumerate}
\end{lemma}

The last inequality in Lemma
\ref{propertiesNonlinCoeffs-xy}(\ref{propertiesNonlinCoeffs-xy-i4})
is proved as follows. Note that
$$
(\hat A_{12}, \hat A_{2})(x,y)=(O_2, O_5)(D\Kphi(x,y), \Kphi(x,y),
x),
$$
where $O_2$ and $O_5$ are given by (\ref{erTerms-xy-nontrunc}).
Then, by $\Kphi\in\setK$ and (\ref{condConst-00}), we find that, for
$(x,y)\in \DomS$, i.e., $x\in (0, 2\varepsilon)$,
\begin{eqnarray*}
|D(\hat A_{12},\hat A_{21},\hat A_{2})(x,y)|
&\le& C(1+M_1\varepsilon)|D\Kphi_y(x,y)|
+(1+M_1)|\Kphi_y(x,y)|\\
&\le& C(1+M_1\varepsilon)M_1x^{1/2}+C(1+M_1)M_1x^{3/2}\le Cx^{1/2};
\end{eqnarray*}
and, for $(x,y)\in \Dom\cap\{\varepsilon\le x\le
12\varepsilon\}\subset\DomU$, we have $\dist(x,\Sigma_0)\ge
c_2/2\ge\bar c_2/4$ so that
$$
|D(\hat A_{12},\hat A_{21},\hat A_{2})(x,y)|
\le C(1+M_2\epsP)M_2\epsP
\le C\varepsilon\le C x.
$$

The next lemma follows directly from both (\ref{defZeta-1}) and the
definition of $A_{ij}$.

\begin{lemma}\label{cutOffEqIsOriginalEq}
Let $\Omega\subset\Dom$, $\psi\in C^2(\Omega)$, and $\psi$ satisfy
equation {\rm (\ref{iterationEquation})} with $\Kphi=\psi$ in
$\Omega$. Assume also that $\psi$, written in the
$(x,y)$--coordinates, satisfies $|\psi_x|\le
4x/\big(3(\gamma+1)\big)$ in
$\Omega'\defd\Omega\cap\{c_2-r<4\varepsilon\}$. Then $\psi$
satisfies \eqref{potent-flow-nondiv-psi-1} in $\Omega$.
\end{lemma}

\Subsec{\bf\large  The iteration procedure and choice of the
constants} \label{overViewProcedureSubsection} With the previous
analysis, our iteration procedure will consist of the following ten
steps, in which Steps 2--9 will be carried out in detail in
\S\ref{unifElliptApproxSection}--\S\ref{removeCutoffSection} and the
main theorem is completed in \S\ref{proofSection}.

\medskip
{\em Step 1.} Fix $\Kphi\in\setK$. This determines the domain
$\Omega^+(\Kphi)$, equation (\ref{iterationEquation}), and condition
(\ref{iterationRH}) on $\shock(\Kphi)$, as described in \S
\ref{Constr-iter-section}--\S\ref{eqForIterationSection} above.

\medskip
{\em Step 2.} In \S\ref{unifElliptApproxSection}, using the
vanishing viscosity approximation of equation
(\ref{iterationEquation}) via a uniformly elliptic equation
$$
\Nl(\psi)+\delta\Delta\psi=0\,\,\qquad\mbox{for }\;\delta\in(0, 1)
$$
and sending $\delta\to 0$, we establish the existence of a solution
$\psi\in C(\overline{\Omega^+(\Kphi)})\cap
C^1(\overline{\Omega^+(\Kphi)}\setminus\overline\sonic) \cap
C^2(\Omega^+(\Kphi))$
%$\psi\in C^2(\Omega^+(\Kphi))\cap
%C^{1}(\overline{\Omega^+(\Kphi)})$
to problem
(\ref{iterationEquation})--(\ref{iterationCondOnSymmtryLine}). This
solution satisfies
\begin{equation}
0\le\psi\le C\epsP \qquad \mbox{in }\;\Omega^+(\Kphi),
\label{L-ifty-BdIteratioOverview}
\end{equation}
where $C$ depends only on the data.

\medskip
{\em Step 3.} For every $s\in(0, c_2/2)$, set $\ElDomU_s:=
\Omega^+(\Kphi)\cap\{c_2-r>s\}$. By Lemma
\ref{propertiesNonlinCoeffs}, if {\rm (\ref{condConst-00})} holds
with sufficiently large $\hat C$ depending only on the data, then
equation (\ref{iterationEquation}) is uniformly elliptic in
$\ElDomU_s$ for every $s\in (0, c_2/2)$, the ellipticity constant
depends only on the data and s, and the bounds of coefficients in
the corresponding H\"{o}lder norms also depend only on the data and
$s$. Furthermore, (\ref{iterationEquation}) is linear on
$\{c_2-r>4\varepsilon\}$, which implies that it is also linear near
the corners $\PtLwL$ and $\PtLwR$. Then, by the standard elliptic
estimates in the interior and near the smooth parts of
$\partial\Omega^+(\Kphi)\cap \overline {\ElDomU_s}$ and using
Lieberman's estimates \cite{Lieberman88} for linear equations with
the oblique derivative conditions near the corners $(-u_2,-v_2)$ and
$\shock(\Kphi)\cap\{\eta=-v_2\}$, we have
\begin{equation}\label{ellipticEstimates-unif-halfDomain}
\|\psi\|^{(-1-\alpha,\Sigma_0)}_{2,\alpha,\ElDomU_{s/2}}
\le C(s)(\|\psi\|_{L^\infty(\overline {\ElDomU_s})}
  +|v_2|),
\end{equation}
if $\|\psi\|_{L^\infty(\overline {\ElDomU_s})}+
 |v_2|\le 1$,
where the second term in the right-hand side comes from the boundary
condition (\ref{iterationCondOnSymmtryLine}), and the constant
$C(s)$ depends only on the ellipticity constants, the angles at the
corners
 $\PtLwL=\shock(\Kphi)\cap\{\eta=-v_2\}$ and $\PtLwR=(-u_2,-v_2)$, the
norm of $\shock(\Kphi)$ in $C^{1,\alpha}$, and $s$, which implies
that $C(s)$ depends only on the data and $s$.

 Now, using (\ref{L-ifty-BdIteratioOverview}) and
(\ref{u2-v2-bound}), we obtain
$\|\psi\|_{L^\infty(\overline{\ElDomU_s})}+|v_2|\le 1$ if $\epsP$ is
sufficiently small, which is achieved by choosing $\hat C$ in
(\ref{condConst-00}) sufficiently large. Then, from
(\ref{ellipticEstimates-unif-halfDomain}),
we obtain
\begin{equation}\label{ellipticEstimates-unif-halfDomain-2}
\|\psi\|^{(-1-\alpha,\Sigma_0)}_{2,\alpha,\ElDomU_{s/2}}\le
C(s)\epsP
\end{equation}
for every $s\in (0, c_2/2)$, where $C$ depends only on the data and $s$.

\medskip
{\em Step 4. Estimates of $\psi$ in
$\hElDomS(\Kphi)\defd\Omega^+(\Kphi)\cap \{c_2-r<\varepsilon\}$.} We
work in the $(x, y)$--coordinates, and then equation
(\ref{iterationEquation}) is equation
(\ref{iteration-equation-sonicStruct}) in $\ElDomS$.

\medskip
{\em Step 4.1. $L^\infty$ estimates of $\psi$ in
$\Omega^+(\Kphi)\cap \DomS$.} Since $\Kphi\in\setK$, the estimates
in (\ref{estSmallterms-iter}) hold for large $\hat C$ in
(\ref{condConst-00}) depending only on the data. We also rewrite the
boundary condition (\ref{iterationRH}) in the $(x, y)$--coordinates
and obtain (\ref{RH-psi-3-xy}) with $\hat E_i$ replaced by $\hat
E_i^\Kphi(x,y)\defd \hat E_i(D\Kphi(x,y),\Kphi(x,y),x,y)$. Using
$\Kphi\in\setK$,  (\ref{RH-psi-2-error-term-xy-1}),
(\ref{RH-psi-2-error-term-xy-2}), and (\ref{holder-hat-f-S}) with
$\hat f_\Kphi(0)=\hat f_0(0)=y_1$, we obtain
\begin{equation}
|\hat E_i^\Kphi(x,y)|\le C(M_1\varepsilon+M_2\epsP)\le C/\hat C,\qquad
i=1,2,
\label{RH-psi-2-error-term-iter}
\end{equation}
for $(x,y)\in \shock(\Kphi)\cap\{0<x<2\varepsilon\}$. Then, if $\hat
C$ in (\ref{condConst-00}) is large, we find that the function
$$
w(x,y)=\frac{3x^2}{5(\gamma+1)}
%w(x,y)=\frac{1+\frac{1}{50}}{2(\gamma+1)}x^2
$$
is a supersolution of equation
(\ref{iteration-equation-sonicStruct}) in $\ElDomS(\Kphi)$ with the
boundary condition (\ref{iterationRH}) on $\shock(\Kphi)\cap
\{0<x<2\varepsilon\}$. That is, the right-hand sides of
(\ref{iterationRH}) and (\ref{iteration-equation-sonicStruct}) are
negative on $w(x,y)$ in the domains given above. Also, $w(x,y)$
satisfies the boundary conditions
(\ref{iterationCondOnSonicLine})--(\ref{iterationCondOnWedge})
within $\ElDomS(\Kphi)$. Thus,
\begin{equation}
0\le\psi(x,y)\le \frac{3x^2}{5(\gamma+1)}
%0\le\psi(x,y)\le \frac{1+\frac{1}{50}}{2(\gamma+1)}x^2
\qquad\mbox{in }\;\ElDomS(\Kphi),
\label{L-ifty-BdIteratin-Sonic-Overview}
\end{equation}
if $w\ge \psi$ on $x=\varepsilon$. By (\ref{L-ifty-BdIteratioOverview}),
$w\ge \psi$ on $x=\varepsilon$ if
$$
C\epsP\le \varepsilon^2,
$$
where $C$ is a large constant depending only on the data, i.e., if
(\ref{condConst-00}) is satisfied with large $\hat C$. The details
of the argument of  Step 4.1 are in Lemma
\ref{quadraticGrowthPsi-Lemma}.

\medskip
{\em Step 4.2. Estimates of the norm
$\|\psi\|_{2,\alpha,\hElDomS(\Kphi)}^{(par)}$}. We use the parabolic
rescaling in the rectangle $R_{z}$ defined by
(\ref{parabRectangles}) in which $\ElDomS$ is replaced by
$\ElDomS(\Kphi)$. Note that $R_{z}\subset \ElDomS$ for every
$z=(x,y)\in\hElDomS(\Kphi)$. Thus, $\psi$ satisfies
(\ref{iteration-equation-sonicStruct}) in $R_{z}$. For every
$z\in\hElDomS(\Kphi)$, we define the functions $\psi^{(z)}$ and
$\Kphi^{(z)}$ by (\ref{parabRescaling}) in the domain $Q_1^{(z)}$
defined by (\ref{rescaled-parabRectangles}). Then equation
(\ref{iteration-equation-sonicStruct}) for $\psi$ yields the
following equation for $\psi^{(z)}(S,T)$ in  $Q_1^{(z)}$:
\begin{eqnarray}
&&\quad\big((1+\frac{S}{4})\big(2-
(\gamma+1)\zeta_1(\frac{4\psi^{(z)}_S}{1+S/4})\big)
+xO_1^{(\Kphi,z)} \big)\psi^{(z)}_{SS}
+xO_2^{(\Kphi,z)}\psi^{(z)}_{ST}
\label{iteration-equation-sonicStruct-ParabRescaled}\\
&&\qquad\quad
+ \big({1\over c_2}+x O_3^{(\Kphi,z)}\big)\psi^{(z)}_{TT}
-({1\over 4}+xO_4^{(\Kphi,z)})\psi^{(z)}_{S}
+x^2O_5^{(\Kphi,z)}\psi^{(z)}_{T}=0,
 \nonumber
\end{eqnarray}
where the terms $O_k^{(\Kphi,z)}(S,T,p)$, $k=1,\dots, 5$, satisfy
\begin{equation}\label{rescaledErrorTermsEqEstimate}
\|O_k^{(\Kphi,z)}\|_{C^{1,\alpha}(\overline{Q_1^{(z)}}\times\bR^2)}\le C(1+M_1^2).
\end{equation}
Estimate (\ref{rescaledErrorTermsEqEstimate}) follows from the
explicit expressions of $O_k^{(\Kphi,z)}$ obtained from both
(\ref{iteration-equation-sonicStruct-coef}) by rescaling and the
fact that
$$
\|\Kphi^{(z)}\|_{C^{2,\alpha}(\overline{Q_1^{(z)}})}\le CM_1,
$$
which is true since $\|\Kphi\|_{2,\alpha,\ElDomS(\Kphi)}^{(par)}\le
M_1$. Now, since every term $O_k^{(\Kphi,z)}$ in
(\ref{iteration-equation-sonicStruct-ParabRescaled}) is multiplied
by $x^{\beta_k}$ with $\beta_k\ge 1$ and $x\in (0, \varepsilon)$,
condition (\ref{condConst-00}) (possibly after increasing $\hat C$
depending only on the data) implies that equation
(\ref{iteration-equation-sonicStruct-ParabRescaled}) is uniformly
elliptic in $Q_1^{(z)}$ and has the $C^{1,\alpha}$ bounds on the
coefficients by a constant depending only on the data.

Now, if the rectangle $R_{z}$ does not intersect $\partial
\Omega^+(\Kphi)$, then $Q_1^{(z)}=Q_1$, where $Q_s=(-s, s)^2$ for
$s>0$. Thus, the interior elliptic estimates in Theorem
\ref{locEstElliptEq} in Appendix imply
\begin{equation}\label{ellipticEstimates-interior}
\|\psi^{(z)}\|_{C^{2,\alpha}(\overline{Q_{1/2}})}\le C,
\end{equation}
where $C$ depends only on the data
and $\|\psi^{(z)}\|_{L^\infty(\overline{Q_{1}})}$.
{}From (\ref{L-ifty-BdIteratin-Sonic-Overview}),
we have
$$
\|\psi^{(z)}\|_{L^\infty(\overline{Q_{1}})} \le 1/(\gamma+1).
%\le \frac{1+\frac{1}{50}}{2(\gamma+1)}.
$$
Therefore, we obtain (\ref{ellipticEstimates-interior}) with $C$
depending only on the data.

Now consider the case when the rectangle $R_{z}$ intersects
$\partial \Omega^+(\Kphi)$. {}From its definition, $R_{z}$ does not
intersect $\sonic$. Thus, $R_{z}$ intersects either $\shock$ or the
wedge boundary $\wedgeB$. On these boundaries, we have the
homogeneous oblique derivative conditions (\ref{iterationRH}) and
(\ref{iterationCondOnWedge}).
 In the case when $R_{z}$
intersects  $\wedgeB$, the rescaled condition
(\ref{iterationCondOnWedge}) remains the same form, thus oblique,
and we use the estimates for the oblique derivative problem in
Theorem \ref{locEstElliptEq-Dirichlet} to obtain
\begin{equation}\label{ellipticEstimates-bdry-rescaled}
\|\psi^{(z)}\|_{C^{2,\alpha}(\overline{Q^{(z)}_{1/2}})}\le C,
\end{equation}
where $C$ depends only on the data, since the $L^\infty$ bound of
$\psi^{(y)}$ in $Q^{(z)}_{1}$ follows from
(\ref{L-ifty-BdIteratin-Sonic-Overview}).
 In the case
when $R_{z}$ intersects $\shock$, the obliqueness in the rescaled
condition  (\ref{iterationRH}) is of order $x^{1/2}$, which is small
since $x\in (0, 2\varepsilon)$. Thus we use the estimates for the
``almost tangential derivative" problem in Theorem
\ref{locEstElliptEq-non-oblique} to obtain
(\ref{ellipticEstimates-bdry-rescaled}).

Finally, rescaling back, we have
\begin{equation}\label{ellipticEstimates-par-halfDomain}
 \|\psi\|_{2,\alpha,\hElDomS(\Kphi)}^{(par)}\le C.
\end{equation}
The details of the argument of Step 4.2 are in Lemma
\ref{EstParabolicHolder-Lemma}.

\medskip
{\em Step 5.} In Lemma \ref{extension-Lemma}, we extend $\psi$ from
the domain $\Omega^+(\Kphi)$ to $\Dom$ working in the $(x,
y)$--coordinates (or, equivalently in the polar coordinates) near
the sonic line and in the rest of the domain in the $(\mxx,
\mxy)$--coordinates, by using the procedure of \cite{ChenFeldman1}.
If $\hat C$ is sufficiently large, the  extension of $\psi$
satisfies
\begin{eqnarray}\label{ellipticEstimates-par-Domain}
&& \|\psi\|_{2,\alpha,\DomS}^{(par)}\le C,
 \\
 \label{ellipticEstimates-unif-Domain}
&&\|\psi\|^{(-1-\alpha,\Sigma_0)}_{2,\alpha,\DomU}
\le C(\varepsilon)\epsP,
\end{eqnarray}
%in the domain $\hat\DomS=\DomS\cap\{c_2-r<\varepsilon\}$
with $C$ depending only on the data in
(\ref{ellipticEstimates-par-Domain}) and $C(\varepsilon)$ depending
only on the data and $\varepsilon$ in
(\ref{ellipticEstimates-unif-Domain}). This is obtained by using
(\ref{ellipticEstimates-par-halfDomain}) and
(\ref{ellipticEstimates-unif-halfDomain-2}) with $s>0$ determined by
the data and $\varepsilon$, and by using the estimates of the
functions $f_\Kphi$ and $\hat f_\Kphi$ in (\ref{OmegaPL}),
(\ref{holder-hat-f}), and (\ref{holder-hat-f-S}).

\medskip
{\em Step 6.} We fix $\hat C$ in (\ref{condConst-00}) large
depending only on the data, so that Lemmas
\ref{propertiesNonlinCoeffs}--\ref{propertiesNonlinCoeffs-xy} hold
and the requirements on  $\hat C$ stated in Steps 1--5 above are
satisfied. Set $M_1=\max(2C,1)$ for the constant $C$ in
(\ref{ellipticEstimates-par-Domain}) and choose
\begin{equation}
\varepsilon =\frac{1}{10\max((\hat CM_1)^4, \hat C)}.
\label{condConst-1-1}
\end{equation}
This choice of $\varepsilon$ fixes $C$ in (\ref{ellipticEstimates-unif-Domain})
depending only on the data and $\hat C$.
Now set $M_2=\max(C,1)$ for $C$ from (\ref{ellipticEstimates-unif-Domain}) and let
$$
0<\epsP\le \epsP_0:={({\hat
C}^{-1}-\varepsilon-\varepsilon^{1/4}M_1)\varepsilon^2 \over
2\left(\varepsilon^2\max(M_1, M_2)+M_2\right)},
$$
where $\epsP_0>0$ since $\varepsilon$ is defined by
(\ref{condConst-1-1}). Then (\ref{condConst-00}) holds with constant
$\hat C$  fixed above.

Note that the constants $\epsP_0, \varepsilon, M_1$, and $M_2$
depend only on the data and $\hat C$.

\medskip
{\em Step 7.} With the constants $\epsP, \varepsilon, M_1$, and
$M_2$ chosen in Step 6, estimates
(\ref{ellipticEstimates-par-Domain})--(\ref{ellipticEstimates-unif-Domain})
imply
$$
\|\Kphi\|_{2,\alpha,\DomS}^{(par)} \leq M_1, \qquad
\|\psi\|^{(-1-\alpha,\Sigma_0)}_{2,\alpha,\DomU}
\le M_2\epsP.
$$
Thus, $\psi\in\setK(\epsP, \varepsilon, M_1, M_2)$.
Then the iteration map $J:\setK\to\setK$ is
defined.

\medskip
{\em Step 8.} In Lemma \ref{extension-Lemma} and Proposition
\ref{existenceFixedPt}, by the argument similar to
\cite{ChenFeldman1} and the fact that $\setK$ is a compact and
convex subset of $C^{1,\alpha/2}(\overline\Dom)$, we  show that the
iteration map $J$ is continuous, by uniqueness of the solution
$\psi\in C^{1,\alpha}(\overline\Dom)\cap C^2(\Dom)$ of
(\ref{iterationEquation})--(\ref{iterationCondOnSymmtryLine}). Then,
by the Schauder Fixed Point Theorem, there exists a fixed point
$\psi\in\setK$. This is a solution of the free boundary problem.

\medskip
{\em Step 9.} Removal of the cutoff: By Lemma
\ref{cutOffEqIsOriginalEq}, a fixed point $\psi=\Kphi$ satisfies the
original equation \eqref{potent-flow-nondiv-psi-1} in
$\Omega^+(\psi)$ if $|\psi_x|\le 4x/\big(3(\gamma+1)\big)$ in
$\Omega^+(\psi)\cap \{c_2-r<4\varepsilon\}$. We prove this estimate
in \S\ref{removeCutoffSection} by choosing $\hat C$ sufficiently
large depending only on the data.

\medskip
{\em Step 10.} Since the fixed point $\psi\in \setK$
of the iteration map $J$ is a solution
of (\ref{iterationEquation})--(\ref{iterationCondOnSymmtryLine})
for $\Kphi=\psi$, we conclude

\begin{enumerate}\renewcommand{\theenumi}{\roman{enumi}}
\item \label{iterProcItem-first}
$\psi\in C^{1,\alpha}(\overline{\Omega^+(\psi)})\cap
C^{2,\alpha}(\Omega^+(\psi))$;

\item  $\psi=0$ on $\sonic$ by (\ref{iterationCondOnSonicLine}),
and $\psi$ satisfies  the original equation
\eqref{potent-flow-nondiv-psi-1} in $\Omega^+(\psi)$ by Step 9;

\item  $D\psi=0$ on $\sonic$
    since $\|\Kphi\|_{2,\alpha,\DomS}^{(par)}\leq M_1$;

\item \label{iterProcItem}
$\psi=\varphi_1-\varphi_2$ on $\shock(\psi)$
  by (\ref{shockPL})--(\ref{shockIterDef}) since $\Kphi=\psi$;

\item \label{iterProcItem-last}
The Rankine-Hugoniot  gradient jump condition (\ref{RH-psi}) holds
on $\shock(\psi)$. Indeed, as we have showed in (\ref{iterProcItem})
above, the function $\varphi=\psi+\varphi_2$ satisfies
(\ref{cont-accross-shock-mod-phi}) on $\shock(\psi)$.
%From this,
Since $\psi\in\setK$, it follows that $\psi$ satisfies
(\ref{cont-accross-shock-psi-resolved}). Also, $\psi$ on
$\shock(\psi)$ satisfies (\ref{iterationRH}) with $\Kphi=\psi$,
which is (\ref{RH-psi-2}). Since $\psi\in\setK$ satisfies
(\ref{cont-accross-shock-psi-resolved}) and  (\ref{RH-psi-2}), it
has been shown in \S\ref{equationForPsiSection} that $\varphi$
satisfies (\ref{RH-mod-phi}) on $\shock(\psi)$, i.e., $\psi$
satisfies  (\ref{RH-psi}).
\end{enumerate}

Extend the function $\varphi=\psi+\varphi_2$ from
$\Omega:=\Omega^+(\psi)$ to the whole domain $\Lambda$ by using
(\ref{phi-states-0-1-2}) to define $\varphi$ in
$\Lambda\setminus\Omega$. Denote $\Lambda_0:=\{\xi>\xi_0\}\cap
\Lambda$, $\Lambda_1$ the domain with $\xi<\xi_0$ and above the
reflection shock $P_0\PtUpL\PtLwL$, and
$\Lambda_2:=\Lambda\setminus(\overline \Lambda_0\cup \overline
\Lambda_1)$. Set $S_0:=\{\xi=\xi_0\}\cap\Lambda$ the incident shock
and $S_1:=P_0\PtUpL\PtLwL\cap \Lambda$ the reflected shock. We show
in \S\ref{proofSection} that $S_1$ is a $C^2$--curve. Then we
conclude that the domains $\Lambda_0$,  $\Lambda_1$, and $\Lambda_2$
are disjoint, $\partial \Lambda_0\cap \Lambda=S_0$,
 $\partial \Lambda_1\cap \Lambda=S_0\cup S_1$, and $\partial \Lambda_2\cap\Lambda=S_1$.
Properties (\ref{iterProcItem-first})--(\ref{iterProcItem-last})
above and the fact that $\psi$ satisfies
\eqref{potent-flow-nondiv-psi-1}
 in $\Omega$ imply that
$$
\varphi\in W^{1,\infty}_{loc}(\Lambda), \quad \varphi\in
C^1(\overline{\Lambda_i}) \cap C^{1,1}(\Lambda_i)\quad\mbox{for
}i=0,1,2,
$$
$\varphi$ satisfies  equation \eqref{1.1.5} a.e. in $\Lambda$ and
the Rankine-Hugoniot condition (\ref{FBConditionSelfSim-0}) on the
$C^2$-curves $S_0$ and $S_1$, which intersect only at
$P_0\in\partial \Lambda$ and are transversal at the intersection
point. Using this, Definition \ref{def2.1}, and the remarks after
Definition \ref{def2.1}, we conclude that $\varphi$ is a weak
solution of Problem 2, thus of Problem 1. Note that the solution is
obtained for every $\epsP\in(0,\epsP_0]$, i.e., for every
$\theta_w\in [\pi/2-\epsP_0, \pi/2]$ by (\ref{angleCloseToPiOver2}),
and that $\epsP_0$ depends only on the data since $\hat C$ is fixed
in Step 9.

\section{Vanishing Viscosity Approximation and Existence of
Solutions of\\ Problem
(\ref{iterationEquation})--(\ref{iterationCondOnSymmtryLine})}
\label{unifElliptApproxSection}

\numberwithin{equation}{section}

In this section we perform Step 2 of the iteration procedure
described in \S \ref{overViewProcedureSubsection}. Through this
section, we keep $\Kphi\in\setK$ fixed, denote by $\corners:=\{P_1,
P_2, P_3, P_4\}$ the set of the corner points of $\Omega^+(\Kphi)$,
and use $\alpha\in (0, 1/2)$ defined in \S
\ref{Constr-iter-section}.

We regularize equation (\ref{iterationEquation}) by the vanishing
viscosity approximation via the uniformly elliptic equations
$$
\Nl(\psi)+\delta\Delta\psi=0\qquad\mbox{for }\;\delta\in(0, 1).
$$
That is, we consider the equation
\begin{equation}
\Nl_\delta(\psi)\defd (A_{11}+\delta)\psi_{\mxx\mxx}+ 2A_{12}\psi_{\mxx\mxy}
+(A_{22}+\delta)\psi_{\mxy\mxy}=0 \qquad \mbox{ in
}\;\;\Omega^+(\Kphi). \label{unif-ellipt-iterationEquation}
\end{equation}
In the domain $\ElDomS$ in the $(x,y)$--coordinates defined by
(\ref{coordNearSonic}), this equation has the form
\begin{eqnarray}\label{unif-ellipt-iteration-equation-sonicStruct}
&&\qquad\big(
\delta+2x-(\gamma+1)x\zeta_1\big(\frac{\psi_x}{x}\big)+O_1^\Kphi
\big)\psi_{xx}
+O_2^\Kphi\, \psi_{xy}\\
&&\qquad\,\,+ \big({1\over c_2}+{\delta\over (c_2-x)^2}+O_3^\Kphi
\big)\psi_{yy}
-(1-{\delta\over c_2-x}+O_4^\Kphi)\psi_{x} +O_5^\Kphi\,\psi_{y}=0
\nonumber
\end{eqnarray}
by using (\ref{iteration-equation-sonicStruct}) and writing the
Laplacian operator $\Delta$ in the $(x,y)$--coordinates, which is
easily derived from the form of $\Delta$ in the polar coordinates.
The terms $O_k^\Kphi$ in
(\ref{unif-ellipt-iteration-equation-sonicStruct}) are defined by
(\ref{iteration-equation-sonicStruct-coef}).

We now study equation (\ref{unif-ellipt-iterationEquation}) in
$\Omega^+(\Kphi)$ with the boundary conditions
(\ref{iterationRH})--(\ref{iterationCondOnSymmtryLine}).

We first note some properties of the boundary condition
(\ref{iterationRH}). Using Lemma \ref{relatingNorms} with
$\alpha\in(0, 1/2)$ and  (\ref{condConst-00}), we find
$\|\Kphi\|_{2,\alpha,\Dom}^{(-1-\alpha, \Sigma_0\cup\sonic)} \le C$,
where $C$ depends only on the data. Then, writing
(\ref{iterationRH}) as
\begin{equation}
\Ml(\psi)(\mxx,\mxy):= b_1(\mxx,\mxy)\psi_\mxx +
b_2(\mxx,\mxy)\psi_\mxy + b_3(\mxx,\mxy)\psi=0 \qquad\mbox{on
}\;\shock(\Kphi) \label{iterationRH-lf}
\end{equation}
and using (\ref{RH-psi-2-error-term1})--(\ref{obliquenessRH}), we
obtain
\begin{equation}\label{estCoefsIterRH-0}
\|b_i\|^{(-\alpha,\{\PtUpL,\PtLwL\})}_{1,\alpha, \shock(\Kphi)}\le C
\qquad \mbox{for }\;i=1,2,3,
\end{equation}
where $C$ depends only on the data.

Furthermore, $\Kphi\in\setK$ with (\ref{condConst-00}) implies that
$$
\|\Kphi\|_{C^1}\le M_1\varepsilon+M_2\epsP\le \varepsilon^{3/4}/\hat
C.
$$
Then, using (\ref{RH-psi-2-error-term1})--(\ref{obliquenessRH})
and assuming that $\hat C$ in (\ref{condConst-00}) is sufficiently
large, we have
\begin{equation}\label{estCoefsIterRH}
\begin{array}{l}
(b_1(\mxx, \mxy), b_2(\mxx, \mxy))\cdot\nu(\mxx, \mxy)\ge
\frac{1}{4}\rho_2'(c_2^2-\hat{\mxx}^2)>0\qquad\mbox{for any }\; (\mxx,
\mxy)\in\shock(\Kphi), \displaystyle
\\
b_1(\mxx, \mxy)\ge
\frac{1}{2}\rho_2'(c_2^2-\hat{\mxx}^2)>0\qquad\mbox{for any }\; (\mxx,
\mxy)\in\shock(\Kphi), \displaystyle
\\
\left|b_2(\mxx, \mxy)-\mxy\big(\frac{\rho_2-\rho_1}{u_1}-\rho_2'\hat\mxx
\big)\right|\le \varepsilon^{3/4}\qquad\mbox{for any }\; (\mxx,
\mxy)\in\shock(\Kphi), \displaystyle
\\
\left|b_3(\mxx, \mxy)+\big(\frac{\rho_2-\rho_1}{u_1}-\rho_2'\hat\mxx
\big)\right|\le\varepsilon^{3/4}\qquad\mbox{for any }\; (\mxx,
\mxy)\in\shock(\Kphi). \displaystyle
\end{array}
\end{equation}

Now we write condition (\ref{iterationRH}) in the
$(x,y)$--coordinates on $\shock(\Kphi)\cap\overline\DomS$. Then we
obtain the following condition of the form
\begin{equation}
\Ml(\psi)(x,y)=\hat b_1(x,y)\psi_x + \hat b_2(x,y)\psi_y + \hat
b_3(x,y)\psi=0 \qquad\mbox{on }\;\shock(\Kphi)\cap\overline\DomS,
\label{iterationRH-lf-flattened}
\end{equation}
where $\hat b_1(x,y)=b_1(\xi,\eta){\partial x\over\partial\xi}
  +b_2(\xi,\eta){\partial x\over\partial\eta},
\hat b_2(x,y)=b_1(\xi,\eta){\partial y\over\partial\xi}
  +b_2(\xi,\eta){\partial y\over\partial\eta},$
and $\hat b_3(x,y)=b_3(\xi,\eta)$.
Condition (\ref{iterationRH}) is oblique, by the first inequality in
(\ref{estCoefsIterRH}). Then, since transformation
(\ref{coordNearSonic}) is smooth on $\{0<c_2-r<2\varepsilon\}$ and
has nonzero Jacobian, it follows that
(\ref{iterationRH-lf-flattened}) is oblique, that is,
%we have
\begin{equation}\label{obliqInxy-1}
(\hat b_1(x,y),\hat b_2(x,y))\cdot\nu_s(x,y)\ge C^{-1}>0 \qquad
\mbox{ on }\; \shock(\Kphi)\cap\overline\DomS,
\end{equation}
where $\hat\nu_s=\hat\nu_s(x,y)$ is the interior unit normal at
$(x,y)\in\shock(\Kphi)\cap \DomS$ to $\ElDom(\Kphi)$.

As we have showed in \S \ref{sectionEqNearSonicLine}, writing the
left-hand side of (\ref{RH-psi-2}) in the $(x,y)$--coordinates, we
obtain the left-hand side of (\ref{RH-psi-3-xy}). Thus,
(\ref{iterationRH-lf-flattened}) is obtained from
(\ref{RH-psi-3-xy}) by substituting $\Kphi(x,y)$ into $\hat E_1$ and
$\hat E_2$. Also, from (\ref{holder-hat-f-S}) with $\hat
f_\Kphi(0)=\hat f_0(0)=y_1$, we estimate $|y-y_1|=|\hat
f_\Kphi(x)-\hat f_\Kphi(0)|\le CM_1\varepsilon$ on
$\shock\cap\{x<2\varepsilon\}$. Then, using
(\ref{RH-psi-3-xy})--(\ref{RH-psi-2-error-term-xy-2}) and
$\mxx_1<0$, we find that, if $\hat C$ in (\ref{condConst-00}) is
sufficiently large depending only on the data, then
\begin{equation}\label{estCoefsIterRH-flattened}
\begin{array}{l}
\|\hat b_i\|^{(-1,\{\PtUpL\})}_{1,\alpha, \shock(\Kphi)\cap\overline\DomS}\le
CM_1 \qquad \mbox{for }\;i=1,2,3, \displaystyle
\\
\hat b_1(x,y)\le -\frac{1}{2}\frac{\rho_2-\rho_1}{u_1}
\frac{\mxy^2_1}{c_2}<0\qquad\mbox{for }\;
(x,y)\in\shock(\Kphi)\cap\overline\DomS, \displaystyle
\\
\hat b_2(x,y)\le -\frac{1}{2}\mxy_1 \big(\rho_2'
+\frac{\rho_2-\rho_1}{u_1 c_2^2}|\mxx_1| \big)<0\qquad\mbox{for }\;
(x,y)\in\shock(\Kphi)\cap\overline\DomS, \displaystyle
\\
\hat
b_3(x,y)\le-\frac{1}{2}\big(\rho_2'|\mxx_1|+\frac{\rho_2-\rho_1}{u_1}
\big)<0\qquad\mbox{for }\; (x,y)\in\shock(\Kphi)\cap\overline\DomS,
\displaystyle
\end{array}
\end{equation}
where $C$ depends only on the data.

\medskip
Now we state the main existence result for the regularized problem.
%*****************BEGIN PROPOSITION*************
\begin{proposition}\label{existSolUnifEllipt}
There exist $\hat C, C, \delta_0>0$ depending only
on the data such that, if $\epsP, \varepsilon>0$ and $M_1, M_2\ge 1$
in {\rm (\ref{defSetK_R})} satisfy {\rm (\ref{condConst-00})}, then,
for every $\delta\in(0,\delta_0)$, there exists a unique solution
$\psi\in C^{(-1-\alpha, \corners)}_{2,\alpha,\Omega^+(\Kphi)}$ of
{\rm (\ref{unif-ellipt-iterationEquation})} and {\rm
(\ref{iterationRH})}--{\rm (\ref{iterationCondOnSymmtryLine})}, and
this solution satisfies
\begin{eqnarray}
&&0\le\psi(\mxx, \mxy)\le C\epsP \qquad\mbox{for }\;(\mxx,
\mxy)\in\Omega^+(\Kphi), \label{L-infty-for-unif-ellipt}
\\
&&|\psi(x,y)|\le C{\epsP\over \varepsilon} x\qquad \mbox{for
}\;(x,y)\in\ElDomS, \label{barier-for-unif-ellipt}
\end{eqnarray}
where we have used coordinates {\rm (\ref{coordNearSonic})} in {\rm
(\ref{barier-for-unif-ellipt})}. Moreover, for any $s\in (0,
c_2/4)$, there exists $C(s)>0$ depending only on the data and $s$,
but independent of $\delta\in (0, \delta_0)$, such that
\begin{equation}
\|\psi\|^{(-1-\alpha,
\{\PtLwL,\PtLwR\})}_{2,\alpha,\Omega^+_s(\Kphi)}\le C(s)\epsP,
\label{Hoder-est-for-unif-ellipt}
\end{equation}
where $\Omega^+_s(\Kphi):=\Omega^+(\Kphi)\cap\{c_2-r>s\}$.
\end{proposition}
%***************END PROPOSITION***************

%*************** BEGIN PROOF ***************
\Proof
Note that equation (\ref{unif-ellipt-iterationEquation}) is
nonlinear and the boundary conditions
(\ref{iterationRH})--(\ref{iterationCondOnSymmtryLine}) are linear.
We find a solution of
(\ref{iterationRH})--(\ref{iterationCondOnSymmtryLine}) and
(\ref{unif-ellipt-iterationEquation}) as a fixed point of the map
\begin{equation}\label{6.9a}
\hat J: C^{1,\alpha/2}(\overline{\Omega^+(\Kphi)}) \to
C^{1,\alpha/2}(\overline{\Omega^+(\Kphi)})
\end{equation}
defined as follows: For $\hat \psi\in
C^{1,\alpha/2}(\overline{\Omega^+(\Kphi)})$, we consider the linear
 elliptic equation obtained by substituting
$\hat\psi$ into the coefficients of equation
(\ref{unif-ellipt-iterationEquation}):
\begin{equation}
a_{11}\psi_{\mxx\mxx}+ 2a_{12}\psi_{\mxx\mxy} +a_{22}\psi_{\mxy\mxy}=0
\qquad \mbox{ in }\;\;\Omega^+(\Kphi),
\label{unif-ellipt-linear-iterationEquation}
\end{equation}
where
\begin{equation}
a_{ij}(\mxx,\mxy)=A_{ij}(D\hat \psi(\mxx,\mxy),\mxx,\mxy) +\delta\,
\delta_{ij}\qquad \mbox{ for }\;(\mxx,\mxy)\in\Omega^+(\Kphi),\,\,
i,j=1,2, \label{coef-linear-iterationEquation}
\end{equation}
with $\delta_{ij}=1$ for $i=j$ and $0$ for $i\ne j$, $i,j=1,2$. We
establish below the existence of a unique solution $\psi\in
C^{(-1-\alpha, \corners)}_{2,\alpha/2,\Omega^+(\Kphi)}$ to the
linear elliptic equation
(\ref{unif-ellipt-linear-iterationEquation}) with the boundary
conditions (\ref{iterationRH})--(\ref{iterationCondOnSymmtryLine}).
Then we define $\hat J(\hat\psi)=\psi$.

We first state some properties of equation
(\ref{unif-ellipt-linear-iterationEquation}).

%*****************BEGIN LEMMA*************
\begin{lemma}\label{unifEllipticityOfLinearEq-Lemma}
There exists $\hat C>0$ depending only on the data such that, if
$\epsP, \varepsilon>0$ and $M_1, M_2\ge 1$ in {\rm
(\ref{defSetK_R})} satisfy {\rm (\ref{condConst-00})}, and
$\delta\in(0,1)$, then, for any  $\hat \psi\in
C^{1,\alpha/2}(\overline{\Omega^+(\Kphi)})$, equation {\rm
(\ref{unif-ellipt-linear-iterationEquation})} is uniformly elliptic
in $\Omega^+(\Kphi)$:
\begin{equation}\label{ellipticityOfIterEq-Linear-0}
\delta|\mu|^2\le \sum_{i,j=1}^2 a_{ij}(\mxx,\mxy)\mu_i\mu_j \le
2{\lambda}^{-1}|\mu|^2 \qquad\mbox{for }(\mxx,
\mxy)\in\Omega^+(\Kphi),\,
 \mu\in\bR^2,
\end{equation}
where $\lambda$ is from Lemma {\rm\ref{propertiesNonlinCoeffs}}.
Moreover, for any $s\in (0, c_2/2)$, the ellipticity constants
depend only on the data and are independent of $\delta$
 in $\Omega^+_s(\Kphi)=\Omega^+(\Kphi)\cap\{c_2-r>s\}$:
\begin{equation}\label{ellipticityOfIterEq-Linear}
\lambda(c_2-s)|\mu|^2 \le\sum_{i,j=1}^2a_{ij}(\mxx,\mxy)\mu_i\mu_j \le
2{\lambda}^{-1}|\mu|^2 \qquad\mbox{for
}z=(\mxx,\mxy)\in\Omega^+_s(\Kphi),\,\mu\in\bR^2.
\end{equation}
Furthermore,
\begin{equation}\label{holder-coef}
a_{ij}\in C^{\alpha/2}(\overline{\Omega^+(\Kphi)}).
\end{equation}
\end{lemma}
%*****************END LEMMA*************
%*************** BEGIN PROOF ***************
\Proof Facts
(\ref{ellipticityOfIterEq-Linear-0})--(\ref{ellipticityOfIterEq-Linear})
directly follow from the definition of $a_{ij}$ and both the
definition and properties of $A_{ij}$ in \S
\ref{eqForIterationSection} and Lemma \ref{propertiesNonlinCoeffs}.

Since $A_{ij}(p,\mxx,\mxy)$ are independent of $p$ in
$\Omega^+(\Kphi)\cap \{c_2-r>4\varepsilon\}$, it follows from
(\ref{iterationUniforDomEquation}),
(\ref{iterationSonicDomEquation}), and $\Kphi\in\setK$ that $a_{ij}
\in C^{(-\alpha,\Sigma_0)}_{1,\alpha/2,\Omega^+(\Kphi)\cap \DomU}
\subset  C^\alpha(\overline{\Omega^+(\Kphi)\cap \DomU})$.

To show $a_{ij}\in C^{\alpha/2}(\overline{\Omega^+(\Kphi)})$, it
remains to prove that $a_{ij}\in
C^{\alpha/2}(\overline{\Omega(\Kphi)\cap \DomS})$. To achieve this,
we note that the nonlinear terms in the coefficients
$A_{ij}(p,\mxx,\mxy)$ are only the terms
$$
(c_2-r)\zeta_1(\frac{\mxx\psi_\mxx+\mxy\psi_\mxy}{r(c_2-r)}).
$$
Since $\zeta_1$ is a bounded and $C^\infty$-smooth function on
$\bR$,
 and $\zeta'_1$ has compact
support, then there exists $C>0$ such that, for any $s>0$, $q\in
\bR$,
\begin{equation}\label{propZeta1}
\left|s\zeta_1(\frac{q}{s})\right|\le
\big(\sup_{t\in\bR}|\zeta_1(t)|\big)s, \qquad
\left|D_{(q,s)}\big(s\zeta_1(\frac{q}{s})\big)\right| \le C.
\end{equation}
Then it follows that the function
$$
F(p,\mxx, \mxy)= (c_2-r)\zeta_1(\frac{\mxx p_1+\mxy p_2}{r(c_2-r)})
$$
satisfies $|F(p,\mxx, \mxy)|\le \|\zeta_1\|_{L^\infty(\bR)}(c_2-r)$
for any $(p,\mxx,\mxy)\in\bR^2\times\DomS$, and $|\grad_{(p,\mxx,
\mxy)}F|$ is bounded on compact subsets of
$\bR^2\times\overline\DomS$. From this and $\hat\psi\in
C^{1,\alpha/2}(\overline{\Omega^+(\Kphi)})$, we have $a_{ij}\in
C^{\alpha/2}(\overline{\Omega^+(\Kphi)})$.
\Endproof
%*************** END PROOF ***************

Now we state some properties of equation
(\ref{unif-ellipt-linear-iterationEquation}) written in the
$(x,y)$--coordinates.

%*****************BEGIN LEMMA*************
\begin{lemma}\label{unifEllipticityOfLinearEq-xy-Lemma}
There exist $\lambda>0$ and $C,\, \hat C>0$ depending only on the
data such that, if $\epsP, \varepsilon>0$ and $M_1, M_2\ge 1$ in
{\rm (\ref{defSetK_R})} satisfy {\rm (\ref{condConst-00})}, and
$\delta\in(0,1)$, then, for any  $\hat \psi\in
C^{1,\alpha/2}(\overline{\Omega^+(\Kphi)})$, equation {\rm
(\ref{unif-ellipt-linear-iterationEquation})} written in the
$(x,y)$--coordinates has the structure
\begin{equation}
\hat a_{11}\psi_{xx}+ 2\hat a_{12}\psi_{xy} +\hat a_{22} \psi_{yy}
+\hat a_1\psi_{x} +\hat a_2\psi_{y} =0 \quad\; \mbox{in
}\;\Omega^+(\Kphi)\cap\DomS_{4\varepsilon},
\label{unif-ellipt-linear-xy-iterationEquation}
\end{equation}
where $\hat a_{ij}=\hat a_{ij}(x,y)$ and $\hat a_i=\hat a_{i}(x,y)$
satisfy
\begin{equation}\label{holder-coef-xy}
\hat a_{ij}, \hat a_{i}\in
C^{\alpha/2}(\overline{\Omega^+(\Kphi)\cap\DomS_{4\varepsilon}})
\qquad\mbox{for }\;i,j=1,2,
\end{equation}
and the ellipticity condition
\begin{equation}\label{ellipticityOfIterEq-Linear-xy-0}
\begin{array}{l}
\displaystyle \delta\lambda|\mu|^2\le \sum_{i,j=1}^2 \hat
a_{ij}(\mxx,\mxy)\mu_i\mu_j \le \lambda^{-1}|\mu|^2 \quad\mbox{ for
any }(x,y)\in\Omega^+(\Kphi)\cap\DomS_{4\varepsilon},\, \mu\in\bR^2.
\end{array}
\end{equation}
Moreover,
\begin{eqnarray}
&& \qquad\qquad \delta\le \hat a_{11}(x,y)\le\delta+2x, \,\,\,
{1\over 2c_2} \le \hat a_{22}(x,y) \le {2\over c_2}, \,\,\, -2\le
\hat a_{1}(x,y)\le
-{1\over 2}, \nonumber\\
&&\qquad\qquad |(\hat a_{12}, \hat a_{21}, \hat a_{2})(x,y)|\le
 C|x|,
 \qquad |D(\hat a_{12}, \hat a_{21}, \hat a_{2})(x,y)| \le  C|x|^{1/2},
\label{estSmallterms-iter-lin}
\\
&& \qquad \qquad|\hat a_{ii}(x,y)-\hat a_{ii}(0, \tilde y)|\le
C\,|(x,y)-(0,\tilde y)|^\alpha
 \qquad\mbox{for }\;i=1,2,\nonumber
\end{eqnarray}
for all $(x,y), (0,\tilde y)\in
\Omega^+(\Kphi)\cap\DomS_{4\varepsilon}$. \nonumber
\end{lemma}
%*****************END LEMMA*************
%*************** BEGIN PROOF ***************
\Proof
By (\ref{condRewritingRH-0}), if $\varepsilon\le \bar c_2/10$, then
 the change of variables from $(\mxx, \mxy)$ to $(x,y)$ in
$\DomS_{4\varepsilon}$ is smooth and smoothly invertible with
Jacobian bounded away from zero, where the norms and lower bound of
the Jacobian depend only on the data. Now
(\ref{ellipticityOfIterEq-Linear-xy-0}) follows from
(\ref{ellipticityOfIterEq-Linear}).

Equation (\ref{unif-ellipt-linear-iterationEquation}) written in the
$(x,y)$--coordinates can be obtained by substituting $\hat \psi$
into the term $\displaystyle x\zeta_1(\frac{\psi_x}{x})$ in the
coefficients of equation
(\ref{unif-ellipt-iteration-equation-sonicStruct}). Using
(\ref{propZeta1}), the assertions in
 (\ref{holder-coef-xy}) and (\ref{estSmallterms-iter-lin}),
except the last inequality, follow directly from
(\ref{unif-ellipt-iteration-equation-sonicStruct}) with
(\ref{iteration-equation-sonicStruct-coef}) and
(\ref{erTerms-xy-nontrunc}), $\Kphi\in\setK$ with
(\ref{condConst-00}), and $\hat \psi\in
C^{1,\alpha/2}(\overline{\Omega^+(\Kphi)})$.

Then we prove the last inequality in (\ref{estSmallterms-iter-lin}).
We note that, from
(\ref{unif-ellipt-iteration-equation-sonicStruct}) and
(\ref{iteration-equation-sonicStruct-coef}), it follows that $\hat
a_{ii}(x,y)=F_{ii}(\grad\Kphi,\Kphi,x,y)+G_{ii}(x)x\zeta_1(\frac{\hat
\psi_x}{x})$, where $F_{ii}$ and $G_{ii}$ are smooth functions, and
$\Kphi$ and $\hat \psi$ are evaluated at $(x,y)$. In particular,
since $\zeta_1(\cdot)$ is bounded, $\hat
a_{ii}(0,y)=F_{ii}(\grad\Kphi(0,y),\Kphi(0,y),0,y)$. Thus, assuming
$x>0$, we use the boundedness of $\zeta_1$ and $G_{ii}$, smoothness
of $F_{ii}$, and $\Kphi\in\setK$ with Lemma \ref{relatingNorms} to
obtain
\begin{eqnarray*}
&&|\hat a_{ii}(x,y)-\hat a_{ii}(0, \tilde y)|\\
&&\le |F_{ii}(\grad\Kphi(x,y),\Kphi(x,y),x,y)-
F_{ii}(\grad\Kphi(0,\tilde y),\Kphi(0,\tilde y),0,\tilde y)|\\
&&\quad + x|G_{ii}(x)\zeta_1(\frac{\hat \psi_x(x,y)}{x})|\\
&&\le Cx+C(M_1\epsilon^{1-\alpha}+M_2\epsP)|(x,y)-(0,\tilde
y)|^\alpha \le C|(x,y)-(0,\tilde y)|^\alpha,
\end{eqnarray*}
where the last inequality holds since $\alpha\in(0, 1/2)$ and
(\ref{condConst-00}). If $x=0$, the only difference is that the
first term is dropped in the estimates.
\Endproof
%*************** END PROOF *****LEMMA**********

%*****************BEGIN LEMMA*************
\begin{lemma}[Comparison Principle]\label{comparisonPrincipleOfUnifEllipt-Lemma}
There exists $\hat C>0$ depending only on the data such that, if
$\epsP, \varepsilon>0$ and $M_1, M_2\ge 1$ in {\rm
(\ref{defSetK_R})} satisfy {\rm (\ref{condConst-00})}, and
$\delta\in(0,1)$, the following comparison principle holds: Let
$\psi\in C(\overline{\Omega^+(\Kphi)}) \cap
C^{1}(\overline{\Omega^+(\Kphi)}\setminus\overline\sonic) \cap
C^2(\Omega^+(\Kphi))$, let the left-hand sides of {\rm
(\ref{unif-ellipt-linear-iterationEquation})}, {\rm
(\ref{iterationRH})}, and {\rm (\ref{iterationCondOnWedge})}--{\rm
(\ref{iterationCondOnSymmtryLine})} are nonpositive for $\psi$, and
let $\psi\ge 0$ on $\sonic$. Then
$$
\psi\ge 0 \qquad \mbox{in }\, \Omega^+(\Kphi).
$$
\end{lemma}
%*****************END LEMMA*************
%*************** BEGIN PROOF ***************
\Proof
We assume that $\hat C$ is large so that
(\ref{q2-u1})--(\ref{OmegaPL}) hold.

We first note that the boundary condition {\rm (\ref{iterationRH})}
on
 $\shock(\Kphi)$, written as (\ref{iterationRH-lf}),
 satisfies
$$
(b_1, b_2)\cdot\nu>0, \qquad b_3<0 \qquad \mbox{on $\shock(\Kphi)$,}
$$
by (\ref{estCoefsIterRH}) combined with $\hat\mxx<0$ and
$\rho_2>\rho_1$. Thus, if  $\psi$ is not a constant in
$\Omega^+(\Kphi)$, a negative minimum of $\psi$ over
$\overline{\Omega^+(\Kphi)}$ cannot be achieved:
\begin{enumerate}
\item[(i)] In the interior of $\Omega^+(\Kphi)$, by the strong maximum principle for linear
elliptic equations;
\item[(ii)] In the relative interiors of
$\shock(\Kphi)$, $\wedgeB$, and $\partial \Omega^+(\Kphi)\cap
\{\mxy=-v_2\}$, by Hopf's Lemma and the oblique derivative conditions
{\rm (\ref{iterationRH})} and {\rm
(\ref{iterationCondOnWedge})}--{\rm
(\ref{iterationCondOnSymmtryLine})};
\item[(iii)] In the corners $\PtLwL$ and
$\PtLwR$, by the result in Lieberman \cite[Lemma 2.2]{Lieberman87},
via a standard argument as in \cite[Theorem 8.19]{GilbargTrudinger}.
Note that we have to flatten the curve $\shock$ in order to apply
\cite[Lemma 2.2]{Lieberman87} near $\PtLwL$, and this flattening can
be done by using the $C^{1,\alpha}$ regularity of  $\shock$.
\end{enumerate}
Using that  $\psi\ge 0$ on $\sonic$, we conclude the proof.
\Endproof
%*************** END PROOF ***************

%*****************BEGIN LEMMA*************
\begin{lemma}\label{unifEstOfUnifEllipt-Lemma}
There exists $\hat C>0$ depending only on the data such that, if
$\epsP, \varepsilon>0$ and $M_1, M_2\ge 1$ in {\rm
(\ref{defSetK_R})} satisfy {\rm (\ref{condConst-00})}, and
$\delta\in(0,1)$, then any solution $\psi\in
C(\overline{\Omega^+(\Kphi)}) \cap
C^{1}(\overline{\Omega^+(\Kphi)}\setminus\overline\sonic) \cap
C^2(\Omega^+(\Kphi))$ of
 {\rm (\ref{unif-ellipt-linear-iterationEquation})}
and {\rm (\ref{iterationRH})}--{\rm
(\ref{iterationCondOnSymmtryLine})} satisfies {\rm
(\ref{L-infty-for-unif-ellipt})}--{\rm
(\ref{barier-for-unif-ellipt})} with the constant $C$ depending only
on the data.
\end{lemma}
%*****************END LEMMA*************
%*************** BEGIN PROOF ***************
\Proof
First we note that, since $\Omega^+(\Kphi)\subset\{\mxy<c_2\}$,  the
function
$$
w(\mxx,\mxy)=-v_2(\mxy-c_2)
$$
is a nonnegative supersolution of
 (\ref{unif-ellipt-linear-iterationEquation})
and (\ref{iterationRH})--(\ref{iterationCondOnSymmtryLine}): Indeed,
\begin{enumerate}
\item[(i)] $w$ satisfies (\ref{unif-ellipt-linear-iterationEquation}) and
  (\ref{iterationCondOnSymmtryLine});
\item[(ii)]
$w$ is a supersolution of (\ref{iterationRH}). This can be seen by
using (\ref{iterationRH-lf}), (\ref{estCoefsIterRH}), $\rho_2>\rho_1$,
$u_1>0$, $\rho'_2>0$
$\hat\mxx<0$,
and $|\mxy|\le c_2$ to compute on $\shock$:
$$
\Ml(w)=-b_2v_2 -b_3v_2(\mxy-c_2)\le
-v_2\Big(\rho_2'|\hat\mxx|
+\frac{\rho_2-\rho_1}{u_1}-\varepsilon^{3/4}(1+2c_2) \Big)<0
$$
if $\varepsilon$ is small depending on the data, which is achieved
by the choice of $\hat C$ in (\ref{condConst-00});
\item[(iii)]
$w$ is a supersolution of (\ref{iterationCondOnWedge}). This follows
from $\grad w\cdot\nu=-c_2\cos\theta_w<0$ since the interior unit
normal on $\wedgeB$ is $\nu=(-\sin\theta_w, \cos\theta_w)$;

\item[(iv)]
$w\ge 0$ on $\sonic$.
\end{enumerate}
Similarly, $\tilde w\equiv 0$ is a subsolution of
(\ref{unif-ellipt-linear-iterationEquation}) and
(\ref{iterationRH})--(\ref{iterationCondOnSymmtryLine}). Thus, by
the Comparison Principle (Lemma
\ref{comparisonPrincipleOfUnifEllipt-Lemma}), any solution $\psi\in
C(\overline{\Omega^+(\Kphi)}) \cap
C^{1}(\overline{\Omega^+(\Kphi)}\setminus\overline\sonic) \cap
C^2(\Omega^+(\Kphi))$ satisfies
$$
0\le\psi(\mxx, \mxy)\le w(\mxx, \mxy)\qquad\mbox{for any }\;(\mxx,
\mxy)\in\Omega^+(\Kphi).
$$
Since $|v_2|\le C\epsP$, then (\ref{L-infty-for-unif-ellipt})
follows.

To prove (\ref{barier-for-unif-ellipt}), we work in the
$(x,y)$--coordinates in $\DomS\cap \Omega^+(\Kphi)$ and assume that
$\hat C$ in (\ref{condConst-00}) is sufficiently large so that the
assertions of Lemma \ref{unifEllipticityOfLinearEq-xy-Lemma} hold.
Let $v(x,y)=L\epsP x$ for $L>0$. Then

(i) $v$ is a supersolution of equation
(\ref{unif-ellipt-linear-xy-iterationEquation}) in
$\ElDomS\cap\{x<\varepsilon\}$: Indeed, the left-hand side of
(\ref{unif-ellipt-linear-xy-iterationEquation}) on $v(x,y)=L\epsP x$
is $\hat a_1(x,y) L\epsP$, which is negative in $\DomS\cap
\Omega^+(\Kphi)$ by (\ref{estSmallterms-iter-lin});

(ii) $v$ satisfies the boundary conditions
(\ref{condOnSonicLine-Psi-xy}) on $\partial\Omega^+(\Kphi)\cap
\{x=0\}$ and (\ref{condOnWedge-Psi-xy}) on $\partial
\Omega^+(\Kphi)\cap \{y=0\}$;

(iii) The left-hand side of (\ref{iterationRH-lf-flattened}) is
negative for $v$ on $\shock\cap\{x<\varepsilon\}$: Indeed,
$\Ml(v)(x,y)=L\epsP(\hat b_1 +\hat b_3x)<0$ by
(\ref{estCoefsIterRH-flattened}) and since $x\ge0$ in $\overline\ElDomS$.

Now, choosing $L$ large so that $ L\varepsilon>C$ where $C$ is the
constant in (\ref{L-infty-for-unif-ellipt}), we have by
(\ref{L-infty-for-unif-ellipt}) that $v\ge\psi$ on
$\{x=\varepsilon\}$. By the Comparison Principle, which holds since
equation (\ref{unif-ellipt-linear-xy-iterationEquation})
 is elliptic and condition (\ref{iterationRH-lf-flattened}) satisfies
 (\ref{obliqInxy-1}) and $\hat b_3<0$
where the last inequality follows from (\ref{estCoefsIterRH-flattened}),
we obtain $v\ge\psi$
in $\Omega^+(\Kphi)\cap\{x<\varepsilon\}$. Similarly, $-\psi\ge -v$
in $\Omega^+(\Kphi)\cap\{x<\varepsilon\}$. Then
(\ref{barier-for-unif-ellipt}) follows.
\Endproof
%*************** END PROOF ***************

%*****************BEGIN LEMMA*************
\begin{lemma}\label{C2alpha-Lin-UnifEllipt-Lemma}
There exists $\hat C>0$ depending only on the data such that, if
$\epsP, \varepsilon>0$ and $M_1, M_2\ge 1$ in {\rm
(\ref{defSetK_R})} satisfy {\rm (\ref{condConst-00})}, and
$\delta\in(0,1)$, any solution $\psi\in
C(\overline{\Omega^+(\Kphi)}) \cap
C^{1}(\overline{\Omega^+(\Kphi)}\setminus\overline\sonic) \cap
C^2(\Omega^+(\Kphi))$ of
 {\rm (\ref{unif-ellipt-linear-iterationEquation})}
and {\rm (\ref{iterationRH})}--{\rm
(\ref{iterationCondOnSymmtryLine})} satisfies
\begin{equation}
\|\psi\|^{(-1-\alpha,
\{\PtLwL,\PtLwR\})}_{2,\alpha/2,\Omega^+_s(\Kphi)} \le
C(s,\hat\psi)\epsP \label{Hoder-est-for-lin-unif-ellipt}
\end{equation}
for any $s\in (0, c_2/2)$, where the constant $C(s, \hat\psi)$
depends only on the data,
$\|\hat\psi\|_{C^{1,\alpha/2}(\overline{\Omega^+(\Kphi)})}$, and
$s$.
%%but is independent of $\delta\in (0, \delta_0)$.
\end{lemma}
%*****************END LEMMA*************
%*************** BEGIN PROOF ***************
\Proof {}From (\ref{OmegaPL}), (\ref{OmegaPL-f-higher}),
(\ref{estCoefsIterRH-0})--(\ref{estCoefsIterRH}),
(\ref{ellipticityOfIterEq-Linear})--(\ref{holder-coef}), and the
choice of $\alpha$ in \S \ref{Constr-iter-section}, it follows by
\cite[Lemma 1.3]{Lieberman88} that
\begin{equation}\label{weakerLinearObliqueEstimate}
\|\psi\|^{(-1-\alpha, \Sigma_0\cup\shock(\Kphi)\cup\wedgeB )}_{2,\alpha/2,\Omega^+_s(\Kphi)}\le
 C(s,\hat\psi)(\|\psi\|_{C(\Omega^+(\Kphi))}+|v_2|)\le  C(s,\hat\psi)\epsP,
\end{equation}
where
we have used (\ref{u2-v2-bound}) and Lemma
\ref{unifEstOfUnifEllipt-Lemma} in the second inequality.

In deriving (\ref{weakerLinearObliqueEstimate}), we have used
(\ref{OmegaPL-f-higher}) and (\ref{estCoefsIterRH-0}) only to infer
that $\shock(\Kphi)$ is a $C^{1,\alpha}$--curve and $b_i\in
C^\alpha(\overline{\shock(\Kphi)})$. To improve
(\ref{weakerLinearObliqueEstimate}) to
(\ref{Hoder-est-for-lin-unif-ellipt}), we use the higher regularity
of $\shock(\Kphi)$ and $b_i$, given by  (\ref{OmegaPL-f-higher}) and
(\ref{estCoefsIterRH-0}) (and a similar regularity for the boundary
conditions
(\ref{iterationCondOnWedge})--(\ref{iterationCondOnSymmtryLine}),
which are given on the flat segments and have constant
coefficients), combined with  rescaling from the balls
$B_{d/2}(z)\cap \Omega^+(\Kphi)$ for any $z\in
\overline{\Omega^+_s(\Kphi)}\setminus\{\PtLwL, \PtLwR\}$ (with
 $d=\dist(z, \{\PtLwL,\PtLwR\}\cup\Sigma_0)$)
into the unit ball and the standard estimates for the oblique
derivative problems for linear elliptic equations.
\Endproof
%*************** END PROOF *****LEMMA**********

Now we show that the solution $\psi$ is $C^{2,\alpha/2}$ near the
corner $\PtUpR=\sonic\cap\wedgeB(\Kphi)$. We  work  in $\DomS$ in
the $(x,y)$--coordinates.

%*****************BEGIN LEMMA*************
\begin{lemma}\label{wedge-sonic-Lin-UnifEllipt-Lemma}
There exists $\hat C>0$ depending only on the data such that, if
$\epsP, \varepsilon>0$ and $M_1, M_2\ge 1$ in {\rm
(\ref{defSetK_R})} satisfy {\rm (\ref{condConst-00})}, and
$\delta\in(0,1)$, any solution $\psi\in
C(\overline{\Omega^+(\Kphi)}) \cap
C^{1}(\overline{\Omega^+(\Kphi)}\setminus\overline\sonic) \cap
C^2(\Omega^+(\Kphi))$ of {\rm
(\ref{unif-ellipt-linear-iterationEquation})} and {\rm
(\ref{iterationRH})}--{\rm (\ref{iterationCondOnSymmtryLine})} is in
$C^{2,\alpha/2}(\overline{B_\vr(\PtUpR)\cap\Omega^+(\Kphi)})$ for
sufficiently small $\vr>0$.
\end{lemma}
%*****************END LEMMA*************
%*************** BEGIN PROOF ***************

\Proof
In this proof, the constant $C$ depends only on the data, $\delta$,
and $\|(\hat a_{ij}, \hat a_{i})\|_{
C^{\alpha/2}(\overline{\Omega^+(\Kphi)})}$ for $i, j=1,2$, i.e., $C$
is independent of $\vr$.

{\em Step 1.} We work in the $(x,y)$--coordinates. Then
$\PtUpR=(0,0)$ and $\Omega^+(\Kphi)\cap B_{2\vr}= \{x>0,y>0\})\cap
B_{2\vr}$ for $\vr\in (0, \varepsilon)$. Denote
$$
B_\vr^{+}:=B_\vr(0)\cap\{x>0\},\qquad B_\vr^{++}:=B_\vr(0)\cap\{x>0,
y>0\}.
$$
Then $\psi$ satisfies equation
(\ref{unif-ellipt-linear-xy-iterationEquation}) in $B_{2\vr}^{++}$
and
\begin{eqnarray}
&&\psi=0  \qquad\mbox{on }\;\sonic\cap B_{2\vr} =B_{2\vr}\cap \{x=0,
y>0\}, \label{condOnSonicLine-Psi-xy-loc}
\\
&&\psi_\nu\equiv\psi_y=0 \qquad\mbox{on }\;\wedgeB\cap
B_{2\vr}=B_{2\vr}\cap \{y=0, x>0\}. \label{condOnWedge-Psi-xy-loc}
\end{eqnarray}
Rescale $\psi$ by
$$
v(z)=\psi(\vr z)\qquad \mbox{ for } z=(x,y)\in B_2^{++}.
$$
Then $v\in C(\overline{B_2^{++}}) \cap
C^{1}(\overline{B_2^{++}}\setminus\overline{\{x=0\}}) \cap
C^2(B_2^{++})$ satisfies
\begin{equation}
\|v\|_{L^\infty(B_2^{++})}= \|\psi\|_{L^\infty(B_{2\vr}^{++})},
\label{L-infty-for-unif-ellipt-v}
\end{equation}
and $v$ is a solution of
\begin{eqnarray}
&&\hat a_{11}^{(\vr)}v_{xx}+ 2\hat a_{12}^{(\vr)}v_{xy} +\hat
a_{22}^{(\vr)}v_{yy} +\hat a_1^{(\vr)}v_{x} +\hat a_2^{(\vr)}v_{y}
=0 \qquad\mbox{in }\;B_2^{++}, \label{unif-ellipt-linear-xy-v}
\\
&&v=0\qquad\mbox{on }\;\partial B_2^{++}\cap \{x=0\},
\label{condOnSonicLine-v-xy-loc}
\\
&&v_\nu\equiv v_y=0\qquad\mbox{on }\;\partial B_2^{++}\cap \{y=0\},
\label{condOnWedge-v-xy-loc}
\end{eqnarray}
where
\begin{equation}\label{def-rescaled-coef-v}
\hat a_{ij}^{(\vr)}(x,y)=\hat a_{ij}(\vr x,\vr y),\,\,\, \hat
a_{i}^{(\vr)}(x,y)=\vr\,\hat a_{i}(\vr x,\vr y)\quad\mbox{for }\;
(x,y)\in B_2^{++},\;i,j=1,2.
\end{equation}
Thus, $\hat a_{ij}^{(\vr)}$ satisfy
(\ref{ellipticityOfIterEq-Linear-xy-0}) with the unchanged constant
$\lambda>0$ and, since $\vr\le 1$,
\begin{equation}\label{holder-coef-xy-v}
\|(\hat a_{ij}^{(\vr)}, \hat a_{i}^{(\vr)})\|_{
C^{\alpha/2}(\overline{B_2^{++}})} \le \|(\hat a_{ij}, \hat
a_{i})\|_{ C^{\alpha/2}(\overline{\Omega^+(\Kphi)})} \qquad\mbox{for
}\;i,j=1,2.
\end{equation}
Denote $Q:=\{z\in  B_2^{++}\;:\;\dist(z, \partial B_2^{++})>1/50\}$.
The interior estimates for the elliptic equation
(\ref{unif-ellipt-linear-xy-v}) imply
$\|v\|_{C^{2,\alpha/2}(\overline{Q})}\le
C\|v\|_{L^\infty(B_2^{++})}$. The local estimates for the Dirichlet
problem
(\ref{unif-ellipt-linear-xy-v})--(\ref{condOnSonicLine-v-xy-loc})
imply
\begin{equation}\label{localEstNearBdry}
\|v\|_{C^{2,\alpha/2}(\overline{B_{1/10}(z)\cap B_2^{++}})} \le
C\|v\|_{L^\infty(B_2^{++})}
\end{equation}
for every $z=(x,y)\in\{x=0, 1/2\le y\le 3/2\}$. The local estimates
for the oblique derivative problem (\ref{unif-ellipt-linear-xy-v})
and (\ref{condOnWedge-v-xy-loc}) imply (\ref{localEstNearBdry}) for
every $z\in\{1/2\le x\le 3/2, y=0\}$. Then we have
\begin{equation}\label{Hoder-est-near-Br}
\|v\|_{C^{2,\alpha/2}(\overline{B_{3/2}^{++}\setminus
B_{1/2}^{++}})} \le C\|v\|_{L^\infty(B_2^{++})}.
\end{equation}

\medskip
{\em Step 2.} We modify the domain $B_1^{++}$ by mollifying the
corner at $(0,1)$ and denote the resulting domain by $D^{++}$. That
is, $D^{++}$ denotes an open domain satisfying
$$
D^{++}\subset B_{1}^{++},\qquad D^{++}\setminus B_{1/10}(0,1)=
B_1^{++}\setminus B_{1/10}(0,1),
$$
and
$$
\partial D^{++}\cap B_{1/5}(0,1)\qquad
\mbox{is a  $C^{2,{\alpha/2}}$--curve}.
$$
Then we prove the following fact: For any $g\in
C^{\alpha/2}(\overline{D^{++}})$, there exists a unique solution
$w\in C^{2,\alpha/2}(\overline{D^{++}})$ of the problem:
\begin{equation}
\begin{array}{ll}
&\hat a_{11}^{(\vr)}w_{xx} +\hat a_{22}^{(\vr)}w_{yy} +\hat
a_1^{(\vr)}w_{x} =g \qquad\mbox{in }\;D^{++},
\\
&w=0\qquad\mbox{on }\;\partial D^{++}\cap \{x=0, y>0\},
\\
&w_\nu\equiv w_y=0\qquad\mbox{on }\;\partial D^{++}\cap \{x>0,
y=0\}, %\label{condOnWedge-w-xy-loc}
\\
&w=v\qquad\mbox{on }\;\partial D^{++}\cap \{x>0, y>0\},
\end{array}
\label{unif-ellipt-linear-xy-w}
\end{equation}
with
\begin{equation}\label{Hoder-est-near-origin}
\|w\|_{C^{2,{\alpha/2}}(\overline{D^{++}})}\le
C(\|v\|_{L^\infty(B_{2}^{++})}+
\|g\|_{C^{\alpha/2}(\overline{D^{++}})}).
\end{equation}

\medskip
This can be seen as follows.
Denote by $D^+$ the even extension of $D^{++}$ from $\{x,y >0\}$
into  $\{x >0\}$, i.e.,
$$
D^+:=D^{++}\cup\{(x,0)\; : \;x\in(0, 1)\}\cup D^{+-},
$$
where $D^{+-}:=\{(x, y)\; : \;(x,-y)\in D^{++}\}$. Then
$B^+_{7/8}\subset D^+\subset B^+_{1}$ and $\partial D^+$ is a
$C^{2,{\alpha/2}}$--curve. Extend $F=(v,g,\hat a_{11}^{(\vr)},\hat
a_{22}^{(\vr)},\hat a_{1}^{(\vr)})$ from $\overline{B_2^{++}}$ to
$\overline{B_2^+}$ by setting
$$
F(x,-y)=F(x,y) \qquad\,\,\mbox{for }\;(x,y)\in \overline{B_2^{++}}.
$$
Then it follows from
(\ref{condOnSonicLine-v-xy-loc})--(\ref{condOnWedge-v-xy-loc}) and
(\ref{Hoder-est-near-Br}) that, denoting by $\hat v$ the restriction
of (extended) $v$ to $\partial D^+$, we have $\hat v\in
C^{2,{\alpha/2}}(\partial D^+)$ with
\begin{equation}\label{Hoder-est-bdry-funct}
\|\hat v\|_{C^{2,{\alpha/2}}(\partial D^+)}\le
C\|v\|_{L^\infty(B_{2}^{++})}.
\end{equation}
Also, the extended $g$ satisfies $g\in C^{\alpha/2}(\overline{D^+})$
with $\|g\|_{C^{\alpha/2}(\overline{D^+})}
=\|g\|_{C^{\alpha/2}(\overline{D^{++}})}$. The extended $(\hat
a_{11}^{(\vr)}, \hat a_{22}^{(\vr)}, \hat a_{1}^{(\vr)})$ satisfy
(\ref{ellipticityOfIterEq-Linear-xy-0}) and
\begin{eqnarray*}
\|(\hat a_{11}^{(\vr)}, \hat a_{22}^{(\vr)}, \hat a_{1}^{(\vr)})\|_{
C^{\alpha/2}(\overline{B_2^{+}})}&=&\|(\hat a_{11}^{(\vr)}, \hat
a_{22}^{(\vr)}, \hat a_{1}^{(\vr)})\|_{
C^{\alpha/2}(\overline{B_2^{++}})}\\
&\le& \sum_{i,j=1}^2\|(\hat a_{ij}, \hat a_{i})\|_{
C^{\alpha/2}(\overline{\Omega^+(\Kphi)})}.
\end{eqnarray*}
Then, by
\cite[Theorem 6.8]{GilbargTrudinger}, there exists a unique solution
$w\in C^{2,{\alpha/2}}(D^+)$ of the Dirichlet problem
\begin{eqnarray}
&&\hat a_{11}^{(\vr)}w_{xx} +\hat a_{22}^{(\vr)}w_{yy} +\hat
a_1^{(\vr)}w_{x} =g \qquad\mbox{in }\;D^{+},
\label{unif-ellipt-linear-xy-v-ext}
\\
&& w=\hat v\qquad\mbox{on }\;\partial D^{+}, \label{dirichlet-v-ext}
\end{eqnarray}
and $w$ satisfies
\begin{equation}\label{Hoder-est-near-Br-w}
\|w\|_{C^{2,{\alpha/2}}(\overline{D^{+}})}\le C(\|\hat
v\|_{C^{2,{\alpha/2}}(\partial D^+)}+
\|g\|_{C^{\alpha/2}(\overline{D^+})}).
\end{equation}
{}From the structure of equation
(\ref{unif-ellipt-linear-xy-v-ext}) and the symmetry of the domain
and the coefficients and right-hand sides obtained by the even
extension, it follows that $\hat w$, defined by $\hat w(x,y)=w(x,
-y)$ in $D^+$, is also a solution of
(\ref{unif-ellipt-linear-xy-v-ext})--(\ref{dirichlet-v-ext}). By
uniqueness for
(\ref{unif-ellipt-linear-xy-v-ext})--(\ref{dirichlet-v-ext}), we
find
$$
w(x,y)=w(x, -y) \qquad  \mbox{in }  D^+.
$$
Thus, $w$ restricted to $D^{++}$ is a solution of
(\ref{unif-ellipt-linear-xy-w}), where we use
(\ref{condOnSonicLine-v-xy-loc}) to see that $w=0$ on $\partial
D^{++}\cap \{x=0, y>0\}$. Moreover, (\ref{Hoder-est-bdry-funct}) and
(\ref{Hoder-est-near-Br-w}) imply (\ref{Hoder-est-near-origin}). The
uniqueness of the solution $w\in
C^{2,{\alpha/2}}(\overline{D^{++}})$ of
(\ref{unif-ellipt-linear-xy-w}) follows from the Comparison
Principle (Lemma \ref{comparisonPrincipleOfUnifEllipt-Lemma}).

\medskip
{\em Step 3.} Now we prove the existence of a solution $w\in
C^{2,{\alpha/2}}(\overline{D^{++}})$ of the problem:
\begin{equation}
\begin{array}{ll}
&\hat a_{11}^{(\vr)}w_{xx}+2\hat a_{12}^{(\vr)}w_{xy} +\hat
a_{22}^{(\vr)}w_{yy} +\hat a_1^{(\vr)}w_{x} +\hat a_2^{(\vr)}w_{y}
=0
\qquad\mbox{in }\;D^{++},\\
&w=0\qquad\mbox{on }\;\partial D^{++}\cap \{x=0, y>0\},
\\
&w_\nu\equiv w_y=0\qquad\mbox{on }\;\partial D^{++}\cap \{y=0,
x>0\}, %\label{condOnWedge-w-xy-loc-fullEq}
\\
&w=v\qquad\mbox{on }\;\partial D^{++}\cap \{x>0, y>0\}.
\end{array}
\label{unif-ellipt-linear-xy-w-fullEq}
\end{equation}
Moreover, we prove that $w$ satisfies
\begin{equation}\label{Hoder-est-near-Br-fullEq}
\|w\|_{C^{2,{\alpha/2}}(\overline{D^{++}})}\le
C\|v\|_{L^\infty(B_{2}^{++})}.
\end{equation}

We obtain such $w$ as a fixed point of map
$K:C^{2,{\alpha/2}}(\overline{D^{++}}) \to
C^{2,{\alpha/2}}(\overline{D^{++}})$ defined as follows. Let $W\in
C^{2,{\alpha/2}}(\overline{D^{++}})$. Define
\begin{equation}\label{defG-corner}
g=-2\hat a_{12}^{(\vr)}W_{xy}-\hat a_2^{(\vr)}W_{y}.
\end{equation}
By (\ref{estSmallterms-iter-lin}) and
(\ref{def-rescaled-coef-v}) with $\vr\in(0,1)$, we find
\begin{equation}\label{estimateCoefR}
\|(a_{12}^{(\vr)}, a_2^{(\vr)})\|_{C^{\alpha/2}(\overline{D^{++}})}
\le C\vr^{1/2},
\end{equation}
which implies
$$
g\in C^{\alpha/2}(\overline{D^{++}}).
$$
Then, by the results of Step 2, there exists a unique solution
$w\in C^{2,{\alpha/2}}(\overline{D^{++}})$ of
(\ref{unif-ellipt-linear-xy-w})
with $g$
defined by (\ref{defG-corner}). We set $K[W]=w$.

Now we prove that, if $\vr>0$ is sufficiently small, the map $K$ is
a contraction map. Let $W^{(i)}\in
C^{2,{\alpha/2}}(\overline{D^{++}})$ and $w^{(i)}:=K[W^{(i)}]$ for
$i=1,2$. Then $w\defd w^{(1)}-w^{(2)}$ is a solution of
(\ref{unif-ellipt-linear-xy-w})
%--(\ref{dirichlet-w})
with
\begin{eqnarray*}
&&g=-2\hat a_{12}^{(\vr)}(W^{(1)}_{xy}-W^{(2)}_{xy})-
\hat a_2^{(\vr)}(W^{(1)}_{y}-W^{(2)}_{y}), \\
&& v\equiv 0.
\end{eqnarray*}
Then  $g\in C^{\alpha/2}(\overline{D^{++}})$ and, by
(\ref{estimateCoefR}),
$$
\|g\|_{C^{\alpha/2}(\overline{D^{++}})} \le
C\vr^{1/2}\|W^{(1)}-W^{(2)}\|_{C^{2,{\alpha/2}}(\overline{D^{++}})}.
$$
Since $v\equiv 0$ satisfies
(\ref{condOnSonicLine-v-xy-loc})--(\ref{condOnWedge-v-xy-loc}), we
can apply both (\ref{Hoder-est-near-origin}) and the results of Step
2 to obtain
\begin{eqnarray*}
\|w^{(1)}-w^{(2)}\|_{C^{2,{\alpha/2}}(\overline{D^{++}})}&\le&
C\vr^{1/2}\|W^{(1)}-W^{(2)}\|_{C^{2,{\alpha/2}}(\overline{D^{++}})}\\
&\le& {1\over
2}\|W^{(1)}-W^{(2)}\|_{C^{2,{\alpha/2}}(\overline{D^{++}})},
\end{eqnarray*}
where the last inequality holds if $\vr>0$ is sufficiently
small. We fix such $\vr$. Then the map $K$ has a fixed point $w\in
C^{2,{\alpha/2}}(\overline{D^{++}})$ which is a solution of
(\ref{unif-ellipt-linear-xy-w-fullEq}).

\medskip
{\em Step 4.} Since $v$ satisfies
(\ref{unif-ellipt-linear-xy-v})--(\ref{condOnWedge-v-xy-loc}), it
follows from the uniqueness of solutions in $C(\overline{D^{++}})
\cap C^{1}(\overline{D^{++}}\setminus\overline{\{x=0\}}) \cap
C^2(D^{++})$ of problem (\ref{unif-ellipt-linear-xy-w-fullEq}) that
$w=v$ in $D^{++}$. Thus $v\in C^{2, {\alpha/2}}(\overline{D^{++}})$
so that $\psi\in C^{2,
{\alpha/2}}(\overline{B_{\vr/2}(\PtUpR)\cap\Omega^+(\Kphi)})$.
\Endproof
%*************** END PROOF *****LEMMA**********

Now we prove that the solution $\psi$ is $C^{1,\alpha}$ near the
corner $\PtUpL=\sonic\cap\shock(\Kphi)$ if $\delta$ is small.
%*****************BEGIN LEMMA*************
\begin{lemma}\label{C2alpha-near-P4-Lin-UnifEllipt-Lemma}
There exist $\hat C>0$ and $\delta_0\in (0,1)$ depending only on the
data such that, if $\epsP, \varepsilon>0$ and $M_1, M_2\ge 1$ in
{\rm (\ref{defSetK_R})} satisfy {\rm (\ref{condConst-00})}, and
$\delta\in (0,\delta_0)$, then any solution $\psi\in
C(\overline{\Omega^+(\Kphi)}) \cap
C^1(\overline{\Omega^+(\Kphi)}\setminus\overline\sonic) \cap
C^2(\Omega^+(\Kphi))$ of {\rm
(\ref{unif-ellipt-linear-iterationEquation})} and {\rm
(\ref{iterationRH})}--{\rm (\ref{iterationCondOnSymmtryLine})} is in
$C^{1,{\alpha}}(\overline{B_\vr(\PtUpL)\cap\Omega^+(\Kphi)}) \cap
C^{2,{\alpha/2}}(B_\vr(\PtUpL)\cap\Omega^+(\Kphi))$, for
sufficiently small $\vr>0$ depending only on the data and $\delta$,
and satisfies
\begin{equation}
\|\psi\|^{(-1-\alpha, \{\PtUpL\})}_{2,\alpha/2,\Omega^+(\Kphi)}\le
 C(\delta,\hat\psi)\epsP,
\label{Hoder-est-P4-for-lin-unif-ellipt}
\end{equation}
where $C$ depends only on the data, $\delta$, and
$\|\hat\psi\|_{C^{1,\alpha/2}(\overline{\Omega^+(\Kphi)})}$.
Moreover, for $\delta$ as above,
\begin{equation}
|\psi(x)|\le \tilde C(\delta)(\dist(x, \PtUpL))^{1+\alpha}
\qquad\mbox{for any }\; x\in \Omega^+(\Kphi),
\label{growth-est-P4-for-lin-unif-ellipt}
\end{equation}
where $\tilde C$ depends only on the data and $\delta$,
and is independent of $\hat\psi$.
\end{lemma}
%*****************END LEMMA*************
%*************** BEGIN PROOF ***************
\Proof
In Steps 1--3 of this proof below, the positive constants $C$ and
$L_i, 1\le i\le 4,$ depend only on the data.

\medskip
{\em Step 1.} We work in the $(x,y)$--coordinates. Then the point
$\PtUpL$ has the coordinates $(0, y_{\PtUpL})$ with $y_{\PtUpL}
=\pi/2+\arctan{(|\mxx_1|/\mxy_1)}-\theta_w>0$. {}From
(\ref{domain-in-rescaled-lemma})--(\ref{holder-hat-f}), we have
$$
\Omega^+(\Kphi)\cap B_\kappa(\PtUpL) =\{x>0,y<\hat
f_\Kphi(x)\}\cap B_\varepsilon(\PtUpL),
$$
where $\hat f_\Kphi(0)=y_{\PtUpL}$, $\hat f_\Kphi'(0)>0$, and $\hat
f_\Kphi>y_{\PtUpL}$ on $\bR_+$ by (\ref{domain-in-xy-funct-0}) and
(\ref{holder-hat-f}).

\medskip
{\em Step 2.} We change the variables in such a way that $\PtUpL$
becomes the origin and the second-order part of equation
(\ref{unif-ellipt-linear-iterationEquation}) at $\PtUpL$ becomes the
Laplacian. Denote
\begin{equation}
\Dnu=\sqrt{\hat a_{11}(\PtUpL)/\hat a_{22}(\PtUpL)}.
\label{defLambdaDistortion}
\end{equation}
Then, using (\ref{estSmallterms-iter-lin}) and $x_{\PtUpL}=0$, we
have
\begin{equation}
\sqrt{c_2\delta/2}\le \Dnu\le \sqrt{2c_2\delta}.
\label{estElliptDistortion}
\end{equation}
Now we introduce the variables
$$
(X, Y):=(x/\Dnu,y_{\PtUpL}-y).
$$
Then, for $\vr=\varepsilon$, we have
\begin{equation} \label{domainResceled-distortion}
\Omega^+(\Kphi)\cap B_\vr=\{X>0,\;Y> F(X)\}\cap B_\vr,
\end{equation}
where $F(X)=y_{\PtUpL}-\hat f_\Kphi(\Dnu X)$. By
(\ref{holder-hat-f}), we have $0<\hat f_\Kphi'(X)\le C$ for all
$X\in[0,2\varepsilon]$ if $\hat C$ is sufficiently large in
(\ref{condConst-00}) so that $2\varepsilon\le \kappa$. With this, we
use $\hat f_\Kphi(0)=y_{\PtUpL}$ and (\ref{estElliptDistortion}) to
obtain
\begin{eqnarray}
F(0)=0, \qquad \quad-L_1\sqrt{\delta}\le F'(X)<0
 \quad  \mbox{ for }\;X\in [0,\vr].
\label{properties-F-cap}
\end{eqnarray}

We now write $\psi$ in the $(X,Y)$--coordinates. Introduce the
function
$$
v(X,Y):=\psi(x,y)=\psi(\Dnu X, y_{\PtUpL}-Y).
$$
Since $\psi$ satisfies equation (\ref{iterationRH-lf-flattened}) and
the boundary conditions (\ref{iterationCondOnWedge}) and
(\ref{unif-ellipt-linear-xy-iterationEquation}), then $v$ satisfies
\begin{eqnarray}
&&\qquad\,\, Av\defd{1\over\Dnu^2}\tilde a_{11}v_{XX}-
{2\over\Dnu}\tilde a_{12}v_{XY} +\tilde a_{22} v_{YY}
+{1\over\Dnu}\tilde a_1v_{X} -\tilde a_2v_{Y}
=0\label{unif-ellipt-linear-XY-for-v}
\\
&&\qquad\qquad\qquad\qquad\qquad \qquad\qquad\qquad\quad \mbox{ in
}\;\;\{X>0,\;Y> F(X)\}\cap B_\vr,\nonumber
\\
&&\qquad\,\, Bv\defd{1\over\Dnu}\tilde b_1 v_X - \tilde b_2 v_Y +
\tilde b_3 v=0 \qquad\mbox{on }\;\{X>0,\;Y=F(X)\}\cap B_\vr,
\label{RH-XY-for-v}
\\
&&\qquad\,\, v=0 \qquad\mbox{on }\;\{X=0,\;Y>0\}\cap B_\vr,
\label{cond-on-sonic-for-v}
\end{eqnarray}
where
\begin{eqnarray*}
&&\tilde a_{ij}(X,Y)=\hat a_{ij}(\Dnu X, y_{\PtUpL}-Y), \,\, \tilde
a_{i}(X,Y)=\hat a_{i}(\Dnu X, y_{\PtUpL}-Y), \\
&&\tilde b_{i}(X,Y)=\hat b_{i}(\Dnu X, y_{\PtUpL}-Y).
\end{eqnarray*}
In particular, from (\ref{holder-coef-xy}),
(\ref{estSmallterms-iter-lin}), and (\ref{defLambdaDistortion}), we
have
\begin{eqnarray}
&&\qquad \tilde a_{ij}, \tilde a_{i}\in
C^{\alpha/2}(\overline{\{X>0,\;Y> F(X)\}\cap B_\vr}),\\
&& \qquad \tilde a_{22}(0,0)={1\over\Dnu^2}\tilde a_{11}(0,0),\quad
\tilde a_{12}(0,0)=\tilde a_{2}(0,0)=0,
\label{estSmallterms-iter-lin-XY}
\\
&&\qquad |\tilde a_{ii}(X,Y)-\tilde  a_{ii}(0, 0)|\le C|(X,
Y)|^\alpha \qquad\mbox{for }\;i=1,2,
\label{estSmallterms-iter-lin-contII}
\\
&&\qquad |\tilde a_{12}(X,Y)|+|\tilde a_{21}(X,Y)|+ |\tilde
a_{2}(X,Y)|\le C|X|^{1/2}, \quad |\tilde a_{1}(X,Y)|\le C.
\label{estSmallterms-iter-lin-contIJ}
\end{eqnarray}
{}From (\ref{estCoefsIterRH-flattened}), there exists $L_2>0$ such that
\begin{equation}
\label{estCoefsIterRH-flattened-XY} -L_2^{-1}\le \tilde
b_i(X,Y)\le -L_2 \qquad\mbox{ for any }\;
(X,Y)\in\{X>0,\;Y=F(X)\}\cap B_\vr.
\end{equation}
Moreover, (\ref{obliqInxy-1}) implies
\begin{equation}\label{obliqueness-blowup}
(\tilde b_1, \tilde b_2)\cdot\nu_F>0\qquad\mbox{ on
}\;\{X>0,\;Y=F(X)\}\cap B_\vr,
\end{equation}
where $\nu_F=\nu_F(X,Y)$ is the interior unit normal at
$(X,Y)\in\{X>0,\;Y=F(X)\}\cap B_\vr$. Thus condition
(\ref{RH-XY-for-v}) is oblique.

\medskip
{\em Step 3.} We use the polar coordinates $(r, \theta)$ on the
$(X,Y)$--plane, i.e.,
$$
(X,Y)=(r\cos\theta, r\sin\theta).
$$
{}From (\ref{properties-F-cap}), we have $F, F'<0$ on $(0,\vr)$,
which implies that $(X^2+F(X)^2)'>0$ on $(0,\vr)$. Then it follows
from (\ref{properties-F-cap}) that, if $\delta>0$ is a small
constant depending only on the data and $\vr$ is a small constant
depending only on the data and $\delta$, there exist a function
$\theta_F\in C^1(\R_+)$ and a constant $L_3>0$ such that
\begin{equation} \label{domain-XY-polar}
\{X>0,\;Y> F(X)\}\cap B_\vr=\{0<r<\vr,\; \theta_F(r)<\theta<\pi/2\}
\end{equation}
with
\begin{equation} \label{bdry-XY-polar}
-L_3\sqrt\delta\le\theta_F(r)\le 0.
\end{equation}

Choosing sufficiently small $\delta_0>0$, we show that, for any
$\delta\in(0,\delta_0)$, a function
\begin{equation}
\label{supersolution-P4} w(r,\theta)=r^{1+\alpha}\cos
G(\theta),\qquad\mbox{with } G(\theta)={3+\alpha\over
2}(\theta-{\pi\over 4}),
\end{equation}
is a positive supersolution of
(\ref{unif-ellipt-linear-XY-for-v})--(\ref{cond-on-sonic-for-v}) in
$\{X>0,\;Y> F(X)\}\cap B_\vr$.

By (\ref{domainResceled-distortion}) and
(\ref{domain-XY-polar})--(\ref{bdry-XY-polar}), we find that, when $
0<\delta\le
\delta_0\le\big(\frac{(1-\alpha)\pi}{8(3+\alpha)L_3}\big)^2,
$
$$
-{\pi\over 2}+{1-\alpha\over 16}\pi\le G(\theta)\le {\pi\over
2}-{1-\alpha\over 8}\pi\qquad \mbox{for all }(r,\theta)\in
\Omega^+(\Kphi)\cap B_\vr.
$$
In particular,
\begin{equation}\label{polarAngleBounds}
\cos(G(\theta))\ge \sin\big({1-\alpha\over 16}\pi\big)>0\qquad
\mbox{for all }(r,\theta)\in \overline{\Omega^+(\Kphi)\cap
B_\vr}\setminus\{X=Y=0\},
\end{equation}
which implies
$$
w> 0 \qquad\mbox{ in }\,\, \{X>0,\;Y> F(X)\}\cap B_\vr.
$$
By (\ref{domain-XY-polar})--(\ref{bdry-XY-polar}), we
find that, for all $r\in(0, \vr)$
and $\delta\in (0,\delta_0)$ with  small  $\delta_0>0$,
$$
\cos(\theta_F(r))\ge 1-C\delta_0>0,\qquad |\sin(\theta_F(r))|\le
C\sqrt\delta_0.
$$

Now, possibly further reducing $\delta_0$, we show that $w$ is a
supersolution of (\ref{RH-XY-for-v}). Using
(\ref{estElliptDistortion}), (\ref{RH-XY-for-v}),
(\ref{estCoefsIterRH-flattened-XY}), the above estimates of
$(\theta_F, G(\theta_F))$ derived above, and the fact that
$\theta=\theta_F$ on $\{X>0,\;Y=F(X)\}\cap B_\vr$, we have
\begin{eqnarray*}
Bw&\le&{\tilde b_1\over
\Dnu}r^\alpha\Big((\alpha+1)\cos(\theta_F)\cos(G(\theta_F))
+{3+\alpha\over 2}\sin(\theta_F)\sin(G(\theta_F)) \Big)\\
&& +Cr^\alpha|\tilde b_2|+Cr^{\alpha+1}|\tilde b_3|\\
&\le&-r^\alpha \Big((1-C\delta_0)({L_2\sin({1-\alpha\over
16}\pi)\over C\sqrt\delta_0}-{C\over L_2})
 -C \Big)<0,
\end{eqnarray*}
if $\delta_0$ is sufficiently small. We now fix  $\delta_0$ that
satisfies all the smallness assumptions made above.

Finally, we show that $w$ is a supersolution of equation
(\ref{unif-ellipt-linear-XY-for-v}) in $(X,Y)\in \{X>0,\;Y>
F(X)\}\cap B_\vr$ if $\vr$ is small. Denote by $A_0$ the operator
obtained by fixing the coefficients of $A$ in
(\ref{unif-ellipt-linear-XY-for-v}) at $(X,Y)=(0,0)$. Then
$A_0=\tilde a_{22}(0,0)\Delta$ by (\ref{estSmallterms-iter-lin-XY}). By
(\ref{estSmallterms-iter-lin}), we obtain
$\tilde a_{22}(0,0)=\hat a_{22}(0, y_{\PtUpL})\ge 1/(4\bar{c}_2)>0$.
Now, by an explicit
calculation and using (\ref{estElliptDistortion}),
(\ref{estSmallterms-iter-lin-XY})--(\ref{estSmallterms-iter-lin-contIJ}),
(\ref{domain-XY-polar}), and (\ref{polarAngleBounds}), we find that,
for $\delta\in(0,\delta_0)$ and $(X,Y)\in \{X>0,\;Y> F(X)\}\cap
B_\vr$,
\begin{eqnarray*}
Aw(r,\theta) &=& a_2(0,0)\Delta w(r,\theta) +(A-A_0)w(r,\theta)
\\
&\le&\tilde a_{22}(0,0)r^{\alpha-1}\big((\alpha+1)^2
-({3+\alpha\over 2})^2\big)\cos(G(\theta))
\\
&&+Cr^{\alpha-1}\bigg( {1\over\Dnu^2}|\tilde a_{11}(X,Y)-\tilde
a_{11}(0,0)|
+|\tilde a_{22}(X,Y)-\tilde a_{22}(0,0)|
\bigg)\\
&&+{C\over \Dnu}r^{\alpha-1}|\tilde a_{12}(X,Y)| +
{C\over \Dnu}r^{\alpha}|\tilde a_{1}(X,Y)|+
Cr^{\alpha}|\tilde a_{2}(X,Y)|
\\
&\le& r^{\alpha-1}\left( - {(1-\alpha)(5+3\alpha)\over 8\bar{c}_2}
\sin\big({1-\alpha\over 16}\pi\big)+C{\vr^{\alpha/ 2}\over\sqrt{\delta}} \right)<0
\end{eqnarray*}
for sufficiently small $\vr>0$ depending only on the data and $\delta$.

Thus, all the estimates above hold for small $\delta_0>0$  and
$\vr>0$ depending only on the data.

Now, since
$$
\displaystyle\min_{\{X\ge 0,\;Y\ge F(X)\}\cap\partial
B_\vr}w(X,Y)=L_4>0,
$$
we use the Comparison Principle (Lemma
\ref{comparisonPrincipleOfUnifEllipt-Lemma}) (which holds since
condition (\ref{RH-XY-for-v}) satisfies (\ref{obliqueness-blowup})
and $\tilde b_3<0$ by (\ref{estCoefsIterRH-flattened-XY})) to obtain
$$
{\|\psi\|_{L^\infty(\Omega^+(\Kphi))}\over L_4} w \ge v\qquad\mbox{ in
}\;\{X>0,\;Y> F(X)\}\cap B_\vr.
$$
Similar estimate can be obtained for $-v$.
Thus, using (\ref{L-infty-for-unif-ellipt}), we obtain
(\ref{growth-est-P4-for-lin-unif-ellipt}) in $B_\vr$. Since $\vr$
depends only on the data and $\delta>0$, then we use
(\ref{L-infty-for-unif-ellipt}) to obtain the full estimate
(\ref{growth-est-P4-for-lin-unif-ellipt}).

\medskip
{\em Step 4.} Estimate
(\ref{Hoder-est-P4-for-lin-unif-ellipt}) can be obtained from
(\ref{estCoefsIterRH-flattened}), (\ref{holder-coef-xy}), and
(\ref{growth-est-P4-for-lin-unif-ellipt}), combined with  rescaling
from the balls $B_{d_z/L}(z)\cap \Omega^+(\Kphi)$ for $z\in
\overline{\Omega^+_s(\Kphi)}\setminus\{\PtUpL\}$ (with $d_z=\dist(z,
\PtUpL)$ and $L$ sufficiently large depending only on the data) into
the unit ball and the standard interior estimates for the linear
elliptic equations and the local estimates for the linear Dirichlet
and oblique derivative problems in smooth domains. Specifically,
from the definition of sets $\setK$ and $\Omega^+(\Kphi)$ and by
(\ref{condConst-00}), there exists $L\ge 1$ depending only on the
data such that
$$
B_{d/L}(z)\cap (\partial \Omega^+(\Kphi)\setminus \shock)=\emptyset
\qquad\mbox{for any }\,\, z\in \shock\cap\Omega_\vr,
$$
and
$$
B_{d/L}(z)\cap (\partial \Omega^+(\Kphi)\setminus \sonic)=\emptyset
\qquad\mbox{for any } z\in \sonic\cap\Omega_\vr.
$$
Then, for any $z\in \Omega^+(\Kphi)\cap B_\vr(\PtUpL)$, we have at
least one of the following three cases:

\begin{enumerate}\renewcommand{\theenumi}{\arabic{enumi}}

\item  $B_{d\over 10L}(z)\subset \Omega^+(\Kphi)$;

\item $z\in B_{d_{z_1}\over 2L}(z_1)$
and ${d_z\over d_{z_1}}\in({1\over 2}, 2)$ for some $z_1\in \sonic$;

\item $z\in B_{d_{z_1}\over 2L}(z_1)$
and ${d_z\over d_{z_1}}\in({1\over 2}, 2)$ for some $z_1\in \shock$.
\end{enumerate}

Thus, it suffices to make the $C^{2,\alpha}$--estimates of $\psi$ in
the following subdomains for $z_0=(x_0, y_0)$:

\begin{enumerate}\renewcommand{\theenumi}{\roman{enumi}}
\item \label{cases-HolderCorner-1}
$B_{d_{z_0}\over 20L}({z_0})$ when $B_{d_{z_0}\over 10L}({z_0})
\subset \Omega^+(\Kphi)$;

\item \label{cases-HolderCorner-2}
$B_{d_{z_0}\over 2L}({z_0})\cap\Omega^+(\Kphi)$ for  ${z_0}\in
\sonic\cap B_\vr(\PtUpL)$;

\item \label{cases-HolderCorner-3}
$B_{d_{z_0}\over 2L}({z_0})\cap\Omega^+(\Kphi)$ for  ${z_0}\in
\shock\cap B_\vr(\PtUpL)$.
\end{enumerate}

We discuss only case (\ref{cases-HolderCorner-3}), since the other
cases are simpler and can be handled similarly.

Let ${z_0}\in \shock\cap B_\vr(\PtUpL)$. Denote $\hat
d={d_{z_0}\over 2L}>0$. Without loss of generality, we can assume
that $\hat d\le 1$.

We rescale $z=(x,y)$ near $z_0$:
$$
Z=(X, Y):={1\over \hat d}(x-x_0, y-y_0).
$$
Since $B_{\hat d}(z_0)\cap (\partial \Omega^+(\Kphi)\setminus
\shock)=\emptyset$, then, for $\rho\in (0, 1)$, the domain obtained
by rescaling $\Omega^+(\Kphi)\cap B_{\rho\hat d}(z_0)$ is
$$
\hat\Omega^{z_0}_\rho\defd B_\rho\cap \big\{Y<\hat F(X)\defd {\hat
f_\Kphi(x_0+ \hat d X)-\hat f_\Kphi(x_0)\over \hat d}\big\},
$$
where $\hat f_\Kphi$ is the function in
(\ref{domain-in-rescaled-lemma}). Note that $y_0=\hat f_\Kphi(x_0)$ since
$(x_0, y_0)\in\shock$. Since $L\ge 1$, we have
$$
\|\hat F\|_{C^{2,\alpha}([-1,1])}\le \|\hat f_\Kphi\|^{(-1-\alpha,
\{0\})}_{2,\alpha,\bR_+}
$$
and $\|\hat f_\Kphi\|^{(-1-\alpha, \{0\})}_{2,\alpha,\bR_+}$ is
estimated in terms of the data by (\ref{holder-hat-f}).

\medskip
Define
\begin{equation}\label{defV-mixed-corner}
v(Z)={1\over \hat d^{1+\alpha}}\psi(z_0+\hat dZ)\qquad \mbox{for }\;
Z\in \hat\Omega^{z_0}_1.
\end{equation}
Then
\begin{equation}\label{V-mixed-corner-Linfty}
\|v\|_{L^\infty(\hat\Omega^{z_0}_1)}\le C
\end{equation}
by (\ref{growth-est-P4-for-lin-unif-ellipt}) with $C$ depending only
on the data.

Since $\psi$ satisfies equation
(\ref{unif-ellipt-linear-xy-iterationEquation}) in
$\Omega^+(\Kphi)\cap\DomS_{4\varepsilon}$ and the oblique derivative
condition (\ref{iterationRH-lf-flattened}) on
$\shock\cap\overline{\DomS_{4\varepsilon}}$, then $v$ satisfies an
equation and an oblique derivative condition of the similar form in
$\hat\Omega^{z_0}_1$ and on $\partial\hat\Omega^{z_0}_1\cap\{Y=\hat
F(X)\}$, respectively, whose coefficients satisfy properties
(\ref{estCoefsIterRH-flattened}) and
(\ref{ellipticityOfIterEq-Linear-xy-0}) with the same constants as
for the original equations, where we have used $\hat d\le 1$ and the
$C^{\alpha/2}$--estimates of the coefficients of the equation
depending only on the data, $\delta$, and $\hat\psi$. Then, from the
standard local estimates for linear oblique derivative problems, we
have
$$
\|v\|_{C^{2,{\alpha/2}}(\overline{\hat{\Omega}^{z_0}_{1/2}})}\le C,
$$
with $C$ depending only on the data, $\delta$, and $\hat\psi$.

We obtain similar estimates for cases
(\ref{cases-HolderCorner-1})--(\ref{cases-HolderCorner-2}), by using
the interior estimates for elliptic equations for case
(\ref{cases-HolderCorner-1}) and the local estimates for the
Dirichlet problem for linear elliptic equations for case
(\ref{cases-HolderCorner-2}).

Writing the above estimates in terms of $\psi$ and using the fact
that the whole domain $\Omega^+(\Kphi)\cap B_\vr(\PtUpL)$ is covered
by the subdomains in
(\ref{cases-HolderCorner-1})--(\ref{cases-HolderCorner-3}), we
obtain (\ref{Hoder-est-P4-for-lin-unif-ellipt}) by an argument
similar to the proof of  \cite[Theorem 4.8]{GilbargTrudinger} (see
also the proof of Lemma \ref{partIntSeminorm-est-lemma} below).
\Endproof
%*************** END PROOF *****LEMMA**********

%*****************BEGIN LEMMA*************
\begin{lemma}\label{existence-Lin-UnifEllipt-Lemma}
There exist $\hat C>0$ and $\delta_0\in (0,1)$ depending only on the
data such that, if $\epsP, \varepsilon>0$ and $M_1, M_2\ge 1$ in
{\rm (\ref{defSetK_R})} satisfy {\rm (\ref{condConst-00})}, and
$\delta\in (0,\delta_0)$, there exists a unique solution $\psi\in
C^{(-1-\alpha, \corners)}_{2,{\alpha/2},\Omega^+(\Kphi)}$ of
{\rm (\ref{unif-ellipt-linear-iterationEquation})} and {\rm
(\ref{iterationRH})}--{\rm (\ref{iterationCondOnSymmtryLine})}. The
solution $\psi$ satisfies {\rm
(\ref{L-infty-for-unif-ellipt})}--{\rm
(\ref{barier-for-unif-ellipt})}.
\end{lemma}
%*****************END LEMMA*************

%*************** BEGIN PROOF ***************
\Proof
In this proof, for simplicity, we write $\Omega^+$ for
$\Omega^+(\Kphi)$ and denote by $\Gamma_1$, $\Gamma_2$, $\Gamma_3$,
and $\Gamma_D$ the relative interiors of the curves $\shock(\Kphi)$,
$\Sigma_0(\Kphi)$, $\wedgeB$, and $\sonic$ respectively.

We first prove the existence of a solution for a general problem
${\mathcal P}$ of the form
$$
\sum_{i,j=1}^2 a_{ij}D^2_{ij} \psi=f\;\;\mbox{in }\Omega^+;\quad
\sum_{i=1}^2 b^{(k)}_{i}D_i \psi=g_i\;\;\mbox{on }\Gamma_k,\;
k=1,2,3;\quad \psi=0\;\;\mbox{on }\Gamma_D,
$$
where the equation is uniformly elliptic in $\Omega^+$ and the
boundary conditions on $\Gamma_k$, $k=1,2,3,$ are uniformly oblique,
i.e., there exist constants $\lambda_1,\lambda_2, \lambda_3>0$ such
that \begin{eqnarray*}
&& \lambda_1|\mu|^2 \le
 \sum_{i,j=1}^2
a_{ij}(\mxx,\mxy)\mu_i\mu_j\le \lambda_1^{-1}|\mu|^2\qquad \mbox{for
all} \,\, (\mxx,\mxy)\in\Omega^+, \mu\in\bR^2,\\
&& \sum_{i=1}^2 b^{(k)}_{i}(\mxx,\mxy)\nu_i(\mxx,\mxy)\ge\lambda_2,\\
&&\displaystyle \left|{(b^{(k)}_1, b^{(k)}_2)\over |(b^{(k)}_1,
b^{(k)}_2)|} (P_k) -{(b^{(k-1)}_1, b^{(k-1)}_2)\over |(b^{(k-1)}_1,
b^{(k-1)}_2)|} (P_k)\right|\ge \lambda_3 \qquad\mbox{for }\;k=2,3,
\end{eqnarray*}
and  $\|a_{ij}\|_{C^\alpha(\overline{\Omega^+})}+
\|b^{(k)}_{i}\|_{C^{1,\alpha}(\overline{\Gamma_k})}\le L$ for some
$L>0$.

First we derive an apriori estimate of a solution of problem
${\mathcal P}$. For that, we define the following norm for $\psi\in
C^{k,\beta}(\Omega^+)$, $k=0,1,2,\dots$, and $\beta\in (0,1)$:
$$
\|\psi\|_{*,k,\beta}:=\sum_{i=2}^3 \|\psi\|^{-k+1-\beta,
\{P_i\}}_{k,\beta, B_{2\vr}(P_i)\cap \Omega^+} + \sum_{i=1,4}
\|\psi\|^{-k+2-\beta, \{P_i\}}_{k,\beta, B_{2\vr}(P_i)\cap
\Omega^+}+\|\psi\|_{C^{k,\beta}(\overline{
\Omega^+\setminus(\cup_{i=1}^4 B_{\vr}(P_i))})},
$$
where $\vr>0$ is chosen small so that the balls $B_{2\vr}(P_i)$ for
$i=1,\dots,4$ are disjoint. Denote $C^{*,k,\beta}:=\{\psi\in
C^{*,k,\beta}\; : \;\|\psi\|_{*,k,\beta}<\infty\}$. Then
$C^{*,k,\beta}$ with norm $\|\cdot\|_{*,k,\beta}$ is a Banach space.
Similarly, define
$$
\|g_k\|_{*,\beta}:=\sum_{i=2}^3 \|g_k\|^{-\beta, \{P_i\}}_{k,\beta,
B_{2\vr}(P_i)\cap \Gamma_k} + \sum_{i=1,4} \|g_k\|^{1-\beta,
\{P_i\}}_{k,\beta, B_{2\vr}(P_i)\cap
\Gamma_k}+\|g_k\|_{C^{1,\beta}(\overline{
\Gamma_k\setminus(\cup_{i=1}^4 B_{\vr}(P_i))})},
$$
where the respective terms are zero if $B_{2\vr}(P_i)\cap
\Gamma_k=\emptyset$. Using the regularity of boundary of $\Omega^+$,
from the localized version of the estimates of \cite[Theorem
2]{Lieberman} applied in $B_{2r}(P_i)\cap \Omega^+$, $i=1,4$, and of
the estimates of \cite[Lemma 1.3]{Lieberman88} applied in
$B_{2r}(P_i)\cap \Omega^+$, $i=2,3$, and the standard local
estimates for the  Dirichlet and oblique derivative problems of
elliptic equations in smooth domains applied similarly to Step 4 in
the proof of Lemma \ref{C2alpha-near-P4-Lin-UnifEllipt-Lemma}, we
obtain that there exists $\beta=\beta(\Omega^+, \lambda_2,
\lambda_3)\in (0,1)$ such that any solution $\psi\in
C^{\beta}(\overline{\Omega^+}) \cap C^{1,\beta}(\overline{\Omega^+}
\setminus\overline\Gamma_D) \cap C^2(\Omega^+)$ of problem
${\mathcal P}$ satisfies
\begin{equation}\label{LiebermanEst-1}
\|\psi\|_{*,2,\beta}\le
C\big(\|f\|_{*,0,\beta}+\sum_{k=1}^3\|g_k\|_{*,\beta}
+\|\psi\|_{0,\Omega^+}\big)
\end{equation}
for $C=C(\Omega^+, \lambda_1,  \lambda_2,  \lambda_3, L)$. Next, we
show that $\psi$ satisfies
\begin{equation}\label{LiebermanEst-2}
\|\psi\|_{*,2,\beta}\le
C(\|f\|_{*,0,\beta}+\sum_{k=1}^3\|g_k\|_{*,\beta})
\end{equation}
for $C=C(\Omega^+, \lambda_1,  \lambda_2,  \lambda_3, L)$. By
(\ref{LiebermanEst-1}), it suffices to estimate
$\|\psi\|_{0,\Omega^+}$ by the right-hand side of
(\ref{LiebermanEst-2}). Suppose such an estimate is false. Then
there exists a sequence of problems ${\mathcal P}^m$ for
$m=1,2,\dots$ with coefficients $a_{ij}^m$ and $b_i^{(k),m}$, the
right-hand sides $f^m$ and $g_k^m$, and solutions $\psi^m\in
C^{*,2,\beta}$, where the assumptions on  $a_{ij}^m$ and
$b_i^{(k),m}$ stated above are satisfied with uniform constants
$\lambda_1, \lambda_2,  \lambda_3$, and $L$, and
$\|f^m\|_{*,0,\beta}+\sum_{k=1}^3\|g^m_k\|_{*,\beta}\to 0$ as
$m\to\infty$, but $\|\psi^m\|_{0,\Omega^+}=1$ for $m=1,2,\dots$.
Then, from (\ref{LiebermanEst-1}), we obtain
$\|\psi^m\|_{*,2,\beta}\le C$ with $C$ independent of $m$. Thus,
passing to a subsequence (without change of notation), we find
$a_{ij}^m\to a_{ij}^0$ in $C^{\beta/2}(\overline{\Omega^+})$,
$b_i^{(k),m}\to b_i^{(k),0}$ in
$C^{1,\beta/2}(\overline{\Gamma_k})$, and $\psi^m\to\psi^0$ in
$C^{*,2,\beta/2}$, where $\|\psi^0\|_{0,\Omega^+}=1$, and $a_{ij}^0$
and $b_i^{(k),0}$ satisfy the same ellipticity, obliqueness, and
regularity conditions as $a_{ij}^m$ and $b_i^{(k),m}$. Moreover,
$\psi^0$ is a solution of the homogeneous Problem ${\mathcal P}$
with coefficients  $a_{ij}^0$ and $b_i^{(k),0}$. Since
$\|\psi^0\|_{0,\Omega^+}=1$, this contradicts the uniqueness of a
solution in $C^{*,2,\beta}$ of problem ${\mathcal P}$ (the
uniqueness for problem ${\mathcal P}$ follows by the same argument
as in Lemma \ref{comparisonPrincipleOfUnifEllipt-Lemma}). Thus
(\ref{LiebermanEst-2}) is proved.

Now we show the existence of a solution for  problem ${\mathcal P}$
if $\hat C$ in (\ref{condConst-00}) is sufficiently large. We first
consider problem ${\mathcal P}_0$ defined as follows:
$$
\Delta \psi=f\;\;\mbox{in }\Omega^+;\quad D_\nu\psi=g_k\;\;\mbox{on
}\Gamma_k,\; k=1,2,3;\quad \psi=0\;\;\mbox{on }\Gamma_D.
$$
Using the fact that $\Gamma_2$ and $\Gamma_3$ lie on $\eta=0$ and
$\eta=\xi\tan\theta_w$ respectively, and using
(\ref{angleCloseToPiOver2}) and (\ref{OmegaPL-f-higher}), it is easy
to construct a diffemorphism
$$
F: \,\, \Omega^+\to Q:=\{(X, Y)\in(0,1)^2\}
$$
satisfying
\begin{eqnarray*}
&&\|F\|_{C^{1,\alpha}(\overline\Omega^+)}\le C,\qquad
\|F^{-1}\|_{C^{1,\alpha}(\overline Q)}\le C,\\
&&F(\Gamma_D)=\Sigma_D\defd\{X=1, Y\in(0,1)\},
\end{eqnarray*}
and
\begin{equation}\label{diffenorm}
\|DF^{-1}-Id\|_{C^\alpha(Q\cap\{X<\eta_1/2\})}\le
C\varepsilon^{1/4},
\end{equation}
where $C$ depends only on the data, and $(\xi_1,\eta_1)$ are the
coordinates of $\PtUpL$ defined by (\ref{coord-P4}) with $\eta_1>0$.
The mapping $F$ transforms problem ${\mathcal P}_0$ into the
following problem $\tilde{\mathcal P}_0$:
\begin{eqnarray*}
&&\sum_{i,j=1}^2 D_i(\tilde a_{ij}D_j u)=\tilde f\;\qquad\mbox{in
}Q;\\
&&\sum_{i,j=1}^2 \tilde a_{ij}D_ju\,\nu_i=\tilde g_k\;\qquad\mbox{on
}I_k,\; k=1,2,3;\\
&&u=0\;\qquad \mbox{on }\Sigma_D,
\end{eqnarray*}
where $I_k=F(G_k)$ are the respective sides of $\partial Q$, $\nu$
is the unit normal on $I_k$, $\|\tilde a_{ij}\|_{C^\alpha(\overline
Q)}\le C$, and $\tilde a_{ij}$ satisfy the uniform ellipticity in
$\overline Q$ with elliptic constant $\tilde\lambda>0$. Using
(\ref{diffenorm}), we obtain
\begin{equation}\label{coeffsConorm}
\|\tilde a_{ij}-\delta_i^j\|_{C^\alpha(Q\cap\{X<\eta_1/2\})}\le
C\varepsilon^{1/4},
\end{equation}
where $\delta_i^i=1$ and $\delta_i^j=0$ for $i\ne j$, and $C$
depends only on the data. If  $\varepsilon>0$ is sufficiently small
depending only on the data, then, by \cite[Theorem 3.2, Proposition
3.3]{ChenFeldman3},
 there exists
$\beta\in(0,1)$ such that, for any $\tilde f\in
C^\beta(\overline{Q})$ and $\tilde g_k\in C^{\beta}(\overline{I_k})$
with $k=1,2,3$, there exists a unique weak solution $u\in H^1(Q)$ of
problem $\tilde {\mathcal P}_0$, and this solution satisfies $u\in
C^{\beta}(\overline Q) \cap C^{1,\beta}(\overline Q
\setminus\overline\Sigma_D)$. We note that, in \cite[Theorem 3.2,
Proposition 3.3]{ChenFeldman3}, condition (\ref{coeffsConorm}) is
stated in the whole $Q$, but in fact this condition was used only in
a neighborhood of $I_2=\{0\}\times(0,1)$, i.e., the results can be
applied to the present case. We can assume that $\beta\le \alpha$.
Then, mapping back to $\Omega^+$, we obtain the existence of a
solution $\psi\in C^{\beta}(\overline{\Omega^+}) \cap
C^{1,\beta}(\overline{\Omega^+} \setminus\overline\Gamma_D) \cap
C^2(\Omega^+)$ of problem ${\mathcal P}_0$ for any $f\in
C^\beta(\overline{\Omega^+})$ and $g_k\in
C^{\beta}(\overline{\Gamma_k})$, $k=1,2,3$. Now, reducing $\beta$ if
necessary and using (\ref{LiebermanEst-2}), we conclude that, for
any $(f, g_1, g_2, g_3)\in{\mathcal Y}^\beta\defd \{(f, g_1, g_2,
g_3)\;: \;
\|f\|_{*,0,\beta}+\sum_{k=1}^3\|g_k\|_{*,\beta}<\infty\}$, there
exists a unique solution $\psi\in C^{*,2,\beta}$ of problem
${\mathcal P}_0$, and $\psi$ satisfies (\ref{LiebermanEst-2}).

Now the existence of a unique solution $\psi\in C^{*,2,\beta}$ of
 problem ${\mathcal P}$, for any $(f, g_1, g_2, g_3)\in{\mathcal Y}^\beta$
with  sufficiently small $\beta\in (0,1)$, follows by the method of
continuity, applied to the family of problems $t{\mathcal
P}+(1-t){\mathcal P}_0$ for $t\in[0,1]$. This proves the existence
of a solution $\psi\in C^{*,2,\beta}$ of problem {\rm
(\ref{unif-ellipt-linear-iterationEquation})} and {\rm
(\ref{iterationRH})}--{\rm (\ref{iterationCondOnSymmtryLine})}.

%%%%%%%%%%%%%%%%%%%%%%%%%%%%%%%%%%%%%%%%%%%%%%%%%%%%%%%%%%%%%%%%%%%%%%%%%%%%

Estimates
(\ref{L-infty-for-unif-ellipt})--(\ref{barier-for-unif-ellipt}) then
follow from Lemma \ref{unifEstOfUnifEllipt-Lemma}. The higher
regularity $\psi\in C^{(-1-\alpha,
\corners)}_{2,{\alpha/2},\Omega^+(\Kphi)}$ follows from Lemmas
\ref{C2alpha-Lin-UnifEllipt-Lemma}--\ref{C2alpha-near-P4-Lin-UnifEllipt-Lemma}
and the standard estimates for the Dirichlet problem near the flat
boundary, applied in a neighborhood of
$\sonic\setminus(B_{\vr/2}(\PtUpL)\cup B_{\vr/2}(\PtUpR))$ in the
$(x,y)$--coordinates, where $\vr>0$ may be smaller than the constant
$\vr$ in Lemmas
\ref{wedge-sonic-Lin-UnifEllipt-Lemma}--\ref{C2alpha-near-P4-Lin-UnifEllipt-Lemma}.
In fact, from Lemma \ref{wedge-sonic-Lin-UnifEllipt-Lemma},  we
obtain even a higher regularity than that in the statement of Lemma
\ref{existence-Lin-UnifEllipt-Lemma}: $\psi\in
C^{(-1-\alpha,\{\PtLwL,\PtLwR,\PtUpR\})}_{2,{\alpha/2},\Omega^+(\Kphi)}$.
The uniqueness of solutions follows from the Comparison Principle
(Lemma \ref{comparisonPrincipleOfUnifEllipt-Lemma}).
\Endproof
%*************** END PROOF *****LEMMA**********

Lemma \ref{existence-Lin-UnifEllipt-Lemma} justifies the definition
of map $\hat{J}$ in \eqref{6.9a} defined by
$\hat{J}(\hat\psi):=\psi$. In order to apply the Leray-Schauder
Theorem, we make the following apriori estimates for solutions of
the nonlinear equation.

%*****************BEGIN LEMMA*************
\begin{lemma}\label{estimates-nonlin-UnifEllipt}
There exist $\hat C>0$ and $\delta_0\in (0,1)$ depending only on the
data such that the following holds. Let $\epsP, \varepsilon>0$ and
$M_1, M_2\ge 1$ in {\rm (\ref{defSetK_R})} satisfy {\rm
(\ref{condConst-00})}. Let $\delta\in (0,\delta_0)$ and
$\mu\in[0,1]$. Let $\psi\in C^{(-1-\alpha,
\corners)}_{2,{\alpha/2},\Omega^+(\Kphi)}$ be a solution of
{\rm (\ref{unif-ellipt-iterationEquation})}, {\rm
(\ref{iterationRH})}--{\rm (\ref{iterationCondOnWedge})},
 and
\begin{equation}
\psi_\mxy=-\mu v_2\qquad\mbox{on }\;
\Sigma_0(\Kphi)\defd\partial \Omega^+(\Kphi)\cap\{\mxy=-v_2\}.
\label{iterationCondOnSymmtryLine-LS}
\end{equation}
Then

\begin{enumerate}\renewcommand{\theenumi}{\roman{enumi}}

\item \label{estimates-nonlin-UnifEllipt-i1}
There exists $C>0$ independent of $\psi$ and $\mu$ such that
\begin{equation*}
\|\psi\|_{C^{1,\alpha}(\overline{\Omega^+(\Kphi)})}\le C;
%\label{Hoder-est-for-unif-ellipt-fixed-pt}
\end{equation*}

\item \label{estimates-nonlin-UnifEllipt-i2}
$\psi$ satisfies {\rm (\ref{L-infty-for-unif-ellipt})}--{\rm
(\ref{barier-for-unif-ellipt})}  with constant $C$ depending only on
the data;

\item \label{estimates-nonlin-UnifEllipt-i2.1}
$\psi\in C^{(-1-\alpha, \corners)}_{2,\alpha,\Omega^+(\Kphi)}$.
Moreover, for every  $s\in (0, {c_2/2})$, estimate {\rm
(\ref{Hoder-est-for-unif-ellipt})} holds with constant $C$ depending
only on the data and $s$;

\item \label{estimates-nonlin-UnifEllipt-i3}
Solutions  of problem {\rm (\ref{unif-ellipt-iterationEquation})},
{\rm (\ref{iterationRH})}--{\rm (\ref{iterationCondOnWedge})}, and
{\rm (\ref{iterationCondOnSymmtryLine-LS})} satisfy the following
comparison principle: Denote by $\Nl_\delta(\psi)$, $B_1(\psi)$,
$B_2(\psi)$, and $B_3(\psi)$ the left-hand sides of {\rm
(\ref{unif-ellipt-iterationEquation})}, {\rm (\ref{iterationRH})},
{\rm (\ref{iterationCondOnWedge})}, and {\rm
(\ref{iterationCondOnSymmtryLine-LS})} respectively. If
$\psi_1,\psi_2\in
C^{(-1-\alpha,\corners)}_{2,\alpha,\Omega^+(\Kphi)}$ satisfy
\begin{eqnarray*}
&&\Nl_\delta(\psi_1)\le\Nl_\delta(\psi_2)
  \qquad\mbox{ in } \Omega^+(\Kphi),\\
&&B_k(\psi_1)\le B_k(\psi_2)
  \qquad \mbox{on } \shock(\Kphi),\,
   \wedgeB, \,\mbox{and }\,\Sigma_0(\Kphi)
   \mbox{ for } k=1,2,3,\\
&&\psi_1\ge \psi_2\qquad\mbox{ on } \,
  \sonic,
\end{eqnarray*}
then
$$
\psi_1\ge \psi_2 \qquad\mbox{ in }\, \Omega^+(\Kphi).
$$
In particular, problem {\rm (\ref{unif-ellipt-iterationEquation})},
{\rm (\ref{iterationRH})}--{\rm (\ref{iterationCondOnWedge})}, and
{\rm (\ref{iterationCondOnSymmtryLine-LS})} has at most one solution
$\psi\in C^{(-1-\alpha, \corners)}_{2,\alpha,\Omega^+(\Kphi)}$.
\end{enumerate}
\end{lemma}
%***************END LEMMA***************
%*************** BEGIN PROOF ***************
\Proof The proof consists of six steps.

{\em Step 1.} Since a solution $\psi\in C^{(-1-\alpha,
\corners)}_{2,\alpha,\Omega^+(\Kphi)}$ of {\rm
(\ref{unif-ellipt-iterationEquation})}, {\rm
(\ref{iterationRH})}--{\rm (\ref{iterationCondOnWedge})}, and
(\ref{iterationCondOnSymmtryLine-LS}) with $\mu\in[0,1]$ is the
solution of the linear problem for equation
(\ref{unif-ellipt-linear-iterationEquation}) with
$\hat\psi\defd\psi$ and boundary conditions
(\ref{iterationRH})--(\ref{iterationCondOnWedge}) and
(\ref{iterationCondOnSymmtryLine-LS}). Thus, estimates
(\ref{L-infty-for-unif-ellipt})--(\ref{barier-for-unif-ellipt}) with
constant $C$ depending only on the data follow directly from Lemma
\ref{unifEstOfUnifEllipt-Lemma}.

\medskip
{\em Step 2.} Now, from Lemma
\ref{propertiesNonlinCoeffs}(\ref{propertiesNonlinCoeffs-i2}),
equation (\ref{unif-ellipt-iterationEquation}) is linear in
$\Omega^+(\Kphi)\cap \{c_2-r>4\varepsilon\}$, i.e.,
(\ref{unif-ellipt-iterationEquation}) is
(\ref{unif-ellipt-linear-iterationEquation}) in $\Omega^+(\Kphi)\cap
\{c_2-r>4\varepsilon\}$, with coefficients
$a_{ij}(\mxx,\mxy)=A^1_{ij}(\mxx,\mxy)+\delta\delta_{ij}$ for $A^1_{ij}$
defined by (\ref{iterationUniforDomEquation}). Then, by Lemma
\ref{propertiesNonlinCoeffs}(\ref{propertiesNonlinCoeffs-i2}),
$a_{ij}\in C^\alpha(\overline{\Omega^+(\Kphi)\cap
\{c_2-r>4\varepsilon\}})$ with the norm estimated in terms of the
data. Also, $\shock(\Kphi)$ and the coefficients $b_i$ of
(\ref{iterationRH-lf}) satisfy (\ref{OmegaPL-f-higher}) and
(\ref{estCoefsIterRH-0})--(\ref{estCoefsIterRH}). Then, repeating
the proof of Lemma \ref{C2alpha-Lin-UnifEllipt-Lemma} with the use
of the $L^\infty$ estimates of $\psi$ obtained in Step 1 of the
present proof, we conclude that $\psi\in C^{(-1-\alpha,
\{\PtLwL,\PtLwR\})}_{2,\alpha,\Omega^+(\Kphi) \cap
\{c_2-r>6\varepsilon\}}$ with
\begin{equation}
\|\psi\|^{(-1-\alpha, \{\PtLwL,\PtLwR\})}_{2,\alpha,\Omega^+(\Kphi)
\cap \{c_2-r>6\varepsilon\}}\le C\epsP
\label{Hoder-est-for-nonlin-unif-ellipt}
\end{equation}
for $C$ depending only on the data.

\medskip
{\em Step 3.} Now we prove (\ref{Hoder-est-for-unif-ellipt}) for all
$s\in (0, {c_2/2})$. If $s\ge 6\varepsilon$, then
(\ref{Hoder-est-for-unif-ellipt}) follows from
(\ref{Hoder-est-for-nonlin-unif-ellipt}). Thus, it suffices to
consider the case $s\in (0,  6\varepsilon)$ and show that
\begin{equation}
\|\psi\|_{C^{2,\alpha}(\overline{\Omega^+(\Kphi)\cap
\{{s/2}<c_2-r<6\varepsilon+{s/4}\}})}\le C(s)\epsP,
\label{Hoder-est-for-nonlin-unif-ellipt-s}
\end{equation}
with $C$ depending only on the data and $s$. Indeed,
(\ref{Hoder-est-for-nonlin-unif-ellipt})--(\ref{Hoder-est-for-nonlin-unif-ellipt-s})
imply (\ref{Hoder-est-for-unif-ellipt}).

In order to prove (\ref{Hoder-est-for-nonlin-unif-ellipt-s}), it
suffices to prove the existence of $C(s)$ depending only on the data
and $s$ such that
\begin{equation}\label{nonlin-est-unif-balls-inter}
\|\psi\|_{C^{2,\alpha}(\overline{B_{s/16}(z)})}\le
 C(s)\|\psi\|_{L^\infty(B_{s/8}(z))}
\end{equation}
for all $z:=(\mxx,\mxy)\in \Omega^+(\Kphi)\cap
\{{s/2}<c_2-r<6\varepsilon+{s/4}\}$ with $\dist(z,
\partial\Omega^+(\Kphi))>{s/8}$ and that
\begin{equation}\label{nonlin-est-unif-balls-oblique}
\|\psi\|_{C^{2,\alpha}(\overline{B_{s/8}(z)\cap\Omega^+(\Kphi)})}
\le C(s)\|\psi\|_{L^\infty(B_{s/4}(z)\cap\Omega^+(\Kphi))}
\end{equation}
for all
$z\in(\shock(\Kphi)\cup\wedgeB)\cap\{{s/2}<c_2-r<6\varepsilon+{s/4}\}$.
Note that all the domains in (\ref{nonlin-est-unif-balls-inter}) and
(\ref{nonlin-est-unif-balls-oblique}) lie within
$\Omega^+(\Kphi)\cap \{{s/4}<c_2-r<12\varepsilon\}$. We can assume
that $\varepsilon < c_2/24$. Since equation
(\ref{unif-ellipt-iterationEquation}) is uniformly elliptic in
$\Omega^+(\Kphi)\cap \{{s/4}<c_2-r<12\varepsilon\}$ by Lemma
\ref{propertiesNonlinCoeffs}(\ref{propertiesNonlinCoeffs-i1}), and
the boundary conditions (\ref{iterationRH}) and
(\ref{iterationCondOnWedge}) are linear and oblique with
$C^{1,\alpha}$--coefficients estimated in terms of the data, then
(\ref{nonlin-est-unif-balls-inter}) follows from Theorem
\ref{locEstElliptEq}
and (\ref{nonlin-est-unif-balls-oblique}) follows from Theorem
\ref{locEstElliptEq-oblique} (in Appendix \ref{append-1-section}).
Since $\|\psi\|_{L^\infty(\Omega^+(\varphi))}\le 1$ by
(\ref{L-infty-for-unif-ellipt}), the constants in the local
estimates depend only on the ellipticity, the constants in Lemma
\ref{propertiesNonlinCoeffs}(\ref{propertiesNonlinCoeffs-i3}), and,
for the case of (\ref{nonlin-est-unif-balls-oblique}), also on the
$C^{2,\alpha}$--norms of the boundary curves and the obliqueness and
$C^{1,\alpha}$--bounds of the coefficients in the boundary
conditions (which, for condition (\ref{iterationRH}), follow from
(\ref{OmegaPL-f-higher}) and (\ref{estCoefsIterRH-0}) since our
domain is away from the points $\PtUpL$ and $\PtLwL$). All these
quantities depend only on the data and $s$. Thus, the constant
$C(s)$ in
(\ref{nonlin-est-unif-balls-inter})--(\ref{nonlin-est-unif-balls-oblique})
depends only on the data and $s$.

\medskip
{\em Step 4.} In this step, the universal constant $C$ depends only
on the data and $\delta$, unless specified otherwise. We prove that
$\psi\in C^{2,\alpha}(\overline{B_\vr(\PtUpR)\cap\Omega^+(\Kphi)})$
for sufficiently small $\vr>0$, depending only on the data and
$\delta$, and
%the estimate holds
\begin{equation}
\|\psi\|_{C^{2,\alpha}(\overline{B_\vr(\PtUpR)\cap\Omega^+(\Kphi)})}\le
C. \label{Hoder-est-wedge-sonic-nonLin-unif-ellipt}
\end{equation}

We follow the proof of Lemma \ref{wedge-sonic-Lin-UnifEllipt-Lemma}.
Since $B_\vr(\PtUpR)\cap\Omega^+(\Kphi)\subset \DomS$ for small
$\vr$, we work in the $(x,y)$--coordinates. We use the notations
$B^+_\vr$ and $B^{++}_\vr$, introduced in Step 1 of Lemma
\ref{wedge-sonic-Lin-UnifEllipt-Lemma}, and consider the function
$$
v(x,y)={1\over \vr}\psi(\vr x, \vr y).
$$
Then, by (\ref{barier-for-unif-ellipt}), $v$ satisfies
\begin{equation}\label{L-infty-for-unif-ellipt-nonlin-v}
\|v\|_{L^\infty(B_2^{++})}\le 2C{\epsP\over \varepsilon}\le 1,
\end{equation}
where the last inequality holds if $\hat C$ in
(\ref{condConst-00}) is sufficiently large. Moreover, $v$ is a solution of
\begin{eqnarray}
&&\qquad\hat A_{11}^{(\vr)}v_{xx}+ 2\hat A_{12}^{(\vr)}v_{xy} +\hat
A_{22}^{(\vr)}v_{yy} +\hat A_1^{(\vr)}v_{x} +\hat A_2^{(\vr)}v_{y}
=0 \qquad\mbox{in }\;B_2^{++}, \label{unif-ellipt-NONlinear-xy-v}
\\
&&\qquad v=0\qquad\mbox{on }\;B_{2}\cap \{x=0, y>0\},
\label{condOnSonicLine-v-xy-loc-NONlin}
\\
&&\qquad v_\nu\equiv v_y=0\qquad\mbox{on }\;B_{2}\cap \{y=0, x>0\},
\label{condOnWedge-v-xy-loc-NONlin}
\end{eqnarray}
with $(A_{ij}^{(\vr)}, A_{i}^{(\vr)})=(A_{ij}^{(\vr)},
A_{i}^{(\vr)})(Dv, x, y)$, where we use
(\ref{unif-ellipt-iteration-equation-sonicStruct}) to find that, for
$(x,y)\in B_2^{++}$, $p\in\bR^2$, $i,j=1,2,$
\begin{equation}\label{def-rescaled-coef-v-NONlinear}
\begin{array}{l}
\displaystyle
\hat A_{11}^{(\vr)}(p,x,y)=\hat A_{11}(p,\vr x,\vr y)+\delta, \\
\displaystyle
\hat A_{12}^{(\vr)}(p,x,y)=\hat A_{21}^{(\vr)}(p,x,y)=\hat A_{12}(p,\vr x,\vr y),\\
\displaystyle
\hat A_{22}^{(\vr)}(p,x,y)=\hat A_{22}(p,\vr x,\vr y)+{\delta\over (c_2-\vr x)^2},
\\
\displaystyle\hat A_{1}^{(\vr)}(p, x,y)= \vr\hat A_{1}(p, \vr x,\vr
y)+{\delta\over (c_2-\vr x)},\quad \hat A_{2}^{(\vr)}(p, x,y)=
\vr\hat A_{2}(p, \vr x,\vr y),
\end{array}
\end{equation}
with $\hat A_{ij}$ and $\hat A_{i}$ as in Lemma
(\ref{propertiesNonlinCoeffs-xy}). Since $\vr\le 1$, $\hat
A_{ij}^{(\vr)}$ and $\hat A_{i}^{(\vr)}$ satisfy the assertions of
Lemma
\ref{propertiesNonlinCoeffs-xy}(\ref{propertiesNonlinCoeffs-xy-i1})--(\ref{propertiesNonlinCoeffs-xy-i3})
with the unchanged constants. Moreover, $\hat A_{11}^{(\vr)}$, $\hat
A_{22}^{(\vr)}$, and $\hat A_{1}^{(\vr)}$ satisfy the property in
Lemma
\ref{propertiesNonlinCoeffs-xy}(\ref{propertiesNonlinCoeffs-xy-i4-0}).
The property in Lemma
\ref{propertiesNonlinCoeffs-xy}(\ref{propertiesNonlinCoeffs-xy-i4})
is now improved to
\begin{equation}\label{estRescaledA-nonlin}
|(\hat A_{12}^{(\vr)},\hat A_{21}^{(\vr)},\hat
A_{2}^{(\vr)})(x,y)|\le C\vr |x|, \,\quad |D(\hat
A_{12}^{(\vr)},\hat A_{21}^{(\vr)}, \hat A_{2}^{(\vr)})(x,y)|\le
C|\vr x|^{1/2}.
\end{equation}

Combining the estimates in Theorems \ref{locEstElliptEq} and
\ref{locEstElliptEq-Dirichlet}--\ref{locEstElliptEq-oblique} with
the argument that has led to (\ref{Hoder-est-near-Br}), we have
\begin{equation}\label{Hoder-est-near-Br-nonlin}
\|v\|_{C^{2,\alpha}(\overline{B_{3/2}^{++}\setminus B_{1/2}^{++}})}
\le C,
\end{equation}
where $C$ depends only on the data and $\delta>0$ by
(\ref{L-infty-for-unif-ellipt-nonlin-v}), since $\hat
A_{ij}^{(\vr)}$ and $\hat A_{i}^{(\vr)}$ satisfy
(\ref{locEstElliptEq-i1-0})--(\ref{locEstElliptEq-i2-0}) with the
constants depending only on the data and $\delta$. In particular,
$C$ in (\ref{Hoder-est-near-Br-nonlin}) is independent of $\vr$.

We now use the domain $D^{++}$ introduced in Step 2 of the proof of
Lemma \ref{wedge-sonic-Lin-UnifEllipt-Lemma}. We prove that, for
 any $g\in C^{\alpha}(\overline{D^{++}})$ with
$
\|g\|_{C^{\alpha}(\overline{D^{++}})}\le 1,
$
there exists a unique solution $w\in
C^{2,\alpha}(\overline{D^{++}})$ of the problem:
\begin{eqnarray}
&&\hat A_{11}^{(\vr)}w_{xx} +\hat A_{22}^{(\vr)}w_{yy} +\hat
A_1^{(\vr)}w_{x} =g \qquad\mbox{in }\;D^{++},
\label{unif-ellipt-NONlinear-xy-w}
\\
&&w=0\qquad\mbox{on }\;\partial D^{++}\cap \{x=0, y>0\},
\label{condOnSonicLine-w-xy-loc-nonlin}
\\
&&w_\nu\equiv w_y=0\qquad\mbox{on }\;\partial D^{++}\cap \{x>0,
y=0\}, \label{condOnWedge-w-xy-loc-nonlin}
\\
&& w=v\qquad\mbox{on }\;\partial D^{++}\cap \{x>0, y>0\},
\label{dirichlet-w-nonlin}
\end{eqnarray}
with $(A_{ii}^{(\vr)}, A_{1}^{(\vr)})=(A_{ii}^{(\vr)},
A_{1}^{(\vr)})(Dw, x, y)$. Moreover, we show
\begin{equation}\label{Hoder-est-near-origin-nonlin}
\|w\|_{C^{2,{\alpha}}(\overline{D^{++}})}\le C,
\end{equation}
where $C$ depends only on the data and is independent of $\vr$. For
that, similar to  Step 2 of the proof of Lemma
\ref{wedge-sonic-Lin-UnifEllipt-Lemma}, we consider the even
reflection $D^+$ of the set $D^{++}$, and the even reflection of
$(v,g,\hat A_{11}^{(\vr)},\hat A_{22}^{(\vr)},\hat A_{1}^{(\vr)})$
from $\overline{B_2^{++}}$ to $\overline{B_2^+}$, without change of
notation, where the even reflection of $(\hat A_{11}^{(\vr)},\hat
A_{22}^{(\vr)}, \hat A_{1}^{(\vr)})$, which depends on $(p,x,y)$, is
defined by
$$
\hat A_{ii}^{(\vr)}(p,x,-y)=\hat A_{ii}^{(\vr)}(p,x,y), \,\,\, \hat
A_{1}^{(\vr)}(p,x,-y)=\hat A_{1}^{(\vr)}(p,x,y) \quad\,\,\mbox{for
}\;(x,y)\in \overline{B_2^{++}}.
$$

Also, denote by  $\hat v$ the restriction of (the extended) $v$ to
$\partial D^+$. It follows from
(\ref{condOnSonicLine-v-xy-loc-NONlin})--(\ref{condOnWedge-v-xy-loc-NONlin})
and (\ref{Hoder-est-near-Br-nonlin}) that $\hat v\in
C^{2,{\alpha}}(\partial D^+)$ with
\begin{equation}\label{Hoder-est-bdry-funct-nonlin}
\|\hat v\|_{C^{2,\alpha}(\partial D^+)}\le C,
\end{equation}
depending only on the data and $\delta$. Furthermore, the extended
$g$ satisfies $g\in C^{\alpha}(\overline{D^+})$ with
$\|g\|_{C^{\alpha}(\overline{D^+})}
=\|g\|_{C^{\alpha/2}(\overline{D^{++}})}\le 1$. The extended $\hat
A_{11}^{(\vr)}, \hat A_{22}^{(\vr)}$, and $\hat A_{1}^{(\vr)}$
satisfy (\ref{locEstElliptEq-i1-0})--(\ref{locEstElliptEq-i2-0}) in
$D^+$ with the same constants as the estimates satisfied by $A_{ii}$
and $A_i$ in $\Omega^+(\Kphi)$. We consider the Dirichlet problem
\begin{eqnarray}
&&\hat A_{11}^{(\vr)}w_{xx} +\hat A_{22}^{(\vr)}w_{yy} +\hat
A_1^{(\vr)}w_{x} =g \qquad\mbox{in }\;D^{+},
\label{unif-ellipt-NONlinear-xy-v-ext}
\\
&& w=\hat v\qquad\mbox{on }\;\partial D^{+},
\label{dirichlet-v-ext-nonlin}
\end{eqnarray}
with $(A_{ii}^{(\vr)}, A_{1}^{(\vr)}):=(A_{ii}^{(\vr)},
A_{1}^{(\vr)})(Dw, x, y)$. By the Maximum Principle,
$$
\|w\|_{L^\infty(D^+)}\le \|\hat v\|_{L^\infty(D^+)}.
$$
Thus, using
(\ref{Hoder-est-bdry-funct-nonlin}), we obtain an estimate of
$\|w\|_{L^\infty(D^+)}$. Now, using Theorems \ref{locEstElliptEq}
and \ref{locEstElliptEq-Dirichlet} and the estimates of
$\|g\|_{C^{\alpha}(\overline{D^+})}$ and $\|\hat
v\|_{C^{2,\alpha}(\partial D^+)}$ discussed above, we obtain the
a-priori estimate for the $C^{2,\alpha}$--solution $w$ of
(\ref{unif-ellipt-NONlinear-xy-v-ext})--(\ref{dirichlet-v-ext-nonlin}):
\begin{equation}\label{Hoder-est-near-Br-w-nonlin}
\|w\|_{C^{2,\alpha}(\overline{D^{+}})}\le C,
\end{equation}
where $C$ depends only on the data and $\delta$. Moreover, for every
$\hat w\in C^{1,{\alpha}}(\overline{D^+})$, the existence of a
unique solution $w\in C^{2,{\alpha}}(\overline{D^+})$ of the linear
Dirichlet problem, obtained by substituting $\hat w$ into the
coefficients of (\ref{unif-ellipt-NONlinear-xy-v-ext}), follows from
\cite[Theorem 6.8]{GilbargTrudinger}.
 Now, by a standard
application of the Leray-Schauder Theorem, there exists a unique
solution $w\in C^{2,{\alpha}}(\overline{D^+})$ of the Dirichlet
problem
(\ref{unif-ellipt-NONlinear-xy-v-ext})--(\ref{dirichlet-v-ext-nonlin})
which satisfies (\ref{Hoder-est-near-Br-w-nonlin}).

{}From the structure of equation
(\ref{unif-ellipt-NONlinear-xy-v-ext}), especially the fact that
$\hat A_{11}^{(\vr)}$, $\hat A_{22}^{(\vr)}$, and $\hat A_1^{(\vr)}$
are independent of $p_2$ by Lemma \ref{propertiesNonlinCoeffs-xy}
(\ref{propertiesNonlinCoeffs-xy-i4-0}), and from the symmetry of the
domain and the coefficients and right-hand sides obtained by the
even extension, it follows that $\hat w$, defined by $\hat
w(x,y)=w(x, -y)$, is also a solution of
(\ref{unif-ellipt-NONlinear-xy-v-ext})--(\ref{dirichlet-v-ext-nonlin}).
By uniqueness for problem
(\ref{unif-ellipt-NONlinear-xy-v-ext})--(\ref{dirichlet-v-ext-nonlin}),
we find $ w(x,y)=w(x, -y) $ in $D^+$. Thus, $w$ restricted to
$D^{++}$ is a solution of
(\ref{unif-ellipt-NONlinear-xy-w})--(\ref{dirichlet-w-nonlin}),
where (\ref{condOnSonicLine-w-xy-loc-nonlin}) follows from
(\ref{condOnSonicLine-v-xy-loc-NONlin}) and
(\ref{dirichlet-v-ext-nonlin}). Moreover,
(\ref{Hoder-est-near-Br-w-nonlin}) implies
(\ref{Hoder-est-near-origin-nonlin}).

The uniqueness of a solution $w\in
C^{2,{\alpha}}(\overline{D^{++}})$ of
(\ref{unif-ellipt-NONlinear-xy-w})--(\ref{dirichlet-w-nonlin})
follows from the Comparison Principle (Lemma
\ref{comparisonPrincipleOfUnifEllipt-Lemma}).

\medskip
Now we prove the existence of a solution $w\in
C^{2,{\alpha}}(\overline{D^{++}})$ of the problem:
\begin{equation}
\begin{array}{ll}
&\hat A_{11}^{(\vr)}w_{xx}+2\hat A_{12}^{(\vr)}w_{xy} +\hat
A_{22}^{(\vr)}w_{yy} +\hat A_1^{(\vr)}w_{x} +\hat A_2^{(\vr)}w_{y}
=0 \qquad\mbox{in }\;D^{++},
\\
&w=0\qquad\mbox{on }\;\partial D^{++}\cap \{x=0, y>0\},
%\label{condOnSonicLine-w-xy-loc-fullEq-nonlin}
\\
&w_\nu\equiv w_y=0\qquad\mbox{on }\;\partial D^{++}\cap \{y=0,
x>0\}, %\label{condOnWedge-w-xy-loc-fullEq-nonlin}
\\
& w=v\qquad\mbox{on }\;\partial D^{++}\cap \{x>0, y>0\},
%\label{dirichlet-w-fullEq-nonlin}
\end{array}
\label{unif-ellipt-NONlinear-xy-w-fullEq}
\end{equation}
where $(A_{ij}^{(\vr)}, A_{i}^{(\vr)}):=(A_{ij}^{(\vr)},
A_{i}^{(\vr)})(Dw, x, y)$. Moreover, we prove that $w$ satisfies
\begin{equation}\label{Hoder-est-near-Br-fullEq-nonlin}
\|w\|_{C^{2,{\alpha}}(\overline{D^{++}})}\le C
\end{equation}
for $C>0$ depending only on the data and $\delta$.

Let $N$ be chosen below. Define
\begin{equation}\label{defSetM}
\setS:=\left\{W\in C^{2,\alpha}(\overline{D^{++}})\,\,:\,\,
\|W\|_{C^{2,\alpha}(\overline{D^{++}})}\le N\right\}.
\end{equation}

We obtain such $w$ as a fixed point of the map $K:\setS \to \setS$
defined as follows (if $R$ is small and $N$ is large, as specified
below). For $W\in \setS$, define
\begin{equation}\label{defG-corner-nonlin}
g=-2\hat A_{12}^{(\vr)}(x,y)W_{xy}-\hat A_2^{(\vr)}(x,y)W_{y}.
\end{equation}
By (\ref{estRescaledA-nonlin}),
$$
\|g\|_{C^{\alpha}(\overline{D^{++}})}\le CN \sqrt{\vr}\le 1,
$$
if $\vr\le \vr_0$ with $\vr_0= {1\over CN^2}$, for $C$ depending
only on the data and $\delta$. Then, as we have proved above, there
exists a unique solution $w\in C^{2,\alpha}(\overline{D^{++}})$ of
(\ref{unif-ellipt-NONlinear-xy-w})--(\ref{dirichlet-w-nonlin}) with
$g$ defined by (\ref{defG-corner-nonlin}). Moreover, $w$ satisfies
(\ref{Hoder-est-near-origin-nonlin}). Then, if we choose $N$ to be
the constant $C$ in (\ref{Hoder-est-near-origin-nonlin}), we get
$w\in\setS$. Thus, $N$ is chosen depending only on the data and
$\delta$. Now our choice of $\vr_0= {1\over CN^2}$ and $\vr\le
\vr_0$ (and the other smallness conditions stated above) determines
$\vr$ in terms of the data and $\delta$. We define $K[W]:=w$ and
thus obtain $K:\setS \to \setS$.

\medskip
Now the existence of a fixed point of $K$ follows from the Schauder
Fixed Point Theorem in the following setting: From its definition,
$\setS$ is a compact and  convex subset in
$C^{2,{\alpha/2}}(\overline{D^{++}})$. The map $K:\setS \to \setS$
is continuous in $C^{2,{\alpha/2}}(\overline{D^{++}})$: Indeed, if
$W_k\in \setS$ for $k=1,\dots$, and $W_k\to W$ in
$C^{2,{\alpha/2}}(\overline{D^{++}})$, then it is easy to see that
$W\in \setS$. Define $g_k$ and $g$ by  (\ref{defG-corner-nonlin})
for $W_k$ and $W$ respectively. Then $g_k\to g$ in
$C^{{\alpha/2}}(\overline{D^{++}})$ since $(\hat A_{12}, \hat A_{2})
=(\hat A_{12}, \hat A_{2})(x,y)$  by Lemma
\ref{propertiesNonlinCoeffs-xy}(\ref{propertiesNonlinCoeffs-xy-i4}).
Let $w_k=K[W_k]$. Then $w_k\in \setS$, and $\setS$ is bounded in
$C^{2,\alpha}(\overline{D^{++}})$. Thus, for any subsequence
$w_{k_l}$, there exists a further subsequence $ w_{k_{l_m}}$
converging in $C^{2,{\alpha/2}}(\overline{D^{++}})$. Then the limit
$\tilde w$ is a solution of
(\ref{unif-ellipt-NONlinear-xy-w})--(\ref{dirichlet-w-nonlin}) with
the limiting function $g$ in the right-hand side of
(\ref{unif-ellipt-NONlinear-xy-w}). By uniqueness of solutions in
$\setS$ to
(\ref{unif-ellipt-NONlinear-xy-w})--(\ref{dirichlet-w-nonlin}), we
have $\tilde w=K[W]$. Then it follows that the whole sequence
$K[W_k]$ converges to $K[W]$. Thus $K:\setS \to \setS$ is continuous
in $C^{2,{\alpha/2}}(\overline{D^{++}})$. Therefore, there exists
$w\in \setS$ which is a fixed point of $K$. This function $w$ is a
solution of (\ref{unif-ellipt-NONlinear-xy-w-fullEq}).

\medskip
Since $v$ satisfies
(\ref{unif-ellipt-NONlinear-xy-v})--(\ref{condOnWedge-v-xy-loc-NONlin}),
it follows from the uniqueness of solutions in $C(\overline{D^{++}})
\cap C^{1}(\overline{D^{++}}\setminus\overline{\{x=0\}}) \cap
C^2(D^{++})$ of problem (\ref{unif-ellipt-NONlinear-xy-w-fullEq})
%--(\ref{dirichlet-w-fullEq-nonlin})
that $w=v$ in $D^{++}$. Thus, $v\in C^{2,\alpha}(\overline{D^{++}})$
and satisfies (\ref{Hoder-est-wedge-sonic-nonLin-unif-ellipt}).

\medskip
{\em Step 5.} It remains to make the following estimate near the
corner $\PtUpL$:
\begin{equation}
\|\psi\|^{(-1-\alpha, \{\PtUpL\})}_{2,\alpha,\Omega^+(\Kphi)}\le
 C,
\label{Hoder-est-P4-for-NONlin-unif-ellipt}
\end{equation}
where $C$ depends only on the data, $\epsP$, and $\delta$.

Since $\psi$ is a solution of the linear equation {\rm
(\ref{unif-ellipt-linear-iterationEquation})} for $\hat\psi=\psi$
and satisfies the boundary conditions {\rm
(\ref{iterationRH})}--{\rm (\ref{iterationCondOnSymmtryLine})}, it
follows from Lemma \ref{C2alpha-near-P4-Lin-UnifEllipt-Lemma} that
$\psi$ satisfies (\ref{growth-est-P4-for-lin-unif-ellipt}) with
constant $\hat C$ depending only on the data and $\delta$.

Now we follow the argument of Lemma
\ref{C2alpha-near-P4-Lin-UnifEllipt-Lemma} (Step 4): We consider
cases (\ref{cases-HolderCorner-1})--(\ref{cases-HolderCorner-3}) and
define the function $v(X,Y)$ by (\ref{defV-mixed-corner}). Then
$\psi$ is a solution of the nonlinear equation
(\ref{unif-ellipt-iteration-equation-sonicStruct}). We apply the
estimates in Appendix \ref{append-1-section}.
 From Lemma
\ref{propertiesNonlinCoeffs-xy} and the properties of the Laplacian
in the polar coordinates, the coefficients of
(\ref{unif-ellipt-iteration-equation-sonicStruct})
 satisfy (\ref{locEstElliptEq-i1-0})--(\ref{locEstElliptEq-i2-0})
with $\lambda$ depending only on the data and $\delta$. It is easy
to see that $v$ defined by (\ref{defV-mixed-corner}) satisfies an
equation of the similar structure and properties
(\ref{locEstElliptEq-i1-0})--(\ref{locEstElliptEq-i2-0}) with the
same $\lambda$, where we use that $0\le \hat d\le 1$. Also, $v$
satisfies the same boundary conditions as in the proof of Lemma
\ref{C2alpha-near-P4-Lin-UnifEllipt-Lemma} (Step 4). Furthermore,
since $\psi$ satisfies (\ref{growth-est-P4-for-lin-unif-ellipt}), we
obtain the $L^\infty$ estimates of $v$ in terms of the data and
$\delta$, e.g., $v$ satisfies (\ref{V-mixed-corner-Linfty}) in case
(\ref{cases-HolderCorner-3}). Now we obtain the
$C^{2,\alpha}$--estimates of $v$ by using Theorem
\ref{locEstElliptEq} for case (\ref{cases-HolderCorner-1}), Theorem
\ref{locEstElliptEq-Dirichlet} for case
(\ref{cases-HolderCorner-2}), and Theorem
\ref{locEstElliptEq-oblique} for case (\ref{cases-HolderCorner-3}).
Writing these estimates in terms of $\psi$, we obtain
(\ref{Hoder-est-P4-for-NONlin-unif-ellipt}), similar to the proof of
Lemma \ref{C2alpha-near-P4-Lin-UnifEllipt-Lemma} (Step 4).

\medskip
{\em Step 6.} Finally, we prove the comparison principle, assertion
(\ref{estimates-nonlin-UnifEllipt-i3}). The function
$u=\psi_1-\psi_2$ is a solution of a linear problem of form
(\ref{unif-ellipt-linear-iterationEquation}), (\ref{iterationRH}),
(\ref{iterationCondOnWedge}), and (\ref{iterationCondOnSymmtryLine})
with right-hand sides $\Nl_\delta(\psi_1)-\Nl_\delta(\psi_2)$ and
$B_k(\psi_1)- B_k(\psi_2)$ for $k=1,2,3$, respectively, and $u\ge 0$
on $\sonic$. Now the comparison principle follows from Lemma
\ref{comparisonPrincipleOfUnifEllipt-Lemma}.
\Endproof
%*************** END PROOF *****LEMMA**********

Using Lemma \ref{existence-Lin-UnifEllipt-Lemma} and  the definition
of map $\hat{J}$ in \eqref{6.9a}, and using Lemma
\ref{estimates-nonlin-UnifEllipt} and the Leray-Schauder Theorem, we
conclude the proof of Proposition \ref{existSolUnifEllipt}.
\Endproof
%*************** END PROOF *****PROPOSITION**********

Using Proposition \ref{existSolUnifEllipt} and sending $\delta\to
0$, we establish the existence of a solution
of problem
(\ref{iterationEquation})--(\ref{iterationCondOnSymmtryLine}).

%*****************BEGIN PROPOSITION*************
\begin{proposition}\label{existSolDegenEllipt}
Let $\epsP, \varepsilon, M_1,$ and $M_2$ be as in Proposition {\rm
\ref{existSolUnifEllipt}}. Then there exists a solution $\psi\in
C(\overline{\Omega^+(\Kphi)})\cap
C^1(\overline{\Omega^+(\Kphi)}\setminus\overline\sonic) \cap
C^2(\Omega^+(\Kphi))$ of problem {\rm
(\ref{iterationEquation})}--{\rm (\ref{iterationCondOnSymmtryLine})}
so that the solution $\psi$ satisfies {\rm
(\ref{L-infty-for-unif-ellipt})}--{\rm
(\ref{Hoder-est-for-unif-ellipt})}.
\end{proposition}
%***************END PROPOSITION***************
\Proof
%*************** BEGIN PROOF *****PROPOSITION**********
Let $\delta\in(0,\delta_0)$. Let $\psi_\delta$ be a solution of
(\ref{unif-ellipt-iterationEquation}) and
(\ref{iterationRH})--(\ref{iterationCondOnSymmtryLine}) obtained in
Proposition \ref{existSolUnifEllipt}. Using
(\ref{Hoder-est-for-unif-ellipt}), we can find a sequence $\delta_j$
for $j=1,\dots$ and $\psi\in
C^1(\overline{\Omega^+(\Kphi)}\setminus\overline\sonic)\cap
C^2(\Omega^+(\Kphi))$
 such that, as $j\to\infty$,
we have
%stop
\begin{enumerate}\renewcommand{\theenumi}{\roman{enumi}}
\item $\delta_j \to 0$;
\item  $\psi_{\delta_j}\to\psi$ in $C^1(\overline{\Omega^+_s(\Kphi)})$
for every $s\in (0,c_2/2)$, where
$\Omega^+_s(\Kphi)=\Omega^+(\Kphi)\cap\{c_2-r>s\}$;
\item $\psi_{\delta_j}\to\psi$ in $C^2(K)$
for every compact $K\subset\Omega^+(\Kphi)$.
\end{enumerate}
Then, since each $\psi_{\delta_j}$ satisfies
(\ref{unif-ellipt-iterationEquation}), (\ref{iterationRH}), and
(\ref{iterationCondOnWedge})--(\ref{iterationCondOnSymmtryLine}), it
follows that $\psi$ satisfies
(\ref{iterationEquation})--(\ref{iterationRH}) and
(\ref{iterationCondOnWedge})--(\ref{iterationCondOnSymmtryLine}).
Also, since each $\psi_{\delta_j}$ satisfies
(\ref{L-infty-for-unif-ellipt})--(\ref{Hoder-est-for-unif-ellipt}),
$\psi$ also satisfies these estimates. {}From
(\ref{barier-for-unif-ellipt}), we conclude that $\psi\in
C(\overline{\Omega^+(\Kphi)})$ and satisfies
(\ref{iterationCondOnSonicLine}). \Endproof
%*************** END PROOF *****PROPOSITION**********

%******************* END SECTION ***************************

\section{Existence of the Iteration Map and Its Fixed Point}
\label{fixedPtSection} In this section we perform Steps 4--8 of the
procedure described in \S \ref{overViewProcedureSubsection}. In the
proofs of this section, the universal constant $C$ depends only on
the data.

We assume that $\Kphi\in\setK$ and the coefficients in problem
(\ref{iterationEquation})--(\ref{iterationCondOnSymmtryLine}) are
determined by $\Kphi$. Then the existence of a solution $\psi\in
C(\overline{\Omega^+(\Kphi)})\cap
C^1(\overline{\Omega^+(\Kphi)}\setminus\overline\sonic) \cap
C^2(\Omega^+(\Kphi))$ of
(\ref{iterationEquation})--(\ref{iterationCondOnSymmtryLine})
follows from Proposition \ref{existSolDegenEllipt}.

We first show that a comparison principle holds for
 (\ref{iterationEquation})--(\ref{iterationCondOnSymmtryLine}).
We use the operators $\Nl$ and $\Ml$ introduced in
(\ref{iterationEquation}) and (\ref{iterationRH}). Also, for
$\mu>0$, we denote
\begin{eqnarray*}
&&\Omega^{+,\mu}(\Kphi):=\Omega^+(\Kphi)\cap\{c_2-r<\mu\}, \quad
\shock^\mu(\Kphi):=\shock(\Kphi)\cap\{c_2-r<\mu\},
\\
&& \wedgeB^\mu:=\wedgeB\cap\{c_2-r<\mu\}.
\end{eqnarray*}
%*****************BEGIN LEMMA*************
\begin{lemma}\label{comparisonPrincipleOfDegenEllipt-Lemma}
Let $\epsP, \varepsilon, M_1$, and $M_2$ be as in Proposition {\rm
\ref{existSolDegenEllipt}}, and $\mu\in (0, \kappa)$, where $\kappa$
is defined in \S\,{\rm \ref{iteration-dom-subsect}}. Then the
following comparison principle holds:  If $\psi_1, \psi_2\in
C(\overline{\Omega^{+,\mu}(\Kphi)}) \cap
C^{1}(\overline{\Omega^{+,\mu}(\Kphi)}\setminus\overline\sonic) \cap
C^2(\Omega^{+,\mu}(\Kphi))$ satisfy that
\begin{eqnarray*}
&&\Nl(\psi_1)\le\Nl(\psi_2) \qquad \text{in}\,\,
\Omega^{+,\mu}(\Kphi),\\
&&\Ml(\psi_1)\le\Ml(\psi_2)\qquad \text{on}\,\,
\shock^\mu(\Kphi),\\
&&\partial_\nu \psi_1\le\partial_\nu \psi_2\qquad \mbox{on}\,\,
\wedgeB^\mu, \\
&&\psi_1\ge\psi_2\qquad \mbox{on} \,\, \sonic\,\mbox{and }\,\,
\Omega^+(\Kphi)\cap\{c_2-r=\mu\},
\end{eqnarray*}
then
$$
\psi_1\ge \psi_2 \qquad \mbox{in }\, \Omega^{+,\mu}.
$$
\end{lemma}
%*****************END LEMMA*************
\Proof
%*************** BEGIN PROOF *****LEMMA**********
Denote $\Sigma_\mu:=\Omega^+(\Kphi)\cap\{c_2-r=\mu\}$. If $\mu\in
(0, \kappa)$, then $\partial \Omega^{+,\mu}(\Kphi)=
\shock^\mu(\Kphi)\cup\wedgeB^\mu\cup\overline\sonic\cup\overline{\Sigma_\mu}$.

{}From $\Nl(\psi_1)\le\Nl(\psi_2)$, the difference $\psi_1-\psi_2$
is a supersolution of a linear equation of form
(\ref{unif-ellipt-linear-iterationEquation}) in
$\Omega^{+,\mu}(\Kphi)$ and, by Lemma \ref{propertiesNonlinCoeffs}
(\ref{propertiesNonlinCoeffs-i1}), this equation is uniformly
elliptic in $\Omega^{+,\mu}(\Kphi)\cap\{c_2-r>s\}$ for any $s\in (0,
\mu)$. Then the argument of Steps (i)--(ii) in the proof of Lemma
\ref{comparisonPrincipleOfUnifEllipt-Lemma} implies that
$\psi_1-\psi_2$ cannot achieve a negative minimum in the interior of
$\Omega^{+,\mu}(\Kphi)\cap\{c_2-r>s\}$ and in the relative interiors
of $\shock^\mu(\Kphi)\cap\{c_2-r>s\}$ and
$\wedgeB^\mu\cap\{c_2-r>s\}$. Sending $s\to 0+$, we conclude the
proof.
\Endproof
%*************** END PROOF *****LEMMA**********

%*****************BEGIN LEMMA*************
\begin{lemma}\label{uniquenessCor}
A solution $\psi\in C(\overline{\Omega^{+}(\Kphi)}) \cap
C^{1}(\overline{\Omega^{+}(\Kphi)}\setminus\overline\sonic) \cap
C^2(\Omega^{+}(\Kphi))$ of {\rm (\ref{iterationEquation})}--{\rm
(\ref{iterationCondOnSymmtryLine})} is unique.
\end{lemma}
%*************** END LEMMA**********
\Proof
%*************** BEGIN PROOF *****LEMMA**********
If $\psi_1$ and $\psi_2$ are two solutions, then we repeat the proof
of Lemma \ref{comparisonPrincipleOfDegenEllipt-Lemma} to show that
$\psi_1-\psi_2$ cannot achieve a negative minimum in
$\Omega^+(\Kphi)$ and in the relative interiors of $\shock(\Kphi)$
and $\wedgeB$. Now equation (\ref{iterationEquation}) is linear,
uniformly elliptic near $\Sigma_0$ (by Lemma
\ref{propertiesNonlinCoeffs}), and the function $\psi_1-\psi_2$ is
$C^1$ up to the boundary in a neighborhood of $\Sigma_0$. Then the
boundary condition  {\rm (\ref{iterationCondOnSymmtryLine})}
combined with Hopf's Lemma yields that $\psi_1-\psi_2$ cannot
achieve a minimum in the relative interior of $\Sigma_0$. By the
argument of Step (iii) in the proof of Lemma
\ref{comparisonPrincipleOfUnifEllipt-Lemma}, $\psi_1-\psi_2$ cannot
achieve a negative minimum at the points $\PtLwL$ and $\PtLwR$.
Thus, $\psi_1\ge \psi_2$ in $\Omega^+(\Kphi)$ and, by symmetry, the
opposite is also true.
\Endproof
%*************** END PROOF *****LEMMA**********

%*****************BEGIN LEMMA*************
\begin{lemma}\label{quadraticGrowthPsi-Lemma}
There exists $\hat C>0$ depending only on the data such that, if
$\epsP, \varepsilon, M_1$, and $M_2$ satisfy {\rm
(\ref{condConst-00})},  the solution $\psi\in
C(\overline{\Omega^+(\Kphi)})\cap
C^1(\overline{\Omega^+(\Kphi)}\setminus\overline\sonic) \cap
C^2(\Omega^+(\Kphi))$ of {\rm (\ref{iterationEquation})}--{\rm
(\ref{iterationCondOnSymmtryLine})} satisfies
\begin{equation}
0\le\psi(x,y)\le \frac{3}{5(\gamma+1)}x^2 \qquad\mbox{in
}\;\ElDomS(\Kphi):=\Omega^{+,2\varepsilon}(\Kphi).
\label{L-ifty-BdIteratin-Sonic}
\end{equation}
\end{lemma}
%*****************END LEMMA*************

\Proof
%*************** BEGIN PROOF *****LEMMA**********
We first notice that $\psi\ge 0$ in $\Omega^+(\Kphi)$ by Proposition
\ref{existSolDegenEllipt}. Now we make estimate
\eqref{L-ifty-BdIteratin-Sonic}. Set
$$
w(x,y):=\frac{3}{5(\gamma+1)}x^2.
$$

We first show that $w$ is a supersolution of equation
(\ref{iterationEquation}). Since (\ref{iterationEquation}) rewritten
in the $(x,y)$--coordinates in $\ElDomS(\Kphi)$ has form
(\ref{iteration-equation-sonicStruct}), we write it as
$$
\Nl_1(\psi)+\Nl_2(\psi)=0,
$$
where
\begin{eqnarray*}
&&\Nl_1(\psi)= \big(
2x-(\gamma+1)x\zeta_1(\frac{\psi_x}{x})\big)\psi_{xx} +{1\over
c_2}\psi_{yy}
-\psi_{x},\\
&&\Nl_2(\psi)=O_1^\Kphi\psi_{xx}+O_2^\Kphi\psi_{xy}
+O_3^\Kphi\psi_{yy}-O_4^\Kphi\psi_{x}+O_5^\Kphi\psi_{y}.
\end{eqnarray*}
Now we substitute $w(x,y)$.
By (\ref{defZeta-1}),
$$
\zeta_1\big(\frac{w_x}{x}\big)=
\zeta_1\big(\frac{6}{5(\gamma+1)}\big)=\frac{6}{5(\gamma+1)},
$$
thus
$$
\Nl_1(w)=-\frac{6}{25(\gamma+1)}x.
$$
Using (\ref{estSmallterms-iter}), we have
$$
|\Nl_2(w)| = \Big|\frac{6}{5(\gamma+1)}O_1^\Kphi(Dw, x, y)
+\frac{6x}{5(\gamma+1)} O_4^\Kphi(Dw, x, y)\Big| \le Cx^{3/2}\le
C\varepsilon^{1/2}x,
$$
where the last inequality holds since $x\in(0, 2\varepsilon)$ in $\ElDomS(\Kphi)$. Thus, if
$\varepsilon$ is small, we find
$$
\Nl(w)<0\qquad\mbox{ in }\;\ElDomS(\Kphi).
$$
The required smallness of $\varepsilon$ is achieved if
(\ref{condConst-00}) is satisfied with large $\hat C$.

Also, $w$ is a supersolution of (\ref{iterationRH}): Indeed, since
(\ref{iterationRH}) rewritten in the $(x,y)$--coordinates has form
(\ref{iterationRH-lf-flattened}), estimates
(\ref{estCoefsIterRH-flattened}) hold, and $x>0$, we find
$$
\Ml(w)=\hat b_1(x,y)\frac{6}{5(\gamma+1)}x+\hat
b_3(x,y)\frac{3}{5(\gamma+1)}x^2<0 \qquad \mbox{on }\, \,
\shock(\Kphi)\cap\overline\DomS.
$$

Moreover, on $\wedgeB$, $w_\nu\equiv w_y=0=\psi_\nu$. Furthermore,
$w=0=\psi$ on $\sonic$ and, by (\ref{L-infty-for-unif-ellipt}),
$\psi\le w$ on $\{x=2\varepsilon\}$ if
$$
C\epsP\le \varepsilon^2,
$$
where $C$ is a large constant depending only on the data, i.e., if
(\ref{condConst-00}) is satisfied with large $\hat C$. Thus,
$\psi\le w$ in $\ElDomS(\Kphi)$ by Lemma
\ref{comparisonPrincipleOfDegenEllipt-Lemma}.
\Endproof
%*************** END PROOF *****LEMMA**********

We now estimate the norm
$\|\psi\|_{2,\alpha,\hElDomS(\Kphi)}^{(par)}$ in the  subdomain
$\hElDomS(\Kphi)\defd\Omega^+(\Kphi)\cap \{c_2-r<\varepsilon\}$ of
$\ElDomS(\Kphi)\defd\Omega^+(\Kphi)\cap \{c_2-r<2\varepsilon\}$.
%*****************BEGIN LEMMA*************
\begin{lemma}\label{EstParabolicHolder-Lemma}
There exist $\hat C, C>0$ depending only on the data such that, if
$\epsP, \varepsilon, M_1$, and $M_2$ satisfy {\rm
(\ref{condConst-00})}, the solution $\psi\in
C(\overline{\Omega^+(\Kphi)})\cap
C^1(\overline{\Omega^+(\Kphi)}\setminus\overline\sonic) \cap
C^2(\Omega^+(\Kphi))$ of {\rm (\ref{iterationEquation})}--{\rm
(\ref{iterationCondOnSymmtryLine})} satisfies
\begin{equation}
\|\psi\|_{2,\alpha,\hElDomS(\Kphi)}^{(par)} \leq C.
\label{ParabolicHolder-BdIteratin-Sonic}
\end{equation}
\end{lemma}

%*****************END LEMMA*************
\Proof
%*************** BEGIN PROOF *****LEMMA**********
We assume $\hat C$ in (\ref{condConst-00}) is sufficiently large so
that $\epsP, \varepsilon, M_1$, and $M_2$ satisfy the conditions of
Lemma {\rm \ref{quadraticGrowthPsi-Lemma}}.

\medskip
{\em Step 1.} We work in the $(x,y)$--coordinates and, in
particular, we use
(\ref{domain-in-rescaled-lemma})--(\ref{holder-hat-f}). We can
assume $\varepsilon<\kappa/20$, which can be achieved by increasing
$\hat C$ in (\ref{condConst-00}).

For $z:=(x,y)\in \hElDomS(\Kphi)$ and $\rho\in(0,1)$, define
\begin{equation}\label{parabRectangles-1}
\quad \tilde R_{z,\rho}:=\left\{(s,t)\;\;:\;\;
|s-x|<\frac{\rho}{4}x, |t-y|<\frac{\rho}{4}\sqrt{x}\right\}, \quad
R_{z,\rho}:=\tilde R_{z,\rho} \cap \Omega^+(\Kphi).
\end{equation}
Since $\ElDomS(\Kphi)=\Omega^+(\Kphi)\cap \{c_2-r<2\varepsilon\}$,
then, for any $z\in \hElDomS(\Kphi)$ and $\rho\in(0,1)$,
\begin{equation}\label{localizeRectangle}
R_{z, \rho}\subset \Omega^+(\Kphi)\cap \{(s,t)\;\;:\;\;
\frac{3}{4}x<s<\frac{5}{4}x\}
\subset\ElDomS(\Kphi).
%\quad\mbox{for any } z=(x,y)\in\hElDomS(\Kphi),\;\rho\in(0,1).
\end{equation}

For any $z\in \hElDomS(\Kphi)$, we have at least one of the
following three cases:

\begin{enumerate}\renewcommand{\theenumi}{\roman{enumi}}

\item  $R_{z, 1/10}=\tilde R_{z, 1/10}$;
\item $z\in R_{z_w, 1/2}$ for $z_w=(x,0)\in  \wedgeB$;
\item $z\in R_{z_s, 1/2}$ for $z_s=(x,\hat f_\Kphi(x))\in  \shock(\Kphi)$.
\end{enumerate}
Thus, it suffices to make the local estimates of $D\psi$ and
$D^2\psi$ in the following rectangles with $z_0:=(x_0, y_0)$:
\begin{enumerate}\renewcommand{\theenumi}{\roman{enumi}}
\item \label{cases-ParabolicHolder-1}
$R_{z_0, {1/20}}$ for $z_0\in\hElDomS(\Kphi)$ and $R_{z_0,
1/10}=\tilde R_{z_0, 1/10}$;

\item \label{cases-ParabolicHolder-2}
$R_{z_0, {1/2}}$ for  ${z_0}\in \wedgeB\cap \{x<\varepsilon\}$;

\item \label{cases-ParabolicHolder-3}
$R_{z_0, {1/2}}$ for  ${z_0}\in \shock(\Kphi)\cap
\{x<\varepsilon\}$.
\end{enumerate}

\medskip
{\em Step 2.} We first consider case (\ref{cases-ParabolicHolder-1})
in Step 1. Then
$$
R_{z_0, {1/10}}=\big\{(x_0+\frac{x_0}{4}S,
 y_0+\frac{\sqrt{x_0}}{4}T)\; : \;
(S,T)\in Q_{1/10}\big\},
$$
where $Q_\rho\defd (-\rho, \rho)^2$ for $\rho>0$.

Rescale $\psi$ in $R_{z_0, {1/10}}$ by defining
\begin{equation}\label{parabRescaling-1}
\psi^{(z_0)}(S, T):=\frac{1}{x_0^2}\psi(x_0+\frac{x_0}{4}S,
y_0+\frac{\sqrt{x_0}}{4}T) \qquad\mbox{for }\,\,(S, T)\in Q_{1/10}.
\end{equation}
Then, by (\ref{L-ifty-BdIteratin-Sonic}) and
(\ref{localizeRectangle}),
\begin{equation}\label{estInterLinfty-Rescaled-pf}
\|\psi^{(z_0)}\|_{C(\overline{Q_{1/10}})}\le 1/(\gamma+1).
\end{equation}
Moreover, since $\psi$ satisfies equation
(\ref{iteration-equation-sonicStruct})--(\ref{iteration-equation-sonicStruct-coef})
 in $R_{z_0, {1/10}}$, then
$\psi^{(z_0)}$ satisfies
\begin{eqnarray}\label{iteration-equation-sonicStruct-ParabRescaled-1}
&&\qquad \Big((1+\frac{1}{4}S)\big(2-
(\gamma+1)\zeta_1(\frac{4\psi^{(z_0)}_S}{1+S/4})\big) +x_0
O_1^{(\Kphi,z_0)} \Big)\psi^{(z_0)}_{SS}
+x_0 O_2^{(\Kphi,z_0)}\psi^{(z_0)}_{ST}
 \\
&&\qquad\qquad
+ \big({1\over c_2}+x_0 O_3^{(\Kphi,z_0)}\big)\psi^{(z_0)}_{TT}
-({1\over 4}+x_0O_4^{(\Kphi,z_0)})\psi^{(z_0)}_{S} +x_0^2
O_5^{(\Kphi,z_0)}\psi^{(z_0)}_{T}=0 \nonumber
\end{eqnarray}
in $Q_{1/10}$, where
\begin{equation}\label{iteration-equation-sonicStruct-coef-rescaled}
\begin{array}{ll}
\tilde O_1^{\Kphi, z_0}(p,S,T)=
-\frac{(1+{S/4})^2}{c_2}+{\gamma+1\over 2c_2}
\Big(2(1+{S/4})^2\zeta_1\big(\frac{4p_1}{1+{S/4}}\big)
-16|\Kphi^{(z_0)}_S|^2\Big)
%\nonumber
\\
\quad \qquad\qquad\qquad-{\gamma-1\over c_2}\Big(\Kphi^{(z_0)}+
{8x_0\over
(c_2-x_0(1+{S/4}))^2}|\Kphi^{(z_0)}_T|^2\Big), %\nonumber
\\
\tilde O_2^{\Kphi, z_0}(p,S,T)=-{8\over c_2(c_2-x_0(1+{S/4}))^2}
\big(4x_0\Kphi^{(z_0)}_S+c_2-x_0(1+{S/4})\big)\Kphi^{(z_0)}_T,
%\nonumber
\\
\tilde O_3^{\Kphi, z_0}(p,S,T)\\
\quad ={1\over c_2(c_2-x_0(1+{S/4}))^2}
\bigg\{(1+{S/4})\big(2c_2-x_0(1+{S/4})\big)
%\nonumber
\\
  \qquad-(\gamma-1) \Big( x_0\Kphi^{(z_0)}+ (c_2-x_0(1+{S/4}))(1+{S/4})
\zeta_1\big(\frac{4p_1}{1+{S/4}}\big)+8x_0|\Kphi^{(z_0)}_S|^2 \Big)
%\nonumber
\\
\qquad -\frac{8(\gamma+1)}{(c_2-x_0(1+{S/4}))^2}
x_0^2|\Kphi^{(z_0)}_T|^2\bigg\}, %\nonumber
\\
\tilde O_4^{\Kphi,z_0}(p,S,T) =
\frac{1}{c_2-x_0(1+{S/4})}\bigg\{1+{S/4} - {\gamma-1\over
c_2}\bigg(x_0\Kphi^{(z_0)} +8x_0|\Kphi^{(z_0)}_S|^2 %\nonumber
\\
 \qquad\qquad \qquad\quad +\big(c_2-x_0(1+{S/4})\big) (1+{S/4})
\zeta_1\big(\frac{4p_1}{1+{S/4}}\big) +\frac{
8|x_0\Kphi^{(z_0)}_T|^2} {(c_2-x_0(1+{S/4}))^2} \bigg)
\bigg\}, %\nonumber
\\
\tilde O_5^{\Kphi, z_0}(p,S,T)={8\over c_2(c_2-x_0(1+{S/4}))^2}
\big(4x_0\Kphi^{(z_0)}_S+2c_2-2x_0(1+{S/4})\big) \Kphi^{(z_0)}_T,
%\nonumber
\end{array}
\end{equation}
where $\Kphi^{(z_0)}$ is the rescaled $\Kphi$ as in
(\ref{parabRescaling-1}). By (\ref{localizeRectangle}) and
$\Kphi\in\setK$, we have
$$
\|\Kphi^{(z_0)}\|_{C^{2,\alpha}(\overline{Q_{1/10}})}\le CM_1,
$$
and thus
\begin{equation}\label{errorTermsRescaled-pf}
\|\tilde O_k^{\Kphi,
z_0}\|_{C^1(\overline{Q_{1/10}^{(z)}}\times\bR^2)} \le C(1+M_1^2),
\quad k=1,\dots, 5.
\end{equation}
Now, since every term $O_k^{(\Kphi,z_0)}$ in
(\ref{iteration-equation-sonicStruct-ParabRescaled-1}) is multiplied
by $x_0^{\beta_k}$ with $\beta_k\ge 1$ and $x_0\in (0,
\varepsilon)$, condition (\ref{condConst-00}) (possibly after
increasing $\hat C$) depending only on the data implies that
equation (\ref{iteration-equation-sonicStruct-ParabRescaled-1})
satisfies conditions
(\ref{locEstElliptEq-i1-0})--(\ref{locEstElliptEq-i2-0}) in
$Q_{1/10}$ with $\lambda>0$
depending only on
$c_2$, i.e., on the data by (\ref{condRewritingRH-0}). Then, using
Theorem \ref{locEstElliptEq} and (\ref{estInterLinfty-Rescaled-pf}),
we find
\begin{equation}\label{estInterHolder-Rescaled-1-pf}
\|\psi^{(z_0)}\|_{C^{2,\alpha}(\overline{Q_{1/20}})}\le C.
\end{equation}

\medskip
{\em Step 3.} We then consider case (\ref{cases-ParabolicHolder-2})
in Step 1. Let ${z_0}\in \wedgeB\cap \{x<\varepsilon\}$. Using
(\ref{domain-in-rescaled-lemma}) and assuming that $\epsP$ and
$\varepsilon$ are sufficiently small depending only on the data, we
have $\overline{R_{z_0,
1}}\cap\partial\Omega^+(\Kphi)\subset\wedgeB$ and thus, for any
$\rho\in(0,1]$,
$$
R_{z_0, \rho}=\bigg\{(x_0+\frac{x_0}{4}S,
 y_0+\frac{\sqrt{x_0}}{4}T)\; : \;
(S,T)\in Q_\rho\cap\{T>0\}\bigg\}.
$$
The choice of parameters for that can be made as follows: First
choose $\epsP$ small so that $|\bar\mxx-\mxx_1|\le |\bar\mxx|/10$,
where $\bar\mxx$ is defined by (\ref{3.14}), which is possible since
$\mxx_1\to \bar\mxx$ as $\theta_w\to \pi/2$, and then choose
$\varepsilon<(|\bar\mxx|/10)^2$.

Define $\psi^{(z_0)}(S,T)$ by (\ref{parabRescaling-1}) for $(S,T)\in
Q_1\cap\{T>0\}$. Then, by (\ref{L-ifty-BdIteratin-Sonic}) and
(\ref{localizeRectangle}),
\begin{equation}\label{estInterLinfty-2-Rescaled-pf}
\|\psi^{(z_0)}\|_{C(\overline{Q_1}\cap\{T\ge 0\})}\le 1/(\gamma+1).
\end{equation}
Moreover,  similar to Step 2, $\psi^{(z_0)}$ satisfies equation
(\ref{iteration-equation-sonicStruct-ParabRescaled-1}) in
$Q_1\cap\{T>0\}$, and the terms $\tilde O_k^{\Kphi, z_0}$ satisfy
estimate (\ref{errorTermsRescaled-pf}) in $Q_1\cap\{T>0\}$. Then, as
in Step 2, we conclude that
(\ref{iteration-equation-sonicStruct-ParabRescaled-1}) satisfies
conditions (\ref{locEstElliptEq-i1-0})--(\ref{locEstElliptEq-i2-0})
in $Q_1\cap\{T>0\}$ if (\ref{condConst-00}) holds with sufficiently
large $\hat C$. Moreover, since $\psi$ satisfies
(\ref{iterationCondOnWedge}), it follows that
$$
\partial_T\psi^{(z_0)}=0\qquad\mbox{on }
\{T=0\}\cap Q_{1/2}.
$$
Then, from Theorem \ref{locEstElliptEq-oblique},
\begin{equation}\label{estInterHolder-Rescaled-2-pf}
\|\psi^{(z_0)}\|_{C^{2,\alpha}(\overline{Q_{1/2}}\cap\{T\ge 0\})}\le
C.
\end{equation}

\medskip
{\em Step 4.} We now consider case (\ref{cases-ParabolicHolder-3})
in Step 1. Let ${z_0}\in \shock(\Kphi)\cap \{x<\varepsilon\}$. Using
(\ref{domain-in-rescaled-lemma}) and the fact that $y_0=\hat
f_\Kphi(x_0)$ for $z_0\in\shock(\Kphi)\cap\{x<\varepsilon\}$, and
assuming that $\epsP$ and $\varepsilon$ are small as in Step 3, we
have $\overline{R_{z_0,
1}}\cap\partial\Omega^+(\Kphi)\subset\shock(\Kphi)$ and thus, for
any $\rho\in(0,1]$,
$$
R_{z_0, \rho }=\bigg\{(x_0+\frac{x_0}{4} S,
 y_0+\frac{\sqrt{x_0}}{4} T)\; : \;
(S,T)\in Q_\rho\cap\{T<\varepsilon^{1/4}F_{(z_0)}(S)\}\bigg\}
$$
with
$$
F_{(z_0)}(S)=4\frac{\hat f_\Kphi(x_0+\frac{x_0}{4}S)- \hat
f_\Kphi(x_0)}{\varepsilon^{1/4}\sqrt{ x_0}}.
$$
Then we use (\ref{holder-hat-f-S}) and $x_0\in (0, 2\varepsilon)$ to
obtain
\begin{eqnarray*}
&&F_{(z_0)}(0)=0,\\
&&\|F_{(z_0)}\|_{C^1([-{1/2}, {1/2}])} \le {\|\hat
f'_\Kphi\|_{L^\infty([0,2\varepsilon])}x_0\over\varepsilon^{1/4}
\sqrt{x_0}}\le C(1+M_1\varepsilon)\varepsilon^{1/4},\\
&& \|F_{(z_0)}''\|_{C^\alpha([-{1/2}, {1/2}])}
\le {\|\hat
f''_\Kphi\|_{L^\infty([0,2\varepsilon])}x_0^2
      +[\hat f_\Kphi'']_{\alpha, (x_0/2, \varepsilon)}x_0^{2+\alpha}
\over 4\varepsilon^{1/4}\sqrt{x_0}}\le C(1+M_1)\varepsilon^{5/4},
\end{eqnarray*}
and thus, from (\ref{condConst-00}),
\begin{equation}\label{norm-rescaled-bdry-func}
\|F_{(z_0)}\|_{C^{2,\alpha}([-{1/2},{1/2}])}\le C/\hat C\le 1
\end{equation}
if $\hat C$ is large. Define $\psi^{(z_0)}(S, T)$ by
(\ref{parabRescaling-1}) for $(S,T)\in
Q_1\cap\{T<\varepsilon^{1/4}F_{(z_0)}(S)\}$. Then, by
(\ref{L-ifty-BdIteratin-Sonic}) and (\ref{localizeRectangle}),
\begin{equation}\label{estInterLinfty-3-Rescaled-pf}
\|\psi^{(z_0)}\|_{C(\overline{Q_1}\cap\{T\le F_{(z_0)}(S)\})} \le
1/(\gamma+1).
\end{equation}
Similar to Steps 2--3, $\psi^{(z_0)}$ satisfies equation
(\ref{iteration-equation-sonicStruct-ParabRescaled-1}) in
$Q_1\cap\{T<\varepsilon^{1/4}F_{(z_0)}(S)\}$ and the terms $\tilde
O_k^{\Kphi, z_0}$ satisfy estimate (\ref{errorTermsRescaled-pf}) in
$Q_1\cap\{T<\varepsilon^{1/4}F_{(z_0)}(S)\}$. Then, as in Steps
2--3, we conclude that
(\ref{iteration-equation-sonicStruct-ParabRescaled-1}) satisfies
conditions (\ref{locEstElliptEq-i1-0})--(\ref{locEstElliptEq-i2-0})
in $Q_1\cap\{T<\varepsilon^{1/4}F_{(z_0)}(S)\}$ if
(\ref{condConst-00}) holds with sufficiently large $\hat C$.
Moreover, $\psi$ satisfies (\ref{iterationRH}) on $\shock(\Kphi)$,
which can be written in form (\ref{iterationRH-lf-flattened}) on
$\shock(\Kphi)\cap \DomS$. This implies that $\psi^{(z_0)}$
satisfies
$$
\partial_S\psi^{(z_0)}
=\varepsilon^{1/4}\left(B_2\partial_T\psi^{(z_0)}+B_3\psi^{(z_0)}\right)
\qquad\mbox{on } \{T=\varepsilon^{1/4}F_{(z_0)}(S)\}\cap Q_{1/2},
$$
where
\begin{eqnarray*}
&&B_2(S,T)=-{\sqrt{x_0}\over\varepsilon^{1/4}} {\hat b_2\over\hat
b_1}(x_0+\frac{x_0}{4}S,
 y_0+\frac{\sqrt{x_0}}{4}T),\\
&&B_3(S,T)=-{x_0\over 4\varepsilon^{1/4}} {\hat b_3\over\hat
b_1}(x_0+\frac{x_0}{4}S,
 y_0+\frac{\sqrt{x_0}}{4}T).
\end{eqnarray*}
{}From (\ref{estCoefsIterRH-flattened}),
$$
\|(B_2, B_3)\|_{1,\alpha, \overline{Q_1}\cap
\{T\le\varepsilon^{1/4}F_{(z_0)}(S)\}}\le C\varepsilon^{1/4}M_1\le
C/\hat C\le 1.
$$

Now, if $\varepsilon$ is sufficiently small,
it follows from Theorem \ref{locEstElliptEq-non-oblique} that
\begin{equation}\label{estInterHolder-Rescaled-3-pf}
\|\psi^{(z_0)}\|_{C^{2,\alpha}(\overline{Q_{1/2}}\cap
\{T\le\varepsilon^{1/4}F_{(z_0)}(S)\})}\le C.
\end{equation}
The required smallness of $\varepsilon$ is achieved by choosing large
$\hat C$ in (\ref{condConst-00}).

\medskip
{\em Step 5.} Combining (\ref{estInterHolder-Rescaled-1-pf}),
(\ref{estInterHolder-Rescaled-2-pf}), and
(\ref{estInterHolder-Rescaled-3-pf}) with an argument similar to the
proof of  \cite[Theorem 4.8]{GilbargTrudinger} (see also the proof
of Lemma \ref{partIntSeminorm-est-lemma} below), we obtain
(\ref{ParabolicHolder-BdIteratin-Sonic}).
\Endproof
%*************** END PROOF *****LEMMA**********

Now we define the extension of solution $\psi$ from the domain
$\Omega^+(\Kphi)$ to the domain $\Dom$.

%*****************BEGIN LEMMA*************
\begin{lemma}\label{extension-Lemma}
There exist $\hat C, C_1>0$ depending only on the data such that, if
$\epsP, \varepsilon, M_1$, and $M_2$ satisfy {\rm
(\ref{condConst-00})}, there exists $C_2(\varepsilon)$ depending
only on the data and $\varepsilon$ and, for any  $\Kphi\in \setK$,
there exists an extension operator
$$
\ExCFn_\Kphi: C^{1, \alpha}(\overline{\Omega^+(\Kphi)})\cap C^{2,
\alpha}(\overline{\Omega^+(\Kphi)}\setminus
\overline{\sonic\cup\Sigma_0}) \to C^{1, \alpha}(\overline\Dom)\cap
C^{2, \alpha}(\Dom)
$$
satisfying the following two properties:

\begin{enumerate}\renewcommand{\theenumi}{\roman{enumi}}
\item  \label{Uextended_Est}
If  $\psi\in  C^{1, \alpha}(\overline{\Omega^+(\Kphi)})\cap
C^{2,\alpha}(\overline{\Omega^+(\Kphi)} \setminus
\overline{\sonic\cup\Sigma_0})$ is a solution of problem {\rm
(\ref{iterationEquation})}--{\rm
 (\ref{iterationCondOnSymmtryLine})}, then
\begin{eqnarray}
&&\|\ExCFn_\Kphi\psi\|_{2,\alpha, \DomS}^{(par)}
\leq C_1,
\label{Hold-Ext-Sonic}\\
&&\|\ExCFn_\Kphi\psi\|^{(-1-\alpha, \Sigma_0)}_{2,\alpha,
\DomU} \le C_2(\varepsilon)\epsP;
\label{Hold-Ext-Unif}
\end{eqnarray}

\item \label{Uextended_Cont}
Let $\beta \in (0, \alpha)$. If a sequence $\Kphi_k \in \setK$
converges to $\Kphi$ in $C^{1, \beta}(\overline{\Dom})$, then
$\Kphi\in \setK$. Furthermore, if $\psi_k\in C^{1,
\alpha}(\overline{\Omega^+(\Kphi_k)})\cap C^{2,
\alpha}(\overline{\Omega^+(\Kphi_k)}\setminus\overline{\sonic\cup\Sigma_0})$
and $\psi \in C^{1, \alpha}(\overline{\Omega^+(\Kphi)})\cap C^{2,
\alpha}(\overline{\Omega^+(\Kphi)}\setminus
\overline{\sonic\cup\Sigma_0})$ are the solutions of problems {\rm
(\ref{iterationEquation})}--{\rm (\ref{iterationCondOnSymmtryLine})}
for $\Kphi_k$ and $\Kphi$ respectively, then
$\ExCFn_{\Kphi_k}\psi_k\rightarrow \ExCFn_\Kphi\psi$ in $C^{1,
\beta}(\overline{\Dom})$.
\end{enumerate}
\end{lemma}
%*****************END LEMMA*************
\Proof
%*************** BEGIN PROOF *****LEMMA**********
Let $\kappa>0$ be the constant in (\ref{domain-in-rescaled-lemma})
and $\varepsilon<\kappa/20$. For any $\Kphi\in\setK$, we first
define the extension operator separately on the domains
$\Omega_1\defd\Omega^+(\Kphi)\cap\{c_2-r<\kappa\}$ and
$\Omega_2\defd\Omega^+(\Kphi)\cap\{c_2-r>\kappa/2\}$ and then
combine them to obtain the operator $\ExCFn_\Kphi$ globally.

In the argument below, we will state various smallness requirements
on $\epsP$ and $\varepsilon$, which will depend only on the data,
and can be achieved by choosing $\hat C$ sufficiently large in
(\ref{condConst-00}). Also, the constant $C$ in this proof depends
only on the data.

\medskip
{\em Step 1.} First we discuss some properties on the domains
$\Omega^+(\Kphi)$ and $\Dom$ to be used below. Recall $\bar\mxx<0$
defined by (\ref{3.14}), and the coordinates $(\mxx_1,\mxy_1)$ of the
point $\PtUpL$ defined by (\ref{coord-P4}). We assume $\epsP$ small
so that $|\bar\mxx-\mxx_1|\le |\bar\mxx|/10$, which is possible since
$\mxx_1\to \bar\mxx$ as $\theta_w\to \pi/2$. Then $\mxx_1<0$. By
(\ref{OmegaPL-f-higher}) and $\PtUpL\in \shock(\Kphi)$, it follows
that
\begin{equation}\label{domainExtension-0'}
\shock(\Kphi)\subset \overline\Dom \cap \{\mxx<\mxx_1+\varepsilon^{1/4}\}.
\end{equation}
Also, choosing $\varepsilon^{1/4}<|\bar\mxx|/10$, we have
\begin{equation}\label{domainExtension-0}
\mxx_1+\varepsilon^{1/4}<\bar\mxx/2<0.
\end{equation}
Furthermore, when $\epsP$ is sufficiently small,
\begin{equation}\label{domainExtension}
\mbox{\em if $(\mxx, \mxy)\in\Dom\cap\{\mxx<\mxx_1+\varepsilon^{1/4}\}$,
$(\mxx', \mxy)\in\Dom$, and $\mxx'>\mxx$, then
$|\mxx'|<|\mxx|$. }
\end{equation}
Indeed, from the conditions in (\ref{domainExtension}), we have
$$
-c_2<\mxx<\mxx_1+\varepsilon^{1/4}<\bar\mxx/2<0.
$$
Thus, $|\mxx'|<|\mxx|$ if $\mxx'<0$. It remains to consider the case
$\mxx'>0$. Since $\Dom\subset
B_{c_2}(0)\cap\{\mxx<\mxy\cot\theta_w\}$, it follows that
$|\mxx'|\le c_2\cos \theta_w$. Thus $|\mxx'|<|\mxx|$ if $c_2\cos
\theta_w\le |\bar\mxx|/2$. Using (\ref{condRewritingRH-0}) and
(\ref{angleCloseToPiOver2}), we see that the last inequality holds
 if $\epsP>0$ is small depending only on the data.
Then (\ref{domainExtension}) is proved.

Now we define the extensions.

\medskip
{\em Step 2.} First, on $\Omega_1$, we work in the
$(x,y)$--coordinates. Then $\Omega_1=\{0< x< \kappa,\, 0<y<\hat
f_\Kphi(x)\}$ by (\ref{domain-in-rescaled-lemma}).
% and $\hat f_\Kphi$ satisfies (\ref{holder-hat-f})--(\ref{holder-hat-f-S}).
 Denote
$Q_{(a,b)}\defd (0,\kappa) \times (a,b)$. Define the mapping
$\Phi:Q_{(-\infty,\infty)}\to Q_{(-\infty,\infty)}$ by
$$
\Phi(x,y)=(x,1-{y/\hat f_\Kphi(x)}).
$$
The mapping $\Phi$ is invertible with
 the inverse $\Phi^{-1}(x,y)=(x, \hat f_\Kphi(x)(1-y))$.
By definition of $\Phi$,
\begin{eqnarray}
&&\Phi(\Omega_1)=Q_{(0,1)},\qquad
\Phi(\shock(\Kphi)\cap\{0<x<\kappa\})=(0,\kappa)\times\{0\},
\nonumber \\
&&\Phi(\Dom\cap\{0<x<\kappa\})\subset Q_{(-1,1)},
\label{holder-hat-mapF}
\end{eqnarray}
where the last property can be seen as follows: First we note that
$\hat f_\Kphi(x)\ge {\hat f_{0,0}(0)\over 2}>0$ for $x\in(0,\kappa)$
by (\ref{fbFUnctionCloseToNormalPolar}) and (\ref{holder-hat-f}),
then we use  $\Dom\cap\{0<x<\kappa\}=\{0< x< \kappa,\, 0<y<\hat
f_0(x)\}$ and (\ref{holder-hat-f-S}) to obtain ${y\over \hat
f_\Kphi(x)}>0$ on $\Dom\cap\{0<x<\kappa\}$ and
\begin{eqnarray*}
\sup_{(x,y)\in\Dom\cap\{0<x<\kappa\}}{y\over \hat
f_\Kphi(x)}&=&\sup_{x\in(0,\kappa)}{\hat f_0(x)\over \hat
f_\Kphi(x)} \le 1+{2\over  \hat f_{0,0}(0)}\| \hat f_\Kphi-\hat
f_0\|_{C(0,\kappa)}\\
&<& 1+C(M_1\varepsilon+M_2\epsP)<2,
\end{eqnarray*}
if
$M_1\varepsilon$ and $M_2\epsP$ are small, which can be achieved by
choosing $\hat C$ in (\ref{condConst-00}) sufficiently large.

We first define the extension operator:
$$
\Ex: C^{1,\beta}(\overline{Q_{(0,1)}})\cap
C^{2,\beta}(\overline{Q_{(0,1)}}\setminus\{x=0\}) \to
C^{1,\beta}(\overline{Q_{(-1,1)}})\cap
C^{2,\beta}(\overline{Q_{(-1,1)}}\setminus\{x=0\})
$$
for any $\beta\in (0, 1]$. Let $v\in
C^{1,\beta}(\overline{Q_{(0,1)}})\cap
C^{2,\beta}(\overline{Q_{(0,1)}}\setminus\{x=0\})$. Define $\Ex v=v$
in $Q_{(0,1)}$. For $(x,y) \in Q_{(-1,0)}$, define
\begin{equation}
\label{genExtension}
\Ex v(x,y)=
   \sum_{i=1}^3 a_i v(x,-\frac{y}{i}),
\end{equation}
where $a_1=6$, $a_2=-32$, and $a_3=27$, which are determined by
$\sum_{i=1}^3 a_i\left(-\frac{1}{i} \right)^m=1$ for $ m=0, 1, 2$.

Now let $\psi\in  C^{1, \alpha}(\overline{\Omega^+(\Kphi)})\cap
C^{2, \alpha}(\overline{\Omega^+(\Kphi)}\setminus
\overline{\sonic\cup\Sigma_0})$. Let
$$
v=\psi|_{\Omega_1}\circ \Phi^{-1}.
$$
Then $v\in C^{1,\beta}(\overline{Q_{(0,1)}})\cap
C^{2,\beta}(\overline{Q_{(0,1)}}\setminus\{x=0\})$. By
(\ref{holder-hat-mapF}), we have
$\Dom\cap\{c_2-r<\kappa\}\subset\Phi^{-1}(Q_{(-1,1)})$. Thus, we
define an extension operator on $\Omega_1$ by
$$
\ExCFn_\Kphi^1 \psi=(\Ex v)\circ\Phi \qquad \mbox{ on }\;
\Dom\cap\{c_2-r<\kappa\}.
$$
Then $\ExCFn_\Kphi^1 \psi\in C^{1, \alpha}(\overline{\Dom_1})\cap
C^{2, \alpha}(\overline{\Dom_1}\setminus\overline{\sonic})$ with
$\Dom_1:=\Dom\cap\{c_2-r<\kappa\}$.

Next we estimate $\ExCFn_\Kphi^1$ separately on the domains
$\DomS=\Dom\cap\{c_2-r<2\varepsilon\}$ and
$\Dom_1\cap\{c_2-r>\varepsilon/2\}$.

In order to estimate the H\"{o}lder norms of $\ExCFn_\Kphi^1$ on
$\DomS$, we note that
$\Phi(\ElDomS(\Kphi))=(0,2\varepsilon)\times(0,1)$
and $\DomS\subset\Phi^{-1}((0,2\varepsilon)\times(-1,1))$
 in the
$(x,y)$--coordinates. We first show the following estimates, in
which the sets are defined in the $(x,y)$--coordinates:
\begin{eqnarray}
&&\qquad\quad \|\psi\circ \Phi^{-1}\|_{2,\alpha,
(0,2\varepsilon)\times(0, 1)}^{(par)} \le C\|\psi\|_{2,\alpha,
\ElDomS(\Kphi)}^{(par)}\quad\mbox{for any } \psi\in C_{2,\alpha,
(0,2\varepsilon)\times(0, 1)}^{(par)}, \label{normOfSubstitution}
\\
&&\qquad\quad \|w\circ \Phi\|_{2,\alpha, \DomS}^{(par)} \le
C\|w\|_{2,\alpha, (0,2\varepsilon)\times(-1,
1)}^{(par)}\quad\mbox{for any } w\in C_{2,\alpha,
(0,\varepsilon)\times(-1, 1)}^{(par)},\label{normOfSubstitution-bk}
\\
&&\qquad\quad \|\Ex v\|_{2,\alpha,(0,2\varepsilon)\times(-1,
1)}^{(par)} \le C\|v\|_{2,\alpha,(0,2\varepsilon)\times(0,
1)}^{(par)}\,\,\,\mbox{for any } v\in
C_{2,\alpha,(0,2\varepsilon)\times(-1, 1)}^{(par)}.
\label{simpleExtInHolderNorms}
\end{eqnarray}
To show (\ref{normOfSubstitution}), we denote  $v=\psi\circ
\Phi^{-1}$ and estimate every term in definition
(\ref{parabNormsApp}) for $v$. Note that $v(x,y)=\psi(x, \hat
f_\Kphi(x)(1-y))$. In the calculations below, we denote $(v, Dv,
D^2v)=(v, Dv, D^2v)(x,y)$, $(\psi, D\psi, D^2\psi)=(\psi, D\psi,
D^2\psi)(x, \hat f_\Kphi(x)(1-y))$, and $(\hat f_\Kphi, \hat
f_\Kphi', \hat f_\Kphi'')= (\hat f_\Kphi, \hat f_\Kphi', \hat
f_\Kphi'')(x)$. We use that, for $x\in (0,2\varepsilon)$,
$0<M_1x<2M_1\varepsilon<{2/\hat C}$ by  (\ref{condConst-00}). Then,
for any $(x,y)\in (0,2\varepsilon)\times(0,1)$, we have
\begin{eqnarray*}
|v|&=&|\psi|\le \|\psi\|_{2,\alpha,
\ElDomS(\Kphi)}^{(par)}x^2,\\
|v_x|&=&|\psi_x+(1-y)\psi_y \hat f_\Kphi'| \le \|\psi\|_{2,\alpha,
\ElDomS(\Kphi)}^{(par)}\left(x+x^{3/2}(1+M_1x)\right)\le
C\|\psi\|_{2,\alpha,
\ElDomS(\Kphi)}^{(par)}x,\\
|v_{xx}|&=&|\psi_{xx}+2
(1-y)\psi_{xy}\hat f_\Kphi'
+(1-y)^2\psi_{yy}(\hat f_\Kphi')^2
+(1-y)\psi_y\hat f_\Kphi''|\\
&\le& \|\psi\|_{2,\alpha,
\ElDomS(\Kphi)}^{(par)}\left(1+x^{1/2}(1+M_1x)+x(1+M_1x)^2+M_1x^{3/2}\right)\\
&\le& C\|\psi\|_{2,\alpha, \ElDomS(\Kphi)}^{(par)}.
\end{eqnarray*}
The estimates of the other terms in (\ref{parabNormsApp}) for $v$
follow from similar straightforward (but lengthy) calculations.
Thus, (\ref{normOfSubstitution}) is proved. The proof of
(\ref{normOfSubstitution-bk}) is similar by using that $\hat
f_\Kphi(x)\ge {\hat f_{0,0}(0)/2}>0$ for $x\in(0,\kappa)$ from
(\ref{fbFUnctionCloseToNormalPolar}) and (\ref{holder-hat-f}) and
that $\hat f_{0,0}(0)$ depends only on the data. Finally, estimate
(\ref{simpleExtInHolderNorms}) follows readily from
(\ref{genExtension}).

Now, let $\psi\in  C^{1, \alpha}(\overline{\Omega^+(\Kphi)})\cap
C^{2, \alpha}(\overline{\Omega^+(\Kphi)}\setminus
\overline{\sonic\cup\Sigma_0})$ be a solution of
(\ref{iterationEquation})--(\ref{iterationCondOnSymmtryLine}). Then
\begin{eqnarray*}
\|\ExCFn_\Kphi^1 \psi\|_{2,\alpha, \DomS}^{(par)}
&=&\|\Ex (\psi|_{\Omega_1}\circ \Phi^{-1})\circ\Phi\|_{2,\alpha, \DomS}^{(par)}
\;\le\; C\|\Ex (\psi|_{\Omega_1}\circ \Phi^{-1})\|_{2,\alpha,
(0,2\varepsilon)\times(-1, 1)}^{(par)}
\\
&\le& C\|\psi|_{\Omega_1}\circ \Phi^{-1}\|_{2,\alpha, (0,2\varepsilon)\times(0, 1)}^{(par)}
\; \le\; C\|\psi\|_{2,\alpha, \ElDomS(\Kphi)}^{(par)}
\; \le\; C,
\end{eqnarray*}
where the first inequality is obtained from
(\ref{normOfSubstitution-bk}), the second inequality  from
(\ref{simpleExtInHolderNorms}), the third inequality  from
(\ref{normOfSubstitution}), and the last inequality  from
(\ref{ParabolicHolder-BdIteratin-Sonic}). Thus,
(\ref{Hold-Ext-Sonic}) holds for $\ExCFn_\Kphi^1$.

Furthermore, using the second estimate in (\ref{holder-hat-f-S}),
noting that $M_2\epsP\le 1$ by (\ref{condConst-00}), and using the
definition of $\ExCFn_\Kphi^1$ and the fact that the change of
coordinates $(x,y)\to(\mxx, \mxy)$ is smooth and invertible in
$\Dom\cap\{{\varepsilon/2}<x<\kappa\}$,  we find that, in the $(\mxx,
\mxy)$--coordinates,
\begin{equation}
\label{EstExt-1} \|\ExCFn_\Kphi^1
\psi\|_{C^{2,\alpha}(\overline\Dom\cap\{{\varepsilon/2}\le
c_2-r\le\kappa\})} \le
C\|\psi\|_{C^{2,\alpha}(\overline{\ElDom^+(\Kphi)}\cap\{{\varepsilon/2}\le
c_2-r \le \kappa\})}.
\end{equation}

\medskip
{\em Step 3.} Now we define an extension operator in the $(\mxx,
\mxy)$--coordinates. Let
\begin{eqnarray*}
\tilde\Ex &\, :\, & C^1( [0,1]\times [-v_2, \mxy_1])\cap
C^2([0,1]\times
(-v_2, \mxy_1] )\\
&& \to C^1( [-1,1]\times [-v_2, \mxy_1])\cap C^2( [-1,1]\times(-v_2,
\mxy_1])
\end{eqnarray*}
be defined by
$$
\tilde\Ex v(X,Y):=
   \sum_{i=1}^3 a_i v(-\frac{X}{i}, Y)\qquad
\mbox{for } \;(X,Y)\in(-1,0)\times (-v_2, \mxy_1),
$$
where $a_1, a_2$, and $a_3$ are the same as in (\ref{genExtension}).

Let $\hat\Omega_2\defd\Omega^+(\Kphi)\cap \{0\le\mxy\le\mxy_1\}$. Define
the mapping $\Psi:\hat\Omega_2\to (0,1)\times (-v_2, \mxy_1)$ by
$$
\Psi(\mxx,\mxy):=({\mxx- f_\Kphi(\mxy)\over
\mxy\cot\theta_w-f_\Kphi(\mxy)}, \mxy),
$$
where $f_\Kphi(\cdot)$ is the function from
(\ref{shockPL})--(\ref{OmegaPL}). Then the inverse of $\Psi$ is
$$
\Psi^{-1}(X,Y)=(f_\Kphi(Y)+X(Y\cot\theta_w-f_\Kphi(Y)), Y),
$$
and thus, from (\ref{OmegaPL-f-higher}), % and (\ref{condConst-00}),
\begin{equation}\label{holder-hat-mapF-xx}
\|\Psi\|^{(-1-\alpha, [0, 1]\times\{-v_2,
\mxy_1\})}_{2,\alpha, \hat\Omega_2}+
\|\Psi^{-1}\|^{(-1-\alpha, [0, 1]\times\{-v_2,
\mxy_1\})}_{2,\alpha, (0,1)\times (-v_2, \mxy_1)}\le C.
\end{equation}
Moreover, by (\ref{OmegaPL-f-higher}), for sufficiently small
$\varepsilon$ and $\epsP$ (which are achieved by choosing large
$\hat C$ in (\ref{condConst-00})), we have
$\Dom\cap\{-v_2<\mxy<\mxy_1\}\subset\Psi^{-1}([-1,1]\times [-v_2,
\mxy_1])$. Define
$$
\ExCFn_\Kphi^2 \psi:=\tilde\Ex (\psi\circ \Psi^{-1})\circ \Psi\qquad
\mbox{ on }\; \Dom\cap\{-v_2<\mxy<\mxy_1\}.
$$
Then $\ExCFn_\Kphi^2 \psi\in C^{1, \alpha}(\overline{\Dom})\cap
C^{2,\alpha}(\overline{\Dom}
\setminus\overline{\sonic\cup\Sigma_0})$ since $\Dom\setminus
\Omega^+(\Kphi)\subset \Dom\cap\{-v_2<\mxy<\mxy_1\}$. Furthermore,
using (\ref{holder-hat-mapF-xx}) and the definition of
$\ExCFn_\Kphi^2$, we find that, for any $s\in(-v_2, \mxy_1]$,
\begin{equation}
\label{EstExt-2'}
\|\ExCFn_\Kphi^2 \psi\|^{(-1-\alpha,
\Sigma_0)}_{2,\alpha, \Dom\cap\{\mxy\le s\}}
\le C(\mxy_1-s)\|\psi\|^{(-1-\alpha,
\{\PtLwL,\PtLwR\})}_{2,\alpha, \ElDom^+(\Kphi)\cap\{\mxy\le s\}},
\end{equation}
where $C(\mxy_1-s)$ depends only on the data and $\mxy_1-s>0$.

Choosing $\hat C$ large in  (\ref{condConst-00}), we have
$\varepsilon<\kappa/100$. Then (\ref{domain-in-rescaled-lemma})
implies that there exists a unique point
$P'=\shock(\Kphi)\cap\{c_2-r=\kappa/8\}$. Let $P'=(\mxx', \mxy')$ in
the $(\mxx, \mxy)$--coordinates. Then $\mxy'>0$. Using
(\ref{domainExtension-0'}) and (\ref{domainExtension}), we find
\begin{eqnarray*}
&&(\Dom\setminus \Omega^+(\Kphi))\cap\{c_2-r>\kappa/8\}\subset
\Dom\cap\{\mxy\le \mxy'\}, \\
&&\Omega^+(\Kphi)\cap\{\mxy\le \mxy'\}\subset
\Omega^+(\Kphi)\cap\{c_2-r>\kappa/8\}.
\end{eqnarray*}
Also, $\kappa/C\le \mxy_1-\mxy'\le C\kappa$ by (\ref{OmegaPL}),
(\ref{OmegaPL-f-higher}), and (\ref{reflected-shock-s2}). These
facts and (\ref{EstExt-2'}) with $s=\mxy'$ imply
\begin{equation}
\label{EstExt-2} \|\ExCFn_\Kphi^2 \psi\|^{(-1-\alpha,
\Sigma_0)}_{2,\alpha, \Dom\cap\{c_2-r>\kappa/8\}} \le
C\|\psi\|^{(-1-\alpha, \{\PtLwL,\PtLwR\})}_{2,\alpha,
\ElDom^+(\Kphi)\cap\{c_2-r>\kappa/8\}}.
\end{equation}

\medskip
{\em Step 4.} Finally, we choose a cutoff function $\zeta\in
C^\infty(\bR)$ satisfying
$$
\zeta\equiv 1 \mbox{ on }(-\infty, \kappa/4),\quad \zeta\equiv 0
\mbox{ on }(3\kappa/4, \infty),\quad \zeta'\le 0\,\, \mbox{ on }\bR,
$$
and define
$$
\ExCFn_\Kphi\psi:=\zeta(c_2-r)\ExCFn_\Kphi^1\psi+(1-\zeta(c_2-r))\ExCFn_\Kphi^2\psi
\qquad \mbox{in }\Dom.
$$
Since $\ExCFn_\Kphi^k\psi=\psi$ on $\Omega^+(\Kphi)$ for $k=1,2$, so is
$\ExCFn_\Kphi\psi$. Also, from the properties of $\ExCFn_\Kphi^k$
above, $\ExCFn_\Kphi\psi\in C^{1, \alpha}(\overline\Dom)\cap C^{2,
\alpha}(\Dom)$ if $\psi\in C^{1,
\alpha}(\overline{\Omega^+(\Kphi)})\cap C^{2,
\alpha}(\overline{\Omega^+(\Kphi)}\setminus\overline{\sonic\cup\Sigma_0})$.
If such $\psi$ is a solution of
(\ref{iterationEquation})--(\ref{iterationCondOnSymmtryLine}), then
we prove  (\ref{Hold-Ext-Sonic})--(\ref{Hold-Ext-Unif}):
 $\ExCFn_\Kphi\psi\equiv \ExCFn_\Kphi^1\psi$
on $\DomS$ by the definition of $\zeta$ and by
$\varepsilon<\kappa/100$.  Thus, since (\ref{Hold-Ext-Sonic}) has
been proved  in Step 2 for $\ExCFn_\Kphi^1\psi$, we obtain
(\ref{Hold-Ext-Sonic}) for $\ExCFn_\Kphi\psi$. Also, $\psi$
satisfies (\ref{Hoder-est-for-unif-ellipt}) by Proposition
\ref{existSolDegenEllipt}. Using (\ref{Hoder-est-for-unif-ellipt})
with $s={\varepsilon/2}$, (\ref{EstExt-1}), and (\ref{EstExt-2}),
%, and (\ref{domainExtension-2}),
we obtain (\ref{Hold-Ext-Unif}). Assertion (\ref{Uextended_Est}) is
then proved.

\medskip
{\em Step 5.} Finally we prove assertion (\ref{Uextended_Cont}). Let
$\Kphi_k \in \setK$ converge to $\Kphi$ in $C^{1,
\beta}(\overline{\Dom})$. Then obviously $\Kphi\in \setK$. By
(\ref{nondegeneracy})--(\ref{OmegaPL}), it follows that
\begin{equation}\label{f-converges}
f_{\Kphi_k}\to f_\Kphi\qquad\mbox{in } C^{1,\beta}([-v_2, \mxy_1]),
\end{equation}
where $f_{\Kphi_k},f_\Kphi\in C^{(-1-\alpha, \{-v_2,
\mxy_1\})}_{2,\alpha,(-v_2, \mxy_1)}$ are the functions from
 (\ref{shockPL}) corresponding to $\Kphi_k,\Kphi$, respectively.
Let $\psi_k, \psi\in C^{1, \alpha}(\overline{\Omega^+(\Kphi_k)})\cap
C^{2,
\alpha}(\overline{\Omega^+(\Kphi_k)}\setminus\overline{\sonic\cup\Sigma_0})$
be the solutions of problems
 (\ref{iterationEquation})--(\ref{iterationCondOnSymmtryLine})
for $\Kphi_k, \Kphi$. Let $\{\psi_{k_m}\}$ be any subsequence of
$\{\psi_k\}$. By (\ref{Hold-Ext-Sonic})--(\ref{Hold-Ext-Unif}), it
follows that there exist a further subsequence $\{\Kphi_{k_{m_n}}\}$
and a function $\bar\psi\in C^{1, \alpha}(\overline\Dom)\cap C^{2,
\alpha}(\Dom)$ such that
$$
\ExCFn_{\Kphi_{k_{m_n}}}\psi_{k_{m_n}}\to \bar\psi \qquad\mbox{in
$C^{2,{\alpha/2}}$ on compact subsets of $\Dom$
 and in
$C^{1,{\alpha/2}}(\overline\Dom)$}.
$$
Then, using (\ref{f-converges}) and the convergence $\Kphi_k \to
\Kphi$ in $C^{1, \beta}(\overline{\Dom})$, we prove (by the argument
as in \cite[page 479]{ChenFeldman1}) that $\bar\psi$ is a solution
of problem
(\ref{iterationEquation})--(\ref{iterationCondOnSymmtryLine}) for
$\Kphi$. By uniqueness in Lemma \ref{uniquenessCor}, $\bar\psi=\psi$
in $\Omega^+(\Kphi)$. Now, using (\ref{f-converges}) and the
explicit definitions of extensions $\ExCFn_\Kphi^1$ and
$\ExCFn_\Kphi^2$, it follows by the argument as in \cite[pp.
477--478]{ChenFeldman1} that
\begin{eqnarray*}
&&\zeta\ExCFn_{\Kphi_{k_{m_n}}}^1(\psi_{k_{m_n}})\to
\zeta\ExCFn_\Kphi^1(\bar\psi_{|\Omega^+(\Kphi)})\qquad \mbox{ in
} C^{1,\beta}(\overline\Dom),\\
&&(1-\zeta)\ExCFn_{\Kphi_{k_{m_n}}}^2(\psi_{k_{m_n}})\to
(1-\zeta)\ExCFn_\Kphi^2(\bar\psi_{|\Omega^+(\Kphi)})\qquad \mbox{ in
} C^{1,\beta}(\overline\Dom).
\end{eqnarray*}
Therefore,
$\bar\psi=\psi$ in $\Dom$. Since a convergent subsequence
$\{\psi_{k_{m_n}}\}$ can be extracted from any subsequence
$\{\psi_{k_{m}}\}$ of $\{\psi_k\}$ and the limit $\bar\psi=\psi$ is
independent of the choice of subsequences $\{\psi_{k_{m}}\}$ and
$\{\psi_{k_{m_n}}\}$, it follows that the whole sequence $\psi_k$
converges to $\psi$ in $C^{1,\beta}(\overline\Dom)$. This completes
the proof. \Endproof
%*************** END PROOF *****LEMMA**********

Now we denote by $\hat C_0$ the constant in (\ref{condConst-00})
sufficiently large to satisfy the conditions of Proposition
\ref{existSolDegenEllipt} and Lemma \ref{extension-Lemma}. Fix $\hat
C\ge \hat C_0$. Choose $M_1=\max(2C_1, 1)$ for the constant $C_1$ in
(\ref{Hold-Ext-Sonic}) and define $\varepsilon$ by
(\ref{condConst-1-1}). This choice of  $\varepsilon$ fixes the
constant $C_2(\varepsilon)$ in (\ref{Hold-Ext-Unif}). Define
$M_2=\max(C_2(\varepsilon),1)$. Finally, let
$$
\epsP_0={{\hat C}^{-1}-\varepsilon-\varepsilon^{1/4}M_1 \over
2\left(M_2^2+ \varepsilon^2\max(M_1,M_2)\right)}\varepsilon^2.
$$
Then $\epsP_0>0$, since $\varepsilon$ is defined by
(\ref{condConst-1-1}). Moreover, $\epsP_0$, $\varepsilon$, $M_1$,
and $M_2$ depend only on the data and $\hat C$. Furthermore, for any
$\epsP\in[0, \epsP_0]$, the constants  $\epsP$, $\varepsilon$,
$M_1$, and $M_2$ satisfy (\ref{condConst-00}) with  $\hat C$ fixed
above. Also, $\psi\ge 0$ on $\Omega^+(\Kphi)$ by
(\ref{L-infty-for-unif-ellipt}) and thus
\begin{equation}\label{Ext-positive}
\ExCFn_\Kphi\psi\ge 0
\qquad\mbox{on }\Dom
\end{equation}
by the explicit definitions of $\ExCFn_\Kphi^1, \ExCFn_\Kphi^2$, and
$\ExCFn_\Kphi$. Now we define the iteration map $J$ by
$J(\Kphi)=\ExCFn_\Kphi\psi$. By
(\ref{Hold-Ext-Sonic})--(\ref{Hold-Ext-Unif}) and
(\ref{Ext-positive})
 and the choice
of  $\epsP$, $\varepsilon$, $M_1$, and $M_2$, we find that
$J:\setK\to\setK$. Now, $\setK$ is a compact and convex subset of
$C^{1,{\alpha/2}}(\overline\Dom)$. The map $J:\setK\to\setK$ is
continuous in $C^{1,{\alpha/2}}(\overline\Dom)$ by Lemma
\ref{extension-Lemma}(\ref{Uextended_Cont}). Thus, by the Schauder
Fixed Point Theorem, there exists a fixed point $\Kphi\in\setK$ of
the map $J$. By definition of $J$, such $\psi$ is a solution of
(\ref{iterationEquation})--(\ref{iterationCondOnSymmtryLine}) with
$\Kphi=\psi$. Therefore, we have

%*****************BEGIN PROPOSITION*************
\begin{proposition}\label{existenceFixedPt}
There exists $\hat C_0\ge 1$ depending only on the data such that,
for any $\hat C\ge \hat C_0$, there exist $\epsP_0, \varepsilon>0$
and $M_1, M_2\ge 1$ satisfying {\rm (\ref{condConst-00})} so that,
for any $\epsP\in (0, \epsP_0]$,  there exists a solution
$\psi\in\setK(\epsP, \varepsilon, M_1, M_2)$
  of problem {\rm (\ref{iterationEquation})}--{\rm
(\ref{iterationCondOnSymmtryLine})} with $\Kphi=\psi$ (i.e., $\psi$
is a ``fixed point'' solution). Moreover, $\psi$ satisfies {\rm
(\ref{Hoder-est-for-unif-ellipt})} for all $s\in (0, c_2/2)$ with
$C(s)$ depending only on the data and $s$.
\end{proposition}
%*****************END PROPOSITION*************

%*****************BEGIN SECTION*************
%*****************************************888
\section{Removal of the Ellipticity Cutoff}
\label{removeCutoffSection} In this section we assume that $\hat
C_0\ge 1$ is as in Proposition \ref{existenceFixedPt} which depends
only on the data, $\hat C\ge \hat C_0$, and assume that $\epsP_0,
\varepsilon>0$ and $M_1,M_2\ge 1$ are defined by $\hat C$ as in
Proposition \ref{existenceFixedPt} and $\epsP\in (0, \epsP_0]$. We
fix a ``fixed point'' solution $\psi$ of problem
(\ref{iterationEquation})--(\ref{iterationCondOnSymmtryLine}), that
is, $\psi\in\setK(\epsP, \varepsilon, M_1, M_2)$ satisfying
(\ref{iterationEquation})--(\ref{iterationCondOnSymmtryLine}) with
 $\Kphi=\psi$. Its existence is established in  Proposition \ref{existenceFixedPt}.
To simplify notations, in this section we write $\Omega^+$,
$\shock$, and $\Sigma_0$ for $\Omega^+(\psi)$, $\shock(\psi)$, and
$\Sigma_0(\psi)$, respectively, and the universal constant $C$
depends only on the data.

We now prove that the ``fixed point'' solution $\psi$ satisfies
$|\psi_x|\le 4x/\big(3(\gamma+1)\big)$ in $\Omega^+\cap
\{c_2-r<4\varepsilon\}$ for sufficiently large $\hat C$, depending
only on the data, so that $\psi$ is a solution of the regular
reflection problem; see Step 10 of \S
\ref{overViewProcedureSubsection}.

We also note the higher regularity of $\psi$ away from the corners
and the sonic circle. Since equation (\ref{iterationEquation}) is
uniformly elliptic in every compact subset of $\Omega^+$ (by Lemma
\ref{propertiesNonlinCoeffs}) and the coefficients
$A_{ij}(p,\mxx,\mxy)$ of (\ref{iterationEquation}) are $C^{1,\alpha}$
functions of $(p,\mxx,\mxy)$ in every compact subset of
$\bR^2\times\Omega^+$ (which follows from the explicit expressions
of $A_{ij}(p,\mxx,\mxy)$ given by (\ref{iterationUniforDomEquation}),
(\ref{iterationSonicDomEquation}), and (\ref{combineCoeffs})), then
substituting $p=\grad \psi(\mxx,\mxy)$ with $\psi\in\setK$ into
$A_{ij}(p,\mxx,\mxy)$, rewriting (\ref{iterationEquation}) as a linear
equation with coefficients being $C^{1,\alpha}$  in compact subsets
of $\Omega^+$, and using the interior regularity results for linear,
uniformly elliptic equations yield
\begin{equation}\label{higherRegOfPsi}
\psi\in C^{3,\alpha}\left( \Omega^+\right).
\end{equation}

First we  bound $\psi_x$ from above. We work in the
$(x,y)$--coordinates in $\Omega^+\cap \{c_2-r<4\varepsilon\}$. By
(\ref{domain-in-rescaled-lemma}),
\begin{equation}\label{domain-in-rescaled-lemma-remCutoff}
\Omega^+(\Kphi)\cap\{c_2-r<4\varepsilon\} =\{0< x< \kappa,\,\,
0<y<\hat f_\Kphi(x)\},
\end{equation}
where $\hat f_\Kphi$ satisfies (\ref{holder-hat-f}).
%*****************BEGIN PROPOSITION*************
\begin{proposition}
\label{boundPsiXfromAbove-Prop}
For sufficiently large $\hat C$
depending only on the data,
\begin{equation}\label{boundPsiXfromAbove}
\psi_x\le {4\over 3(\gamma+1)}x \qquad\mbox{ in }\;
\Omega^+\cap\{x\le 4\varepsilon\}.
\end{equation}
\end{proposition}
%*****************END PROPOSITION*************
%*************** BEGIN PROOF ***************
\Proof
To simplify notations, we denote $A=\frac{4}{3(\gamma+1)}$ and
$$
\Omega^+_s\defd \Omega^+\cap\{x\le s\} \qquad\mbox{for}\,\, s>0.
$$
Define a function
\begin{equation}\label{defV-remove-1}
v(x,y):=Ax-\psi_x(x,y) \qquad\mbox{on }\Omega^+_{4\varepsilon}.
\end{equation}
{}From $\psi\in\setK$ and (\ref{higherRegOfPsi}), it follows that
\begin{equation}\label{v-regular-remCutoff}
v\in C^{0,1}\big(\overline{\Omega^+_{4\varepsilon}}\big) \cap
C^1\big(\overline{\Omega^+_{4\varepsilon}}\setminus\{x=0\}\big) \cap
C^2\left(\Omega^+_{4\varepsilon}\right).
\end{equation}
Since $\psi\in\setK$, we have $|\psi_x(x, y)|\le M_1x$ in
$\Omega^+_{4\varepsilon}$. Thus
\begin{equation}\label{bdryCondV-sonic}
v=0\qquad\mbox{on }\, \partial\Omega^+_{4\varepsilon}\cap\{x=0\}.
\end{equation}
We now use the fact that $\psi$ satisfies (\ref{iterationRH}), which
can be written as (\ref{iterationRH-lf-flattened}) in the
$(x,y)$--coordinates, and (\ref{estCoefsIterRH-flattened}) holds.
Since $\psi\in\setK$ implies that
$$
|\psi(x,y)|\le M_1x^2, \qquad |\psi_y(x,y)|\le M_1x^{3/2},
$$
it
follows  from (\ref{iterationRH-lf-flattened}) and
(\ref{estCoefsIterRH-flattened}) that
$$
|\psi_x|\le C(|\psi_y|+|\psi|)\le CM_1x^{3/2}\qquad\mbox{ on } \shock\cap\{x<2\varepsilon\},
$$
and hence, by (\ref{condConst-00}), if $\hat C$ is large depending
only on the data, then
$$
|\psi_x|<A x\qquad\mbox{ on } \shock\cap\{0<x<2\varepsilon\}.
$$
Thus we have
\begin{equation}\label{bdryCondV-shock}
v\ge 0\qquad\mbox{on }\shock\cap\{0<x<2\varepsilon\}.
\end{equation}
Furthermore, condition (\ref{iterationCondOnWedge}) on   $\wedgeB$
in the $(x,y)$--coordinates is
$$
\psi_y=0\qquad\mbox{ on }\; \{0<x<2\varepsilon, \, y=0\}.
$$
Since $\psi\in\setK$ implies that $\psi$ is $C^2$ up to $\wedgeB$,
then differentiating the condition on  $\wedgeB$ with respect to
$x$, i.e., in the tangential direction to  $\wedgeB$, yields
$\psi_{xy}=0$ on $\{0<x<2\varepsilon,\, y=0\}$, which implies
\begin{equation}\label{bdryCondV-wedge}
v_y= 0\qquad\mbox{on }\wedgeB\cap\{0<x<2\varepsilon\}.
\end{equation}
Furthermore, since $\psi\in\setK$,
\begin{equation}\label{sigmaSmall}
|\psi_x|\le M_2\epsP\le A\varepsilon \qquad\mbox{on
}\; \Omega^+\cap\{\varepsilon/2\le x\le 4\varepsilon\},
\end{equation}
where the second inequality holds by (\ref{condConst-00}) if $\hat
C$ is large, depending only on the data. Thus, for such  $\hat C$,
\begin{equation}\label{bdryCondV-upbdry}
v\ge 0\qquad\mbox{on
}\Omega^+_{4\varepsilon}\cap\{x=2\varepsilon\}.
\end{equation}

Now we show that, for large $\hat C$,
 $v$ is a supersolution of a linear homogeneous elliptic
equation on  $\Omega^+_{2\varepsilon}$.
Since $\psi$ satisfies
equation (\ref{iteration-equation-sonicStruct}) with
(\ref{iteration-equation-sonicStruct-coef}) in
 $\Omega^+_{4\varepsilon}$,
we differentiate the equation with respect to $x$ and use the
regularity of $\psi$ in  (\ref{higherRegOfPsi}) and the definition
$v$ in (\ref{defV-remove-1}) to obtain
\begin{equation}\label{eqnForV-diff-1}
\begin{array}{l}
\displaystyle a_{11}v_{xx}+a_{12}v_{xx}+ a_{22}v_{yy}\\
\quad +\left(A-v_x\right) \big(-1 +(\gamma+1) \big(\zeta_1(A-{v\over
x}) + \zeta_1'(A-{v\over x})({v\over x}-v_x) \big)\big)=E(x,y),
\end{array}
\end{equation}
where
\begin{eqnarray}
&&a_{11}=2x-(\gamma+1)x\zeta_1\big({\psi_x\over x}\big)+\hat O_1,
\quad a_{12}=\hat O_2,
\quad a_{22}={1\over c_2}+\hat O_3, \label{defA22-rem_cutoff-1}
\\
&& E(x,y)=\psi_{xx}\partial_x\hat O_1 +\psi_{xy}\partial_x\hat O_2
+\psi_{yy}\partial_x\hat O_3 -\psi_{xx}\hat O_4
-\psi_{x}\partial_x\hat O_4 \label{defE-rem_cutoff-1}
\\
&& \qquad\qquad +\psi_{xy}\hat O_5 +\psi_{y}\partial_x\hat O_5,
\nonumber
\end{eqnarray}
with
\begin{equation}
\label{defO-rem_cutoff} \hat O_k(x,y)=O_k^\psi(D\psi(x,y),x,y)\qquad
\mbox{ for }\;k=1,\dots, 5,
\end{equation}
for $O_k^\psi$ defined by
(\ref{iteration-equation-sonicStruct-coef}) with $\Kphi=\psi$.
{}From (\ref{defZeta-1}), we have
$$
\zeta_1\left(A\right) = A.
% \qquad\mbox{for }\;A=\frac{4}{3(\gamma+1)}.
$$
Thus we can rewrite (\ref{eqnForV-diff-1}) in the form
\begin{equation}\label{eqnForV-diff-3}
a_{11}v_{xx}+a_{12}v_{xx}+ a_{22}v_{yy}+
bv_x+cv=-A\big((\gamma+1)A-1\big)+E(x,y),
\end{equation}
with
\begin{eqnarray}
&&b(x,y)=1 -(\gamma+1) \big(\zeta_1(A-{v\over x})+
\zeta_1'(A-{v\over x}) ({v\over x}-v_x-A) \big),
\label{def-b-rem_cutoff-1}
\\
&&c(x,y)=(\gamma+1){A\over x}\big(\zeta_1'(A-{v\over x})-
\int_0^1\zeta_1'(A-s{v\over x})ds \big), \label{def-c-rem_cutoff-1}
\end{eqnarray}
where $v$ and $v_x$ are evaluated at the point $(x,y)$.

Since $\psi\in\setK$ and $v$ is defined by (\ref{defV-remove-1}), we
have
$$
a_{ij}, b, c \in
C\big(\overline{\Omega^+_{4\varepsilon}}\setminus\{x=0\}\big).
$$

Combining (\ref{defA22-rem_cutoff-1}) with (\ref{condConst-00}),
(\ref{defZeta-1}), (\ref{estSmallterms-iter-S}), and
(\ref{defO-rem_cutoff}), we obtain that, for sufficiently large
$\hat C$  depending only on the data,
$$
a_{11}\ge{1\over 6}x, \qquad a_{22}\ge{1\over 2c_2}, \qquad
|a_{12}|\le{1\over 3\sqrt{c_2}}x^{1/2} \qquad\mbox{on }\,
\Omega^+_{2\varepsilon}.
$$
Thus, $4a_{11}a_{22}-(a_{12})^2\ge{2\over 9c_2}x$ on
$\Omega^+_{2\varepsilon}$, which implies that equation
(\ref{eqnForV-diff-3}) is elliptic on
 $\Omega^+_{2\varepsilon}$ and uniformly elliptic on every compact subset
of $\overline{\Omega^+_{2\varepsilon}}\setminus\{x=0\}$.

Furthermore, using (\ref{zeta-1-concave}) and
(\ref{def-c-rem_cutoff-1}) and noting $A>0$ and $x>0$, we have
\begin{equation}
\label{cNegativeAtVneg} c(x,y)\le 0\qquad\mbox{for every
}\;(x,y)\in\Omega^+_{2\varepsilon}\;\mbox{ such that } v(x,y)\le 0.
\end{equation}

Now we estimate $E(x,y)$. Using (\ref{defO-rem_cutoff}),
(\ref{iteration-equation-sonicStruct-coef}),
(\ref{erTerms-xy-nontrunc}), and $\psi\in\setK$, we find that, on
$\Omega^+_{2\varepsilon}$,
\begin{eqnarray*}
 |\partial_x\hat O_1|&\le&
C\big(x+|\psi|+|D\psi|+x|\psi_{xx}|+|\psi_x\psi_{xx}|+|\psi_y\psi_{xy}|+|D\psi|^2\big)
\le CM_1^2x,
\\
|\partial_x\hat O_{2,5}|&\le&
C\big(|D\psi|+|D\psi|^2+|\psi_y\psi_{xx}|+(1+|\psi_x|)|\psi_{xy}|\big)
\le CM_1x^{1/2}(1+M_1x),
\\
|\partial_x\hat O_{3,4}|&\le& C\big(1+|\psi|+ \big|{\psi_x\over
x}\zeta_1'\big({\psi_x\over x}\big)\big| +(1+|D\psi|)|D^2\psi|
+|D\psi|^2\big) \\
&\le& CM_1(1+M_1x),
\end{eqnarray*}
where we have used the fact that $|s\zeta_1'(s)|\le C$ on $\bR$.
Combining these estimates with
(\ref{defE-rem_cutoff-1})--(\ref{defO-rem_cutoff}),
(\ref{estSmallterms-iter}), and $\psi\in\setK$, we obtain from
(\ref{defE-rem_cutoff-1}) that
$$
|E(x,y)|\le CM_1^2x(1+M_1x)\le {C/\hat C}\qquad \mbox{ on
}\Omega^+_{2\varepsilon}.
$$
{}From this and $(\gamma+1)A>1$, we conclude that the right-hand
side of (\ref{eqnForV-diff-3}) is strictly negative in
$\Omega^+_{2\varepsilon}$ if $\hat C$ is sufficiently large,
depending only on the data.

We fix $\hat C$ satisfying all the requirements above (thus
depending only on the data). Then we have
 \begin{equation}\label{eqnForV-diff-4}
a_{11}v_{xx}+a_{12}v_{xx}+ a_{22}v_{yy}+ bv_x+cv<0 \qquad\mbox{on
}\Omega^+_{2\varepsilon},
\end{equation}
the equation is elliptic in $\Omega^+_{2\varepsilon}$ and uniformly
elliptic on compact subsets of
$\overline{\Omega^+_{2\varepsilon}}\setminus \{x=0\}$, and
(\ref{cNegativeAtVneg}) holds. Moreover, $v$ satisfies
(\ref{v-regular-remCutoff}) and the boundary conditions
(\ref{bdryCondV-sonic})--(\ref{bdryCondV-wedge}) and
(\ref{bdryCondV-upbdry}). Then it follows that
$$
v\ge 0\qquad \mbox{ on }\, \Omega^+_{2\varepsilon}.
$$
Indeed, let $z_0:=(x_0, y_0)\in \overline{\Omega^+_{2\varepsilon}}$
be a minimum point of $v$ over $\overline{\Omega^+_{2\varepsilon}}$
and $v(z_0)<0$. Then, by
(\ref{bdryCondV-sonic})--(\ref{bdryCondV-shock}) and
(\ref{bdryCondV-upbdry}), either $z_0$ is an interior point of
$\Omega^+_{2\varepsilon}$ or $z_0\in
\wedgeB\cap\{0<x<2\varepsilon\}$. If $z_0$ is an interior point of
$\Omega^+_{2\varepsilon}$, then (\ref{eqnForV-diff-4}) is violated
since  (\ref{eqnForV-diff-4}) is elliptic, $v(z_0)<0$, and
$c(z_0)\le 0$ by (\ref{cNegativeAtVneg}). Thus, the only possibility
is $z_0\in \wedgeB\cap\{0<x<2\varepsilon\}$, i.e., $z_0=(x_0, 0)$
with $x_0>0$. Then, by (\ref{domain-in-rescaled-lemma-remCutoff}),
there exists $\rho>0$ such that
 $B_{\rho}(z_0)\cap \Omega^+_{2\varepsilon}=B_{\rho}(z_0)\cap\{y>0\}$.
Equation (\ref{eqnForV-diff-4}) is uniformly elliptic in
$\overline{B_{\rho/2}(z_0)\cap\{y\ge0\}}$, with the coefficients
$a_{ij}, b, c\in C(\overline{B_{\rho/2}(z_0)\cap\{y\ge0\}})$. Since
$v(z_0)<0$ and $v$ satisfies (\ref{v-regular-remCutoff}), then,
reducing $\rho>0$ if necessary, we have $v<0$ in
$B_{\rho}(z_0)\cap\{y>0\}$. Thus, $c\le 0$ on
$B_{\rho}(z_0)\cap\{y>0\}$ by (\ref{cNegativeAtVneg}). Moreover,
$v(x,y)$ is not a constant in
$\overline{B_{x_0/2}(x_0)\cap\{y\ge0\}}$ since its negative minimum
is achieved at $(x_0,0)$ and cannot be achieved in any interior
point, as we have showed above. Thus, $\partial_y v(z_0)>0$ by
Hopf's Lemma, which contradicts (\ref{bdryCondV-wedge}). Therefore,
$v\ge 0$ on $\Omega^+_{2\varepsilon}$ so that
(\ref{boundPsiXfromAbove}) holds on $\Omega^+_{2\varepsilon}$. Then,
using (\ref{sigmaSmall}), we obtain (\ref{boundPsiXfromAbove}) on
$\Omega^+_{4\varepsilon}$.
\Endproof
%*************** END PROOF ***************

Now we  bound $\psi_x$ from below. We first prove the following
lemma in the $(\mxx, \mxy)$--coordinates.

%*****************BEGIN LEMMA*************
\begin{lemma}\label{negativeDerivPsiLemma}
If $\hat C$ in {\rm (\ref{condConst-00})} is sufficiently large,
depending only on the data, then
\begin{equation}\label{negativeDerivPsi}
\psi_{\mxy}\le 0 \qquad\mbox{in } \Omega^+.
\end{equation}
\end{lemma}
%*****************END LEMMA*************
\Proof We divide the proof into six steps.
%*************** BEGIN PROOF *****LEMMA**********

{\em Step 1}. Set
$w=\psi_{\mxy}.$
{}From $\psi\in\setK$ and (\ref{higherRegOfPsi}),
\begin{equation}\label{w-regular-negativeDeriv}
w\in C^{0,\alpha}\big(\overline{\Omega^+}\big) \cap
C^1\big(\overline{\Omega^+}\setminus\overline{\sonic\cup
\Sigma_0}\big) \cap C^2\left(\Omega^+\right).
\end{equation}

In the next steps, we derive the equation and boundary conditions
for $w$ in $\Omega^+$. To achieve this, we use the following facts:

(i) If $\hat C$ in (\ref{condConst-00}) is sufficiently large, then
the coefficient $A_{11}$ of  (\ref{iterationEquation}) satisfies
\begin{equation}\label{A11-bdd-away-fromZero}
|A_{11}\left(\grad\psi(\mxx,\mxy), \mxx,\mxy \right)|\ge
{{\bar{c}_2}^2-\bar\mxx^2\over 2}>0 \qquad\mbox{in } \Omega^+,
\end{equation}
where $\bar{c}_2$ and $\bar\xi$ are defined in \S\ref{section:4}.
Indeed, since $\bar{c}_2>|\bar\xi|$ by
(\ref{sonic-intersect-shock-normal}) and $(c_2,\tilde\mxx)\to (
\bar{c}_2,\bar\mxx)$ as $\theta_w\to \pi/2$ by \S \ref{section:3.3},
we have $\displaystyle c_2^2-\tilde\mxx^2\ge
{9({\bar{c}_2}^2-\bar\mxx^2)/10}>0$ if $\epsP$ is small.
Furthermore, for any $(\mxx, \mxy)\in\Dom$, we have
$c_2\cos\theta_w\ge\mxx\ge \tilde\mxx$ and thus, assuming that
$\sigma$ is small so that $|\tilde\mxx|\le 2|\bar\xi|$ and $c_2\le
2\bar{c}_2$, we obtain $|\mxx|\le C$. Now, since  $\psi\in\setK$, it
follows that, if $\hat C$ in (\ref{condConst-00}) is sufficiently
large, then $A_{11}^1$ defined in (\ref{iterationUniforDomEquation})
 with $\Kphi=\psi$ implies
$A_{11}^1\ge {({\bar{c}_2}^2-\bar\mxx^2)/ 2}$   on $\Dom$, and
$A_{11}^2$ in (\ref{iterationSonicDomEquation}) with $\Kphi=\psi$
implies $A_{11}^2\ge {({\bar{c}_2}^2-\bar\mxx^2)/ 2}$ on $\Dom\cap
\{c_2-r<4\varepsilon\}$. Then (\ref{A11-bdd-away-fromZero}) follows
from (\ref{combineCoeffs}).

(ii) Since $\psi$ satisfies  equation (\ref{iterationEquation}) in
$\Omega^+$ with (\ref{A11-bdd-away-fromZero}), we have
\begin{equation}\label{express2ndDeriv}
\psi_{\mxx\mxx}=-{2\hat A_{12}\psi_{\mxx\mxy}+\hat
A_{22}\psi_{\mxy\mxy}\over \hat A_{11}} \qquad\mbox{ in }\Omega^+,
\end{equation}
where $\hat A_{ij}(\mxx, \mxy)=A_{ij}\left(\grad\psi(\mxx,\mxy), \mxx,\mxy
\right)$ in $\Omega^+$.

\medskip
{\em Step 2.} We differentiate equation (\ref{iterationEquation})
with respect to $\mxy$ and substitute the right-hand side of
(\ref{express2ndDeriv}) for $\psi_{\mxx\mxx}$ to obtain the following
equation for $w$:
\begin{equation}\label{negativeDerivPsi-eqnForW}
\hat A_{11}w_{\mxx\mxx}+2\hat A_{12}w_{\mxx\mxy}+\hat A_{22}w_{\mxy\mxy}
+ 2\big(\partial_\mxy \hat A_{12}-{\partial_\mxy \hat A_{11}\over \hat
A_{11}}
 \hat A_{12}\big) w_\mxx
+\big(\partial_\mxy \hat A_{22}-{\partial_\mxy \hat A_{11}\over \hat
A_{11}} \hat A_{22}\big) w_\mxy =0.
%\label{negativeDerivPsi-eqnForW}
\end{equation}
By Lemma \ref{propertiesNonlinCoeffs},
(\ref{A11-bdd-away-fromZero}), and $\psi\in\setK$, the coefficients
of (\ref{negativeDerivPsi-eqnForW}) are  continuous in
$\overline{\Omega^+}\setminus \overline{\sonic\cup \Sigma_0}$, and
 the equation is uniformly elliptic on compact subsets of
$\overline{\Omega^+}\setminus \overline{\sonic}$.

\medskip
{\em Step 3.} By (\ref{iterationCondOnSymmtryLine}), we have
\begin{equation}\label{negativeDerivPsi-CondSymmetr-ForW}
w=-v_2\qquad\mbox{on }\;\Sigma_0:=\partial \Omega^+\cap
\{\mxy=-v_2\}.
\end{equation}

Since $\psi\in\setK$, it follows that
$|\grad\psi(\mxx, \mxy)|\le CM_1(c_2-r)$ for
all $(\mxx, \mxy)\in\Omega^+\cap\{c_2-r\le 2\varepsilon\}$.
Thus,
\begin{equation}\label{negativeDerivPsi-CondSonic-ForW}
w=0\qquad\mbox{on }\;\sonic.
\end{equation}

\medskip
{\em Step 4.} We derive the boundary condition for $\psi$ on
$\wedgeB$. Then $\psi$ satisfies (\ref{iterationCondOnWedge}), which
can be written as
\begin{equation}\label{iterationCondOnWedge-2}
-\sin\theta_w\,\psi_\mxx+ \cos\theta_w\,\psi_\mxy=0\qquad\,\, \mbox{on
}\;\wedgeB.
\end{equation}
Since $\psi\in\setK$, we have $\psi\in C^2(\overline{\Omega^+}
\setminus\overline{\sonic\cup\Sigma_0})$. Thus we can differentiate
(\ref{iterationCondOnWedge-2}) in the direction tangential to
$\wedgeB$, i.e., apply
%the differential operator
$\partial_\tau:=\cos\theta_w\,\partial_\mxx+
\sin\theta_w\,\partial_\mxy$ to (\ref{iterationCondOnWedge-2}).
Differentiating and substituting the right-hand side of
(\ref{express2ndDeriv}) for $\psi_{\mxx\mxx}$, we have
\begin{equation}\label{negativeDerivPsi-CondWedge-ForW}
\displaystyle \big(\cos(2\theta_w) +{\hat A_{12}\over \hat
A_{11}}\sin(2\theta_w)\big)w_\mxx +\frac{1}{2}\sin(2\theta_w)\big(
1+{\hat A_{22}\over \hat A_{11}}\big)w_\mxy =0 \qquad\mbox{on
}\;\wedgeB.
\end{equation}
This condition is oblique if $\epsP$ is small: Indeed, since the
unit normal on $\wedgeB$ is $(-\sin\theta_w, \cos\theta_w)$, we use
(\ref{angleCloseToPiOver2}) and (\ref{A11-bdd-away-fromZero}) to
find
$$
(\cos(2\theta_w) +{\hat A_{12}\over \hat A_{11}}\sin(2\theta_w),
\frac{1}{2}\sin(2\theta_w)(1+{\hat A_{22}\over \hat A_{11}}))\cdot (
-\sin\theta_w, \cos\theta_w)) \ge 1-C\epsP\ge {1\over 2}.
$$

\medskip
{\em Step 5.} In this step, we derive the condition for $w$ on
$\shock$. Since $\psi$ is a solution of
(\ref{iterationEquation})--(\ref{iterationCondOnSymmtryLine}) for
$\Kphi=\psi$, the Rankine-Hugoniot conditions hold on $\shock$:
 Indeed, the continuous matching of
$\psi$ with $\varphi_1-\varphi_2$ across $\shock$ holds by
(\ref{shockPL})--(\ref{shockIterDef}) since $\Kphi=\psi$. Then
(\ref{cont-accross-shock-psi-resolved}) holds and the
 gradient jump condition (\ref{RH-psi}) can be written in form
(\ref{RH-psi-2}). On the other hand, $\psi$ on $\shock$ satisfies
(\ref{iterationRH}) with $\Kphi=\psi$, which is (\ref{RH-psi-2}).
Thus, $\psi$ satisfies (\ref{RH-psi}).

Since $\psi\in\setK$ which implies $\psi\in C^2(\overline\Omega^+
\setminus\overline{\sonic\cup\Sigma_0})$,  we can differentiate
(\ref{RH-psi})
 in the direction tangential to $\shock$.
The unit normal $\nu_s$ on  $\shock$ is given by
(\ref{norm-to-Shock}). Then the vector
\begin{equation}\label{tangent-to-Shock-2}
\tau_s\equiv (\tau_s^1,\tau_s^2) :=({v_2+\psi_\mxy\over
u_1-u_2},\; 1-{\psi_\mxx\over u_1-u_2})
\end{equation}
is tangential to $\shock$. Note that $\tau_s\ne 0$ if $\hat C$ in
(\ref{condConst-00}) is sufficiently large, since
\begin{equation}\label{smallnessPsi-u-v}
|\grad\psi|\le C(\epsP+\varepsilon)\quad \mbox{ in }\;
\overline{\Omega^+},\qquad |u_2|+|v_2|\le C\epsP,
\end{equation}
and $u_1>0$ from $\psi\in\setK$ and \S \ref{section:3.3}. Thus, we
can apply the differential operator
$\partial_{\tau_s}=\tau_s^1\partial_\mxx + \tau_s^2\partial_\mxy$ to
(\ref{RH-psi}).

In the calculation below, we use the notation in \S
\ref{equationForPsiSection}. We have showed in \S
\ref{equationForPsiSection} that
 condition (\ref{RH-psi}) can be written in form (\ref{RH-psi-int1}),
where $F(p, z, u_2, v_2, \mxx, \mxy)$ is defined by
(\ref{RH-psi-func1})--(\ref{RH-psi-func3}) and
satisfies (\ref{rewriteRH-reg-1}). Also, we denote
\begin{equation}\label{tangent-to-Shock-3}
\hat\tau(p, u_2, v_2)\equiv (\hat\tau^1,\hat\tau^2)(p, u_2, v_2)
:=({v_2+p_2\over
u_1-u_2},\; 1-{p_1\over u_1-u_2}),
\end{equation}
where  $p=(p_1,p_2)\in\bR^2$ and $z\in\bR$. Then $\hat\tau\in
C^\infty(\overline{B_{\delta^*}(0)\times B_{u_1/50}(0)})$. Now,
applying the differential operator $\partial_{\tau_s}$, we obtain
that $\psi$ satisfies
\begin{equation}\label{diff-RH-1}
\Phi(D^2\psi, \grad\psi, \psi, u_2, v_2, \mxx, \mxy)=0 \qquad\mbox{on}\;\shock,
\end{equation}
where
\begin{equation}\label{Diff-RH-psi-func1}
\Phi(R, p, z, u_2, v_2, \mxx, \mxy)= \sum_{i,j=1}^2\hat
\tau^iF_{p_j}R_{ij} +\sum_{i=1}^2\hat\tau^i(F_{z}p_i+F_{\xi_i})
\qquad \mbox{for}\,\, R=(R_{ij})_{i,j=1}^2,
\end{equation}
and, in both (\ref{Diff-RH-psi-func1}) and  the calculation below,
$D_{(\xi_1, \xi_2)}F$ denotes as $D_{(\xi, \eta)}F$, $\;(F_{p_j},
F_{z}, F_{\xi_i})$ as $(F_{p_j}, F_{z}, F_{\xi_i})(p, z, u_2, v_2,
\mxx, \mxy)$, $\;(\hat\tau, \hat\nu)$ as $(\hat\tau, \hat\nu)(p,
u_2, v_2)$, and $\tilde\rho$ as $\tilde\rho(p, z, \mxx, \mxy)$, with
$\tilde\rho(\cdot)$ and $\hat\nu(\cdot)$ defined by
(\ref{RH-psi-func2}) and (\ref{RH-psi-func3}), respectively. By
explicit calculation, we apply
(\ref{RH-psi-func1})--(\ref{RH-psi-func3}) and
(\ref{tangent-to-Shock-3}) to obtain that, for every $(p, z, u_2,
v_2, \mxx, \mxy)$,
\begin{equation}\label{Diff-RH-lower-cancel}
\sum_{i=1}^2\hat\tau^i\, (F_{z}p_i+F_{\xi_i})
=(\rho_1-\tilde\rho)\,\hat\tau\cdot \hat\nu =0.
\end{equation}

We note that (\ref{cont-accross-shock-psi-resolved}) holds on
$\shock$. Using (\ref{diff-RH-1}) and (\ref{Diff-RH-lower-cancel})
and expressing $\xi$ from (\ref{cont-accross-shock-psi-resolved}),
we see that $\psi$ satisfies
\begin{equation}\label{diff-RH-2}
\tilde\Phi(D^2\psi, \grad\psi, \psi, u_2, v_2, \mxy)=0 \qquad\mbox{on}\;\shock,
\end{equation}
where
\begin{equation}\label{Diff-RH-psi-func2}
\tilde\Phi(R, p, z, u_2, v_2,  \mxy)= \sum_{i,j=1}^2\hat
\tau^i\Psi_{p_j}(p, z, u_2, v_2,  \mxy)R_{ij},
\end{equation}
$\Psi$ is defined by (\ref{RH-psi-func4}) and satisfies
 $\Psi\in
C^\infty(\overline{\mathcal A})$ with
$\|\Psi\|_{C^k(\overline{\mathcal A})}$ depending only on the data
and $k\in{\bf N}$, and ${\mathcal A}=B_{\delta^*}(0)\times
(-\delta^*, \delta^*)\times B_{u_1/50}(0)\times (-6\bar c_2/5, 6\bar
c_2/5)$.

Now, from (\ref{RH-psi-func1})--(\ref{RH-psi-func3}),
(\ref{RH-psi-func4}), and (\ref{tangent-to-Shock-3}), we find
$$
\hat\tau((0, 0),0,0)=(0,\,1),\quad D_{p}\Psi((0,0), 0, 0,0, \mxy)=
\big(\rho_2'(c_2^2-{\hat\mxx}^2),\;
\big(\frac{\rho_2-\rho_1}{u_1}-\rho_2'\hat\mxx \big)\mxy\big).
$$
Thus, by (\ref{Diff-RH-psi-func2}), we obtain that, on $\bR^{2\times
2}\times{\mathcal A}$,
\begin{eqnarray}\label{Diff-RH-psi-func-at-zero}
&&\qquad \tilde\Phi(R, p, z, u_2, v_2,  \mxy)\\
&&\qquad =\rho_2'(c_2^2-{\hat\mxx}^2) R_{21} +
\big(\frac{\rho_2-\rho_1}{u_1}-\rho_2'\hat\mxx \big)\mxy R_{22}
+\sum_{i.j=1}^2 \hat E_{ij}(p, z, u_2, v_2,  \mxy) R_{ij},\nonumber
\end{eqnarray}
where $\hat E_{ij}\in C^\infty(\overline{\mathcal A})$ and
\begin{eqnarray*}
|\hat E_{ij}(p,z, u_2, v_2, \mxy)| \le C(|p|+|z|+|u_2|+|v_2|)
\qquad\mbox{for any }(p, z, u_2, v_2, \mxy)\in {\mathcal A},
\end{eqnarray*}
with $C$ depending only on $\|D^2\Psi\|_{C^0(\overline{\mathcal
A})}$.

{}From now on, we fix $(u_2, v_2)$ to be equal to the velocity of
state (2) obtained in \S \ref{section:3.3} and write $E_{ij}(p, z,
\mxy)$ for $\hat E_{ij}(p,z,u_2, v_2, \mxy)$. Then, from
(\ref{diff-RH-2}) and (\ref{Diff-RH-psi-func-at-zero}), we conclude
that $\psi$ satisfies
\begin{equation}
\label{differentiateRH-3} \rho_2'(c_2^2-{\hat\mxx}^2) \psi_{\mxx\mxy}+
\big(\frac{\rho_2-\rho_1}{u_1}-\rho_2'\hat\mxx \big)\mxy\psi_{\mxy\mxy}
+\sum_{i,j=1}^2E_{ij}(D\psi, \psi, \mxy)D_{ij}\psi=0 \qquad\mbox{on
}\;\shock,
\end{equation}
and $E_{ij}=E_{ij}(p, z, \mxy), i,j=1,2,$ are smooth on
$$
{\mathcal
B}\defd\overline{B_{\delta^*}(0)\times(-\delta^*, \delta^*) \times
(-6\bar c_2/5, 6\bar c_2/5)}
$$
and satisfy
(\ref{RH-psi-2-error-term1}) with $C$ depending only on the data.
Note that
$$
(\grad\psi(\mxx, \mxy), \,\psi(\mxx, \mxy),\,\mxy)\in {\mathcal B}
\qquad\mbox{on } \shock,
$$
since $\psi\in\setK$ and
(\ref{condConst-00}) holds with sufficiently large $\hat C$.
Expressing
 $\psi_{\mxx\mxx}$ from (\ref{express2ndDeriv}) and using (\ref{A11-bdd-away-fromZero}),
we can rewrite  (\ref{differentiateRH-3}) in the form
$$
\big(\rho_2'(c_2^2-{\hat\mxx}^2)+ E_1(D\psi, \psi, \mxy)\big)
\psi_{\mxx\mxy}+ \big(\big(\frac{\rho_2-\rho_1}{u_1}-\rho_2'\hat\mxx
\big)\mxy+ E_2(D\psi, \psi, \mxy)\big)\psi_{\mxy\mxy}=0
$$
on $\shock$, where the functions $E_i=E_i(p, z, \mxy), i=1,2,$ are
smooth on ${\mathcal B}$ and satisfy (\ref{RH-psi-2-error-term1}).
Thus, $w$ satisfies
\begin{equation}\label{negativeDerivPsi-CondRH-ForW}
\big(\rho_2'(c_2^2-{\hat\mxx}^2)+ E_1(D\psi, \psi, \mxy)\big)
w_{\mxx}+ \big((\frac{\rho_2-\rho_1}{u_1}-\rho_2'\hat\mxx)\mxy +
E_2(D\psi, \psi, \mxy)\big)w_{\mxy}=0
\end{equation}
on $\shock$. Condition (\ref{negativeDerivPsi-CondRH-ForW}) is
oblique if  $\hat C$ is sufficiently large in (\ref{condConst-00}).
Indeed, we have $c_2\ge {9\over 10}\bar{c}_2$, which implies
$c_2^2-|\hat\mxx|^2\ge \bar{c}_2\frac{\bar{c}_2-|\bar{\xi}|}{4}>0$
by using (\ref{inSonicRegion-in-shifed}). Now, combining
(\ref{norm-to-Shock}) and (\ref{RH-psi-2-error-term1}) with
$\psi\in\setK$ and (\ref{u2-v2-bound}), we find that, on $\shock$,
\begin{eqnarray*}
&&(\rho_2'(c_2^2-{\hat\mxx}^2)+ E_1(D\psi, \psi, \mxy),
\big(\frac{\rho_2-\rho_1}{u_1}-\rho_2'\hat\mxx \big)\mxy
+ E_2(D\psi, \psi, \mxy) )\cdot\nu_s \\
&&\quad \ge\rho_2'\bar{c}_2\frac{\bar{c}_2-|\bar{\xi}|}{4}
-C(M_1\varepsilon+M_2\epsP) \ge
\rho_2'\bar{c}_2\frac{\bar{c}_2-|\bar{\xi}|}{8}>0.
\end{eqnarray*}
Also, the coefficients of  (\ref{negativeDerivPsi-CondRH-ForW})
are continuous with respect to $(\mxx, \mxy)\in\shock$.

\medskip
{\em Step 6.} Both the regularity of $w$ in
(\ref{w-regular-negativeDeriv}) and the fact that $w$ satisfies
equation (\ref{negativeDerivPsi-eqnForW}) that is  uniformly
elliptic on compact subsets of $\overline{\Omega^+}\setminus
\overline{\sonic}$ imply that the maximum of $w$ cannot be achieved
in the interior of $\Omega^+$, unless $w$ is constant on $\Omega^+$,
by the Strong Maximum Principle. Since $w$ satisfies the oblique
derivative conditions (\ref{negativeDerivPsi-CondWedge-ForW}) and
(\ref{negativeDerivPsi-CondRH-ForW}) on the straight segment
$\wedgeB$ and on the curve $\shock$ that is $C^{2,\alpha}$ in its
relative interior, and since equation
(\ref{negativeDerivPsi-eqnForW})
 is  uniformly elliptic in a neighborhood of any point from the relative interiors
 of $\wedgeB$ and $\shock$, it follows from Hopf's Lemma that the maximum of
 $w$ cannot be achieved in the relative interiors of $\wedgeB$ and
 $\shock$,
 unless $w$ is constant on $\Omega^+$. Now conditions
 (\ref{negativeDerivPsi-CondSymmetr-ForW})--(\ref{negativeDerivPsi-CondSonic-ForW})
 imply that $w\le 0$ on $\Omega^+$.
This completes the proof.
\Endproof
%*************** END PROOF *****LEMMA**********

Using Lemma \ref{negativeDerivPsiLemma} and working in the
$(x,y)$--coordinates, we have
%*****************BEGIN PROPOSITION*************
\begin{proposition}
\label{boundPsiXfromBelow-Prop} If $\hat C$ in {\rm
(\ref{condConst-00})} is sufficiently large, depending only on the
data, then
\begin{equation}\label{boundPsiXfromBelow}
\psi_x\ge -{4\over 3(\gamma+1)}x \quad\mbox{ in }\;
\Omega^+\cap\{x\le 4\varepsilon\}.
\end{equation}
\end{proposition}
%*****************END PROPOSITION*************
%*************** BEGIN PROOF *****PROPOSITION***********
\Proof
{}By definition of the $(x,y)$--coordinates in
(\ref{coordNearSonic}), we have
\begin{equation}\label{psi-eta}
\psi_\eta=-\sin\theta\, \psi_x+{\cos\theta\over r}\psi_y,
\end{equation}
where $(r,\theta)$ are the polar coordinates in the $(\mxx,
\mxy)$--plane.

{}From (\ref{domainExtension}), it follows that, for sufficiently
small $\epsP$ and $\varepsilon$, depending only on the data,
$$
\mxy\ge\mxy^*\qquad \mbox{for
all } \,\, (\mxx, \mxy)\in\Dom\cap\{c_2-r<4\varepsilon\},
$$
where $(l(\mxy^*), \mxy^*)$ is the unique intersection point of the
segment $\{(l(\mxy),\mxy)\ :\, \mxy\in(0, \mxy_1]\}$ with the circle
$\partial B_{c_2-4\varepsilon}(0)$. Let $\bar\mxy^*$ be the
corresponding point for the case of normal reflection, i.e.,
$\bar\mxy^*=\sqrt{(\bar c_2-4\varepsilon)^2-\bar\mxx^2}$. By
(\ref{sonic-intersect-shock-normal}), $\bar\mxy^*\ge \sqrt{\bar
c_2^2-\bar\mxx^2}/2>0$ if $\varepsilon$ is sufficiently small. Also,
from (\ref{reflected-shock-s2})--(\ref{x1-in-shifed}) and
(\ref{u2-v2-bound}), and using the convergence $(\theta_s, c_2,
\tilde\xi)\to(\pi/2, \bar c_2,\bar\xi)$ as $\theta_w\to \pi/2$, we
obtain $\mxy^*\ge \bar\mxy^*/2$ and $c_2\le 2\bar c_2$ if $\epsP$ and
$\varepsilon$ are sufficiently small. Thus, we conclude that, if
$\hat C$ in (\ref{condConst-00}) is sufficiently large depending
only on the data, then, for every $(\mxx,
\mxy)\in\Dom\cap\{c_2-r<4\varepsilon\}$, the polar angle $\theta$
satisfies
\begin{equation}\label{psi-eta-polar-angle-ineq}
\sin\theta=\frac{\mxy}{\sqrt{\mxx^2+\mxy^2}}>0, \qquad |\cot\theta|=
\left|\frac{\mxx}{\mxy}\right|\le {8\bar c_2\over \sqrt{\bar
c_2^2-\bar\mxx^2}}\le C.
\end{equation}
{}From (\ref{psi-eta})--(\ref{psi-eta-polar-angle-ineq}) and Lemma
\ref{negativeDerivPsiLemma}, we find that, on
$\Omega^+\cap\{c_2-r<4\varepsilon\}$,
\begin{equation}\label{estFromBelowDer}
\psi_x=-{1\over \sin\theta}\psi_\eta+{\cot\theta\over r}\psi_y
\ge{\cot\theta\over r}\psi_y\ge -C|\psi_y|.
\end{equation}
Note that $\psi\in\setK$ implies $|\psi_y(x,y)|\le M_1x^{3/2}$ for
all $(x,y)\in \Omega^+\cap\{c_2-r<2\varepsilon\}$. Then, using
(\ref{estFromBelowDer}) and (\ref{condConst-00}) and choosing large
$\hat C$, we have
$$
\psi_x\ge -{4\over 3(\gamma+1)}x \qquad\mbox{ in }\;
\Omega^+\cap\{x\le 2\varepsilon\}.
$$
Also, $\psi\in\setK$ implies
$$
|\psi_x|\le M_2\epsP\le
{4\over3(\gamma+1)}(2\varepsilon) \qquad\mbox{on }\;
\Omega^+\cap\{2\varepsilon\le x\le 4\varepsilon\},
$$
where the second inequality holds by (\ref{condConst-00}) if $\hat
C$ is sufficiently large depending only on the data. Thus,
(\ref{boundPsiXfromBelow}) holds on $\Omega^+_{4\varepsilon}$.
\Endproof
%*************** END PROOF ******PROPOSITION**********

\section{Proof of Main Theorem}
\label{proofSection} Let $\hat C$ be sufficiently large to satisfy
the conditions in Propositions \ref{existenceFixedPt} and
\ref{boundPsiXfromAbove-Prop}--\ref{boundPsiXfromBelow-Prop}. Then,
by Proposition \ref{existenceFixedPt}, there exist $\epsP_0,
\varepsilon>0$ and $M_1,M_2\ge 1$ such that, for any $\epsP\in (0,
\epsP_0]$,  there exists a solution $\psi\in\setK(\epsP,
\varepsilon, M_1, M_2)$  of problem {\rm
(\ref{iterationEquation})}--{\rm (\ref{iterationCondOnSymmtryLine})}
with $\Kphi=\psi$. Fix $\epsP\in (0, \epsP_0]$ and the corresponding
``fixed point" solution $\psi$, which, by Propositions
\ref{boundPsiXfromAbove-Prop}--\ref{boundPsiXfromBelow-Prop},
satisfies
$$
|\psi_x|\le {4\over 3(\gamma+1)}x \qquad\mbox{ in }\;
\Omega^+\cap\{x\le 4\varepsilon\}.
$$
Then, by Lemma \ref{cutOffEqIsOriginalEq}, $\psi$ satisfies equation
(\ref{potent-flow-nondiv-psi-1}) in $\Omega^+(\psi)$. Moreover,
$\psi$ satisfies properties
(\ref{iterProcItem-first})--(\ref{iterProcItem-last}) in Step 10 of
\S\ref{overViewProcedureSubsection} by following the argument in
Step 10 of \S\ref{overViewProcedureSubsection}. Then, extending the
function $\varphi=\psi+\varphi_2$ from $\Omega:=\Omega^+(\psi)$ to
the whole domain $\Lambda$ by using (\ref{phi-states-0-1-2}) to
define $\varphi$ in $\Lambda\setminus\Omega$, we obtain
$$
\varphi\in W^{1,\infty}_{loc}(\Lambda)\cap \left(\cup_{i=0}^2
C^1(\Lambda_i\cup S)\cap C^{1,1}(\Lambda_i)\right),
$$
where the domains $\Lambda_i$, $i=0,1,2$, are defined in Step 10 of
\S\ref{overViewProcedureSubsection}. {}From the argument in Step 10
of \S\ref{overViewProcedureSubsection}, it follows that $\varphi$ is
a weak solution of Problem 2, provided that the reflected shock
$S_1=P_0\PtUpL\PtLwL\cap \Lambda$ is a $C^2$-curve.

Thus, it remains to show that $S_1=P_0\PtUpL\PtLwL\cap \Lambda$ is a
$C^2$-curve. By definition of $\varphi$ and since
$\psi\in\setK(\epsP, \varepsilon, M_1, M_2)$, the reflected shock
$S_1=P_0\PtUpL\PtLwL\cap\Lambda$ is given by $S_1=\{\mxx=f_{S_1}(\mxy)
\;:\; \mxy_{\PtLwL}<\mxy<\mxy_{P_0}\}$, where $\mxy_{\PtLwL}=  -v_2$,
$\mxy_{P_0}=|\hat\mxx|{\sin\theta_s\sin\theta_w\over
\sin(\theta_w-\theta_s)}>0$, and
\begin{equation}\label{defFulReflShock}
f_{S_1}(\mxy)= \left\{
\begin{array}{ll}
f_\psi(\mxy)\quad & \mbox{if }\mxy\in (\mxy_{\PtLwL}, \mxy_{\PtUpL}),\\
l(\mxy)& \mbox{if }\mxy\in (\mxy_{\PtUpL}, \mxy_{P_0}),
\end{array}
\right.
\end{equation}
where $l(\mxy)$ is defined by (\ref{reflected-shock-s2}),
$\mxy_{\PtUpL}=\mxy_1>0$ is defined by (\ref{coord-P4}), and
$\mxy_{P_0}>\mxy_{\PtUpL}$ if $\epsP$ is sufficiently small, which
follows from the explicit expression of $\mxy_{P_0}$ given above and
the fact that $(\theta_s, c_2, \hat\mxx)\to(\pi/2,\bar{c}_2,
\bar\mxx)$ as $\theta_w\to\pi/2$. The function $f_\psi$ is defined by
(\ref{shockPL}) for $\Kphi=\psi$.

Thus we need to show that $f_{S_1}\in C^2([\mxy_{\PtLwL},
\mxy_{P_0}])$. By (\ref{reflected-shock-s2}) and
(\ref{OmegaPL-f-higher}), it suffices to show that $f_{S_1}$ is
twice differentiable at the points $\mxy_{\PtUpL}$ and
$\mxy_{\PtLwL}$.

First, we consider $f_{S_1}$ near $\mxy_{\PtUpL}$. We change the
coordinates to the $(x,y)$--coordinates in (\ref{coordNearSonic}).
Then, for sufficiently small $\varepsilon_1>0$, the curve
$\{\mxx=f_{S_1}(\mxy)\}\cap \{c_2-\varepsilon_1<r<c_2+\varepsilon_1\}$
has the form $\{y=\hat f_{S_1}(x)\; :
\;-\varepsilon_1<x<\varepsilon_1\}$, where
\begin{equation}\label{defFulReflShock-xy}
\hat f_{S_1}(x)= \left\{
\begin{array}{ll}
\hat f_\psi(x)\quad & \mbox{if }x\in(0, \varepsilon_1),\\
\hat f_0(x)& \mbox{if }x\in (-\varepsilon_1,0),
\end{array}
\right.
\end{equation}
with $\hat f_0$ and $\hat f_\psi$ defined by
(\ref{referenceFB-polar}) and (\ref{domain-in-rescaled-lemma}) for
$\Kphi=\psi$. In order to show that
 $f_{S_1}$ is twice differentiable at $\mxy_{\PtUpL}$,  it suffices
 to show that $\hat f_{S_1}$ is twice differentiable at $x=0$.

{}From (\ref{holder-hat-f})--(\ref{holder-hat-f-S}) and
(\ref{referenceFB-polar}), it follows that $\hat f_{S_1}\in
C^1((-\varepsilon_1, \varepsilon_1))$. Moreover, from
(\ref{nondegenPolar-1}), (\ref{domain-in-xy-0}), (\ref{OmegaPL}),
and (\ref{holder-hat-f-S}), we write $\varphi_1, \varphi_2,$ and
$\psi$ in the $(x,y)$--coordinates to obtain that
\begin{equation}\label{defFulReflShock-deriv-xy}
\hat f_{S_1}'(x)= \left\{
\begin{array}{ll}
\displaystyle -{\partial_x(\varphi_1-\varphi_2-\psi)\over
\partial_y(\varphi_1-\varphi_2-\psi)} (x, \hat f_{S_1}(x))
 \quad & \mbox{if }x\in(0, \varepsilon_1),\\
\displaystyle -{\partial_x(\varphi_1-\varphi_2)\over
\partial_y(\varphi_1-\varphi_2)} (x, \hat f_{S_1}(x)) & \mbox{if
}x\in (-\varepsilon_1,0],
\end{array}
\right.
\end{equation}
and that $\hat f_0'(x)$ is given for $x\in(-\varepsilon_1,
\varepsilon_1)$ by the second line of the right-hand side of
(\ref{defFulReflShock-deriv-xy}). Using (\ref{nondegenPolar-1}) and
$\psi\in\setK$ with (\ref{condConst-00}) for sufficiently large
$\hat C$, we have
\begin{equation}\label{defFulReflShock-deriv-xy-est}
|\hat f_{S_1}'(x)-\hat f_0'(x)|\le C|D_{(x,y)}\psi(x,\hat
f_{\psi}(x))| \qquad\mbox{for all }x\in(0, \varepsilon_1).
\end{equation}
Since $\psi$ satisfies  (\ref{iterationRH}) with $\Kphi=\psi$, it
follows that, in the $(x,y)$--coordinates, $\psi$ satisfies
(\ref{iterationRH-lf-flattened}) on $\{y=\hat f_{\psi}(x)\; :
\;x\in(0,\varepsilon_1)\}$, and (\ref{estCoefsIterRH-flattened})
holds. Then it follows that
$$
|\psi_x(x,\hat f_{\psi}(x))|\le C(|\psi_y(x,\hat f_{\psi}(x))|
+|\psi(x,\hat f_{\psi}(x))|)\le Cx^{3/2},
$$
where the last inequality follows from $\psi\in\setK$. Combining
this with (\ref{defFulReflShock-xy}),
(\ref{defFulReflShock-deriv-xy-est}), and $\hat f_{S_1}, \hat f_0\in
C^1((-\varepsilon_1, \varepsilon_1))$ yields
$$
|\hat f_{S_1}'(x)-\hat f_0'(x)|\le Cx^{3/2}\qquad\mbox{for all }
x\in(-\varepsilon_1, \varepsilon_1).
$$
Then it follows that $\hat f_{S_1}'(x)-\hat f_0'(x)$ is
differentiable at $x=0$. Since $\hat f_0\in
C^\infty((-\varepsilon_1, \varepsilon_1))$, we conclude that $\hat
f_{S_1}$ is twice differentiable at $x=0$. Thus, $f_{S_1}$ is twice
differentiable at $\mxy_{\PtUpL}$.

In order to prove the $C^2$--smoothness of $f_{S_1}$ up to
$\mxy_{\PtLwL}=-v_2$, we extend the solution $\Kphi$ and the free
boundary function $f_{S_1}$ into $\{\mxy<-v_2\}$ by the even
reflection about the line $\Sigma_0\subset\{\mxy=-v_2\}$ so that
$\PtLwL$ becomes an interior point of the shock curve. Note that we
continue to work in the shifted coordinates defined in
\S\ref{shiftCoordSection}, that is, for $(\mxx, \mxy)$ such that
$\mxy<-v_2$ and $(\mxx, -2v_2-\mxy)\in\overline{\Omega^+(\psi)}$, we
define $(\varphi, \varphi_1)(\mxx, \mxy)=(\varphi, \varphi_1)(\mxx,
-2v_2-\mxy)$ and $f_{S_1}(\mxy)=-2v_2-\mxy$ for $\varphi_1$ given by
(\ref{phi-1-shifted}). Denote $\Omega^+_{\varepsilon_1}(\PtLwL)\defd
B_{\varepsilon_1}(\PtLwL) \cap \{\mxx>f_{S_1}(\mxy)\}$ for
sufficiently small $\varepsilon_1 >0$. {}From $\varphi\in
C^{1,\alpha}(\overline{\Omega^+(\psi)})\cap
C^{2,\alpha}(\Omega^+(\psi))$ and (\ref{condOnSymmtryLinePhi}), we
have
$$
\varphi\in C^{1,\alpha}(\overline{\Omega^+_{\varepsilon_1}(\PtLwL)})
\cap C^{2,\alpha}(\Omega^+_{\varepsilon_1}(\PtLwL)).
$$
Also, the extended function $\varphi_1$ is in fact given by
(\ref{phi-1-shifted}). Furthermore, from (\ref{nondegeneracy}) and
(\ref{OmegaPL}), we can see that the same is true for the extended
functions and hence
$$
\{\mxx>f_{S_1}(\mxy)\}\cap B_{\varepsilon_1}(\PtLwL)
=\{\varphi<\varphi_1\}\cap B_{\varepsilon_1}(\PtLwL), \quad
f_{S_1}\in C^{1,\alpha}((-v_2-{\varepsilon_1\over 2},
-v_2+{\varepsilon_1\over 2})).
$$
Furthermore, from (\ref{1.1.5})--(\ref{1.1.6}) and
(\ref{condOnSymmtryLinePhi}), it follows that the extended $\varphi$
satisfies  equation (\ref{1.1.5}) with  (\ref{1.1.6}) in
$\Omega^+_{\varepsilon_1}(\PtLwL)$, where we have used the form of
equation, i.e., the fact that there is no explicit dependence on
$(\mxx, \mxy)$ in the coefficients and that the dependence of
$D\varphi$ is only through $|D\varphi|$. Finally, the boundary
conditions (\ref{cont-accross-shock-mod-phi}) and (\ref{RH-mod-phi})
are satisfied on $\Gamma_{\varepsilon_1}(\PtLwL)
\defd\{\mxx=f_{S_1}(\mxy)\}\cap B_{\varepsilon_1}(\PtLwL)$. Equation (\ref{1.1.5})
is uniformly elliptic in $\Omega^+_{\varepsilon_1}(\PtLwL)$ for
$\varphi$, which follows from $\varphi=\varphi_2+\psi$ and Lemmas
\ref{propertiesNonlinCoeffs} and \ref{cutOffEqIsOriginalEq}.
Condition (\ref{RH-mod-phi}) is uniformly oblique on
$\Gamma_{\varepsilon_1}(\PtLwL)$ for $\varphi$, which follows from
\S\ref{equationForPsiSection}.

Next, we rewrite equation (\ref{1.1.5}) in
$\Omega^+_{\varepsilon_1}(\PtLwL)$
 and the boundary conditions
(\ref{cont-accross-shock-mod-phi})--(\ref{RH-mod-phi}) on
$\Gamma_{\varepsilon_1}(\PtLwL)$ in terms of $u\defd
\varphi_1-\varphi$. Substituting $u+\varphi_1$ for $\varphi$ into
(\ref{1.1.5}) and (\ref{RH-mod-phi}),
we obtain that $u$ satisfies
$$
F(D^2u, Du, u, \mxx, \mxy)=0\quad\mbox{in
}\Omega^+_{\varepsilon_1}(\PtLwL),\qquad\; u=G(Du, u, \mxx, \mxy)=0
\quad\mbox{on } \Gamma_{\varepsilon_1}(\PtLwL),
$$
where the equation is quasilinear and uniformly elliptic, the second
boundary condition is oblique, and the functions $F$ and $G$ are
smooth. Also, from  (\ref{nondegeneracy}) which holds for the even
extensions as well, we find that $\partial_\mxx u>0$ on
$\Gamma_{\varepsilon_1}(\PtLwL)$. Then, applying the hodograph
transform of \cite[\S3]{KinderlehrerNirenberg}, i.e., changing
$(\mxx, \mxy)\to  (X,Y)=(u(\mxx, \mxy), \mxy)$, and denoting the
inverse transform by $(X,Y)\to (\mxx, \mxy)=(v(X,Y), Y)$, we obtain
$$
v\in C^{1,\alpha}(\overline{B^+_\delta((0, -v_2))}) \cap
C^{2,\alpha}(B^+_\delta((0, -v_2))),
$$
where $B_\delta^+((0, -v_2))\defd B_\delta((0, -v_2))\cap\{X>0\}$
for small $\delta>0$,
 $v(X,Y)$ satisfies a uniformly elliptic quasilinear equation
$$
\tilde F(D^2v, Dv, v, X, Y)=0 \qquad\mbox{in } B_\delta^+((0,
-v_2))
$$
and the oblique derivative condition
$$\tilde G(Dv,v,Y)=0 \qquad\mbox{on } \partial B_\delta^+((0, -v_2))\cap \{X=0\},
$$
and the functions
$\tilde F$ and $\tilde G$ are smooth. Then, from the local estimates
near the boundary in the proof of \cite[Theorem 2]{Lieberman86},
$v\in C^{2,\alpha}(\overline{B^+_{\delta/2}((0, -v_2))})$. Since
$f_{S_1}(\eta)=v(0, \eta)$, it follows that $f_{S_1}$ is
$C^{2,\alpha}$ near $\mxy_{\PtLwL}=-v_2$.

It remains to prove the convergence of the solutions to the normal
reflection solution as $\theta_w\to \pi/2$. Let $\theta_w^i\to
\pi/2$ as $i\to\infty$. Denote by $\varphi^i$ and $f^i$ the
corresponding solution and the free-boundary function respectively,
i.e., $P_0\PtUpL\PtLwL\cap \Lambda$ for each $i$ is given by
$\{\mxx=f^i(\mxy) \;:\;\mxy\in(\mxy_{\PtLwL}, \mxy_{P_0})\}$. Denote
by $\varphi^\infty$ and $f^\infty(\eta)=\bar\mxx$ the solution and
the reflected shock for the normal reflection respectively. For each
$i$, we find that $\varphi^i-\varphi_2^i=\psi^i$ in the subsonic
domain $\Omega^+_i$, where $\psi^i$ is the corresponding ``fixed
point solution'' from Proposition \ref{existenceFixedPt} and
$\psi^i\in \setK(\pi/2-\theta_w^i, \varepsilon^i, M_1^i, M_2^i)$
with (\ref{condConst-00}). Moreover,  $f^i$ satisfies
(\ref{OmegaPL-f-higher}).
 We also use the convergence of state (2)
to the corresponding state of the normal reflection obtained in
\S\ref{section:3.3}. Then we conclude that, for a subsequence,
$f^i\to f^\infty$ in $C^1_{loc}$ and $\varphi^i\to\varphi^\infty$ in
$C^1$ on compact subsets of $\{\mxx>\bar\mxx\}$ and
$\{\mxx<\bar\mxx\}$. Also, we obtain $\|(D\varphi^i,
\varphi^i)\|_{L^\infty(K)}\le C(K)$ for every compact set K. Then
$\varphi^i\to \varphi_\infty$ in $W^{1,1}_{loc}(\overline\Lambda)$
by the Dominated Convergence Theorem. Since such a converging
subsequence can be extracted from every sequence $\theta_w^i\to
\pi/2$, it follows that $\varphi_{\theta_w}\to\varphi_\infty$ as
$\theta_w\to \pi/2$.

\appendix
\section{Appendix: Estimates of Solutions to Elliptic Equations}
\label{append-1-section}

\newcommand{\FltBdry}{ \Sigma}

In this appendix, we make some careful estimates of solutions of
boundary value problems for elliptic equations in $\bR^2$, which are
applied in \S\ref{unifElliptApproxSection}--\S\ref{fixedPtSection}.
Throughout the appendix, we denote by $(x,y)$ or $(X,Y)$ the
coordinates in $\bR^2$, by $\bR^2_+:=\{y>0\}$, and, for $z=(x,0)$
and $r>0$, denote $B_r^+(z):=B_r(z)\cap\bR^2_+$ and
$\FltBdry_r(z):=B_r(z)\cap\{y=0\}$. We also denote $B_r:=B_r(0)$,
$B_r^+:=B_r^+(0)$, and $\FltBdry_r:=\FltBdry_r(0)$.

We consider an elliptic equation of the form
\begin{equation}\label{nonlinEq-appdx}
A_{11}u_{xx}+2A_{12}u_{xy}+A_{22}u_{yy}+ A_{1}u_{x}+A_{2}u_{y}=f,
\end{equation}
where $A_{ij}=A_{ij}(Du, x,y)$, $A_{i}=A_{i}(Du, x, y)$, and
$f=f(x,y)$. We study the following three types of boundary
conditions: (i) the Dirichlet condition, (ii) the oblique derivative
condition, (iii) the ``almost tangential derivative" condition.

One of the new ingredients in our estimates below is that we do not
assume that the equation satisfies the ``natural structure
conditions", which are used in the earlier related results; see,
e.g., \cite[Chapter 15]{GilbargTrudinger} for the interior estimates
for the Dirichlet problem and \cite{LiebermanTrudinger} for the
oblique derivative problem. For equation (\ref{nonlinEq-appdx}), the
natural structure conditions include the requirement that
$|p||D_pA_{ij}|\le C$ for all $p\in\bR^2$. Note that equations
(\ref{iteration-equation-sonicStruct}) and
(\ref{nonlinIterEq-xy-lg}) do not satisfy this condition because of
the term $x\zeta_1(\frac{\psi_x}{x})$ in the coefficient of
$\psi_{xx}$. Thus we have to derive the estimates for the equations
without the ``natural structure conditions". We consider only the
two-dimensional case here.

The main point at which the ``natural structure conditions" are
needed is the gradient estimates. The interior gradient estimates
and global
 gradient
estimates for the Dirichlet problem, without requiring the natural
structure conditions, were obtained in the earlier results in the
two-dimensional case; see Trudinger \cite{Trudinger85} and
references therein. However, it is not clear how this approach can
be extended to the oblique and ``almost tangential" derivative
problems. We also note a related result by Lieberman
\cite{Lieberman87-1} for fully nonlinear equations and the boundary
conditions without obliqueness assumption in the two-dimensional
case, in which the H\"{o}lder estimates for the gradient of a
solution depend on both the bounds of the solution and its gradient.

In this appendix, we present the $C^{2,\alpha}$--estimates of the
solution only in terms of its $C$--norm. For simplicity,  we
restrict to the case of quasilinear equation (\ref{nonlinEq-appdx})
and linear boundary conditions, which is the case for the
applications in this paper. Below, we first present the interior
estimate in the form that is used in the other parts of this paper.
Then we give a proof of the $C^{2,\alpha}$--estimates for the
``almost tangential" derivative problem. Since the proofs for the
Dirichlet and oblique derivative problems are similar to that for
the ``almost tangential" derivative problem, we just sketch these
proofs.

%*****************BEGIN THEOREM*************
\begin{theorem}\label{locEstElliptEq}
Let $u\in C^2(B_2)$ be a solution of equation {\rm
(\ref{nonlinEq-appdx})} in $B_2$. Let $A_{ij}(p,x,y)$,
$A_{i}(p,x,y)$, and $f(x,y)$ satisfy that there exist constants
$\lambda>0$
 and $\alpha\in(0,1)$ such that
\begin{eqnarray}
\label{locEstElliptEq-i1-0} &&\qquad\lambda|\mu|^2 \le
\sum_{i,j=1}^nA_{ij}\mu_i\mu_j\le\lambda^{-1}|\mu|^2 \qquad\mbox{for
all } (x,y)\in B_2,\,\, p, \mu\in\bR^2,
\\
\label{locEstElliptEq-i2-0} &&\qquad \|
(A_{ij},\,A_{i})\|_{C^\alpha(\bR^2\times\overline{B_2})}+\| D_p
(A_{ij},\,
A_{i})\|_{C(\bR^2\times\overline{B_2})}+\|f\|_{C^\alpha(\overline{B_2})}
\le \lambda^{-1}.
\end{eqnarray}
Assume that $\|u\|_{C(\overline{B_2})}\le M$. Then there exists
$C>0$ depending only on $(\lambda, M)$ such that
\begin{equation}\label{localEst-appdx-2}
\|u\|_{C^{2,\alpha}(\overline{B_1})}\le C(\|u\|_{C(\overline{B_2})}
+\|f\|_{C^{\alpha}(\overline{B_2})}).
\end{equation}
\end{theorem}
%*****************END THEOREM*************
\Proof
%*****************BEGIN PROOF-THEOREM*************
We use the standard interior H\"{o}lder seminorms and norms as
defined in \cite[Eqs. (4.17), (6.10)]{GilbargTrudinger}. By
\cite[Theorem 12.4]{GilbargTrudinger}, there exists $\beta\in(0,1)$
depending only on $\lambda$ such that
\begin{eqnarray*}
[u]^*_{1,\beta, B_2}&\le& C(\lambda)(\|u\|_{0, B_2}
+\|f-A_1D_1u-A_2D_2u\|^{(2)}_{0, B_2})\\
&\le& C(\lambda, M)(1+\|f\|^{(2)}_{0, B_2}+\|Du\|^{(2)}_{0, B_2}).
\end{eqnarray*}
Then, applying the interpolation inequality
\cite[(6.82)]{GilbargTrudinger} with the argument similar to that
for the proof of \cite[Theorem 12.4]{GilbargTrudinger}, we obtain
$$
\|u\|^*_{1,\beta, B_2}\le C(\lambda, M)(1+\|f\|^{(2)}_{0, B_2}).
$$
Now we consider (\ref{nonlinEq-appdx}) as a linear elliptic equation
$$
\sum_{i,j=1}^na_{ij}(x)u_{x_ix_j}+ \sum_{i=1}^n a_{i}(x)u_{x_i}=f(x)
\qquad\mbox{in }\,\, B_{3/2}
$$
with coefficients $a_{ij}(x)=A_{ij}(Du(x),x)$ and
$a_i=A_{i}(Du(x),x)$ in $C^\beta(\overline{B_{3/2}})$ satisfying
$$
\|(a_{ij}, a_i)\|_{C^\beta(\overline{B_{3/2}})}\le C(\lambda, M).
$$
We can assume $\beta\le\alpha$. Then the local estimates for linear
elliptic equations yield
$$
\|u\|_{C^{2,\beta}(\overline{B_{5/4}})}\le C(\lambda, M)
(\|u\|_{C(\overline{B_{3/2}})}
+\|f\|_{C^{\beta}(\overline{B_{3/2}})}).
$$
With this estimate, we have $ \|(a_{ij},
a_i)\|_{C^\alpha(\overline{B_{5/4}})}\le C(\lambda, M). $ Then the
local estimates for linear elliptic equations in $B_{5/4}$ yield
(\ref{localEst-appdx-2}).
\Endproof
%*****************END PROOF-THEOREM*************

Now we make the estimates for the ``almost tangential derivative"
problem.
%*****************BEGIN THEOREM*************
\begin{theorem}\label{locEstElliptEq-non-oblique}
Let $\lambda>0$, $\alpha\in (0,1)$, and $\varepsilon\ge 0$. Let
$\Phi\in C^{2,\alpha}(\bR)$ satisfy
\begin{equation}\label{bdryNormCond}
\|\Phi\|_{C^{2,\alpha}(\bR)}\le \lambda^{-1},
\end{equation}
and denote $\Omega_R^+:=B_R\cap\{y>\varepsilon\Phi(x)\}$ for $R>0$.
Let $u\in C^2(B_2^+)\cap C^1(\overline{B_2^+})$ satisfy {\rm
(\ref{nonlinEq-appdx})} in $\Omega_2^+$ and
\begin{eqnarray}
&&u_{x}=\varepsilon b(x,y)u_{y}+c(x,y) u\qquad \mbox{ on }\;
\Gamma_\Phi\defd B_2\cap\{y=\varepsilon \Phi(x)\}.
\label{nonlinEq-appdx-nonObl-bc}
\end{eqnarray}
Let $A_{ij}(p,x,y)$,  $A_{i}(p,x,y)$, $a(x,y)$, $b(x,y)$, and
$f(x,y)$ satisfy that there exist constants $\lambda>0$ and
$\alpha\in(0,1)$ such that
\begin{eqnarray}
\label{locEstElliptEq-i1} &&\qquad
\lambda|\mu|^2\le\sum_{i,j=1}^nA_{ij}\mu_i\mu_j\le\lambda^{-1}|\mu|^2
\qquad\mbox{for } (x,y)\in \Omega_2^+,\,\, p, \mu\in\bR^2,
\\
\label{locEstElliptEq-i2}
&&\qquad\|(A_{ij},\,A_{i})\|_{C^\alpha(\overline{\Omega_2^+}\times\bR^2)}+\|D_p
(A_{ij},\, A_{i})\|_{C(\overline{\Omega_2^+}\times\bR^2)}
+\|f\|_{C^{\alpha}(\overline{\Omega_2^+})} \le \lambda^{-1},
\\
\label{bdryCoefBds} &&\qquad
\|(b,c)\|_{C^{1,\alpha}(\overline{\Omega_2^+})}\le \lambda^{-1}.
\end{eqnarray}
Assume that $\|u\|_{C(\overline{\Omega_2^+})}\le M$. Then there
exist $\varepsilon_0(\lambda, M, \alpha)>0$ and $C(\lambda, M,
\alpha)>0$ such that, if $\varepsilon\in (0,\varepsilon_0)$, we have
\begin{equation}\label{localEst-appdx-nonOblique-est}
\|u\|_{C^{2,\alpha}(\overline{\Omega_1^+})}\le
C(\|u\|_{C(\overline{\Omega_2^+})}+\|f\|_{C^{\alpha}(\overline{\Omega_2^+})}).
\end{equation}
\end{theorem}
%*****************END THEOREM*************

%*****************BEGIN PROOF-THEOREM*************
To prove this theorem, we first flatten the boundary part
$\Gamma_\Phi$ by defining the variables $(X,Y)=\Psi(x,y)$ with
$(X,Y)=(x, y-\varepsilon\Phi(x))$. Then $(x,y)=\Psi^{-1}(X,Y)=(X,
Y+\varepsilon\Phi(X))$. {}From (\ref{bdryNormCond}), we have
\begin{equation}\label{defFlatten-Norm}
\|\Psi-Id\|_{C^{2,\alpha}(\overline{\Omega_2^+})}+
\|\Psi^{-1}-Id\|_{C^{2,\alpha}(\overline{B_2^+})}\le
\varepsilon\lambda^{-1}.
\end{equation}
Then, for sufficiently small $\varepsilon$ depending only on
$\lambda$, the transformed domain $\Dom_2^+\defd\Psi(\Omega_2^+)$
satisfies
\begin{equation}\label{rescaledDomain}
\begin{array}{l}
\displaystyle B_{2-2\varepsilon/\lambda}^+\subset\Dom_2^+\subset
B_{2+2\varepsilon/\lambda}^+, \quad
\Dom_2^+\subset\bR^2_+\defd\{Y>0\}, \quad
\partial\Dom_2^+\cap\{Y=0\}=\Psi(\Gamma_\Phi);
\end{array}
\end{equation}
the function
$$
v(X,Y)=u(x,y):=u(\Psi^{-1}(X,Y))
$$
satisfies an equation of form (\ref{nonlinEq-appdx}) in $\Dom_2^+$
with (\ref{locEstElliptEq-i1})--(\ref{locEstElliptEq-i2}) and the
corresponding elliptic constants $\lambda/2$; and the boundary
condition for $v$ by an explicit calculation is
\begin{equation}\label{nonlinEq-appdx-nonObl-bc-fltnd}
v_X=\varepsilon(b(\Psi^{-1}(X,0))+\Phi'(X))v_Y+c(\Psi^{-1}(X,0))v
\qquad \mbox{on }\;\Dom_2^+\cap\{Y=0\},
\end{equation}
i.e., it is of form (\ref{nonlinEq-appdx-nonObl-bc}) with
(\ref{bdryCoefBds})  satisfied on $\overline{\Dom_2^+}$ with
elliptic constant $\lambda/4$. Moreover, by
(\ref{defFlatten-Norm})--(\ref{rescaledDomain}), it suffices for
this theorem to show the following estimate for $v(X,Y)$:
\begin{equation}\label{toProveThem-nontang}
\|v\|_{2,\alpha, B_{6/5}^+}\le C(\lambda,
M,\alpha)\big(\|v\|_{0,B_{2-2\varepsilon/\lambda}^+}+\|f\|_{\alpha,
B_{2-2\varepsilon/\lambda}^+}\big).
\end{equation}
That is, we can consider the equation in
$B_{2-2\varepsilon/\lambda}^+$ and condition
(\ref{nonlinEq-appdx-nonObl-bc-fltnd}) on
$\FltBdry_{2-2\varepsilon/\lambda}$ or, by rescaling, we can simply
consider our equation in $B_2^+$ and condition
(\ref{nonlinEq-appdx-nonObl-bc-fltnd}) on $\FltBdry_2\defd
B_2\cap\{Y=0\}$. In other words, without loss of generality, we can
assume $\Phi\equiv 0$ in the original problem.

For simplicity, we use the original notation $(x,y, u(x,y))$ to
replace the notation $(X, Y, v(X,Y))$. Then we assume that
 $\Phi\equiv 0$. Thus, equation
(\ref{nonlinEq-appdx}) is satisfied in the domain $B_2^+$, the
boundary condition (\ref{nonlinEq-appdx-nonObl-bc}) is prescribed on
$\FltBdry_2=B_2\cap\{y=0\}$, and conditions
(\ref{locEstElliptEq-i1})--(\ref{bdryCoefBds}) hold in $B_2^+$.
Also, we use the partially interior norms \cite[Eq.
4.29]{GilbargTrudinger} in the domain $B_2^+\cup \FltBdry_2$ with
the related distance function
%\begin{equation}\label{semiInterDist}
$d_{z}=\dist(z,\partial B_2^+\setminus \FltBdry_2)$.
%\end{equation}
The universal constant $C$ in the argument below depends only on
$\lambda$ and $M$, unless otherwise specified.

As in \cite[\S13.2]{GilbargTrudinger}, we introduce the functions
$w_i=D_iu$ for $i=1,2$. Then we conclude from equation
(\ref{nonlinEq-appdx}) that $w_1$ and $w_2$ are weak solutions of
the following equations of divergence form:
\begin{equation}
D_1\big({A_{11}\over A_{22}}D_1w_1+{2A_{12}\over
A_{22}}D_2w_1\big) +D_{22}w_1= D_1\big({f\over A_{22}}-{A_{1}\over A_{22}}D_1u
-{A_{2}\over A_{22}}D_2 u\big),
\label{eqnForDerivAppdx-1}
\end{equation}
\begin{equation}
D_{11}w_2+ D_2\big({2A_{12}\over A_{11}}D_1w_2+{A_{22}\over
A_{11}}D_2w_2\big) = D_2\big({f\over A_{11}}-{A_{1}\over A_{11}}D_1
u-{A_{2}\over A_{11}}D_2 u\big). \label{eqnForDerivAppdx-2}
\end{equation}

{}From (\ref{nonlinEq-appdx-nonObl-bc}), we have
\begin{equation}\label{bc-w1}
w_1=g \qquad\mbox{ on } \FltBdry_2,
\end{equation}
where
\begin{equation}\label{bc-w1-defG} g:=\varepsilon
bw_2+cu\qquad\mbox {for }\;B_2^+.
\end{equation}

We first obtain the following H\"{o}lder estimates of $D_1u$.
%*****************BEGIN LEMMA*************
\begin{lemma}\label{appDxnonObl-lemma}
There exist $\beta\in (0,\alpha]$ and $C>0$ depending only on
$\lambda$ such that, for any $z_0\in B_2^+\cup \FltBdry_2$,
\begin{equation}\label{scaledHolderAppdx-claim}
d_{z_0}^\beta[w_1]_{0,\beta, B_{d_{z_0}/16}({z_0})\cap B_2^+}\le
C(\|(Du, f)\|_{0,0,B_{d_{z_0}/2}({z_0})\cap B_2^+} +
d_{z_0}^\beta[g]_{0,\beta, B_{d_{z_0}/2}(z_0)\cap B_2^+}).
\end{equation}
\end{lemma}
%*****************END LEMMA*************
\Proof
%*****************BEGIN PROOF-LEMMA*************
We first prove that, for $z_1\in \FltBdry_2$ and $B_{2R}^+(z_1)\subset
B_2^+$,
\begin{equation}\label{scaledHolderAppdx-claim-1}
R^\beta[w_1]_{0,\beta, B_{R}^+({z_1})}\le
C(\|(Du, Rf)\|_{0,0,B_{2R}^+({z_1})}
+ R^\beta[g]_{0,\beta, B_{2R}^+(z_1)}).
\end{equation}
We rescale $u$, $w_1$, and $f$ in $B_{2R}^+(z_1)$ by defining
\begin{equation}\label{rescaledU}
\hat u(Z)={1\over 2R}u(z_1+2RZ),\quad
\hat f(Z)={2R}f(z_1+2RZ)
\qquad\mbox{ for }\;
Z\in B_1^+,
\end{equation}
and $\hat w_i=D_{Z_i}\hat u$. Then $\hat w_1$ satisfies an equation
of form (\ref{eqnForDerivAppdx-1}) in $B_1^+$ with $u$ replaced by
$\hat u$ whose coefficients $\hat A_{ij}$ and $\hat A_i$ satisfy
(\ref{locEstElliptEq-i1})--(\ref{locEstElliptEq-i2}) with unchanged
constants (this holds for (\ref{locEstElliptEq-i2}) since $R\le 1$).
Then, by the elliptic version of \cite[Theorem 6.33]{LiebermanBook}
stated in the parabolic setting (it can also be obtained by using
\cite[Lemma 4.6]{LiebermanBook} instead of \cite[Lemma
8.23]{GilbargTrudinger} in the proofs of \cite[Theorem 8.27,
8.29]{GilbargTrudinger} to achieve $\alpha=\alpha_0$ in
\cite[Theorem 8.29]{GilbargTrudinger}), we find constants
$\tilde\beta(\lambda)\in (0,1)$ and $C(\lambda)$ such that
$$
[\hat w_1]_{0,\beta, B_{1/2}^+}\le C(\|(D\hat u, \hat f)\|_{0,0,B_1^+} +
[\hat w_1]_{0,\beta, B_1\cap\{y=0\}})
$$
for $\beta=\min(\tilde\beta, \alpha)$. Rescaling back and using
(\ref{bc-w1}), we have (\ref{scaledHolderAppdx-claim-1}).

If $z_1\in B_2^+$  and $B_{2R}(z_1)\subset B_2^+$, then an argument
similar to the proof of (\ref{scaledHolderAppdx-claim-1}) by using
the interior estimates \cite[Theorem 8.24]{GilbargTrudinger} yields
\begin{equation}\label{scaledHolderAppdx-claim-2}
R^\beta[w_1]_{0,\beta, B_{R}({z_1})}\le
C\|(Du, Rf)\|_{0,0,B_{2R}({z_1})}.
\end{equation}

Now let $z_0=(x_0,y_0)\in B_2^+\cup \FltBdry_2$.
When $y_0\le d_{z_0}/8$, then, denoting
$z_0'=(x_0,0)$ and noting that $d_{z_0'}\ge d_{z_0}$,
 it is easy to check that
\begin{eqnarray*}
B_{d_{z_0}/16}(z_0)\cap B_2^+
\subset
B^+_{d_{z_0}/8}(z_0')
\subset
B_2^+,
\quad
B^+_{d_{z_0}/8}(z_0')
\subset
B_{d_{z_0}/2}(z_0)\cap B_2^+,
\end{eqnarray*}
and then applying (\ref{scaledHolderAppdx-claim-1}) with $z_1=z_0'$
and $R={d_{z_0}/8}\le 1$ and using the inclusions stated above yield
(\ref{scaledHolderAppdx-claim}). When $y_0\ge d_{z_0}/8$,
$B_{d_{z_0}/8}(z_0)\subset B_2^+$. Then applying
(\ref{scaledHolderAppdx-claim-2}) with $z_1=z_0$ and
$R={d_{z_0}/16}\le 1$ yields (\ref{scaledHolderAppdx-claim}).
\Endproof
%*****************END PROOF-LEMMA*************

Next, we make the H\"{o}lder estimates for $Du$. We first note that,
by (\ref{bdryCoefBds}) and (\ref{bc-w1-defG}), $g$ satisfies
\begin{eqnarray}
\label{estimateG-appdx-1} &&\qquad |Dg|\le
C(\varepsilon|D^2u|+|Du|+|u|)\qquad \mbox{in }\;B_2^+,
\\
&&\qquad [g]_{0,\beta, B_{d_z/2}(z)\cap B_2^+} \le
C\left(\varepsilon [Du]_{0,\beta, B_{d_z/2}(z)\cap B_2^+} +
\|u\|_{1,0,B_{d_z/2}(z)\cap B_2^+} \right).\label{estimateG-appdx-2}
\end{eqnarray}

%*****************BEGIN LEMMA*************
\begin{lemma}\label{appDxnonObl-DuHolder-lemma}
Let $\beta$ be as in Lemma {\rm \ref{appDxnonObl-lemma}}. Then there
exist $\varepsilon_0(\lambda)>0$ and $C(\lambda)>0$ such that, if
$0\le\varepsilon\le \varepsilon_0$,
\begin{eqnarray}
d_{z_0}^\beta[Du]_{0,\beta, B_{d_{z_0}/ 32}({z_0})\cap B_2^+}
&\le&
C(\|u\|_{1,0,B_{d_{z_0}/ 2}({z_0})\cap B_2^+} +\varepsilon
d_{z_0}^\beta[Du]_{0,\beta, B_{d_{z_0}/ 2}(z_0)\cap B_2^+}
\nonumber
\\
&&\qquad+\|f\|_{0,0,B_{d_{z_0}/ 2}({z_0})\cap B_2^+} )
\label{scaledHolderAppdx-fullDu-claim}
\end{eqnarray}
for any $z_0\in B_2^+\cup\FltBdry_2$.
\end{lemma}
%*****************END LEMMA*************

\Proof
%*****************BEGIN PROOF-LEMMA*************
The H\"{o}lder norm of $D_1u$ has been estimated in Lemma
\ref{appDxnonObl-lemma}. It remains to estimate $D_2u$. We follow
the proof of \cite[Theorem 13.1]{GilbargTrudinger}.

Fix $z_0\in B_2^+\cup\FltBdry_2$. In order to prove
(\ref{scaledHolderAppdx-fullDu-claim}), it suffices to show that,
for every $\hat z\in B_{d_{z_0}/32}({z_0})\cap B_2^+$ and every
$R>0$ such that $B_R(\hat z)\subset B_{d_{z_0}/16}({z_0})$, we have
\begin{equation}\label{MorreyEst}
\int_{B_R(\hat z)\cap B_2^+}|D^2u|^2dz\le {L^2\over
d_{z_0}^{2\beta}}R^{2\beta} ,
\end{equation}
where $L$ is the right-hand side of
(\ref{scaledHolderAppdx-fullDu-claim}) (cf. \cite[Theorem
7.19]{GilbargTrudinger} and \cite[Lemma 4.11]{LiebermanBook}).

In order to prove (\ref{MorreyEst}), we consider separately case
(i) $B_{2R}(\hat z)\cap \FltBdry_2\ne\emptyset$
and case
(ii) $B_{2R}(\hat z)\cap \FltBdry_2=\emptyset$.

We first consider case (i).
Let $B_{2R}(\hat z)\cap \FltBdry_2\ne\emptyset$. Since $B_R(\hat
z)\subset B_{d_{z_0}/32}({z_0})$, then $B_{2R}(\hat z)\subset
B_{d_{z_0}/16}({z_0})$ so that
\begin{equation}\label{R-le-d}
2R\le d_{z_0}.
\end{equation}
Let $\eta\in C^1_0(B_{2R}(\hat z))$ and $\zeta=\eta^2(w_1-g)$. Note
that $\zeta\in W^{1,2}_0(B_{2R}(\hat z)\cap B_2^+)$ by
(\ref{bc-w1}). We use $\zeta$ as a test function in the weak form of
(\ref{eqnForDerivAppdx-1}):
\begin{equation}\label{weakEqForDeriv}
\int_{B_2^+}{1\over A_{22}}\sum_{i,j=1}^2 A_{ij}D_iw_1D_j\zeta dz
=\int_{B_2^+}{1\over A_{22}}\big(\sum_{i=1}^2
A_{i}D_iu+f\big)D_1\zeta dz,
\end{equation}
and apply (\ref{locEstElliptEq-i1})--(\ref{locEstElliptEq-i2}) and
(\ref{estimateG-appdx-1}) to obtain
\begin{eqnarray}\label{intermEstNonObliq}
\\
&&\int_{B_2^+}|Dw_1|^2\eta^2 dz \le
 C\int_{B_2^+}\bigg(
\Big((\delta+\varepsilon)|Dw_1|^2+
\varepsilon|D^2u|^2\Big)\eta^2\nonumber
\\
&&\qquad\qquad\qquad\qquad
+({1\over\delta}+1)\left((|D\eta|^2+|f|\eta^2)(w_1-g)^2
+(|Du|^2+|u|^2)\eta^2\right)\bigg)dz, \nonumber
\end{eqnarray}
where $C$ depends only on $\lambda$, and the sufficiently small
constant $\delta>0$ will be chosen below. Since
\begin{equation}\label{intermEstNonObliq-1}
|Dw_1|^2=(D_{11}u)^2+(D_{12}u)^2,
\end{equation}
it remains to estimate $|D_{22}u|^2$. Using the ellipticity property
(\ref{locEstElliptEq-i1}), we can express $D_{22}u$ from equation
(\ref{nonlinEq-appdx}) to obtain
$$
\int_{B_2^+}|D_{22}u|^2\eta^2 dz\le
 C(\lambda)\int_{B_2^+}(|D_{11}u|^2+|D_{12}u|^2+|Du|^2)\eta^2 dz.
$$
Combining this with
(\ref{intermEstNonObliq})--(\ref{intermEstNonObliq-1}) and using
(\ref{locEstElliptEq-i2}) to estimate $|f|$ yield
\begin{eqnarray}\label{intermEstNonObliq-2}
 \\
 \int_{B_2^+}|D^2u|^2\eta^2 dz
&\le &
 C\int_{B_2^+}\Big(
(\varepsilon+\delta)|D^2u|^2\eta^2\nonumber
\\
\nonumber && \quad\,
+({1\over\delta}+1)\left((|D\eta|^2+\eta^2)(w_1-g)^2
+(|Du|^2+|u|^2)\eta^2\right)\Big)dz.
\end{eqnarray}
Choose $\varepsilon_0=\delta=(4C)^{-1}$.
Then, when $\varepsilon\in
(0,\varepsilon_0)$, we have
\begin{equation}\label{intermEstNonObliq-3}
\int_{B_2^+}|D^2u|^2\eta^2 dz\le C\int_{B_2^+}\left(
(|D\eta|^2+\eta^2)(w_1-g)^2 +(|Du|^2+|u|^2)\eta^2\right)\,dz.
\end{equation}

Now we make a more specific choice of $\eta$: In addition to
$\eta\in C^1_0(B_{2R}(\hat z))$, we assume that $\eta\equiv 1$ on
$B_{R}(\hat z)$, $0\le\eta\le 1$ on $\bR^2$, and $|D\eta|\le
{10/R}$. Also, since $B_{2R}(\hat z)\cap \FltBdry_2\ne\emptyset$,
then, for any fixed $z^*\in B_{2R}(\hat z)\cap \FltBdry_2$, we have
$|z-z^*|\le 2R$ for any $z\in B_{2R}(\hat z)$. Moreover,
$(w_1-g)(z^*)=0$ by (\ref{bc-w1}). Then, since $B_{2R}(\hat
z)\subset B_{d_{z_0}/16}({z_0})$, we find from
(\ref{scaledHolderAppdx-claim}), (\ref{estimateG-appdx-2}), and
(\ref{R-le-d}) that, for any $z\in B_{2R}(\hat z)\cap B_2^+$,
\begin{eqnarray*}
|(w_1-g)(z)|
&=& |(w_1-g)(z) -(w_1-g)(z^*)|
\le
|w_1(z)-w_1(z^*)|+|g(z)-g(z^*)|\\
&\le&
  {C\over d_{z_0}^\beta}\big(\|(Du,f)\|_{0,0,B_{d_{z_0}/2}({z_0})\cap B_2^+}
 +d_{z_0}^\beta[g]_{0,\beta, B_{d_{z_0}/2}(z_0)\cap B_2^+}
 \big)|z-z^*|^\beta
\\
&& +[g]_{0,\beta, B_{d_{z_0}/2}(z_0)\cap B_2^+}|z-z^*|^\beta
\\
&\le& C\big({1\over d_{z_0}^\beta}
\|(Du,f)\|_{0,0,B_{d_{z_0}/2}({z_0})\cap B_2^+} + \varepsilon
[Du]_{0,\beta, B_{d_{z_0}/2}(z_0)\cap B_2^+}
\\
&&\qquad +\|u\|_{0,0,B_{d_{z_0}/2}({z_0})\cap B_2^+} \big)R^\beta.
\end{eqnarray*}
Using this estimate and our choice of $\eta$, we obtain from
(\ref{intermEstNonObliq-3}) that
\begin{eqnarray*}
&&\int_{B_R(\hat z)\cap B_2^+}|D^2u|^2 dz\\
&&\le C\big({1\over d_{z_0}^{2\beta}}\|(Du,f)\|_{0,0,B_{d_{z_0}/
2}({z_0})\cap B_2^+}^2 +\varepsilon^2 [Du]_{0,\beta,
B_{d_{z_0}/2}(z_0)\cap B_2^+}^2\big)R^{2\beta}
\\
&&\quad +C\|u\|_{1,0,B_{d_{z_0}/2}({z_0})\cap
B_2^+}^2(R^{2\beta}+R^2),
\end{eqnarray*}
which implies  (\ref{MorreyEst}) for case (i).

\medskip
Now we consider case (ii): $\hat z\in B_2^+$ and $R>0$ satisfy
$B_R(\hat z)\subset B_{d_{z_0}/32}({z_0})$ and $B_{2R}(\hat z)\cap
\FltBdry_2=\emptyset$. Then $B_{2R}(\hat z)\subset
B_{d_{z_0}/16}({z_0})\cap B_2^+$. Let $\eta\in C^1_0(B_{2R}(\hat
z))$ and $\zeta=\eta^2(w_1-w_1(\hat z))$. Note that $\zeta\in
W^{1,2}_0(B_2^+)$ since $B_{2R}(\hat z)\subset B_2^+$. Thus we can
use $\zeta$ as a test function in (\ref{weakEqForDeriv}). Performing
the estimates similar to those that have been done to obtain
(\ref{intermEstNonObliq-3}), we have
\begin{equation}\label{intermEstNonObliq-3-pr}
\int_{B_2^+}|D^2u|^2\eta^2 dz\le
C(\lambda)\int_{B_2^+}\big((|D\eta|^2+\eta^2)(w_1-w_1(\hat z))^2
+|Du|^2\eta^2\big)\,dz.
\end{equation}
Choose $\eta\in C^1_0(B_{2R}(\hat z))$ so that $\eta\equiv 1$ on
$B_{R}(\hat z)$, $0\le\eta\le 1$ on $\bR^2$, and $|D\eta|\le {10/
R}$. Note that, for any $z\in B_{2R}(\hat z)$,
$$
|w_1(z)-w_1(\hat z)| \le C\big({1\over
d_{z_0}^\beta}\|(Du,f)\|_{0,0,B_{d_{z_0}/2}({z_0})\cap B_2^+}
+\varepsilon [Du]_{0,\beta, B_{d_{z_0}/2}(z_0)\cap
B_2^+}\big)R^\beta
$$
by (\ref{scaledHolderAppdx-claim}) since $B_{2R}(\hat z)\subset
B_{d_{z_0}/16}({z_0})\cap B_2^+$. Now we obtain (\ref{MorreyEst})
from (\ref{intermEstNonObliq-3-pr}) similar to that for case (i).
Then  Lemma \ref{appDxnonObl-DuHolder-lemma} is proved.
\Endproof
%*****************END PROOF-LEMMA*************

%*****************BEGIN LEMMA*************
\begin{lemma}\label{partIntSeminorm-est-lemma}
Let $\beta$ and $\varepsilon_0$ be as in Lemma {\rm
\ref{appDxnonObl-DuHolder-lemma}}. Then, for $\varepsilon\in
(0,\varepsilon_0)$, there exists $C(\lambda)$ such that
\begin{equation}\label{partIntSeminorm-est}
[u]_{1,\beta, B_2^+\cup\FltBdry_2}^*\le
C(\|u\|_{1,0,B_2^+\cup\FltBdry_2}^*
+\varepsilon [u]_{1,\beta, B_2^+\cup\FltBdry_2}^*
+\|f\|_{0,0,B_2^+}),
\end{equation}
where $[\cdot]^*$ and $\|\cdot\|^*$ denote the standard partially
interior seminorms and norms {\rm \cite[Eq.
4.29]{GilbargTrudinger}}.
\end{lemma}
%*****************END LEMMA*************
\Proof
%*****************BEGIN PROOF-LEMMA*************
Estimate (\ref{partIntSeminorm-est}) follows directly from Lemma
\ref{appDxnonObl-DuHolder-lemma} and an argument similar to the
proof of \cite[Theorem 4.8]{GilbargTrudinger}. Let $z_1, z_2\in
B_2^+$ with $d_{z_1}\le d_{z_2}$ (thus $d_{z_1, z_2}=d_{z_1}$) and
let $|z_1-z_2|\le d_{z_1}/64$. Then $z_2\in
B_{d_{z_0}/32}({z_0})\cap B_2^+$ and, by Lemma
\ref{appDxnonObl-DuHolder-lemma} applied to $z_0=z_1$, we find
\begin{eqnarray*}
d_{z_1, z_2}^{1+\beta}{|Du(z_1)-Du(z_2)|\over |z_1-z_2|^\beta} &\le&
C(d_{z_1}\|u\|_{1,0,B_{d_{z_1}/2}({z_1})\cap B_2^+} +\varepsilon
d_{z_1}^{1+\beta}[Du]_{0,\beta, B_{d_{z_1}/2}(z_1)\cap B_2^+}\\
&&\qquad
+\|f\|_{0,0,B_{d_{z_1}/ 2}({z_1})\cap B_2^+})
\\
&\le& C(\|u\|_{1,0,B_2^+\cup\FltBdry_2}^*
+\varepsilon [u]_{1,\beta, B_2^+\cup\FltBdry_2}^*
+\|f\|_{0,0,B_2^+}),
\end{eqnarray*}
where the last inequality holds since $2d_z\ge d_{z_1}$ for all
$z\in B_{d_{z_1}/2}(z_1)\cap B_2^+$. If $z_1, z_2\in B_2^+$
with $d_{z_1}\le d_{z_2}$ and $|z_1-z_2|\ge d_{z_1}/64$, then
\begin{eqnarray*}
d_{z_1, z_2}^{1+\beta}{|Du(z_1)-Du(z_2)|\over |z_1-z_2|^\beta} \le
64(d_{z_1}|Du(z_1)|+d_{z_2}|Du(z_2)|)\le
64\,\|u\|_{1,0,B_2^+\cup\FltBdry_2}^*.
\end{eqnarray*}
This completes the proof.
\Endproof
%*****************END PROOF-LEMMA*************

\medskip
Now we can complete the proof of Theorem
\ref{locEstElliptEq-non-oblique}.  For sufficiently small
$\varepsilon_0>0$ depending only on $\lambda$, when $\varepsilon\in
(0,\varepsilon_0)$, we use Lemma \ref{partIntSeminorm-est-lemma} to
obtain
\begin{equation}\label{partIntSeminorm-est-1}
[u]_{1,\beta, B_2^+\cup\FltBdry_2}^*\le
C(\lambda)(\|u\|_{1,0,B_2^+\cup\FltBdry_2}^* +\|f\|_{0,0,B_2^+}).
\end{equation}
We use the interpolation inequality \cite[Eq.
(6.89)]{GilbargTrudinger} to estimate
$$
 \|u\|_{1,0,B_2^+\cup\FltBdry_2}^*\le
 C(\beta,\delta)\|u\|_{0,B_2^+} +\delta[u]_{1,\beta,
 B_2^+\cup\FltBdry_2}^*
$$
for $\delta>0$. Since $\beta=\beta(\lambda)$, we choose sufficiently
small $\delta(\lambda)>0$ to find
\begin{equation}\label{partIntSeminorm-est-2}
\|u\|_{1,\beta, B_2^+\cup\FltBdry_2}^*\le
C(\lambda)(\|u\|_{0,0,B_2^+} +\|f\|_{0,0,B_2^+})
\end{equation}
from (\ref{partIntSeminorm-est-1}). In particular, we obtain a
global estimate in a smaller half-ball:
\begin{equation}\label{glob-norm-in-smaller-est-2}
\|u\|_{1,\beta, B_{9/5}^+}\le C(\lambda)(\|u\|_{0,0,B_2^+}
+\|f\|_{0,0,B_2^+}).
\end{equation}

We can assume $\beta\le\alpha$. Now we consider
(\ref{eqnForDerivAppdx-1}) as a linear elliptic equation
\begin{equation}\label{linEq-InTangThProof}
\sum_{i,j=1}^2D_i(a_{ij}(x,y)D_jw_1)=D_1F \qquad\mbox{in }
B_{9/5}^+,
\end{equation}
where $a_{ij}(x,y)=(A_{ij}/A_{22})(Du(x,y),x,y)$ for $i+j<4$,
$a_{22}=1$, and $F(x,y)=\big({ A_{1}}D_1 u+{ A_{2}}D_2 u +f\big)/
A_{22}$ with $(A_{ij},A_i)=(A_{ij},A_i)(Du(x,y),x,y)$. Then
(\ref{partIntSeminorm-est-2}), combined with
(\ref{locEstElliptEq-i2}), implies
\begin{equation}\label{normsCoefs-InTangThProof}
\|a_{ij}\|_{0,\beta, B_{9/5}^+}\le C(\lambda, M).
\end{equation}
{}From now on, $d_z$ denotes the distance related to the partially
interior norms in $B_{9/5}^+\cup\FltBdry_{9/5}$, i.e., for $z\in
B_{9/5}^+$, $d_z:=\dist(z,\partial B_{9/5}^+\setminus
\FltBdry_{9/5})$. Now, similar to the proof of Lemma
\ref{appDxnonObl-lemma}, we rescale equation
(\ref{linEq-InTangThProof}) and the Dirichlet condition
(\ref{bc-w1}) from the balls $B^+_R(z'_1)\subset B_{9/5}^+$ and
$B_R(z_1)\subset B_{9/5}^+$ with $R\le 1$ to $B=B_1^+$ or $B=B_1$,
respectively, by defining
$$
(\hat w_1, \hat g, \hat a_{ij})(Z)=(w_1, g, a_{ij})(z_1+RZ),\quad
\hat F(Z)=RF(z_1+RZ)\qquad\mbox{for }Z\in B.
$$
Then $ \sum_{i,j=1}^2D_i(\hat a_{ij}(x,y)D_j\hat w_1)=D_1\hat F $ in
$B$, the ellipticity of this rescaled equation is the same as that
for (\ref{linEq-InTangThProof}), and $\|\hat a_{ij}\|_{0,\beta,
B}\le C$ for $C=C(\lambda, M)$ in (\ref{normsCoefs-InTangThProof}),
where we have used $R\le 1$. This allows us to apply the local
$C^{1,\beta}$ interior and boundary estimates for the Dirichlet
problem \cite[Theorem 8.32, Corollary 8.36]{GilbargTrudinger} to the
rescaled problems  in the balls $B^+_{3d_{z_0}/8}(z_0')$ and
$B_{d_{z_0}/8}(z_0)$ as in Lemma \ref{appDxnonObl-lemma}. Then,
scaling back and multiplying by $d_{z_0}$, applying the covering
argument as in Lemma \ref{appDxnonObl-lemma}, and recalling the
definition of $F$, we obtain that, for any $z_0\in B_{9/5}^+\cup
\FltBdry_{9/5}$,
\begin{equation}\label{scaledHolderAppdx-claim-2nd}
\begin{array}{ll}
d_{z_0}^{2+\beta}[w_1]_{1,\beta, B_{d_{z_0}/16}({z_0})\cap
B_{9/5}^+} +d_{z_0}^2[w_1]_{1,0, B_{d_{z_0}/16}({z_0})\cap
B_{9/5}^+}\\
\displaystyle \quad \le C\Big(
d_{z_0}\|Du\|_{0,0,B_{d_{z_0}/2}({z_0})\cap B_{9/5}^+}
+d_{z_0}^{1+\beta}[u]_{1,\beta, B_{d_{z_0}/2}(z_0)\cap B_{9/5}^+}
+\|f\|_{0,\beta, B_{d_{z_0}/2}(z_0)\cap B_{9/5}^+}
 \\
\displaystyle\qquad\,\,\,\,\,\,\,\,
+d_{z_0}^{2+\beta}[g]_{1,\beta,B_{d_{z_0}/2}(z_0)\cap B_{9/5}^+} +
\sum_{k=0,1} d_{z_0}^{k+1}[g]_{k,0,B_{d_{z_0}/2}(z_0)\cap B_{9/5}^+}
%+d_{z_0}[g]_{0,0,B_{d_{z_0}/2}(z_0)\cap B_{9/5}^+}
\Big),
\end{array}
\end{equation}
where we have used $d_{z_0}<2$. Recall that
$Dw_1=(D_{11}u, D_{12}u)$. Expressing $D_{22} u$ from equation
(\ref{nonlinEq-appdx}) by using
(\ref{locEstElliptEq-i1})--(\ref{locEstElliptEq-i2}) and
(\ref{partIntSeminorm-est-2}) to estimate the H\"{o}lder norms of
$D_{22} u$, in terms of the norms of $D_{11}u, D_{22}u$, and $Du$,
and by using (\ref{bc-w1-defG}) and (\ref{bdryCoefBds}) to estimate
the terms involving $g$ in (\ref{scaledHolderAppdx-claim-2nd}), we
obtain from (\ref{scaledHolderAppdx-claim-2nd}) that, for every
$z_0\in B_{9/5}^+\cup \FltBdry_2$,
$$
\begin{array}{l}\displaystyle
d_{z_0}^{2+\beta}[D^2u]_{0,\beta, B_{d_{z_0}/16}({z_0})\cap B_{9/5}^+}
+d_{z_0}^2[D^2u]_{0,0, B_{d_{z_0}/16}({z_0})\cap B_{9/5}^+}
\\
\displaystyle \le C\Big(d_{z_0}\|Du\|_{C(B_{d_{z_0}/2}({z_0})\cap
B_{9/5}^+)} +d_{z_0}^{1+\beta}[u]_{1,\beta, B_{d_{z_0}/2}(z_0)\cap
B_{9/5}^+}\\
 \qquad\,\, +d_{z_0}\|u\|_{1,0, B_{d_{z_0}/2}(z_0)\cap B_{9/5}^+}
 +\|f\|_{0,\beta, B_{d_{z_0}/2}(z_0)\cap B_{9/5}^+} \\
\qquad\,\, +\varepsilon \big(d_{z_0}^{2+\beta}[D^2u]_{0,\beta,
B_{d_{z_0}/2}(z_0)\cap B_{9/5}^+} +d_{z_0}^{2}[D^2u]_{0,0,
B_{d_{z_0}/2}(z_0)\cap B_{9/5}^+}\big) \Big).
\end{array}
$$
{}From this estimate, the argument of Lemma
\ref{partIntSeminorm-est-lemma} implies
\begin{equation}\label{scaled-C-1-alpha}
\|u\|_{2,\beta, B_{9/5}^+\cup\FltBdry_{9/5}}^*\le
C\big(\|u\|_{1,\beta,B_{9/5}^+\cup\FltBdry_{9/5}}^*
+\varepsilon\|u\|_{2,\beta, B_{9/5}^+\cup\FltBdry_{9/5}}^*
+\|f\|_{0,\beta,B_{9/5}^+}\big).
\end{equation}
Thus, reducing $\varepsilon_0$ if necessary and using
(\ref{glob-norm-in-smaller-est-2}), we  conclude
\begin{equation}\label{partIntSeminorm-est-4}
\|u\|_{2,\beta, B_{9/5}^+\cup\FltBdry_{9/5}}^* \le C(\lambda,
M)(\|u\|_{0,B_2^+}+\|f\|_{0,\beta,B_2^+}).
\end{equation}
Estimate \eqref{partIntSeminorm-est-4} implies a global estimate in
a smaller ball and, in particular,
$ \|u\|_{1,\alpha, B_{8/5}^+} \le C(\lambda,
M)(\|u\|_{0,B_2^+}+\|f\|_{0,\beta,B_2^+}). $
Now we can repeat the argument, which leads from
(\ref{glob-norm-in-smaller-est-2}) to (\ref{partIntSeminorm-est-4})
with $\beta$ replaced by $\alpha$, in $B_{8/5}^+$ (and, in
particular, further reducing $\varepsilon_0$ depending only on
$(\lambda, M,\alpha)$) to obtain
$$
\|u\|_{2,\alpha, B_{8/5}^+\cup\FltBdry_{8/5}}^* \le C(\lambda, M,
\alpha)(\|u\|_{0,B_2^+}+\|f\|_{0,\alpha,B_2^+}),
$$
which implies (\ref{toProveThem-nontang}) and hence
(\ref{localEst-appdx-nonOblique-est}) for the original problem.
Theorem \ref{locEstElliptEq-non-oblique} is proved.
%*****************END PROOF-THEOREM*************

\medskip
Now we show that the estimates also hold for the Dirichlet problem.
%*****************BEGIN THEOREM*************
\begin{theorem}\label{locEstElliptEq-Dirichlet}
Let $\lambda>0$ and $\alpha\in (0,1)$. Let $\Phi\in
C^{2,\alpha}(\bR)$ satisfy {\rm (\ref{bdryNormCond})} and
$\Omega_R^+:=B_R\cap\{y>\Phi(x)\}$ for $R>0$. Let $u\in
C^2(\Omega_2^+)\cap C(\overline{\Omega_2^+})$ satisfy {\rm
(\ref{nonlinEq-appdx})} in $\Omega_2^+$ and
\begin{eqnarray}
&&u=g\qquad
\mbox{ on }\; \Gamma_\Phi\defd B_2\cap\{y=\Phi(x)\},
\label{nonlinEq-appdx-Dirichlet-bc}
\end{eqnarray}
where $A_{ij}=A_{ij}(Du, x,y)$ and $A_{i}=A_{i}(Du, x, y)$,
$i,j=1,2$, and $f=f(x,y)$ satisfy {\rm
(\ref{locEstElliptEq-i1})}--{\rm (\ref{locEstElliptEq-i2})}, and
$g=g(x,y)$ satisfies
\begin{equation}\label{Rhs-Dirichlet}
\|g\|_{C^{2,\alpha}(\overline{\Omega_2^+})}\le \lambda^{-1},
\end{equation}
with $(\lambda, \alpha)$ defined above. Assume that
$\|u\|_{C(\Omega_2^+)}\le M$. Then
\begin{equation}\label{localEst-appdx-Dirichlet-est}
\|u\|_{C^{2,\alpha}(\overline{\Omega_1^+})}\le C(\lambda, M)
(\|u\|_{C(\overline{\Omega_2^+})}+\|f\|_{C^{\alpha}(\overline{\Omega_2^+})}
+\|g\|_{C^{2,\alpha}(\overline{\Omega_2^+})}).
\end{equation}
\end{theorem}
%*****************END THEOREM*************
\Proof
%*****************BEGIN PROOF-THEOREM*************
By replacing $u$ with $u-g$, we can assume without loss of
generality that $g\equiv 0$. Also, by flattening the boundary as in
the proof of Theorem \ref{locEstElliptEq-non-oblique}, we can assume
$\Phi\equiv 0$. That is, we have reduced to the case when
(\ref{nonlinEq-appdx}) holds in $B_2^+$ and $u=0$ on  $\FltBdry_2$.
Thus, $u_x=0$ on $\FltBdry_2$. Then estimate
\eqref{localEst-appdx-Dirichlet-est} follows from Theorem
\ref{locEstElliptEq-non-oblique}. \Endproof
%*****************END PROOF-THEOREM*************

\medskip
We now derive the estimates for the oblique derivative problem.

%*****************BEGIN THEOREM*************
\begin{theorem}\label{locEstElliptEq-oblique}
%Let $B_2\subset\bR^2$,
Let $\lambda>0$ and $\alpha\in (0,1)$. Let $\Phi\in
C^{2,\alpha}(\bR)$ satisfy {\rm (\ref{bdryNormCond})} and
$\Omega_R^+:=B_R\cap\{y>\Phi(x)\}$ for $R>0$. Let $u\in
C^2(\Omega_2^+)\cap C^1(\overline{\Omega_2^+})$ satisfy
%(\ref{nonlinEq-appdx}) in $\Omega_2^+$ and
\begin{eqnarray}
\label{nonlinEq-appdx-Obl-eq}
&&A_{11}u_{xx}+2A_{12}u_{xy}+A_{22}u_{yy}+
A_{1}u_{x}+A_{2}u_{y}=0 \qquad \mbox{ in }\;\Omega_2^+,\\
&&b_1u_{x}+b_2u_{y}+c u=0\qquad \mbox{ on }\; \Gamma_\Phi\defd
B_2\cap\{y=\Phi(x)\}, \label{nonlinEq-appdx-Obl-bc}
\end{eqnarray}
where $A_{ij}=A_{ij}(Du, x,y)$ and $A_{i}=A_{i}(Du, x, y)$,
$i,j=1,2$, satisfy {\rm (\ref{locEstElliptEq-i1})}--{\rm
(\ref{locEstElliptEq-i2})}, and $b_i=b_i(x,y), i=1,2,$ and
$c=c(x,y)$ satisfy the following obliqueness condition and
$C^{1,\alpha}$--bounds:
\begin{eqnarray}\label{obliqueness-appdx}
&&b_2(x,y)\ge\lambda\qquad\mbox{for }\;(x,y)\in \Gamma_\Phi,
\\
\label{bdryCoefBds-obliq} &&\|(b_1, b_2,
c)\|_{C^{1,\alpha}(\overline{\Omega_2^+})}\le \lambda^{-1}.
\end{eqnarray}
Assume that $\|u\|_{C(\overline{\Omega_2^+})}\le M$. Then there
exists $C=C(\lambda, M, \alpha)>0$ such that
\begin{equation}\label{localEst-appdx-Oblique-est}
\|u\|_{C^{2,\alpha}(\overline{\Omega_1^+})}\le
C\|u\|_{C(\overline{\Omega_2^+})}.
%+\|f\|_{C^{\alpha}(\overline{\Omega_2^+})}).
\end{equation}
\end{theorem}
%*****************END THEOREM*************

\Proof
%*****************BEGIN PROOF-THEOREM*************
{\em Step 1.} First, we flatten the boundary $\Gamma_\Phi$ by the
change of coordinates $(X,Y)=\Psi(x,y)=(x, y-\Phi(x))$. Then
$(x,y)=\Psi^{-1}(X,Y)=(X, Y+\Phi(X))$. {}From (\ref{bdryNormCond}),
$
\|\Psi\|_{C^{2,\alpha}(\Omega_2^+)}+\|\Psi^{-1}\|_{C^{2,\alpha}(\Dom_2^+)}\le
C(\lambda),$ where $\Dom_2^+\defd\Psi(\Omega_2^+)$ satisfies
$ \Dom_2^+\subset \bR^2_+\defd\{Y>0\}
$
and
$
\Gamma_0\defd\partial
\Dom_2^+\cap\{Y=0\}=\Psi(\Gamma_\Phi).
$
%\end{equation}
By a standard calculation, $v(X,Y)=u(x,y):=u(\Psi^{-1}(X,Y))$
satisfies the equation of form (\ref{nonlinEq-appdx-Obl-eq}) in
$\Dom_2^+$ and the oblique derivative condition of form
(\ref{nonlinEq-appdx-Obl-bc}) on $\Gamma_0$, where
(\ref{locEstElliptEq-i1})--(\ref{locEstElliptEq-i2}) and
(\ref{obliqueness-appdx})--(\ref{bdryCoefBds-obliq}) are satisfied
with modified constant $\hat \lambda>0$ depending only on $\lambda$.
Also, $\|v\|_{C(\Dom_2^+)}\le M$. Thus,
(\ref{localEst-appdx-Oblique-est}) follows from
\begin{equation}\label{partIntSeminorm-est-obl}
\|v\|_{2,\alpha, \Dom_2^+\cup\Gamma_0}^*\le C(\lambda,
M,\alpha)\|v\|_{0,\Dom_2^+}.
\end{equation}

Next we note that, in order to prove
(\ref{partIntSeminorm-est-obl}), it suffices to prove that there
exist $K$ and $C$ depending only on $(\lambda, M, \alpha)$ such
that, if $v$ satisfies
(\ref{nonlinEq-appdx-Obl-eq})--(\ref{nonlinEq-appdx-Obl-bc}) in
$B^+_1$ and $\FltBdry_1\defd B_1\cap\{y=0\}$ respectively,
(\ref{locEstElliptEq-i1})--(\ref{locEstElliptEq-i2}) and
(\ref{obliqueness-appdx})--(\ref{bdryCoefBds-obliq}) hold in
$B^+_1$, and $|v|\le M$ in $B_1^+$, then
\begin{equation}\label{appdx-Oblique-est-scaled-localized}
\|v\|_{C^{2,\alpha}(\overline{B_{1/K}^+})}\le C\|v\|_{C(B_1^+)}.
\end{equation}
Indeed, if (\ref{appdx-Oblique-est-scaled-localized}) is proved,
then, using also the interior estimates (\ref{localEst-appdx-2}) in
Theorem \ref{locEstElliptEq}  and applying the scaling argument
similar to the proof of Lemma \ref{appDxnonObl-lemma}, we obtain
that, for any $z_0\in \Dom_2^+\cup \FltBdry_2$,
$$
d_{z_0}^{2+\alpha}\|v\|_{C^{2,\alpha}(\overline{B_{d_{z_0}/(16K)}({z_0})\cap
\Dom_2^+})}\le C\|v\|_{C(B_{d_{z_0}/2}({z_0})\cap \Dom_2^+)}.
$$
{}From this, we use the argument of the proof of Lemma
\ref{partIntSeminorm-est-lemma} to obtain
(\ref{partIntSeminorm-est-obl}).

Thus it remains to show (\ref{appdx-Oblique-est-scaled-localized}).
First we make a linear change of variables to normalize the problem
so that
\begin{equation}\label{normalize-obliq-cond}
b_1(0)=0,\quad b_2(0)=1
\end{equation}
for the modified problem. Let
$$(X,Y)=\tilde\Psi(x,y):={1\over
b_2(0)}(b_2(0)x-b_1(0)y,y).
$$
Then
$$
(x,y)=\tilde\Psi^{-1}(X,Y)=(X+b _1(0)Y, b_2(0)Y), \qquad
|D\tilde\Psi|+|D\tilde\Psi^{-1}|\le C(\lambda),
$$
where the estimate follows from
(\ref{obliqueness-appdx})--(\ref{bdryCoefBds-obliq}). Then the
function $w(X,Y)\defd v(x,y)\equiv v(X+b _1(0)Y, b_2(0)Y)$ is a
solution of the equation of form (\ref{nonlinEq-appdx-Obl-eq}) in
the domain $\tilde\Psi(B_1^+)$ and the boundary condition of form
(\ref{nonlinEq-appdx-Obl-bc}) on the boundary part
$\tilde\Psi(\FltBdry_1)$ such that
(\ref{locEstElliptEq-i1})--(\ref{locEstElliptEq-i2}) and
(\ref{obliqueness-appdx})--(\ref{bdryCoefBds-obliq}) are satisfied
with constant $\hat \lambda>0$ depending only on $\lambda$, and
(\ref{normalize-obliq-cond}) holds, which can be verified by a
straightforward calculation. Also, $\|w\|_{C(\tilde\Psi(B_1^+))}\le
M$.

Note that $\tilde\Psi(B_1^+)\subset \bR^2_+:=\{Y>0\}$ and
$\tilde\Psi(\FltBdry_1)=\partial\tilde\Psi(B_1^+)\cap\{Y=0\}$.
Moreover, since $|D\tilde\Psi|+|D\tilde\Psi^{-1}|\le C(\lambda)$,
there exists $K_1=K_1(\lambda)>0$ such that, for any $r>0$, $
B_{r/K_1}\subset \tilde\Psi(B_r)\subset B_{K_1r}$.  Thus it suffices
to prove
$$
\|w\|_{C^{2,\alpha}(\overline{B_{r/2}^+})}\le C\|w\|_{C(B_r^+)}
$$
for some $r\in (0, 1/K_1)$. This estimate implies (\ref{appdx-Oblique-est-scaled-localized})
with $K=2K_1/r$.

\medskip
{\em Step 2.} As a result of the reduction performed in Step 1, it
suffices to prove the following: There exist $\varepsilon\in (0, 1)$
and $C$ depending only on $(\lambda, \alpha, M)$ such that, if $u$
satisfies (\ref{nonlinEq-appdx-Obl-eq}) and
(\ref{nonlinEq-appdx-Obl-bc}) in $B^+_{2\varepsilon}$ and on
$\FltBdry_{2\varepsilon}$ respectively, if
(\ref{locEstElliptEq-i1})--(\ref{locEstElliptEq-i2}) and
(\ref{obliqueness-appdx})--(\ref{bdryCoefBds-obliq}) hold in
$B^+_{2\varepsilon}$, and if (\ref{normalize-obliq-cond}) holds and
$\|u\|_{0,B_{2\varepsilon}^+}\le M$, then
$$
\|u\|_{2,\alpha,B_\varepsilon^+}\le C\|u\|_{0,B_{2\varepsilon}^+}.
$$

We now prove this claim. For $\varepsilon>0$ to be chosen later, we
rescale from $B_{2\varepsilon}^+$ into $B_2^+$ by defining
\begin{equation}\label{def-of-v-obliq}
v(x,y)={1\over \varepsilon}\big(u(\varepsilon x,\varepsilon
y)-u(0,0)\big)\qquad \mbox{for }\;(x,y)\in B_2^+.
\end{equation}
Then $v$ satisfies
\begin{eqnarray}\label{nonlinEq-appdx-Obl-eq-rescaled}
&&\tilde A_{11}v_{xx}+2\tilde A_{12}v_{xy}+\tilde A_{22}v_{yy}+
\tilde A_{1}v_{x}+\tilde A_{2}v_{y}=\tilde f \qquad \mbox{ in }\;B_2^+,\\
&&v_{y}=\tilde b_1 v_{x}+ \tilde b_2 v_{y}+\tilde c v+\tilde c u(0,0)\qquad
\mbox{ on }\; \FltBdry_2,
\label{nonlinEq-appdx-Obl-bc-rescaled}
\end{eqnarray}
where
\begin{eqnarray*}
&&\tilde A_{ij}(p,x,y)=A_{ij}(p,\varepsilon x,\varepsilon y), \quad
  \tilde A_{i}(p,x,y)=\varepsilon A_{i}(p,\varepsilon x,\varepsilon y),\\
&&\tilde b_1(x,y)=-b_1(\varepsilon x,\varepsilon y), \quad \tilde
b_2(x,y)=-b_2(\varepsilon x,\varepsilon y)+1, \quad \tilde
c(x,y)=-\varepsilon c(\varepsilon x,\varepsilon y).
\end{eqnarray*}
Then $\tilde A_{ij}$ and $\tilde A_i$ satisfy
(\ref{locEstElliptEq-i1})--(\ref{locEstElliptEq-i2}) in $B_2^+$ and,
using (\ref{bdryCoefBds-obliq}), (\ref{normalize-obliq-cond}), and
$\varepsilon\le 1$,
\begin{equation}\label{bdryCoefBds-obliq-norm}
\|(\tilde b_1, \tilde b_2, \tilde c)\|_{1,\alpha, B_2^+}\le
C\varepsilon \qquad\text{for some } C=C(\lambda).
\end{equation}

Now we follow the proof of Theorem \ref{locEstElliptEq-non-oblique}.
We use the partially interior norms \cite[Eq.
4.29]{GilbargTrudinger} in the domain $B_2^+\cup \FltBdry_2$ whose
distance function is
$d_{z}=\dist(z,\partial B_2^+\setminus \FltBdry_2).$
We introduce the functions $w_i=D_iv$, $i=1,2$, to conclude from
(\ref{nonlinEq-appdx-Obl-eq-rescaled}) that $w_1$ and $w_2$ are weak
solutions of equations
\begin{eqnarray}
&&\qquad\quad D_1\big({\tilde A_{11}\over \tilde
A_{22}}D_1w_1+{2\tilde A_{12}\over \tilde A_{22}} D_2w_1\big)
+D_{22}w_1= -D_1\big({\tilde A_{1}\over \tilde A_{22}}D_1 v+{\tilde
A_{2}\over \tilde A_{22}}D_2 v\big), \label{eqnForDerivAppdx-1-obl}
\\
&&\qquad\quad D_{11}w_2+ D_2\big({2\tilde A_{12}\over \tilde
A_{11}}D_1w_2+{\tilde A_{22}\over \tilde A_{11}}D_2w_2\big) =
-D_2\big({\tilde A_{1}\over \tilde A_{11}}D_1 v+{\tilde A_{2}\over
\tilde A_{11}}D_2 v\big) \label{eqnForDerivAppdx-2-obl}
\end{eqnarray}
in $B_2^+$, respectively. {}From
(\ref{nonlinEq-appdx-Obl-bc-rescaled}), we have
\begin{equation}\label{bc-w2}
w_2=\tilde g \qquad\mbox{ on } \FltBdry_2,
\end{equation}
where
$\tilde g:=\tilde b_1 v_{x}+ \tilde b_2 v_{y}+\tilde c v+\tilde c
u(0,0)$
in $B_2^+$.

Using equation (\ref{eqnForDerivAppdx-2-obl}) and the Dirichlet
boundary condition (\ref{bc-w2}) for $w_2$ and following the proof
of Lemma \ref{appDxnonObl-lemma}, we can show the existence of
$\beta\in (0,\alpha]$ and C depending only on $\lambda$ such that,
for any $z_0\in B_2^+\cup \FltBdry_2$,
\begin{equation}\label{scaledHolderAppdx-obl-w2}
d_{z_0}^\beta[w_2]_{0,\beta, B_{d_{z_0}/16}({z_0})\cap B_2^+}\le
C\big(\|Dv\|_{0,B_{d_{z_0}/2}({z_0})\cap B_2^+} +
d_{z_0}^\beta[\tilde g]_{0,\beta, B_{d_{z_0}/2}(z_0)\cap
B_2^+}\big).
\end{equation}

Next we obtain the H\"{o}lder estimates of $Dv$ if $\varepsilon$ is
sufficiently small. We first note that, by
(\ref{bdryCoefBds-obliq-norm}), $\tilde g$ satisfies
\begin{eqnarray}
&&|D\tilde g|\le C\varepsilon(|D^2v|+|Dv|+
|v|+
\|u\|_{0,B_{2\varepsilon}^+})\qquad
\mbox{in }\;B_2^+,
\label{estimateG-obl-appdx-1}\\
&&[\tilde g]_{0,\beta, B_{d_z/2}(z)\cap\Dom_2^+}
\le C\varepsilon ( \|v\|_{1,\beta,B_{d_z/2}(z)\cap\Dom_2^+)} +
\|u\|_{0,B_{2\varepsilon}^+})
\qquad\quad\label{estimateG-obl-appdx-2}
\end{eqnarray}
for $C=C(\lambda)$. The term $\varepsilon
\|u\|_{0,B_{2\varepsilon}^+}$ in
(\ref{estimateG-obl-appdx-1})--(\ref{estimateG-obl-appdx-2}) comes
from the term $\tilde c u(0,0)$ in the definition of $\tilde{g}$. We
follow the proof of Lemma \ref{appDxnonObl-DuHolder-lemma}, but we
now use the integral form of equation (\ref{eqnForDerivAppdx-2-obl})
with test functions $\zeta=\eta^2(w_2-\tilde g)$ and
$\zeta=\eta^2(w_2-w_2(\hat z))$ to get an integral estimate of
$|Dw_2|$ and thus of $|D_{ij}v|$ for $i+j>2$, and then use
(\ref{nonlinEq-appdx-Obl-eq-rescaled}) to estimate the remaining
derivative $D_{11}v$. In these estimates, we use
(\ref{scaledHolderAppdx-obl-w2})--(\ref{estimateG-obl-appdx-2}). We
obtain that, for sufficiently small  $\varepsilon$ depending only on
$\lambda$,
\begin{equation}\label{scaledHolderAppdx-fullDu-obliq}
\begin{array}{ll}
\displaystyle d_{z_0}^\beta[Dv]_{0,\beta,
B_{d_{z_0}/32}({z_0})\cap
B_2^+}\\
\le C\big(\|v\|_{C^1(B_{d_{z_0}/2}({z_0})\cap B_2^+)} +\varepsilon
d_{z_0}^\beta[Dv]_{0,\beta, B_{d_{z_0}/2}(z_0)\cap \Dom_2^+}
+\varepsilon d_{z_0}^\beta\|u\|_{0,B_{2\varepsilon}^+}\big)
\end{array}
\end{equation}
for any $z_0\in B_2^+\cup \FltBdry_2$, with $C=C(\lambda)$. Using
(\ref{scaledHolderAppdx-fullDu-obliq}), we follow the proof of Lemma
\ref{partIntSeminorm-est-lemma} to obtain
\begin{equation}\label{partIntSeminorm-est-obliq}
[v]_{1,\beta, B_2^+\cup \FltBdry_2}^*\le C\big(\|v\|_{1,0,B_2^+\cup
\FltBdry_2}^* +\varepsilon [v]_{1,\beta, B_2^+\cup \FltBdry_2}^*
+\varepsilon \|u\|_{0,B_{2\varepsilon}^+}\big).
\end{equation}
%where $C=C(\lambda, \Lambda)$.
Now we choose sufficiently small $\varepsilon>0$ depending only on
$\lambda$ to have
$$
[v]_{1,\beta, B_2^+\cup \FltBdry_2}^*\le
C(\lambda)(\|v\|_{1,0,B_2^+\cup \FltBdry_2}^* +
\|u\|_{0,B_{2\varepsilon}^+}).
$$
Then we use the interpolation inequality, similar to the proof of
(\ref{partIntSeminorm-est-2}), to have
\begin{equation}\label{partIntSeminorm-est-2-obl}
\|v\|_{1,\beta, B_2^+\cup \FltBdry_2}^*\le
C(\lambda)(\|v\|_{0,B_2^+}+\|u\|_{0,B_{2\varepsilon}^+}).
\end{equation}
By (\ref{def-of-v-obliq}) with $\varepsilon=\varepsilon(\lambda)$
chosen above, (\ref{partIntSeminorm-est-2-obl}) implies
\begin{equation}\label{partIntSeminorm-est-3-obl}
\|u\|_{1,\beta, B_{2\varepsilon}^+\cup B_{2\varepsilon}^0}^*\le
C(\lambda)\|u\|_{0,B_{2\varepsilon}^+}.
\end{equation}
Then problem
(\ref{nonlinEq-appdx-Obl-eq})--(\ref{nonlinEq-appdx-Obl-bc}) can be
regarded as a linear oblique derivative problem in
$B_{7\varepsilon/4}^+$ whose coefficients
$a_{ij}(x,y):=A_{ij}(Du(x,y),x,y)$ and $a_{i}(x,y)$
$:=A_{i}(Du(x,y),x,y)$ have the estimate in
$C^{0,\beta}(\overline{B_{7\varepsilon/4}^+})$ by a constant
depending only on $(\lambda, M)$ from
(\ref{partIntSeminorm-est-3-obl}) and (\ref{locEstElliptEq-i2}).
Moreover, we can assume $\beta\le\alpha$ so that
(\ref{bdryCoefBds-obliq}) implies the estimates of $(b_i, c)$ in
$C^{1,\beta}(\overline{B_{7\varepsilon/4}^+})$ with
$\varepsilon=\varepsilon(\lambda)$. Then the standard estimates for
linear oblique derivative problems \cite[Lemma
6.29]{GilbargTrudinger} imply
\begin{equation}\label{partIntSeminorm-est-4-obl}
\|u\|_{2,\beta,B_{3\varepsilon/2}^+}\le
C(\lambda,M)\|u\|_{0,B_{7\varepsilon/4}^+}.
\end{equation}
In particular, the
$C^{0,\alpha}(\overline{B_{3\varepsilon/2}^+})$--norms of the
coefficients $(a_{ij}, a_i)$ of the linear equation
(\ref{nonlinEq-appdx-Obl-eq}) are bounded
 by a constant
depending only on $(\lambda, M)$, which implies
$$
\|u\|_{2,\alpha,B_{\varepsilon}^+}\le C(\lambda,
M)\|u\|_{0,B_{3\varepsilon/2}^+},
$$
by applying again \cite[Lemma 6.29]{GilbargTrudinger}. This implies
the assertion of Step 2, thus Theorem \ref{locEstElliptEq-oblique}.
\Endproof
%*****************END PROOF-THEOREM*************

\bigskip
\medskip
\noindent {\bf Acknowledgments}. Gui-Qiang Chen's research was
supported in part by the National Science Foundation under Grants
DMS-0505473, DMS-0244473, and an Alexander von Humboldt Foundation
Fellowship. Mikhail Feldman's research was supported in part by the
National Science Foundation under Grants DMS-0500722 and
DMS-0354729.

\end{document}